\documentclass[hidelinks,onefignum,onetabnum]{siamart_251216/siamart251216}

\usepackage{amsmath} 
\usepackage{amssymb}  
\usepackage{graphicx}
\usepackage{latexsym}
\usepackage{verbatim}
\usepackage{ifthen}
\usepackage{psfrag}
\usepackage{macros/subfigure}

\newtheorem{nntheorem}{ Theorem}[section]

\newtheorem{nnconjecture}[nntheorem]{ Conjecture}

\newtheorem{nnassumption}[nntheorem]{ Assumption}
\newtheorem{nexample}[nntheorem]{ Example}
\newtheorem{nnremark}[nntheorem]{ Remark}\newtheorem{nnapplication}[nntheorem]{ Application}
\newtheorem{nnproblem}[nntheorem]{ Exercise}

\newtheorem{nntheorem2}{ Theorem}
\newtheorem{app}[nntheorem2]{ Application}

\newenvironment{ttheorem}[1]
{\begin{nntheorem}{\rm\textrm{(#1)}}\sl}
{\end{nntheorem}}

\newenvironment{assumption}[1]
{\begin{nnassumption}{\rm\textrm{
\ifx\hfuzz#1\hfuzz
\else
  (#1)
\fi
}}\sl}
{\end{nnassumption}}

\newenvironment{application}[1]
{\begin{nnapplication}{\rm\textrm{(#1)}}\sl}
{\end{nnapplication}}

{\begin{app}{\rm\textrm{(#1)}}\sl}
{\end{app}}

\newenvironment{remark}[1]
{\begin{nnremark}{\rm\textrm{(#1)}}\sl}
{\end{nnremark}}

\newcommand{\eoe}
           {\hspace*{\fill}{$\vcenter{\hrule height1pt 
                     \hbox{\vrule width1pt height3pt 
            \kern3pt \vrule width1pt} \hrule height1pt}$} }

\newenvironment{example}[1]
{\begin{nexample}{\rm\textrm{(#1)}}\rm}{\eoe
\end{nexample}}

\newcommand{\eop}
           {\hspace*{\fill}{$\vcenter{\hrule height1pt 
                     \hbox{\vrule width1pt height5pt 
            \kern5pt \vrule width1pt} \hrule height1pt}$} }

\newenvironment{proofs}
{\par\noindent\textbf{Proof Sketch.}}{\eop\smallskip\vskip 3 pt}

\newcommand{\dom}{\mathop{\rm dom}\nolimits}

\newcommand{\reals}{{\mathbb R}}

\newcommand{\tto}{\;{\lower 1pt \hbox{$\rightarrow$}}\kern -12pt
           \hbox{\raise 2pt \hbox{$\rightarrow$}}\;}

\newcommand{\A}{\mathcal{A}}

\newcommand{\HS}{\mathcal{H}}

\newcommand{\CUD}{\overline{C}\cup D}

\DeclareMathOperator*{\argmax}{arg\,max}
\DeclareMathOperator*{\argmin}{arg\,min}

\newcommand{\R}{{\mathcal{R}}}
\newcommand{\realsgeq}{{\reals_{\geq 0}}}
\newcommand{\nats}{\mathbb{N}}      

\newcommand{\J}{{\mathcal{J}}}

\usepackage{array}
\usepackage{mathtools}
\usepackage{eqparbox}
\usepackage{dsfont}
\usepackage{pstool}
\usepackage{rotating,graphics,psfrag,epsfig}
\usepackage{macros/scalerel}
\usepackage{balance}
\usepackage[normalem]{ulem}
\usepackage{url}

\usepackage{amsmath}
\usepackage[colorinlistoftodos,prependcaption,textsize=tiny]{todonotes}
\usepackage{varwidth}
\usepackage{parskip}

\usepackage{enumitem}
\usepackage{float}

\usepackage{multicol}
\usepackage{lipsum}
\usepackage{nccmath}

\usepackage{booktabs,arydshln,xcolor}
\newcolumntype{L}{>{\centering\arraybackslash}m{2.1cm}}
\newcolumntype{V}{>{\centering\arraybackslash}m{1.6cm}}
   
\usepackage{amsfonts}              
\usepackage{mathrsfs}
\usepackage{indentfirst}   
\usepackage{setspace}
\usepackage{ifpdf}
\usepackage{subfig}
\allowdisplaybreaks

\newcommand{\EatDot}[1]{} 

\newcommand{\pn}[1]{{\color{black} #1}}
\newcommand{\pnn}[1]{{\color{black} #1}}
\newcommand{\pno}[1]{{\color{black} #1}}
\newcommand{\sj}[1]{{\color{black} #1}}

\usepackage{ifthen}
\newboolean{Conference}
\setboolean{Conference}{false}
\newcommand{\NotConf}[1]{\ifthenelse{\boolean{Conference}}{}{{\color{black}#1}}}
\newcommand{\IfConf}[2]{\ifthenelse{\boolean{Conference}}{#1}{{\color{black}#2}}}   

\newboolean{Inclusion}
\setboolean{Inclusion}{true}
\newcommand{\IfInc}[2]{\ifthenelse{\boolean{Inclusion}}{ \color{cyan}#1\color{black}#2}{#2}}
\newcommand{\IfIncd}[2]{\ifthenelse{\boolean{Inclusion}}{{#1}}{{#2}}}

\newboolean{Personal}
\setboolean{Personal}{false}
\newcommand{\IfPers}[1]{\ifthenelse{\boolean{Personal}}{{\color{purple}#1} }{}}

\newboolean{Infinite-horizon} 
\setboolean{Infinite-horizon}{true} 
\newcommand{\IfIh}[2]{\ifthenelse{\boolean{Infinite-horizon}}{#1}{{\color{purple}#2}}}

\newboolean{Two-Players} 
\setboolean{Two-Players}{true} 
\newcommand{\IfTp}[2]{\ifthenelse{\boolean{Infinite-horizon}}{#1}{{\color{purple}#2}}}

\newboolean{Journal}
\setboolean{Journal}{true}
\newcommand{\IfJ}[1]{\ifthenelse{\boolean{Journal}}{{\color{purple}#1} }{}}

\newboolean{Automatica}
\setboolean{Automatica}{true}
\newcommand{\NotAutom}[1]{\ifthenelse{\boolean{Automatica}}{}{{\color{blue}#1}}}
\newcommand{\IfAutom}[2]{\ifthenelse{\boolean{Automatica}}{#1}{{\color{blue}#2}}}   

\newboolean{SIAM}
\setboolean{SIAM}{false}
\newcommand{\NotSiam}[1]{\ifthenelse{\boolean{SIAM}}{}{{\color{black}#1}}}
\newcommand{\IfSiam}[2]{\ifthenelse{\boolean{SIAM}}{{\color{black}#1}}{{\color{black}#2}}} 

\newboolean{MaxSys}
\setboolean{MaxSys}{true}
\newcommand{\NotMS}[1]{\ifthenelse{\boolean{MaxSys}}{}{{\color{blue}#1}}}
\newcommand{\IfMS}[2]{\ifthenelse{\boolean{MaxSys}}{{#1}}{{\color{blue}#2}}}   

\newboolean{ACC24}
\setboolean{ACC24}{false}
\newcommand{\NotACC}[1]{\ifthenelse{\boolean{ACC24}}{}{{\color{black}#1}}}
\newcommand{\IfACC}[2]{\ifthenelse{\boolean{ACC24}}{#1}{{\color{black}#2}}}

\newboolean{GenGrad}
\setboolean{GenGrad}{true}
\newcommand{\NotGG}[1]{\ifthenelse{\boolean{GenGrad}}{}{{#1}}}
\newcommand{\IfGG}[2]{\ifthenelse{\boolean{GenGrad}}{{\color{black}#1}}{{\color{black}#2}}} 

\newboolean{Comments}
\setboolean{Comments}{false}

\title{Two-Player Zero-Sum Hybrid Games \pnn{with Set-Valued Dynamics}
\thanks{Submitted to the editors DATE.
\funding{
Research partially supported by NSF Grants no. CNS-2039054 and CNS-2111688, by AFOSR Grants no. FA9550-19-1-0169, FA9550-20-1-0238, and FA9550-23-1-0145, by AFRL Grant nos. FA8651-22-1-0017 and FA8651-23-1-0004, and by ARO Grant no. W911NF-20-1-0253. 
}}}
  \author{Santiago J. Leudo\thanks{Electrical and Computer Engineering Department, University of California, Santa Cruz, CA 95064, USA. (
  {sjimen28@ucsc.edu}, 
{ricardo@ucsc.edu}).}
\and Francesco Ferrante\thanks{Department of Applied Mathematics, Fictional University, Boise, ID. (
{francesco.ferrante@unipig.it}).}
\and Ricardo G. Sanfelice\footnotemark[2]
}
\begin{document}
\maketitle
%
\begin{abstract}
In this paper, we formulate a two-player zero-sum game under dynamic constraints defined in terms of a hybrid inclusion. 
The game consists of a min-max problem involving a cost functional that depends on the actions and corresponding (potentially non-unique) solutions, defined as functions of hybrid time and, hence, can flow or jump.
We present sufficient conditions given in terms of Hamilton-Jacobi-Isaacs-like equations to establish a bound and to exactly evaluate 
   the worst-case cost under the optimal strategy. 
Under additional conditions, we show that the proposed optimal state-feedback laws render a set of interest pre-asymptotically stable for the resulting hybrid closed-loop system. The results are illustrated in an application with sporadic jumps and linear quadratic hybrid dynamics.
\IfIncd{}{Applications of this problem, as presented here, include a disturbance rejection scenario for which the effect of the perturbation is upper bounded, and a security scenario in which
we formulate an optimal control problem under the presence of the maximizing adversarial action.}
%
%
\end{abstract}

\IfConf{
\maketitle
}{}

\IfPers{A game is an optimization problem with multiple decision makers (players), a set of constraints (potentially dynamic) that enforces the ``rules'' of the game, and a set of payoff functions to be optimized by selecting decision variables. Constraints on the {state and} decision variables formulated as dynamic relationships lead to {dynamic games}; see \cite{basar1999dynamic} and the references therein. 
In this setting, an interesting scenario arises when the players have different incentives, e.g., to minimize or maximize a cost function. Dynamic noncooperative games focus on the case in which the players select their actions in a competitive manner, such that an individual benefit potentially implies a detriment to the other players  \cite{varaiya1969existence,fudenberg1991game,Owen,isaacs1999differential,Basar2017,Hespanha09}.
This type of dynamic game has been thoroughly studied in the literature when the dynamic constraints are given in terms of difference equations or differential equations, known as differential games.

Challenges arise when the players' dynamics exhibit both continuous and discrete behaviors, {for instance}, due to intermittent information availability, resets, timers in the control algorithms that expire, or non-smooth mechanical behaviors exhibiting impacts, among others. Hybrid {dynamical systems}  conveniently capture this kind of behavior {\cite{65nu,220}.
However,} when designing algorithms {that make optimal choices of the decision variables under constraints given by hybrid dynamics}, relying {only on continuous-time or discrete-time} approaches potentially results in suboptimal solutions. Unfortunately, tools for the design of algorithms for games with such hybrid dynamic constraints, which we refer to as hybrid games, are not fully developed.

Particular classes of dynamic games involving hybrid dynamic constraints have been recently studied in the literature.
A game-theory-based control design approach is presented for timed automata 
in \cite{asarin1995symbolic,asarin1998controller}, for hybrid finite-state automata in \cite{Henzinger1999RectangularHG,tomlin2000game}, and for o-minimal hybrid systems in \cite{bouyer2006control}. In these articles, the specifications to be guaranteed by the system are defined in terms of temporal logic formulae. When the payoff is defined in terms of a terminal cost, this approach allows designing reachability-based controllers through the satisfaction of Hamilton-Jacobi conditions to certify safety of hybrid finite-state automata \cite{ding2011toward}.
{Following an approach that allows for richer dynamics, \cite{vladimerou2009stormed} studies} a class of reachability games between a controller and the environment, for which the dynamics are given in terms of hybrid automata {(STORMED games)} and at each decision step, the players can choose either to have their variables evolve continuously or discretely, following predefined rules.
For {continuous-time systems with state resets,} tools for the computation of the {region of attraction} for hybrid limit cycles under the presence of disturbances are provided in \cite{choi2022computation}. {In \cite{choi2022computation},} only inputs that affect the flow are considered.
%
%

%
%
Efforts pertaining to differential games with impulsive elements include
\cite{platzer2017differential,platzer2015differential}, 
in which the {interaction} between the players is modeled similarly to switched systems,
\cite{dharmatti2006zero}{, which} establishes continuity of bounds on value functions and viscosity solutions, \cite{gromov2017class}{, which} formulates necessary and sufficient conditions for optimality in bimodal linear-quadratic differential games, and \cite{cacace2020stochastic}{, which} studies a class of stochastic two-player differential games in match race problems.
}
\IfPers{
Games involving multiple players with potentially different interests emerge in multi-agent systems, both in benign (or cooperative) and contested (or noncooperative) settings.  A list of examples includes and is not limited to route selection in a road network \cite{Paccagnan}, heavy-duty vehicle platooning \cite{gattami2011establishing}, control of smart grids \cite{tatarenko2014game}, trading modeling in the stock market \cite{firoozi2017mean}, and control of large populations of systems \cite{grammatico2015decentralized}.  Generally speaking, a game is an optimization problem with multiple players, constraints that enforce the ``rules'' of the game, and payoff functions to be optimized through the selection of decision variables. Constraints on the actions played by the players formulated as dynamic relationships (i.e., involving time) lead to 
{dynamic games}. Differential games pertain to the case when these constraints are given in terms of differential equations; see, e.g., \cite{basar1999dynamic} and the references therein.  Of particular interest is the contested setting, which occurs when the players have independent objectives, such as when one player aims at minimizing a cost function and another player aims at maximizing it under dynamic constraints. If the players select their actions seeking their own benefit, a dynamic noncooperative game emerges. These types of dynamic games have been thoroughly studied in the literature, when the dynamic constraints are given in terms of difference equations or differential equations --  in general, referred to as differential games -- including, to just list a few, \cite{varaiya1969existence,fudenberg1991game,Owen,isaacs1999differential,Basar2017,Hespanha09}. 

In recent years, significant progress has been made in the understanding of dynamic games with players sharing information over networks; see, e.g., \cite{marden2015game}. 
Interestingly, the combination of physics, computing, and networks 
leads to dynamic constraints that exhibit both continuous and discrete behavior.
In particular, intermittent information availability, resets of variables, such as expiring timers, and other nonsmooth and instantaneous changes lead to dynamic constraints that can be conveniently captured using hybrid \IfIncd{inclusions}{equations}.
The design of algorithms that guarantee optimality under such hybrid dynamic constraints
requires new tools, since using tools from the differential games literature would most likely lead to  suboptimal solutions.  Unfortunately, tools for the design of algorithms for games with such hybrid dynamic constraints, which we refer to as \IfIncd{set-valued}{} 
hybrid games, are not as developed as those for differential games, as cited above.  
In \cite{Henzinger1999RectangularHG,tomlin2000game}, a control design approach 
using game theory that is applicable to a class of hybrid automata models is presented.  Specifically, the models considered therein are based on finite-state automata, the specifications are defined in terms of  temporal logic formulae, and the payoff is solely given by a terminal cost. 
Decidability for hybrid automata given a winning condition are studied in \cite{vladimerou2009stormed}.
Applications of the approach in \cite{tomlin2000game} include reachability-based controller design \cite{ding2011toward,fisac2015reach}. 
The work in \cite{platzer2017differential,platzer2015differential} pertains to a class of dynamic games 
in which the evolution of the variables associated to each of the players is modeled using differential equations, while the interactions between the players is modeled as switches that occur at isolated time instances, similar to switched systems.
Conditional viability for  impulsive systems with two competing input actions was considered in \cite{aubin2005conditional} and treated as an evolutionary game.
Other efforts pertaining to differential games with impulsive (or discontinuous) elements include establishing continuity of bounds on value functions and (viscosity) solutions \cite{dharmatti2006zero}, formulating necessary and sufficient conditions for optimal strategies for the special case of bimodal linear-quadratic differential games \cite{gromov2017class}, and a class of stochastic two-player differential games in the context of sail boat competitions \cite{cacace2020stochastic}.}

%
\IfIncd{}{\subsection{Contribution and Applications}}

\section{Introduction}
\setlength{\parindent}{3ex}
Optimal control provides powerful tools for the study of multi-agent systems operating in contested scenarios in which {each of the agents (or players) are dynamic and select their} control actions to optimize a cost functional. {When the} constraints are given in terms of differential equations, {such problems are referred to} as differential games \cite{isaacs1999differential}.
The presence of dynamic constraints defined in both continuous and discrete (hybrid) time  
imposes challenges to computing optimal feedback laws and to applying cost evaluation tools.
Such a combination of continuous and discrete constraints can be efficiently captured by hybrid {dynamical} models, yielding {\em hybrid dynamic constraints}.
The classical approaches based on Hamilton-Jacobi-Bellman equations, e.g.,  \cite{bernstein1993nonquadratic,lafflitto2017differential,HaddadLafflittoFeedbacknu}, potentially lead to suboptimal outcomes when employed to compute and evaluate the optimal cost in scenarios with hybrid constraints, which we call {\em hybrid games}.

Following the framework in \cite{65nu}, a hybrid dynamical system, which is denoted $\HS\IfIncd{}{}$, is given in terms of the hybrid \IfIncd{inclusion}{equation}
\begin{equation}
\HS\IfIncd{}{}  \left\{ 
\begin{matrix}
\dot{x} &\IfIncd{\in}{=}& F(x,u_{C1},u_{C2}) \quad  (x,u_{C1},u_{C2}) \in C \\ 
 x^+ &\IfIncd{\in}{=}& G(x,u_{D1},u_{D2}) \quad (x,u_{D1},u_{D2}) \in D 
\end{matrix}
\right. 
\IfIncd{
\label{Hinc}
}
{\label{Heq}}
\end{equation}
where $x \in \mathbb{R}^n$ is the state, $(u_{C1},u_{D1}) \in \mathbb{R}^{m_{C1}} \times \mathbb{R}^{m_{D1}}$ is the input chosen by player $P_1$, $(u_{C2},u_{D2}) \in \mathbb{R}^{m_{C2}} \times \mathbb{R}^{m_{D2}}$ is the input chosen by player $P_2$, and  $(C,F,D,G)$ is the data of $\HS$. 
  Continuous evolution of $\HS$ is allowed when the state {and the input are} in the $\mathit{flow}$ $\mathit{set}$ $C$, and is governed by the $\mathit{flow}$ $\mathit{map}$ $F: \mathbb{R}^n \times \mathbb{R}^{m_C}  \IfIncd{\rightrightarrows}{\rightarrow} \mathbb{R}^n${, with {$m_C:=m_{C1}+m_{C2}$}}.
   Discrete evolution (or a jump) of $\HS$ is allowed when the state {and the input are} in the $\mathit{jump}$ $\mathit{set}$ $D$ and the new value of the state after a jump is captured by the $\mathit{jump}$ $\mathit{map}$ ${G: \mathbb{R}^n  \times \mathbb{R}^{m_D} \IfIncd{\rightrightarrows}{\rightarrow}  \mathbb{R}^n}${, with $m_D:=m_{D1}+m_{D2}$}. 

  For such type of systems  and for the one-player case, a cost functional has been proposed in \cite{ferrante2019certifying}, and optimality is certified via Lyapunov-like conditions. The results in \cite{ferrante2018optimality}  provide cost evaluation tools for {the case of discrete-time set-valued dynamics}.
The work in \cite{GOEBEL2019153nu} provides sufficient conditions to guarantee the existence of optimal solutions. 
A receding-horizon algorithm for optimal control of hybrid systems with single-valued {flow and jump maps} that extends the model predictive paradigm to such setting is presented in \cite{193c}.

{Informally,} a zero-sum two-player {hybrid game} is given as
\begin{equation}\label{min-maxPInt}
{ \min_{(u_{C1},u_{D1})} \max_{(u_{C2},u_{D2})}}\> \mathcal{J}(\xi, u_{C1},u_{C2},u_{D1},u_{D2}).
\end{equation}
where {$\mathcal{J}$} is a cost functional  
associated to the {solutions} to $\HS\IfIncd{}{}$ from the initial condition $\xi$.
The evolution of the game is determined via the selection of the inputs by the players, and its outcome is determined by computing $\J$.

 In \cite{232}, cost evaluation results are established to guarantee optimality and asymptotic stability of a set of interest in a discrete-time setting. A finite-horizon hybrid game, where the dynamical constraints are expressed in terms of hybrid systems, is studied in \cite{leudo2022FHnu}. 
 In \cite{sjlAutom26}, a terminal set is conveniently defined to recast both finite and infinite horizon {hybrid games,} and sufficient conditions are given in terms of Hamilton–Jacobi–Bellman–Isaacs equations to guarantee that a solution to the game is attained.

\subsection{Contribution}

The main contributions of this paper are summarized as follows.
%
\begin{itemize}
  \item In Section 3, we present a 
  framework for the study of two-player zero-sum games with hybrid dynamic {set-valued} constraints.
  \item In Section 4, we present cost evaluation results for hybrid set-valued autonomous systems. {Proposition \ref{prop:Prop1} provides a cost upper bound and Corollary \ref{coro:SingleValued} establishes} conditions to exactly evaluate the cost. 
  \item 
 {We
 present in Theorem \ref{thHJBszihinc}
 sufficient conditions 
  to characterize the solution {to} the min-max problem \eqref{min-maxPInt},  bounding and evaluating the worst-case cost over the set of adversarial strategies, and 
  to characterize the feedback law that attains it. This {result} relies on Propositions \ref{UpperBoundOpenLoop} and \ref{lemma:SolsMaxSys}, in which an upper bound is established on the cost associated to a given control action and the cost is exactly evaluated for a \emph{maximal system}, respectively. 
  }
  \IfIncd{
  }
  {\begin{itemize}
    \item We present in Theorem \ref{thHJBszih} sufficient conditions based on Hamilton-Jacobi-Isaacs-like equations to attain a saddle-point equilibrium and evaluate the game value function without {computing} solutions.}
  {
  \item Connections between optimality and asymptotic stability of a closed set are proposed in Section 6 
  and framed in the game-theoretic approach employed.
  \item We present in Section 7 
  applications to robust and security scenarios by {formulating and solving them as two-player zero-sum set-valued hybrid games}. }
  \IfIncd{}{
\end{itemize}
  }
\end{itemize}

\pnn{A preliminary version of these results appeared in the conference paper\pn{s \cite{ferrante2018optimality} and} \cite{svsjl}. 
\pnn{Compared to \pn{\cite{ferrante2018optimality} and} \cite{svsjl}, this paper provides complete proofs of the results therein and more details for each theorem, example, and application.
  In addition, we include new results and applications, namely, 
  %
  Proposition 4.2 (and Remark 4.1) is adapted from \cite{ferrante2018optimality} to the framework employed herein,
  providing weaker conditions to attain an upper bound in the cost associated to the solutions of a hybrid system with no inputs for a given trajectory (and for an initial condition). Likewise, the case of exact cost evaluation is addressed in Corollary 4.3.
  %
  \IfSiam{\cite[Proposition A.1]{TRHSVGArxiv}(and \cite[Corollary A.2]{TRHSVGArxiv})}{Proposition A.1 (and Corollary A.2)} presents time-dependent conditions to establish an upper (lower) bound in the cost associated to a set-valued hybrid system. These results are instrumental for the proof of Proposition 5.1 and establishing pointwise conditions to characterize the cost bounds and for the proof of Proposition 5.2 and establishing the relation between a maximal system and exact cost evaluation.
  In addition, such a maximal system is studied in \IfSiam{\cite[Lemma A.3]{TRHSVGArxiv}}{Lemma A.3}, where its solutions are shown to be solutions to the associated hybrid system as well. 
  In \IfSiam{\cite[Lemma A.4]{TRHSVGArxiv}}{Lemma A.4}, well-posedness is established, and in \IfSiam{\cite[Proposition A.8]{TRHSVGArxiv}}{Proposition A.8} existence of solutions is built upon the conditions for existence of solutions of a hybrid system (\IfSiam{\cite[Proposition A.7]{TRHSVGArxiv}}{Proposition A.7}).
  These results provide sufficient conditions either to bound or exactly evaluate the cost of the game, which is a key step to characterize a weak saddle-point equilibrium. 
  %
  Finally, in Application 7.1 we study robust hybrid linear quadratic problems with spontaneous jumps, providing Corollary 7.1  giving rise to the hybrid Riccati equation for such problems, and Example 7.2 as a numeric case of study.  
  }

 \pnn{Compared to 
 \cite{sjlAutom26}, which provides a general formulation of hybrid games, we} present a less conservative approach by relaxing the assumption on uniqueness of solutions to the hybrid system defining the constraints of the game. 
 Although, it might not be possible to construct a saddle-point equilibrium as the solution to the game when the dynamics admit nonunique solutions (due to the game being ill-defined by the nonuniqueness of costs associated to a given input), \pno{a weak saddle-point equilibrium} is {defined herein and conditions to characterize it are provided.} } 

  %

%
\IfIncd{}{

This type of hybrid game emerges in scenarios with continuous-time dynamics with logical modes, switching systems,  hybrid automata, impulsive differential equations, and dynamics described by algebraic differential inclusions, and  
 in several applications, as we illustrate next.
\begin{itemize}
  \item
\begin{application}{Robust Control} 
  Given the system $\HS\IfIncd{}{}$ as in (\IfIncd{\ref{Hinc}}{\ref{Heq}}) with state $x$, the disturbance rejection problem consists of establishing conditions such that player $P_1$ selects a control input $(u_{C1},u_{D1})$ that {minimizes} the cost of solutions to $\HS\IfIncd{}{}$ in the presence of a disturbance $(u_{C2},u_{D2})$ chosen by $P_2$. 
\label{RobustApp}
\end{application}

\item
\begin{application}{Security} Given the system $\HS\IfIncd{}{}$ as in (\IfIncd{\ref{Hinc}}{\ref{Heq}}) with state $x$\IfAutom{{, $F(x,u_{C1},u_{C2}) = f_d(x, u_{C1})+f_a(u_{C2})$, and 
$G(x,u_{D1},u_{D2}) = g_d(x, u_{D1})+g_a(u_{D2})$,}}{
  and  
\begin{equation*}
\begin{matrix}
 F(x,u_{C1},u_{C2}) &=& f_d(x, u_{C1})+f_a(u_{C2}) \\ 
 G(x,u_{D1},u_{D2}) &=& g_d(x, u_{D1})+g_a(u_{D2})
\end{matrix}
\label{SecDyn}
\end{equation*}}
the security problem consists of 
finding conditions such that the control input $(u_{C1},u_{D1})$ guarantees the integrity of the system via minimizing a cost functional $\mathcal{J}$ {(representing damage or effect of attacks)} under the action of an attacker $(u_{C2},u_{D2})$ that knows $f_d$ and $g_d$, and is designed to maximize {$\mathcal{J}$}. 
\label{SecurityApp}
\end{application}
\end{itemize}
 
The solution of the game formulated in this paper, known as a saddle-point equilibrium, is given in terms of the action selected by the players. Informally, when a player selects any strategy different from the equilibrium and the adversary plays the equilibrium action, there is no unilateral improvement in its individual outcome. Thus, by formulating the applications above as two-player zero-sum hybrid games, we can synthesize the saddle-point equilibrium and determine the control action that minimizes the cost $\J$ for the maximizing adversarial action. 
}

\pn{ The remainder of this paper is organized as follows. 
 Section 2 presents preliminary definitions that are used throughout the article. 
  In Section 3, we present a formulation of two-player zero-sum set-valued hybrid games and provide the definition of solution.
  Section 4 is dedicated to tools on the cost evaluation for autonomous hybrid inclusions.
 We present in Section 5 the sufficient conditions to solve the type of game studied.
 We study connections between optimality and stability for the proposed type of games in Section 6. 
 Applications and a numerical example covering a type of hybrid systems are presented in Section 7, displaying the versatility of the approach. Section 8 provides conclusions, closing remarks, and future work. }

\subsection{Notation}
Given two vectors, $x$ {and} $y$, we use the equivalent notation $(x,y)=[x^\top y^\top]^\top$ {and $\langle x,y \rangle $ denotes the Euclidean inner product.} The symbol $\mathbb{N}$ denotes the set of natural numbers including zero. The symbol $\mathbb{R}$ denotes the set of real numbers and $\mathbb{R}_{\geq 0}$ denotes the set of nonnegative reals. 
{Given two sets, $A$ and $B$, the Cartesian product $A \times B$ is defined as the set of all ordered pairs $(a, b)$, where $a$ is an element of $A$ and $b$ is an element of $B$.}
Given a nonempty set $\mathcal{A}$, the distance from the vector $x$ to $\mathcal{A}$ is defined as $\left| x \right|_\mathcal{A}:=\inf_{y \in \mathcal{A}}\left | x-y \right |$. We denote by $\mathbb{S}_+^n$ the set of real positive definite matrices of dimension $n$, \IfPers{ i.e., $A \in \mathbb{S}_+^n$ if  $A \in \mathbb{R}^{n \times n}$ $A=A^\top$ and there exists a real nonsingular matrix $M$ such that $A=M M^\top$}
and with $\mathbb{S}_{0+}^n$ the set of real positive semidefinite matrices of dimension $n$. Given a nonempty set $C$, we denote by  $\textup{int} \> C$ its interior and by $\overline{C}$ its closure. 
\IfPers{Given a symmetric matrix $A\in \reals^{n\times n}$, the scalars $\underline{\lambda}(A)$ and $\overline{\lambda}(A)$ denote the \pnn{smallest and largest} eigenvalue of $A$, respectively.}
{We denote by $\mathds{1}_n$ the $n$-dimensional vector of ones.}

\section{Preliminaries}
\IfTp{}{\subsection{Multi-agent network}
The network in which the interactions between players occur is modeled by the directed graph $\Gamma=(\IfTp{\{1,2\}}{\mathcal{V}},\mathcal{E},\mathcal{G})$, with  node set $\IfTp{\{1,2\}}{\mathcal{V}}=\{1,2,\dots,N\}$. The set of edges $\mathcal{E} \subset \IfTp{\{1,2\}}{\mathcal{V}}\times \IfTp{\{1,2\}}{\mathcal{V}}$ models the links between nodes, i.e., every pair $(i,k) \in \mathcal{E}$ implies there exists information sharing from the player $i$ to the player $k$. The $i,k$ entry of the adjacency matrix, $\mathcal{G}$, namely $\mathcal{G}_
{ik}$, is equal to 1 if $(i, k) \in \mathcal{E}$, and takes a zero value otherwise.  The set of indices corresponding to the neighbors that can send
information to the $i$-th player is denoted by $\mathcal{N}(i) := \{k \in \IfTp{\{1,2\}}{\mathcal{V}} : (k, i) \in \mathcal{E}\}$. The number of players is then $N$ and the number of communication links is $|\mathcal{E}|$. 
We consider the case in which there is a limited amount of communication channels. }

\subsection{Hybrid Systems with Inputs}\label{Subsect HySyswIn}

Since solutions to the dynamical system {$\HS\IfIncd{}{} $ as in \IfIncd{\eqref{Hinc}}{\eqref{Heq}}}  can exhibit both continuous and discrete behavior, we use ordinary time $t \in \mathbb{R}_{\geq 0}$ to determine the amount of flow and a counter $j \in \mathbb{N}$ that counts the number of jumps. 
Based on this parametrization of time, the concept of a hybrid time domain, in which solutions are fully described,
{ is proposed.
\begin{definition}[Hybrid time domain]
{ A set $\widetilde{E} \subset \mathbb{R}_{\geq0} \times \mathbb{N}$ is a compact hybrid time domain if 
\pnn{there exists $J \in \nats$ such that}
    $ \widetilde{E} = \bigcup_{j=0}^J ([t_j,t_{j+1}]\times\{j\})$
 for some finite nondecreasing sequence of times $\{t_j\}^{J+1}_{j=0}$ \pn{with $t_{0}=0$}.}
A set $E \subset \mathbb{R}_{\geq0} \times \mathbb{N}$ is a hybrid time domain if
{it is the union of a nondecreasing sequence of compact hybrid time domains, namely, $E$ is the union of compact hybrid time domains $E_j$ with the property that $E_0 \subset E_1 \subset E_2 \subset \ldots E_J \ldots$.}
Each element $(t,j)\in E$ denotes the elapsed hybrid time, which indicates that $t$ seconds of flow time and $j$ jumps have occurred. 
\end{definition}}
A hybrid signal is a function defined on a hybrid time domain. Given a hybrid signal $\phi$ and $j \in \mathbb{N}$, we define $I^j_{\phi}\pnn{:}=\{t:(t,j)\in \dom \phi\}$.

\begin{definition}[Hybrid arc]
A hybrid signal $\phi:\dom \phi \rightarrow \mathbb{R}^n$ is called a hybrid arc if, for each $j \in \mathbb{N}$, the function $t \mapsto \phi(t,j)$ is locally absolutely continuous on the interval $I^j_{\phi}$. A hybrid arc $\phi$ is compact if $\dom \phi$ is compact. 
\label{htd}
\end{definition}

\begin{definition}[Hybrid input] A hybrid signal $u$ is a hybrid input if, for each $j \in \mathbb{N}$, the function $t \mapsto u(t,j)$ is Lebesgue measurable and locally essentially bounded on the interval $I^j_u$.  
  \label{Def:HyInput}
\end{definition}
Let $\mathcal{X}$  be the set of hybrid arcs $\phi: \textup{dom}\> \phi \rightarrow \mathbb{R}^n$ and  $\mathcal{U}=\mathcal{U}_C \times \mathcal{U}_D$ the set of hybrid inputs $u=(u_{C},u_{D}): \textup{dom}\>u \rightarrow \mathbb{R}^{m_C} \times \mathbb{R}^{m_D}$, where $u_C=(u_{C1},u_{C2}\IfTp{}{,\dots, u_{C,N}})$, \IfTp{$m_{C1}+m_{C2}=m_C$,}{} $u_D=(u_{D1},u_{D2}\IfTp{}{,\dots, u_{D,N}})$\IfTp{, and $m_{D1}+m_{D2}=m_D$}{}. {A solution to the hybrid system $\HS\IfIncd{}{} $ with input is defined as follows.}

\begin{definition}[Solution to $\HS\IfIncd{}{}$]
A hybrid signal $(\phi,u)$ {defines a solution pair to (\IfIncd{\ref{Hinc}}{\ref{Heq}})} 
  if $\phi \in \mathcal{X}$, $u=(u_C,u_D) \in \mathcal{U}$,   $\textup{dom}\phi=\textup{dom}u$, and 
  \begin{itemize}
    \item $(\phi(0,0),u_C(0,0)) \in \overline{C}$ or $(\phi(0,0),u_D(0,0)) \in D$,
    \item For each $j\in \mathbb{N}$ such that $I^j_{\phi}$ has a nonempty interior $\textup{int} I^j_{\phi}$, we have, for all
    $t \in \textup{int} I^j_{\phi}$,
    $
        (\phi(t,j),u_C(t,j)) \in C
    $
    and, for almost all $t \in I^j_{\phi}$,
    $
      \frac{d}{dt}\phi(t,j) \in F(\phi(t,j), u_C(t,j)),
    $
    \item For all $(t,j)\in \dom \phi$ such that $(t,j+1)\in \dom \phi$,  
$
	(\phi(t,j),u_D(t,j)) \in 
  D, 
	\phi(t,j+1) 
  \in G(\phi(t,j),u_D(t,j))
.$
\end{itemize}
A solution pair $(\phi, u)$ is a compact solution pair if $\phi$ is a compact hybrid arc
{; see Definition \ref{htd}}.
\label{SolutiontocalH}
\end{definition}
{
Given a solution pair $(\phi,u)$, the component $\phi$ is referred to as the state trajectory.
In this article, the same symbols are used to denote input actions and their values. The context clarifies the meaning of $u$, as follows: ``the function $u$\pnn{,}'' ``the signal $u$,'' or ``the hybrid signal $u$'' that appears in ``the solution pair $(\phi,u)$'' refer to the input action, whereas  ``$u$'' refers to the input value as a point in $\mathbb{R}^{m_C} \times \mathbb{R}^{m_D}$ in any other case. The reader can replace ``the function $u$'' by ``$u_\phi$'', that is the input action yielding the system to a state trajectory $\phi$.}

{A solution pair $(\phi,u)$ to $\HS\IfIncd{}{}$ from $\xi \in \reals^n$ is nontrivial if $\dom (\phi,u)$ contains at least two points. It is complete if $\dom (\phi,u)$ is unbounded. It is maximal if there is no solution $(\psi,w)$ from $\xi$ such that $\phi(t,j)=\psi(t,j)$ and $u(t,j)=w(t,j)$ for all $(t,j)\in \dom (\phi,u)$ and $\dom (\phi,u)$ is a proper subset of $\dom (\psi,w)$.} 
 We denote by $\hat{\mathcal{S}}_{\HS\IfIncd{}{}}(M)$ the set of solution pairs $(\phi,u)$ to $\HS\IfIncd{}{}$ as in (\IfIncd{\ref{Hinc}}{\ref{Heq}}) such that $\phi(0,0) \in M$. The set $\mathcal{S}_{\HS\IfIncd{}{}}(M) \subset \hat{\mathcal{S}}_{\HS\IfIncd{}{}}(M)$ denotes all maximal solution pairs from $M$. 
  \IfPers{For a given $u \in \mathcal{U}$, we denote the set of maximal state trajectories to $\HS\IfIncd{}{} $ from $\xi$ for $u$ by $\mathcal{R}(\xi, u)=\{ \phi :(\phi,u) \in \mathcal{S}_\HS (\xi) \}$. We say $u$ renders a maximal trajectory $\phi$ to $\HS\IfIncd{}{} $ from $\xi$ if $\phi \in \mathcal{R}(\xi, u)$.}
For a given $u \in \mathcal{U}$, we denote the set of maximal state trajectories to $\HS\IfIncd{}{}$ from $\xi$ for $u$ by $\mathcal{R}(\xi, u)=\{ \phi :{\textup{$\phi$ is nontrivial and }}(\phi,u) \in \mathcal{S}_{\HS\IfIncd{}{}} (\xi) \}$. We say $u$ renders a maximal trajectory $\phi$ to $\HS\IfIncd{}{}$ from $\xi$ if $\phi \in \mathcal{R}(\xi, u)$.

We define the projections \pnn{of $\star \in \{C,D\}$ onto $\mathbb{R}^n$
as 
$	\Pi(\star)=\{ \xi \in  \mathbb{R}^n : \exists u_\star \in \mathbb{R}^{m_\star}  \text{ s.t. } 		(\xi,u_\star) \in \star \}. 
$} 
We also define the set-valued maps \pnn{
$	\Pi_u(x,\star)=\{ u_\star \in  \mathbb{R}^{m_\star} : 
		(x,u_\star) \in \star \}, \star \in \{C,D\}, 
$} 
denoting the sets of {allowed input values} for a given state. 
We define $\sup_t \dom \phi :=\sup\{ t \in \mathbb{R}_{\geq 0}: \exists j \in \mathbb{N}_{\geq 0} 	\>\text{s.t.} \> (t,j) \in \dom \phi\}$, $\sup_j \dom \phi :=\sup\{ j \in \mathbb{N}_{\geq 0}: \exists t \in \mathbb{R}_{\geq 0} 	\>\text{s.t.} \> (t,j) \in \dom \phi\}$, and {$\sup \dom \phi = (\sup_t \dom \phi,\sup_j \dom \phi)$.}


\IfIncd{}{
The following conditions guarantee uniqueness of solutions to $\HS\IfIncd{}{}$ as in (\IfIncd{\ref{Hinc}}{\ref{Heq}}) \cite[Proposition 2.11]{65nu}.
\begin{proposition}{Uniqueness of Solutions}
 Consider the hybrid system $\HS\IfIncd{}{} $ as in (\IfIncd{\ref{Hinc}}{\ref{Heq}}). For every $\xi \in \Pi(\overline{C}\cup D)$ and each $u \in \mathcal{U}$  there exists a unique maximal {solution $(\phi,u)$} with $\phi(0,0)=\xi$ provided that {the following holds:}
 \begin{enumerate}
\item   for every $\xi \in \Pi(\overline{C})\setminus \Pi(D)$ and  $T>0$, if two locally absolutely continuous {functions} $z_1,z_2:{I_z} \rightarrow \mathbb{R}^n$ and a {Lebesgue} measurable {function} $u_{z}:{I_z} \rightarrow \mathbb{R}^{m_C}$ {with $I_z$ of the form $I_z=[0,T)$ or $I_z = [0,T]$,} are such that, {for each $i \in \{1,2 \}$,} $\dot{z}_i(t) \IfIncd{\in}{=} F(z_i(t),u_{z}(t))$ for almost all $t \in {I_z}$, $(z_i(t), u_{z}(t)) \in C$ for all $t \in 
\textup{int}I_z$, and $z_i(0)=\xi$,  then $z_1(t)=z_2(t)$ for every $t \in {I_z}$;
%
  \item for every $(\xi,u_D) \in D$, $G(\xi, u_D)$ consists of one point.
  \end{enumerate}
  \label{UniquenessHu}
  \end{proposition}
    \NotAutom{
  \begin{proof}
  Proceeding by contradiction, suppose there exist two maximal solutions to $\HS$, $(\phi_1,{u})$ and $(\phi_2,{u})$, with $u=(u_C,u_D)$, {$\dom\phi_1 = \dom\phi_2 = \dom u$}, and $\phi_1(0,0) = \phi_2(0,0)$ such that $\phi_1$ and $\phi_2$ are not identical over ${\dom u}
 $, {namely, there exists $(t^*,j^*) \in {\dom u}
 $ such that $\phi_1(t^*,j^*) \neq \phi_2(t^*,j^*)$.}  We have the following three cases:
  \begin{enumerate}[label=\alph*)]
	  \item If $(t^*,j^*) \in [t',t'') \times \{{j^*}\} \subset {\dom u}
	  $ for some  $t''>t'\geq 0$, then, the functions $z_1, z_2: I_z
	  \rightarrow \mathbb{R}^n$, and 
	  $u_{z}:I_z \rightarrow \mathbb{R}^{m_C}$, with {$T=t''-t'$} and $I_z$ of the form $I_z = [0,T)$, defined
	  {for each $i \in \{1,2 \}$ as} 
	${z}_i(t) = \phi_i ({t'+t},j)$ for all $t \in I_z$, and as $u_z(t) = u_C({t'+t},j)$ for all $t \in I_z$, 
	satisfy, {by Definition \ref{SolutiontocalH}}, $\dot{z}_i(t) \in F({z}_i(t),u_{z}(t))$ for almost all $t \in I_z$, and $(z_i(t), u_{z}(t)) \in C$ for all $t \in \textup{int}I_z$.
	However, at $t^* \in I_z$, $z_1(t^*)\neq z_2(t^*)$, 
    which contradicts item 1.
	  %
	  %
	    \item If $(t^*,j^*) \in [t',t''] \times \{{j^*}\} \subset {\dom u}
	  $ for some  $t''>t'\geq 0$, then, the functions $z_1, z_2: I_z
	  \rightarrow \mathbb{R}^n$, and 
	  $u_{z}:I_z \rightarrow \mathbb{R}^{m_C}$, with {$T=t''-t'$} and $I_z$ of the form $I_z = [0,T]$, defined
	  {for each $i \in \{1,2 \}$ as} 
	${z}_i(t) = \phi_i ({t'+t},j)$ for all $t \in I_z$, and as $u_z(t) = u_C({t'+t},j)$ for all $t \in I_z$, 
	satisfy, {by Definition \ref{SolutiontocalH}}, $\dot{z}_i(t) \in F({z}_i(t),u_{z}(t))$ for almost all $t \in [0,T)$, and $(z_i(t), u_{z}(t)) \in C$ for all $t \in \textup{int}I_z$.
	However, at $t^* \in I_z$, $z_1(t^*)\neq z_2(t^*)$, 
    which contradicts item 1.
	  %
	  %
	  %
	  %
	  \item If $(t^*,j^*)$ is such that 
	  $(t^*,j^*-1) \in \dom u
	  $,
	  then $(t^*,j^*-1)$ is a jump time of $\phi_1$ and $\phi_2$.
	  Then, by Definition \ref{SolutiontocalH},
for each $i \in \{1,2\},$ 
$\phi_i(t^*, j^*) = G(\phi_i(t^*,j^*-1),u_D(t^*, j^*-1))$.
Since, 
  $\phi_1(t^*,j^*) \neq \phi_2(t^*,j^*)$,
   and $\phi_1(t^*,  j^*-1)= \phi_2(t^*, j^*-1),
  G(\phi_1(t^*,j^*-1),u_D(t^*, j^*-1))$ takes more than one value, which contradicts item 2. 
	  %
  \end{enumerate}
Thus, provided items $1$ and $2$ hold, for every $\xi \in \Pi(\overline{C}\cup D)$ and each $u \in \mathcal{U}$  there exists a unique maximal {solution $(\phi,u)$} with $\phi(0,0)=\xi$. 
  \end{proof}
}
  }
\subsection{{Hybrid Closed-Loop Systems}}
Given 
a hybrid system $\HS\IfIncd{}{}$ and a function $\kappa:=(\kappa_C, \kappa_D)$ with $\kappa:\mathbb{R}^n \rightarrow \mathbb{R}^{m_C} \times \mathbb{R}^{m_D}$, 
 the autonomous hybrid system resulting from assigning $u=\kappa(x)$, namely, the hybrid closed-loop system, is given by
\begin{equation}
\HS_\kappa  \left\{ 
\begin{matrix}
\dot{x} &\IfIncd{\in}{=}& F(x,\kappa_C(x)) \quad  x \in C_\kappa \\ 
 x^+ &\IfIncd{\in}{=}& G(x,\kappa_D(x)) \quad x \in D_\kappa 
\end{matrix}
\right. 
\label{Hkeq}
\end{equation}
where $C_\kappa:=\{x \in \mathbb{R}^n: (x,\kappa_C(x)) \in C\}$ and $D_\kappa:=\{x \in \mathbb{R}^n: (x,\kappa_D(x)) \in D\}$.

A solution to the {hybrid closed-loop system} $\HS_\kappa$ is defined \IfAutom{similarly to Definition \ref{SolutiontocalH}.}{as follows.

\begin{definition}{Solution to $\HS_\kappa$} 
A hybrid arc $\phi$ defines a solution to the hybrid system $\HS_\kappa$ in (\ref{Hkeq}) 
  if 
\begin{itemize}
    \item $\phi(0,0) \in \overline{C_\kappa} \cup D_\kappa$,
    \item For each $j\in \mathbb{N}$ such that $I^J_{\phi}$ has a nonempty interior $ \textup{int} I^J_{\phi}$, we have, for all
    $t \in \textup{int} I^J_{\phi}$,
    \begin{equation*}
        \phi(t,j) \in C_\kappa
    \end{equation*}  
    and, for almost all $t \in I^J_{\phi}$,
    \begin{equation*}
       \frac{d}{dt}\phi(t,j)\in F(\phi(t,j), \kappa_C(\phi(t,j)))
    \end{equation*}
    \item For all $(t,j) \in \dom \phi$ such that $(t,j+1)\in \dom \phi$,  
    \begin{eqnarray*}
	\phi(t,j) &\in& D_\kappa, 
  \hspace{2cm}
	\phi(t,j+1) 
  \in G(\phi(t,j),\kappa_D({\phi}(t,j))) 
\end{eqnarray*}
\end{itemize}
A solution $\phi$ is a compact solution if $\phi$ is a compact hybrid arc.
\end{definition}
}
We denote by $\hat{\mathcal{S}}_{\HS_\kappa}(M)$ the set of solutions $\phi$ to $\HS_\kappa$ as in (\ref{Hkeq}) such that $\phi(0,0) \in M$. 
The set $\mathcal{S}_{\HS_\kappa}(M) \subset \hat{\mathcal{S}}_{\HS_\kappa}(M)$ denotes all maximal solutions from $M$.
\IfIh{}{Given $\mathcal{T}\subset \realsgeq \times \nats$, let us denote by $\hat{\mathcal{S}}^\mathcal{T}_{\HS  }(M)\subset \hat{\mathcal{S}}_{\HS_\kappa}(M)$ the set of compact solutions with 
terminal time in $\mathcal{T}$, i.e., if $(\phi,u) \in \hat{\mathcal{S}}^\mathcal{T}_{\HS  }(M)$ and $\max \dom (\phi,u) =(T_\phi,J_{\phi^*})$, then $(T_\phi,J_{\phi^*})\in \mathcal{T}$.}


\section{ Two-player Zero-sum Hybrid Games}\label{sec: 2pzshg}
\subsection{Game Formulation}\label{Sec: Game form}
Following the formulation in \cite{basar1999dynamic,sjlAutom26}, for each $i\in\{1,2\}$, consider the $\mathit{i}$-th player $P_i$ with dynamics described by $\HS_i$ as in (\IfIncd{\ref{Hinc}}{\ref{Heq}}) with data $(C_i, F_i, D_i, G_i)$, state $x_i \in \mathbb{R}^{n_i}$, and input $u_{ i}=(u_{Ci}, u_{Di})\in \reals^{m_{Ci}} \times \reals^{m_{Di}}$, where {$C_i\subset \mathbb{R}^{n} \times \mathbb{R}^{m_{C}}$, $F_i: \mathbb{R}^{n} \times \mathbb{R}^{m_{C}} \IfIncd{\rightrightarrows}{\rightarrow}  \mathbb{R}^{n_i}$, $D_i \subset \mathbb{R}^{n}  \times \mathbb{R}^{m_{D}}$ and ${G_i: \mathbb{R}^{n}  \times \mathbb{R}^{m_{D}} \IfIncd{\rightrightarrows}{\rightarrow} \mathbb{R}^{n_i} }$}, with \IfTp{$\star_1+\star_2=\star$ for $\star \in \{n, m_C, m_D\}$}{}. 
We denote by $\mathcal{U}_i=\mathcal{U}_{Ci} \times \mathcal{U}_{Di}$ the set of hybrid inputs for $\HS_i$; see Definition 2.3.

Notice that each player's dynamics are described in terms of data with maps and sets defined in the entire state and input space rather than the individual spaces ($\reals^n$ and $\reals^m$ rather than $\reals^{n_i}$ and $\reals^{m_i}$, respectively). This allows to model the ability of each player's state to evolve according to the state variables and input of the other players, or scenarios in which a player does not have a state but only selects an input.

\begin{definition}[Elements of a two-player zero-sum hybrid game] \IfTp{A two}{An $N$}-player zero-sum hybrid game is composed of 
\begin{enumerate}[label=\arabic*)]
\IfTp{}{\item A communication network described by a graph $\Gamma=(\mathcal{V,E,G})$ such that $\sum_{i= 1}^N n_{i} = n$, $\sum_{i= 1}^N m_{Ci} = m_C$ and $\sum_{i= 1}^N m_{Di} = m_D$, with $n_{i},m_{Ci} , m_{Di} \in \mathbb{N}$ for each $i \in \IfTp{\{1,2\}}{\mathcal{V}}$, where N is the number of players in the game.}
\item The state $x=(x_1, x_2\IfTp{}{, \dots, x_N}) \in \reals^n$, where, for each $i \in \IfTp{\{1,2\}}{\mathcal{V}}$, $x_i \in \reals^{n_i}$ is the state of player $P_i$.
\item The set of joint input actions $\mathcal{U}=\mathcal{U}_1 \times \mathcal{U}_2 \IfTp{}{\times \dots \times \mathcal{U}_N}$ with elements $u=(u_{ 1}, u_{ 2} \IfTp{}{,\dots, u_{ N}})$, where, for each $i \in \IfTp{\{1,2\}}{\mathcal{V}}$, $u_i=(u_{Ci},u_{Di})$ is a hybrid input. For each $i \in \IfTp{\{1,2\}}{\mathcal{V}}$, $P_i$ selects $u_i$ independently of $P_{3-i}$, who selects \IfTp{$u_{3-i}$}{any $u_j$, $j\neq i$}, namely, the joint input action $u$ has components $u_i$ that are independently chosen by each player.
\item The dynamics of the game, described as in $(\IfIncd{\ref{Hinc}}{\ref{Heq}}
)$ and denoted by $\HS\IfIncd{}{}$, with data 
 \IfPers{Comment: Before we had
 \begin{equation*} \hspace{0.3cm}
 \begin{array}{r@{\>\>}c@{\>\>}l}
 C&:=&C_1 \times C_2 \IfTp{}{\times \dots \times   C_N}\\
 D&:=&\{(x,u_D) \in \mathbb{R}^n \times \mathbb{R}^{m_D}\IfConf{\hspace{-0.1cm}}{}:
 (x_i,u_{Di}) \in D_i,
  i \in \IfTp{\{1,2\}}{\mathcal{V}}  \}\\
 G(x,u_D)&:=&\{ \hat{G}_i(x,u_{D}):
 (x_i,u_{Di}) \in D_i, i \in \IfTp{\{1,2\}}{\mathcal{V}} \}
  \hspace{2.2cm}\IfConf{}{ \forall (x,u_D) \in D }
  \end{array}
 \end{equation*}
 We modify it, so the data of each player depends on the full state which implies that whether player $P_i$ flows or jumps (potentially) depends on the full state. Every example we present so far (in which player $P_2$ consists only of an input and no states) follows this structure. Now we have:
 }
\begin{equation*} \hspace{-0.2cm}
    \begin{array}{r@{\>\>}c@{\>\>}l}
    C&:=&{C_1 \cap C_2} \IfTp{}{\times \dots \times   C_N}\\
    F(x,u_C)&:=&(F_1(x, u_{C}),F_2(x, u_{C})\IfTp{}{, \dots, F_N(x, u_{C})}) \hspace{1cm} 
    \\
    D&:=&{ D_1 \cup D_2}
    \\
    G(x,u_D)&:=&\{ \hat{G}_i(x,u_{D}):
    {(x,u_{D})} \in D_i, i \in \IfTp{\{1,2\}}{\mathcal{V}} \}
     \end{array}
    \end{equation*}
%
where \IfTp{$\hat{G}_1(x,u_{D})=(G_1(x, u_{D}),\mathds{1}_{n_2})$, $\hat{G}_2(x,u_{D})=(\mathds{1}_{n_1},G_2(x, u_{D}))$}{\pnn{
$\hat{G}_i(x,u_{D})=(\mathds{1}_{n_1},\mathds{1}_{n_2}, \dots,G_i(x, u_{D}) ,\dots,\mathds{1}_{n_N})$}}, $u_C=(u_{C1},u_{C2})$, and $u_D=(u_{D1},u_{D2})$.
\item For each $i  \in \IfTp{\{1,2\}}{\mathcal{V}}$, a strategy space $\mathcal{K}_i$ of $P_i$ defined as a collection of mappings \IfTp{$\kappa_i:\mathbb{R}^{n\IfTp{}{_{in}}} \rightarrow \mathbb{R}^{m_{Ci}} \times \mathbb{R}^{m_{Di}}$}{$\gamma_i:\realsgeq \times \nats \times \mathbb{R}^{n_{in}} \rightarrow \mathbb{R}^{m_{Ci}} \times \mathbb{R}^{m_{Di}}$}\IfTp{}{, where $n_{in}=\sum_{j=1}^N n_j \mathcal{G}(j,i)$ denotes the size of the state information shared with $P_i$}. The strategy space of the game, namely {$\mathcal{K}=\mathcal{K}_1 \times \mathcal{K}_2$}, is the collection of mappings with elements \IfIh{$\kappa=(\kappa_1,\kappa_2\IfTp{}{,\dots,\kappa_N})$}{$\gamma=(\gamma_1,\gamma_2,\dots,\gamma_N)$}, where \IfIh{$\kappa_i\in \mathcal{K}_i$}{$\gamma_i\in \mathcal{K}_i$} for each $i \in \IfTp{\{1,2\}}{\mathcal{V}}$. 
Each \IfIh{$\kappa_i \in \mathcal{K}_i$}{$\gamma_i \in \mathcal{K}_i$} is said to be a permissible pure\footnote{This is in contrast to when $\mathcal{K}_i$ is defined as a probability distribution, namely, when $\kappa_i\in \mathcal{K}_i$ is the space of mixed strategies.} strategy for $P_i$.
 \item  
%
A scalar-valued functional\IfIncd{\footnote{Given that we do not insist on {uniqueness of} solutions, the cost $\mathcal{J}$ measures the highest cost of the solutions yielded to $\HS\IfIncd{}{} $ from $\xi$ by $u$. Thus, its arguments are {hybrid inputs as in Definition \ref{Def:HyInput}} and not solution pairs.}}{} \IfIh{$(\xi, u) \mapsto \mathcal{J}_i(\xi, u)$}{$\mathcal{J}_i:\mathbb{R}^n \times \mathcal{U} \rightarrow \mathbb{R}$} defined for each $i \in \{1,2\}$, and called the cost associated to $P_i$. For each $u\in \mathcal{U}$, we refer to $\mathcal{J}:=\mathcal{J}_1=-\J_2$ as the 
worst-case cost due to nonuniqueness of solutions to $\HS\IfIncd{}{}$ for the hybrid input $u$  from the initial condition $\xi$.
\end{enumerate}
\label{elements}
\end{definition}
%
\IfPers{Sometimes we present the joint input action from the perspective of the $i$-th player as $u=(u_1,\IfTp{u_2}{\dots, u_i,\dots , u_N})=(u_i,u_{-i})$, where $u_{-i}$ is the set of input actions of \IfTp{player $P_{3-i}$}{players other than $P_i$}. In the games studied in this work, player $P_i$ wants to minimize the cost $\mathcal{J}_i$. For the two-player game with $\mathcal{J}_1=-\mathcal{J}_2$, a gain by player $P_1$ implies a loss to player $P_2$, leading to the term ``zero-sum games." For this type of games, we refer to a single cost functional $\mathcal{J}=\mathcal{J}_1=-\mathcal{J}_2$. }
In this type of game, for each $i \in \{1,2\}$, the player $P_i$ aims to minimize the cost $\J_i$, which, thanks to the definition of $\mathcal{J}$, allows us to define a min-max problem in terms of $\J$ alone.

\IfIncd{}{We say that a game formulation is in normal (or matrix) form when it describes only the correspondences between
different strategies and costs. 
On the other hand, 
{we refer to the mathematical description of a game to be in the Kuhn's extensive form if the evolution of the game defined by its dynamics, the decision-making process defined by the strategies, the sharing of information between the players defined by the communication network, 
and their outcome defined by the cost associated to each player, are described in the formulation. }
If a game is formulated in a Kuhn's extensive form, then it admits a solution.
{
Notice that for a given initial condition {$\xi$}, a given strategy $\kappa \in \mathcal{K}$ potentially leads to nonunique solutions $(\phi^1,u^1), (\phi^2,u^2),\dots, (\phi^k,u^k)$ to $\HS\IfIncd{}{}$, where $u^l=\kappa(\phi^l)$ {and $\phi(0,0)=\xi$} for each $l \in \{1,2,\dots k\}$.} 
Thus, for the formulation in Definition \ref{elements} to be in Kuhn's extensive form, 
{an appropriate cost definition is
required such that each strategy $\kappa \in \mathcal{K}$ has a unique cost correspondence, namely, every solution $(\phi^l,u^l)$ to $\HS\IfIncd{}{} $ with $u^l=\kappa(\phi^l)$ and $l \in \{1,2,\dots, k\}$ is assigned the same cost.
}

\begin{remark}{Relation to the literature}
\hspace{0.1cm} Given a discrete-time \IfTp{two-player}{} zero-sum game with final time{\footnote{This corresponds to the hybrid time $(0,J)$ for $\HS\IfIncd{}{}$.}} $``J$'', $f_k$ and $X$ defining the {single-valued} jump map and jump set, respectively, as in \cite{basar1999dynamic}, setting the data of $\HS\IfIncd{}{}$ as $C = \emptyset$, $G=f_k$ for $k \in \nats_{\leq J}$, and $D=X$ reduces Definition \ref{elements} to  \cite[Def. 5.1]{basar1999dynamic} for the case in which
the output of each player is equal to its state and there is a feedback information structure as in \cite[Def. 5.2]{basar1999dynamic}. 
Thus, items $(vi)-(ix)$ in \cite[Def. 5.1]{basar1999dynamic} are omitted in the formulation  herein and items $(i)-(v)$ and $(x)-(xi)$ are covered by Definition \ref{elements}, the definition of the hybrid time domain with final time $(0,J)$, and the set $\mathcal{S}_{\HS\IfIncd{}{}}$.
\\
\indent Given a continuous-time \IfTp{two-player}{} zero-sum game with final time{\footnote{This corresponds to the hybrid time $(T,0)$ for $\HS\IfIncd{}{}$.}} $``T$'', $f$ and $\mathcal{S}^0$ defining the {single-valued} flow map and flow set, respectively, as in \cite{basar1999dynamic}, setting the data of $\HS\IfIncd{}{}$ as $D = \emptyset$, $F=f$, and $C=\mathcal{S}^0$ reduces Definition \ref{elements} to \cite[Def. 5.5]{basar1999dynamic} for the case in which
the output of each player is equal to its state and there is a feedback information structure as in \cite[Def. 5.6]{basar1999dynamic}. 
Thus, items $(vi)-(vii)$ in \cite[Def. 5.5]{basar1999dynamic} are omitted in the formulation herein and items $(i)-(v)$ and $(viii),(ix)$ are covered by Definition \ref{elements}, the definition of the hybrid time domain with final time $(0,T)$, and the set $\mathcal{S}_{\HS\IfIncd{}{}}$.
\end{remark}

\IfPers{Given the elements of a \pn{zero-sum} game in the context of Definition \ref{elements}, the general solution of the game is given by the following definition.
\pn{\begin{definition}{Pareto-dominance solution}
An outcome of a game is Pareto dominated if some other outcome would make at least one player better off without hurting any other player. That is, some other outcome is weakly preferred by all players and strictly preferred by at least one player. If an outcome is not Pareto dominated by any other, then it is Pareto optimal, named after Vilfredo Pareto.
 \end{definition}
\begin{definition}{Pareto-optimal solution}For the case in which the player $i$ aims to minimize its own cost $\mathcal{J}_i: \mathbb{R}^n \times \mathcal{U} \rightarrow \mathbb{R}_{\geq 0}$, $i \in \IfTp{\{1,2\}}{\mathcal{V}}$, (TO DEFINE IF LCi, LDi, Vi)  following structure (\ref{defJ}), we consider the scenario in which they cooperate, acting in full trust of each other, to attain a mutually beneficial solution. Here, the cost $\mathcal{J}_i$ depends on the full state $x$ and the input $u$ of the network $\HS  $.
For this scenario,  we say an action $u^* \in \mathcal{U}$ is a \textit{Pareto-optimal solution} if for each $\xi \in \Pi(\overline{C})\cup \Pi(D)$, there does not exist another action $u^p \in \mathcal{U}$ such that 
\begin{equation*}
\mathcal{J}_i(\xi,u^p) \leq \mathcal{J}_i(\xi,u^*)
\end{equation*}
for all $i \in \IfTp{\{1,2\}}{\mathcal{V}}$, with 
\begin{equation*}
\mathcal{J}_i(\xi,u^p) < \mathcal{J}_i(\xi,u^* )
\end{equation*}
for some $i \in \IfTp{\{1,2\}}{\mathcal{V}}$. 
In other words,  $u^* \in \mathcal{U}$ is a Pareto-optimal solution if for any action $u^p \in \mathcal{U}$, such that 
\begin{equation*}
\mathcal{J}_i(\xi,u^p) < \mathcal{J}_i(\xi,u^* )
\end{equation*}
for some $i \in \IfTp{\{1,2\}}{\mathcal{V}}$, there exists an $i \in \IfTp{\{1,2\}}{\mathcal{V}}$ such that
\begin{equation*}
\mathcal{J}_i(\xi,u^p) > \mathcal{J}_i(\xi,u^* )
\end{equation*}
Hence, if $u^*$ is a Pareto-optimal solution rendering costs $\mathcal{J}^*_i$, $i \in \IfTp{\{1,2\}}{\mathcal{V}}$, there does not exist another solution to the game for which the cost of some players is reduced without increasing the cost of others (\cite{Hespanha09}).
\end{definition}}

 \pn{\begin{definition}{Pure strategy Nash equilibrium} For a \IfTp{two}{N}-player noncooperative game, with \IfTp{}{communication network $\Gamma=(\mathcal{V,E,G})$ and }dynamics $\HS$ as in (\ref{Heq}),  
we say a strategy \IfIh{$\kappa^*=(\kappa_1^*, \IfTp{\kappa_2^*}{\dots, \kappa_N^*})\in \mathcal{K}$}{$\gamma^*=(\gamma_1^*, \dots, \gamma_N^*)\in \mathcal{K}$} 
 is \textit{a pure\footnote{This is in contrast to when it is defined as a probability distribution of selecting certain mappings, referred to as a mixed strategy.} strategy Nash equilibrium} if for any $\xi \in \Pi(\overline{C}\cup D)$,
 every $u^*=(u^*_1,\IfTp{u^*_2}{ \dots, u^*_N})$ rendering a maximal trajectory $\phi^*$ to $\HS  $ from $\xi$, 
with components defined as \IfIh{$\dom \phi^* \ni (t,j) \mapsto u^*_i(t,j)= \kappa_i^*(\IfIh{}{t,j,}\phi^*(t,j))$}{$\dom \phi^* \ni (t,j) \mapsto u^*_i(t,j)= \gamma_i^*(t,j,\phi^*(t,j))$}, for each $i\in \IfTp{\{1,2\}}{\mathcal{V}}$, 
satisfies
\IfConf{
\begin{eqnarray}
\mathcal{J}_1(\xi,u^*) 
&\leq&  \mathcal{J}_1(\xi,(u_{1}, u_{-1}^*)) \nonumber\\&& \hspace{1cm}  \forall u_1:\exists \phi : (\phi, (u_{1}, u_{-1}^*)) \in \mathcal{S}_{\HS  } (\xi),
\nonumber \\
\mathcal{J}_2(\xi,u^*) 
&\leq&  \mathcal{J}_2(\xi,(u_{2}, u_{-2}^*)) \nonumber\\ && \hspace{1cm}  \forall  u_2:\exists \phi :(\phi, (u_{2}, u_{-2}^*)) \in \mathcal{S}_{\HS  } (\xi), \nonumber \\
\IfTp{}{&\vdots&
\nonumber\\
\mathcal{J}(\xi,u^*) 
&\leq&  \mathcal{J}(\xi,(u_{N}, u_{-N}^*))\nonumber \\ && \hspace{1cm}  \forall  u_N:\exists \phi :(\phi, (u_{N}, u_{-N}^*)) \in \mathcal{S}_{\HS  } (\xi) \nonumber\\}
 \label{NashEqNCoopIneq}
\end{eqnarray}
\vspace{-0.9cm}
}
{
\begin{equation}
\begin{array}{r@{}l}
\mathcal{J}(\xi,u^*) 
&{}\leq  \mathcal{J}(\xi,(u_{1}, u_{2}^*\IfTp{}{, \dots, u_{N}^*})) \hspace{1cm}  \forall u_1:\exists \phi : (\phi, (u_{1}, u_{2}^*\IfTp{}{, \dots, u_{N}^*})) \in \mathcal{S}_{\HS  } (\xi),
 \\
\mathcal{J}(\xi,u^*) 
&{}\leq 
 \mathcal{J}(\xi,(u_{1}^*, u_{2}\IfTp{}{, \dots, u_{N}^*})) \hspace{1cm} \forall  u_2:\exists \phi :(\phi, (u_{1}^*, u_{2}\IfTp{}{, \dots, u_{N}^*})) \in \mathcal{S}_{\HS  } (\xi)\IfTp{}{, 
 \\
&{ }\hspace{0.2cm}\vdots
\\
\mathcal{J}(\xi,u^*) 
&{}\leq  \mathcal{J}(\xi,(u_{ 1}^*, u_{2}^*, \dots, u_{N}))\hspace{1cm} \forall  u_N:\exists \phi :(\phi, (u_{ 1}^*, u_{2}^*, \dots, u_{N})) \in \mathcal{S}_{\HS} (\xi)}.
\end{array}
 \label{NashEqNCoopIneq}
\end{equation}
 } 
 \label{NashEqCoop}
 \end{definition}}}
}
 %
 %
%
\subsection{Equilibrium Solution Concept}
Given the formulation of a zero-sum hybrid game in Definition \ref{elements}, 
its solution is defined as follows.

\begin{definition}[Feedback weak saddle-point equilibrium]\label{NashEqNonCoop} Consider a \\two-player zero-sum game, with 
dynamics $\HS\IfIncd{}{}$ as in (\IfIncd{\ref{Hinc}}{\ref{Heq}}) with $\mathcal{J}_1=\mathcal{J}$, $\mathcal{J}_2=-\mathcal{J}$, for 
a given cost functional $\mathcal{J}:\mathbb{R}^n\times \mathcal{U} \rightarrow \mathbb{R}$. 
We say that a strategy $\IfIh{\kappa}{\gamma}=(\IfIh{\kappa}{\gamma}_1, \IfIh{\kappa}{\gamma}_2) \in \mathcal{K}$ 
is \textit{a weak saddle-point equilibrium} if, for each $\xi \in \Pi(\overline{C}\cup D)$,
every\footnote{Notice that a given strategy $\kappa$ can lead to multiple input actions due to $C \cap D$ being nonempty.}
 solution  $(\phi^*,u^*)=(\phi^*,(u_1^*, u_2^*))\in \mathcal{S}_\HS(\xi)$ attaining the worst-case cost and input components defined, for each $i\in \IfTp{\{1,2\}}{\mathcal{V}}$, as $\dom \phi^* \ni (t,j) \mapsto u_i^*(t,j):= \IfIh{\kappa}{\gamma}_i(\IfIh{}{t,j,}\phi^*(t,j))$,
satisfies
\begin{eqnarray}
 \mathcal{J}(\xi,(u_{ 1}^*, u_{ 2})) \leq
\mathcal{J}(\xi, u^*) \leq
 \mathcal{J}(\xi,(u_{ 1}, u_{ 2}^*))
 \label{SaddlePointIneq}
\end{eqnarray}
%
for all hybrid inputs $u_1$ and $u_2$ such that the sets of maximal state trajectories \\$\mathcal{R}(\xi, (u_{ 1}, u_{ 2}^*))$ and $\mathcal{R}(\xi, (u_{ 1}^*, u_{ 2}))$ are nonempty.
 \end{definition}
Definition \ref{NashEqNonCoop} is a generalization of the classical pure strategy Nash equilibrium \cite[(6.3)]{basar1999dynamic} to the case where the players exhibit dynamics expressed in terms of hybrid inclusions and opposite optimization goals.
 In words, we refer to the strategy {\IfIh{$\kappa=(\kappa_1,\kappa_2)$}{$\gamma^*$}} 
 as a saddle-point when 
a player $P_i$ cannot reduce its associated cost $\J_i$ by playing any strategy different from \IfIh{$\kappa_i$}{$\gamma_i^*$} when player $P_{3-i}$ plays \IfIh{$\kappa_{3-i}$}{$\gamma_{3-i}^*$}.     
%
{Condition \eqref{SaddlePointIneq} is to be verified over the set of inputs that define joint input actions $(u_1^*,u_2)$ and $(u_1,u_2^*)$, yielding at least one nontrivial solution to $\HS\IfIncd{}{}$ from $\xi$.}
\IfPers{Notice that the saddle-point, as a solution to the zero-sum \IfTp{two}{$N$}-player game, is a strategy in $\mathcal{K}$,  though the concept of a solution to a hybrid system $\HS\IfIncd{}{}$ 
is a hybrid arc.}

\NotAutom{
 \begin{remark}{Joint Action}
The pure strategy Nash equilibrium in Definition \ref{NashEqNonCoop} implies that for each $i \in \IfTp{\{1,2\}}{\mathcal{V}}$, the $i-$th component of each $u^*=(u^*_1, \IfTp{u^*_2}{\dots, u^*_N})$ 
satisfies
	\begin{equation}
u_i^* \in
	\argmin_{u_i:(\phi, (u_{i}, u_{3-i}^*)) \in \mathcal{S}_{\HS }  (\xi)}  \mathcal{J}_i(\xi,(u_{i}, u_{3-i}^*))
	\label{minprobNash}
	\end{equation}
%
for which the input action of \IfTp{}{any }player $P_{3-i}$ is fixed.
Thus, the joint input action is composed of solutions to the minimization problem (\ref{minprobNash}) for each $i\in \IfTp{\{1,2\}}{\mathcal{V}}$. 
\end{remark}}
\begin{remark}{Equivalent costs}
Given $\xi \in \Pi ( \overline{C}) \cup \Pi(D)$ and a strategy \IfIh{$\kappa^*=(\kappa_1^*,\IfTp{\kappa_2^*}{ \dots, \kappa_N^*})\in \mathcal{K}$,}{$\gamma^*=(\gamma_1^*, \dots, \gamma_N^*)\in \mathcal{K}$,} define the set  \IfIh{$\mathcal{U}^*(\xi, \kappa^*)$}{$\mathcal{U}^*(\xi, \gamma^*)$} $:= \{
u=(u_1, \IfTp{u_2}{\dots, u_N
})\in \mathcal{U}
: \exists \phi \in \mathcal{S}_{\HS_{\kappa^*}}(\xi) \textup{ s.t. } $ $\dom \phi \ni (t,j) \mapsto u_i(t,j)= \IfIh{\kappa}{\gamma}_i^*(\IfIh{}{t,j,}\phi(t,j)) {\textup{ for each}}$ ${i\in \IfTp{\{1,2\}}{\mathcal{V}}}\}$. 
{The worst-case cost associated to the solutions to $\HS\IfIncd{}{}$ from $\xi$ under the strategy $\IfIh{\kappa}{\gamma}^*$ 
 is given by $\hat{\mathcal{J}}(\xi,\IfIh{\kappa}{\gamma}^*):=\sup_{u \in \mathcal{U}^*(\xi, \IfIh{\kappa}{\gamma}^*)} \mathcal{J}(\xi, u)$. 
 When
{$\mathcal{J}(\xi, u)=\hat{\mathcal{J}}(\xi, \IfIh{\kappa}{\gamma}^*)$ for every $u \in \mathcal{U}^*(\xi, \IfIh{\kappa}{\gamma}^*)$,} 
  an equivalent condition to 
 (\ref{SaddlePointIneq})}
is
\IfPers{\begin{eqnarray*}
\hat{\mathcal{J}}_1(\xi,\IfIh{\kappa}{\gamma}^*) 
&\leq&  \hat{\mathcal{J}}_1(\xi,\IfIh{(\kappa_{1}, \kappa_{2}^*\IfTp{}{, \dots, \kappa_{N}^*})}{(\gamma_{1}, \gamma_{2}^*, \dots, \gamma_{N}^*)})
 \\
\hat{\mathcal{J}}_2(\xi,\IfIh{\kappa}{\gamma}^*) 
&\leq& \hat{\mathcal{J}}_2(\xi,\IfIh{(\kappa_{1}^*, \kappa_{2}\IfTp{}{, \dots, \kappa_{N}^*})}{(\gamma_{1}^*, \gamma_{2}, \dots, \gamma_{N}^*)})
\IfTp{}{, \\
&\vdots&
\\
\hat{\mathcal{J}}(\xi,\IfIh{\kappa}{\gamma}^*) 
&\leq&  \hat{\mathcal{J}}(\xi,\IfIh{(\kappa_{ 1}^*, \kappa_{2}^*, \dots, \kappa_{N})}{(\gamma_{ 1}^*, \gamma_{2}^*, \dots, \gamma_{N})})}.
\end{eqnarray*}}
\begin{equation*}
\hat{\mathcal{J}}(\xi,\IfIh{(\kappa_{1}^*, \kappa_{2}\IfTp{}{, \dots, \kappa_{N}^*})}{(\gamma_{1}^*, \gamma_{2}, \dots, \gamma_{N}^*)})\leq\hat{\mathcal{J}}(\xi,\IfIh{\kappa}{\gamma}^*) 
\leq  \hat{\mathcal{J}}(\xi,\IfIh{(\kappa_{1}, \kappa_{2}^*\IfTp{}{, \dots, \kappa_{N}^*})}{(\gamma_{1}, \gamma_{2}^*, \dots, \gamma_{N}^*)})
\end{equation*}
 for {each} $\IfIh{\kappa}{\gamma}_i \in \mathcal{K}_i$ {and each} $i \in \IfTp{\{1,2\}}{\mathcal{V}}$.
\label{EqCosts} 
\end{remark}
\IfIncd{}{\begin{remark}{Relation to the literature}
\indent By considering a discrete-time system with the single-valued function $G$ or by considering a continuous-time system with $F$ \IfIncd{single valued and}{} Lipschitz continuous in $\overline{C}$, and by removing the initial condition as an argument of the cost functionals and specifying it in the state equation, Remark \ref{EqCosts} presents equivalent conditions to those in \cite[(6.3)]{basar1999dynamic}. 
Thus, Definition \ref{NashEqNonCoop} covers the definitions of a pure strategy Nash equilibrium in \cite[Sec. 6.2, 6.5]{basar1999dynamic}  for the zero-sum case.
\end{remark}
}
%
\IfPers{Let $\mathcal{A} \subset \mathbb{R}^n$ be a closed set and consider
\begin{equation*}
	\mathcal{X}_\mathcal{A}:=\{ \phi \in \mathcal{X}: \displaystyle \lim \limits_{\substack{(t,j) \in \dom \phi\\
	(t,j)\rightarrow \sup \dom \phi}} \left | 
	\phi(t,j) \right|_\mathcal{A}=0 \}
\end{equation*}
\begin{equation*}
\mathcal{U}_\mathcal{A}(\xi) :=  \left\{ u\in \mathcal{U}: \exists  \phi \in \mathcal{R} (\xi,u) \cap \mathcal{X}_{\mathcal{A}} \right\}
\end{equation*}
To move to another section
}
%
\IfIh{}{ \section{One-player Infinite-horizon Game}
Consider a one-player hybrid game with dynamics $\HS$ described by (\ref{Heq}) for given $(C,F,D,G)$. 
\begin{assumption}{}
The flow map $F$ is single valued and Lipschitz continuous in $\overline{C}$. The jump map $G$ is single valued.
 \label{AssLips}
\end{assumption}
 \pnn{Under Assumption \ref{AssLips}, conditions in Proposition \ref{UniquenessHu} are satisfied, rendering solutions to $\HS$ for a given $u$ from $\xi$ unique.}

Given $\xi \in \Pi ( \overline{C}) \cup \Pi(D)$, a joint input action $u=(u_C, u_D)\in \mathcal{U}$ such that maximal solutions to $\HS  $ from $\xi$ for $u$ are complete, 
the stage cost for flows $L_C:\mathbb{R}^n  \times \mathbb{R}^{m_C} \rightarrow \mathbb{R}_{\geq 0}$, the stage cost for jumps $L_D:\mathbb{R}^n  \times \mathbb{R}^{m_D} \rightarrow \mathbb{R}_{\geq 0}$, and the terminal cost $q: \mathbb{R}^n \rightarrow \mathbb{R}$, we define the cost associated to 
the solution $(\phi,u)$ to $\HS$ from $\xi$, under Assumption \ref{AssLips}, as
\begin{multline}
\mathcal{J}(\xi, u) := 
 \sum_{j=0}^{\sup_j \dom \phi} \int_{t_{j}}^{t_{j+1}} L_C(\phi(t,j),u_{C}(t,j))dt + \\ 
 \sum_{j=0}^{\sup_j \dom \phi -1} \hspace{-0.4cm}L_D(\phi(t_{j+1},j),u_{D}(t_{j+1},j))
 + \underset{(t,j) \in \textup{dom}\phi}{\limsup_{t+j\rightarrow \infty} } q(\phi(t,j))\IfConf{\\}{}
\label{defJ}
\end{multline}
where $\{t_j\}_{j=0}^{\sup_j \dom \phi}$ is a nondecreasing sequence 
associated to the definition of the hybrid time domain of $\phi$; see Definition \ref{htd}.
%
%
%
%
%
%
%
%

The one-player infinite-horizon version of the game is formalized as the following optimization problem.

\textit{Problem ($\star$):} 
Given $\xi \in \mathbb{R}^n$, under Assumption \ref{AssLips}, solve 
\begin{eqnarray}
\underset{ u \in \mathcal{U}_{\HS  }^\infty(\xi)
}{\text{minimize}}
&& \mathcal{J}(\xi, u)
\label{problemih}
\end{eqnarray}
where $\mathcal{U}_{\HS  }^\infty$ is the set of joint input actions yielding maximal complete solutions to $\HS$ as defined in Section 2.1.
\begin{assumption}{Existence of a solution to the game}
Given $\xi \in \Pi ( \overline{C}) \cup \Pi(D)$, a solution to Problem ($\star$) exists, namely there exists an input action $u^*=(u_C^*,u_D^*) \in \mathcal{U}_{\HS  }^\infty (\xi)$ \pnn{such that $\mathcal{J}(\xi, u^*)<\infty$}, that attains the minimum in (\ref{problemih}), and as a consequence satisfies (\ref{NashEqNCoopIneq}) in Definition \ref{NashEqCoop}. In addition, there exists a strategy $\IfIh{\kappa}{\gamma} \in \mathcal{K}$ such that every complete solution to the closed-loop system \IfIh{$\HS_\kappa$}{$\HS_\gamma$} as in (\ref{Hkeq}), from $\xi$, has a cost \pnn{that is finite and} equal to the minimum in (\ref{problemih}). 
%
%
\label{SolExistence}
\end{assumption}
\begin{remark}{Optimal control solution}
The solution to Problem $(\star)$, when it exists, can be expressed in terms of a pure strategy Nash equilibrium $\IfIh{\kappa}{\gamma}^*$ as in Definition \ref{NashEqCoop}. Each $u^*$ rendering a trajectory $\phi^*$ such that $(\phi^*,u^*) \in \mathcal{S}_{\HS}^\infty (\xi)$ with $\dom \phi^* \ni (t,j) \mapsto u^*_i(t,j)= \IfIh{\kappa}{\gamma}_i^*(\IfIh{}{t,j,}\phi^*(t,j))$  for each $i\in \IfTp{\{1,2\}}{\mathcal{V}}$, satisfies
\begin{eqnarray*}
u^*\in  
\underset{ u \in \mathcal{U}_{\HS  }^\infty(\xi)
}{\arg\min}
&& \mathcal{J}(\xi, u)
\label{Nasheq1pih}
\end{eqnarray*}
and is it referred to as an optimal control. 
 The pure strategy Nash equilibrium $\IfIh{\kappa}{\gamma}^*$ coincides with the solution to the infinite-horizon optimal control problem for hybrid systems, see \cite[Problem 1]{ferrante2019certifying}.
\end{remark}

\begin{definition}{Value function} Given $\xi \in \Pi ( \overline{C}) \cup \Pi(D)$, \pnn{under Assumption \ref{AssLips},} the value function  at $\xi$ is given by 
\begin{equation}
\mathcal{J}^*(\xi):=\min_{
 u\in \mathcal{U}_{\HS  }^\infty(\xi)}
\mathcal{J}(\xi,u)
 \label{cost2goih}
\end{equation}
\label{ValueFunction}
\end{definition}


\IfPers{\textit{Hybrid Dynamic Programming}\\
Let us consider specific points of the solution to study the different behaviors in their evolution. For a given complete trajectory $\phi$, we denote $(T,J)= \sup \dom \phi$. 
Let us also denote by $(\phi^*,u^*)$ the solution that attains the infimum in (\ref{cost2goih}).\\
\textbf{Continuous Evolution}\\
Pick a $\xi \in \Pi(C)$ and an $h \in \mathbb{R}_{>0}$ such that every solution  from $\xi$ flows for at least $h$ seconds. 
Thus, the cost-to-go from $\xi$ is given by
\begin{multline*}
V(\xi):=\inf_{
(\phi, u)\in \mathcal{S}_{\HS  }(\xi)}
\int_{0}^{h} L_C(\phi(t,j),u_{C}(t,j))dt
+\int_{h}^{t_1} L_C(\phi(t,j),u_{C}(t,j))dt\\
+\sum_{j=1}^{J } \int_{t_{j}}^{t_{j+1}} L_C(\phi(t,j),u_{C}(t,j))dt +  
 \sum_{j=0}^{J-1} L_D(\phi(t_{j+1},j),u_{D}(t_{j+1},j))
+ \limsup_{(t,j)\rightarrow \sup \textup{dom} \phi} V(\phi(t,j))\\
= \inf_{(\phi, u)\in \hat{\mathcal{S}}^{(h,0)}_{\HS  }(\xi)} 
\left ( 
\int_{0}^{h} L_C(\phi(t,j),u_{C}(t,j))dt
\right.
+ 
\inf_{(\phi, u)\in \mathcal{S}_{\HS  }(\phi^*(h,0))}
 \\
\left( \left. \sum_{j=0}^{J} \int_{t_{j}}^{t_{j+1}} L_C(\phi(t,j),u_{C}(t,j))dt +  
 \sum_{j=0}^{J-1} L_D(\phi(t_{j+1},j),u_{D}(t_{j+1},j))
+ \limsup_{(t,j)\rightarrow \sup \textup{dom} \phi} V(\phi(t,j))\right) \right)
\end{multline*}
 The inner infimum in the previous expression is the cost-to-go from $\phi^*(h,0)$, so we can rewrite
\begin{equation*}
V(\xi)=
\inf_{(\phi, u)\in \hat{\mathcal{S}}^{(h,0)}_{\HS  }(\xi)} 
\int_{0}^{h} L_C(\phi(t,i),u_{C}(t,i))dt
+ V(\phi^*(h,0))
\end{equation*}
By dividing by $h$ and taking the limit, we obtain
\begin{equation*}
0= \lim_{h \rightarrow 0} 
\inf_{(\phi, u)\in \hat{\mathcal{S}}^{(h,0)}_{\HS  }(\xi)} 
\frac{1}{h}\int_{0}^{h} L_C(\phi(t,i),u_{C}(t,i))dt
+ \frac{V(\phi^*(h,0))-V(\phi^*(0,0))}{h}.
\end{equation*}
Thus, given that $V$ is time independent, one has the following point-wise condition to guarantee optimal flow 
%
\begin{equation}
0=  \inf_{u \in \Pi_u(x,C)}
\left( 
 L_C(x,u)
+\left\langle \nabla V(x),F(x,u_C) \right\rangle \right) 
\hspace{1cm} \forall x \in \Pi(C)  
\label{HJBih}
\end{equation}
\textbf{Discrete Evolution}\\
Pick now a $\xi \in \Pi(D)$. 
The cost-to-go from $\xi$ 
is given by
\begin{multline*}
V(\xi):=\inf_{
(\phi, u)\in \mathcal{S}_{\HS  }(\xi)}
L_D(\xi,u_{D}(0,0))\\ 
+\sum_{j=1}^{J} \int_{t_{j}}^{t_{j+1}} L_C(\phi(t,j),u_{C}(t,j))dt +  
 \sum_{j=1}^{J-1} L_D(\phi(t_{j+1},j),u_{D}(t_{j+1},j))
 + \limsup_{(t,j)\rightarrow \sup \textup{dom} \phi} V(\phi(t,j))
 \end{multline*}
 \begin{multline*}
= \inf_
{u_D(0, 0)}
\left( L_D(\xi,u_{D}(0,0))
+ 
\inf_{(\phi, u)\in 
\mathcal{S}_{\HS  }(\phi^*(0, 1))}
\right. \\
\left. \left( \sum_{j=0}^{J} \int_{t_{j}}^{t_{j+1}} L_C(\phi(t,j),u_{C}(t,j))dt +  
 \sum_{j=0}^{J-1} L_D(\phi(t_{j+1},j),u_{D}(t_{j+1},j))
 + \limsup_{(t,j)\rightarrow \sup \textup{dom} \phi} V(\phi(t,j))\right) \right)
\end{multline*}
The inner infimum is the minimum cost for a game starting at $\phi^*(0,1)=G(\phi^*(0,0),u^*_{D}(0,0))$ so we can rewrite
\begin{equation*}
V(\xi)
 = \inf_
 {u_D(0,0) \in \Pi_u(x,D)}
\left( L_D(\phi(0,0),u_{D}(0,0))
+ V(G(\phi(0,0),u_{D}(0,0)))
 \right)
\end{equation*}
Thus, one has the following point-wise condition to guarantee optimal jumps 
\begin{equation}
V(x)=  \inf_{u \in \Pi_u(x,D)}
\left( 
 L_D(x,u)
+ V(G(x,u)) \right) 
\hspace{1cm} \forall x \in \Pi(D)
\label{Bellmanih}
\end{equation}
}

\begin{remark}{Existence of optimal processes for hybrid infinite-horizon optimal control problems} Sufficient conditions to guarantee the existence of a solution to a hybrid infinite-horizon optimal control problem, as in Problem ($\star$), are provided in  \cite{GOEBEL2019153nu}. Given Problem ($\star$), a joint input action $u=(u_C, u_D)\in \mathcal{U}$ such that maximal solutions to $\HS  $ from $\xi \in \Pi(C \cup D)$ for $u$ are complete,  the stage cost for flows $L_C:\mathbb{R}^n  \times \mathbb{R}^{m_C} \rightarrow \mathbb{R}_{\geq 0}$, the stage cost for jumps $L_D:\mathbb{R}^n  \times \mathbb{R}^{m_D} \rightarrow \mathbb{R}_{\geq 0}$, and the terminal cost $q:\reals^n \rightarrow \reals$, one has that if
\begin{itemize}
\item $F(x,u_C)=f_0(x)+ \sum_{l=1}^{m_c} f_l(x)u_{C,l}$, where $f_l$ is continuous,  and there exists $M>0$ such that $|f_l(x)|\leq M(1 + |x|)$ for all $x \in \Pi(\overline{C})$, $l=0,1,\dots,m_C;$
\item $C=C_x \times U_C$ for a closed set $C_x \subset \reals^n$ and a closed and convex set $U_C \subset \reals^{m_C}$;
\item $L_C$ is lower semicontinuous, convex in $u_C$ for every $x \in C_x$, and either
$U_C$ is bounded or there exist $r>1,\eta>0, \delta \in \reals$ so that $L_C(x,u_C) \geq \eta |u_C|^r +\delta$ for every $(x,u_C) \in C$; 
\item $G$ is locally bounded;
\item $L_D$ is lower semicontinuous, and either $D=D_x \times U_D$ for a set $D_x \subset \reals^n$ and a bounded set $U_D \subset \reals^{m_D}$ or $L_D(x,u_D) \geq \beta(|u_D|)$ for every $(x,u_D) \in D$ for some $\beta:[0,\infty) \rightarrow [0,\infty)$ such that $\lim_{r \to \infty} \beta(r)=\infty$; 
\item \pn{The set of forward complete hybrid time domains starting at $(0, 0)$ is closed with respect to set convergence;}
\item The function $q$ is lower semicontinuous;
 \end{itemize}
 then $\mathcal{J}^*$ is lower semicontinuous, and if $\mathcal{J}^*(\xi)<\infty$, then there exists a solution to Problem $(\star)$. 
\end{remark}

\begin{theorem}{HJB for Problem ($\star$)}
Given a one-player hybrid game with dynamics $\HS$ as in (\ref{Heq}) with data $(C,F,D,G)$,
 and stage costs $L_C:\mathbb{R}^n  \times \mathbb{R}^{m_C} \rightarrow \mathbb{R}_{\geq 0} ,L_D:\mathbb{R}^n  \times \mathbb{R}^{m_D} \rightarrow \mathbb{R}_{\geq 0}$, and terminal cost $q: \mathbb{R}^n \rightarrow \mathbb{R}$\pnn{, suppose that Assumptions \ref{AssLips} and \ref{SolExistence} hold}. If there exists a function $V: \mathbb{R}^n \rightarrow \mathbb{R}$ that is continuously differentiable on a neighborhood of $\Pi(C)$, satisfies the \textit{Hamilton-Jacobi-Bellman hybrid equations}  
\eqs{
0=  \min_{u_C \in \Pi_u(x,C)}
\left\{ 
 L_C(x,u_C)
+\left\langle \nabla V(x),F(x,u_C) \right\rangle \right\} 
\hspace{1cm} }{ \forall x \in \Pi(\overline{C}),  
\label{HJBih}
}
\eqs{
V(x)=  \min_{u_D \in \Pi_u(x,D)}
\left\{ 
 L_D(x,u_D)
+ V(G(x,u_D)) \right\} 
\hspace{2cm} }{  \forall x \in \Pi(D) \NotConf{\hspace{0.7cm}}
\label{Bellmanih}
}
and, for each $\xi \in \Pi ( \overline{C}) \cup \Pi(D)$,  each $(\phi, u)\in \mathcal{S}_{\HS  }^\infty(\xi)$ satisfies
\begin{equation}
\underset{(t,j) \in \textup{dom}\phi}{\limsup_{t+j\rightarrow \infty} } V(\phi(t,j))=\underset{(t,j) \in \textup{dom}\phi}{\limsup_{t+j\rightarrow \infty} } q(\phi(t,j))
\label{TerminalCond}
\end{equation}
then
\begin{equation}
\J^*(\xi)= V(\xi) \hspace{1cm} \forall \xi \in \Pi ( \overline{C}) \cup \Pi(D)
\end{equation}
 \label{thHJBih}
 \end{theorem}
 \begin{corollary}{Stationary feedback law}
Given a one-player hybrid game with dynamics $\HS$ as in (\ref{Heq}) with data $(C,F,D,G)$,
 and stage costs $L_C:\mathbb{R}^n  \times \mathbb{R}^{m_C} \rightarrow \mathbb{R}_{\geq 0} ,L_D:\mathbb{R}^n  \times \mathbb{R}^{m_D} \rightarrow \mathbb{R}_{\geq 0}$, and terminal cost $q: \mathbb{R}^n \rightarrow \mathbb{R}$, \pnn{suppose that Assumptions \ref{AssLips} and \ref{SolExistence} hold}. If a continuously differentiable function $V$ satisfying (\ref{HJBih}), (\ref{Bellmanih}) exists, so a minimum of (\ref{HJBih}) and (\ref{Bellmanih}) exist, respectively in $\Pi_u(x,C)$ and $\Pi_u(x,D)$ for every $x \in \Pi(C \cup D)$, \IfIh{}{and given that $L_C,L_D, F, q,$ and $G$ are time-independent, }then the feedback law
$\kappa:=(\kappa_C,\kappa_D): \mathbb{R}^n \rightarrow \mathbb{R}^{m_C} \times \mathbb{R}^{m_D}$  with values
\eqs{
\kappa_C(x)\in\underset{u_C \in \Pi_u(x,C)}{\arg \min} \left\{ 
 L_C(x,u_C)
+\left\langle \nabla V(x),F(x,u_C) \right\rangle \right\} \hspace{1cm} }{  \forall x \in \Pi(\overline{C})
\label{kHJBeqih}
}
and
\eqs{
\kappa_D(x)\in\underset{u_D \in \Pi_u(x,D)}{\arg \min}  \left\{ 
 L_D(x,u_D)
+ V(G(x,u_D)) \right\} \hspace{2cm} }{  \forall x \in \Pi(D) 
\label{kBeqih}
}
defines an optimal control 
at every $x \in \Pi ( \overline{C}) \cup \Pi(D)$. In addition, if for each $\xi \in \Pi ( \overline{C}) \cup \Pi(D)$,  each $\phi \in \mathcal{S}_{\HS_\kappa  }^\infty(\xi)$ satisfies (\ref{TerminalCond}), then $\kappa$ is a pure strategy Nash equilibrium for the one-player infinite-horizon game.
\label{FBLrem}
\end{corollary}
%
\begin{proof}
To show the claim we proceed as follows:
\begin{enumerate}
\item \pnn{Given the differentiable function $V$, select} a strategy $\IfIh{\kappa}{\gamma}$ attaining the minimum in (\ref{HJBih}) and (\ref{Bellmanih}) for each point in $\Pi ( \overline{C}) \cup \Pi(D)$.
\item Pick an initial condition $\xi$ and evaluate the cost associated to any 
solution yield by $\IfIh{\kappa}{\gamma}$ from $\xi$. Show that this cost coincides with the value of the function $V$ at $\xi$. 
\item Lower bound the cost associated to any other solution from $\xi$ by the value of the function $V$ evaluated at $\xi$. 
\item By showing that the cost of any other solution from $\xi$ is not less than the cost of the solution yield by $\IfIh{\kappa}{\gamma}$ from $\xi$, we show optimality of $\IfIh{\kappa}{\gamma}$ in Problem ($\star$). 
\end{enumerate}
\NotConf{Proceeding as in item 1 above, 
given 
a function $V: \reals^n \to \reals$ that is continuously differentiable on $\Pi(C)$,
let $\IfIh{\kappa}{\gamma}:=(\IfIh{\kappa}{\gamma}_C,\IfIh{\kappa}{\gamma}_D) \in \mathcal{K}$ 
attain the minimum in (\ref{HJBih}) and (\ref{Bellmanih}) for every $x \in \Pi(\overline{C} \cap D)$, i.e.,
\begin{equation*}
\IfIh{\kappa}{\gamma}_C(x)\in\underset{u_C \in \Pi_u(x,C)}{\argmin} \left\{ 
L_C(x,u_C)
+\left\langle \nabla V(x),F(x,u_C) \right\rangle \right\}
\hspace{1cm} \forall x \in \Pi(\overline{C})
\label{kHJBeqihproofth}
\end{equation*}
and
\begin{equation*}
\IfIh{\kappa}{\gamma}_D(x)\in \underset{u_D \in \Pi_u(x,D)}{\argmin} \left\{ 
L_D(x,u_D) + V(G(x,u_D)) \right\}
\hspace{1cm} \forall x \in \Pi(D) 
\label{kBeqihproofth}
\end{equation*}

Continuing as in item 2 above, pick any $\xi \in \Pi ( \overline{C}) \cup \Pi(D)$ and any $(\phi^*,u^*) \in \mathcal{S}_{\HS }^\infty (\xi)$ with $\dom \phi^* \ni (t,j) \mapsto u^*(t,j)= \IfIh{\kappa}{\gamma}(\IfIh{}{t,j,}\phi^*(t,j))$, 
and $(t_{J_{\phi^*}+1},J_{\phi^*}) \in \dom (\phi^*,u^*) $. 
\pnn{We show that there does not exist a complete solution to $\HS$ with cost smaller than the cost of $(\phi^*,u^*)$.}
Given that $V$ satisfies 
(\ref{HJBih}), and $\IfIh{\kappa}{\gamma}_C$ 
attains the minimum in (\ref{HJBih}), for each $j \in \nats$ such that $I_{\phi^*}^j=[t_j,t_{j+1}]$ has a nonempty interior int$I_{\phi^*}^j$, we have that for every $t \in \textup{int}I_{\phi^*}^j$, 
\begin{equation*}
0=L_C(\phi^*(t,j),\kappa_C(\phi^*(t,j)))+\left\langle \nabla V(\phi^*(t,j),\kappa_C(\phi^*(t,j))),F(\phi^*(t,j),\kappa_C(\phi^*(t,j))) \right\rangle
\end{equation*}
Given that  $\frac{d V}{d t}(\phi^*(t,j)):=\left\langle \nabla V(\phi^*(t,j),\kappa_C(\phi^*(t,j))),F(\phi^*(t,j),\kappa_C(\phi^*(t,j))) \right\rangle$ for every $(t,j):t \in \dom_t \phi^*$, 
by integrating over the interval $[t_j,t_{j+1}]$,
we obtain
\begin{eqnarray*}
0=  \int_{t_{j}}^{t_{j+1}} \left( L_C(\phi^*(t,j),\IfIh{\kappa}{\gamma}_{ C}(\phi^*(t,j)))
+ \frac{d V}{d t}(\phi^*(t,j)) \right )dt
\end{eqnarray*}
from where we have
\begin{eqnarray*}
0=  \int_{t_{j}}^{t_{j+1}} L_C(\phi^*(t,j),\IfIh{\kappa}{\gamma}_{ C}(\phi^*(t,j))) dt
+  V(\phi^*(t_{j+1},j))- V(\phi^*(t_j,j))
\end{eqnarray*}
%
%
%
Summing from $j=0$ to $j=J_{\phi^*}$, 
we obtain 
\begin{eqnarray*}
0
=\sum_{j=0}^{J_{\phi^*}}  \int_{t_{j}}^{t_{j+1}} L_C(\phi^*(t,j),\IfIh{\kappa}{\gamma}_{ C}(\phi^*(t,j))) dt
+\sum_{j=0}^{J_{\phi^*}} \left( V(\phi^*(t_{j+1},j))-V(\phi^*(t_j,j)) \right)
\end{eqnarray*}
Then, solving for $V$ at the initial condition $\phi^*(0,0)$, we obtain
%
%
\begin{eqnarray}
V(\phi^*(0,0))
=\sum_{j=0}^{J_{\phi^*}}  \int_{t_{j}}^{t_{j+1}} L_C(\phi^*(t,j),\IfIh{\kappa}{\gamma}_{ C}(\phi^*(t,j))) dt+ V(\phi^*(t_1,0))\nonumber \\
+\sum_{j=1}^{J_{\phi^*}} \left( V(\phi^*(t_{j+1},j))-V(\phi^*(t_j,j)) \right) \label{ContinCost2goih}
\end{eqnarray}
Given that $V$ satisfies 
(\ref{Bellmanih}) and $\IfIh{\kappa}{\gamma}_D$ 
attains the infimum in (\ref{Bellmanih}), for every $(t_{j+1},j) \in \dom \phi^*$ such that $(t_{j+1},j+1) \in \dom \phi^*$,  we have that
\begin{eqnarray*}
V(\phi^*(t_{j+1},j))&=&  \min_{
u_{D}(t_{j+1},j) \in \Pi_u(\phi^*(t_{j+1},j),D)}
\left\{ 
 L_D(\phi^*(t_{j+1},j),u_{ D}(t_{j+1},j))
+ V(G(\phi^*(t_{j+1},j), u_D(t_{j+1},j) \right\}\\ 
&=&  
 L_D(\phi^*(t_{j+1},j),\IfIh{\kappa}{\gamma}_{ D}(\phi^*(t_{j+1},j)))
+ V(G(\phi^*(t_{j+1},j), \IfIh{\kappa}{\gamma}_D(\phi^*(t_{j+1},j)))) \\
&=& L_D(\phi^*(t_{j+1},j),\IfIh{\kappa}{\gamma}_{ D}(\phi^*(t_{j+1},j)))
+ V(\phi^*(t_{j+1},j+1)) 
\end{eqnarray*}
where $\phi^*(t_{j+1},j)\in D_{\IfIh{\kappa}{\gamma}}$ is defined in (\ref{Hkeq}).
Summing both sides from $j=0$ to $j=J_{\phi^*}-1$, 
we obtain
\begin{eqnarray*}
\sum_{j=0}^{J_{\phi^*}-1} V(\phi^*(t_{j+1},j))
= \sum_{j=0}^{J_{\phi^*}-1} L_D(\phi^*(t_{j+1},j),\IfIh{\kappa}{\gamma}_{ D}(\phi^*(t_{j+1},j))) 
+\sum_{j=0}^{J_{\phi^*}-1} V( \phi^*(t_{j+1},j+1))
\end{eqnarray*}
Then, solving for $V$ at the first jump time $(t_1,0)$, we obtain
%
\begin{eqnarray}
V(\phi^*(t_{1},0))
= V( \phi^*(t_{1},1))
+ \sum_{j=0}^{J_{\phi^*}-1} L_D(\phi^*(t_{j+1},j),\IfIh{\kappa}{\gamma}_{ D}(\phi^*(t_{j+1},j))) \nonumber
\\
+\sum_{j=1}^{J_{\phi^*}-1}\left( V( \phi^*(t_{j+1},j+1))-V(\phi^*(t_{j+1},j)) \right)
\label{DiscreteCost2goih}
\end{eqnarray}
Given that $\phi^*(0,0)=\xi$, by substituting the right-hand side of (\ref{DiscreteCost2goih}) in  (\ref{ContinCost2goih}), we obtain
\begin{eqnarray}
V(\xi)
&=&\sum_{j=0}^{J_{\phi^*}}  \int_{t_{j}}^{t_{j+1}} L_C(\phi^*(t,j),\IfIh{\kappa}{\gamma}_{ C}(\phi^*(t,j))) dt+ V(\phi^*(t_1,0))\nonumber 
+\sum_{j=1}^{J_{\phi^*}} \left( V(\phi^*(t_{j+1},j))-V(\phi^*(t_j,j)) \right) \nonumber \\
&=&\sum_{j=0}^{J_{\phi^*}}  \int_{t_{j}}^{t_{j+1}} L_C(\phi^*(t,j),\IfIh{\kappa}{\gamma}_{ C}(\phi^*(t,j))) dt
+ \sum_{j=0}^{J_{\phi^*}-1}L_D(\phi^*(t_{j+1},j),\IfIh{\kappa}{\gamma}_{ D}(\phi^*(t_{j+1},j)))\nonumber \\ &&+ V(\phi^*(t_{1},1))+\sum_{j=1}^{J_{\phi^*}-1}\left( V(\phi^*(t_{j+1},j+1))-V(\phi^*(t_{j+1},j)) \right) \nonumber 
\\&&+\sum_{j=1}^{J_{\phi^*}} \left( V(\phi^*(t_{j+1},j))-V(\phi^*(t_j,j)) \right) \nonumber 
\end{eqnarray}
\pnn{Since 
\begin{eqnarray}
&& V(\phi^*(t_{1},1))+\sum_{j=1}^{J_{\phi^*}-1}\left( V(\phi^*(t_{j+1},j+1))-V(\phi^*(t_{j+1},j)) \right) \nonumber 
+\sum_{j=1}^{J_{\phi^*}} \left( V(\phi^*(t_{j+1},j))-V(\phi^*(t_j,j)) \right) \nonumber
\\ &=& V(\phi^*(t_{J_{\phi^*}+1},J_{\phi^*}))+V(\phi^*(t_{1},1))+\sum_{j=1}^{J_{\phi^*}-1}\left( V(\phi^*(t_{j+1},j+1)) \right) \nonumber 
-\sum_{j=1}^{J_{\phi^*}} \left(V(\phi^*(t_j,j)) \right) \nonumber
\\ &=&V(\phi^*(t_{J_{\phi^*}+1},J_{\phi^*})) \nonumber
\end{eqnarray}
then we have}
\begin{eqnarray}
V(\xi)
&=&\sum_{j=0}^{J_{\phi^*}}  \int_{t_{j}}^{t_{j+1}} L_C(\phi^*(t,j),\IfIh{\kappa}{\gamma}_{ C}(\phi^*(t,j))) dt
+ \sum_{j=0}^{J_{\phi^*}-1}L_D(\phi^*(t_{j+1},j),\IfIh{\kappa}{\gamma}_{ D}(\phi^*(t_{j+1},j))) \nonumber
\\ &&+V(\phi^*(t_{J_{\phi^*}+1},J_{\phi^*}))
\nonumber
\end{eqnarray}

By taking the limit when $(t_{J_{\phi^*}+1},J_{\phi^*}) \rightarrow \sup \dom \phi^*$, and given that (\ref{TerminalCond}) holds, we have 
\begin{eqnarray}
V(\xi)&=& \sum_{j=0}^{\sup_j \dom \phi^*}   \int_{t_{j}}^{t_{j+1}} L_C(\phi^*(t,j),\IfIh{\kappa}{\gamma}_{ C}(\phi^*(t,j))) dt
+ \sum_{j=0}^{\sup_j \dom \phi^* -1} L_D(\phi^*(t_{j+1},j),\IfIh{\kappa}{\gamma}_{ D}(\phi^*(t_{j+1},j))) \nonumber 
\\ &&
+\underset{(t,j) \in \textup{dom}\phi^*}{\limsup_{t+j\rightarrow \infty} } V(\phi^*(t,j))
\nonumber \\ 
&=& \mathcal{J}(\xi,u^*)
\label{OptimalCostih}
\end{eqnarray}
for any $(\phi^*,u^*) \in \mathcal{S}_{\HS  }^\infty (\xi)$ with $\dom \phi^* \ni (t,j) \mapsto u^*(t,j)= \IfIh{\kappa}{\gamma}(\IfIh{}{t,j,}\phi^*(t,j))$.

Continuing with item 3 as above, pick any $(\phi_s,u) \in \mathcal{S}_{\HS  }^\infty (\xi)$ with 
$\dom \phi_s \ni (t,j) \mapsto u(t,j)= \bar{\IfIh{\kappa}{\gamma}}(\IfIh{}{t,j,}\phi_s(t,j))$ for some $\bar{\IfIh{\kappa}{\gamma}} \in \mathcal{K}$ 
  and $(t_{J_{\phi_s}+1},J_{\phi_s}) \in \dom (\phi_s,u)$. Since $\bar{\IfIh{\kappa}{\gamma}}$ does not necessarily attain the minimum in (\ref{HJBih}), then,  for each $j \in \nats$ such that $I_{\phi_s}^j=[t_j,t_{j+1}]$ has a nonempty interior int$I_{\phi_s}^j$, we have that for every $t \in \textup{int}I_{\phi_s}^j$,
\begin{eqnarray*}
0
\leq   L_C(\phi_s(t,j),u_{ C}(t,j))
+\left\langle \nabla V(\phi_s(t,j),u_C(t,j)),F(\phi_s(t,j),u_C(t,j)) \right\rangle
\end{eqnarray*}
Given that  $\frac{d V}{d t}(\phi_s(t,j)):=\left\langle \nabla V(\phi_s(t,j),u_C(t,j)),F(\phi_s(t,j),u_C(t,j)) \right\rangle$ for every $(t,j):t \in \dom_t \phi_s$, 
by integrating over the interval $[t_j,t_{j+1}]$, 
we obtain
\begin{eqnarray*}
0 \leq \int_{t_{j}}^{t_{j+1}}  \left( L_C(\phi_s(t,j),u_{ C}(t,j))
+  \frac{d V}{d t}(\phi_s(t,j)) \right) dt
\end{eqnarray*}
from which we have
\begin{eqnarray*}
 V(\phi_s(t_j,j))
\leq \int_{t_{j}}^{t_{j+1}}   L_C(\phi_s(t,j),u_{ C}(t,j))dt
+V(\phi_s(t_{j+1},j))
\end{eqnarray*}
Summing both sides from  $j=0$ to $j=J_{\phi_s}$, 
we obtain 
\begin{eqnarray}
\sum_{j=0}^{J_{\phi_s}}V(\phi_s(t_j,j)) 
\leq \sum_{j=0}^{J_{\phi_s}}  \int_{t_{j}}^{t_{j+1}}  L_C(\phi_s(t,j),u_{C}(t,j))dt
+\sum_{j=0}^{J_{\phi_s}}V(\phi_s(t_{j+1},j))
\nonumber
\end{eqnarray}
Then, solving for $V$ at the initial time $(0,0)$, we obtain
%
%
\begin{eqnarray}
V(\phi_s(0,0))
\leq \sum_{j=0}^{J_{\phi_s}}  \int_{t_{j}}^{t_{j+1}}  L_C(\phi_s(t,j),u_{C}(t,j))dt+ V(\phi_s(t_1,0))\nonumber \\+\sum_{j=1}^{J_{\phi_s}} \left( V(\phi_s(t_{j+1},j))-V(\phi_s(t_j,j)) \right) 
\label{ContinCost2goNOih}
\end{eqnarray}
In addition, since $u$ defined by $\bar{\IfIh{\kappa}{\gamma}}$ does not necessarily attain the infimum in (\ref{Bellmanih}), then for every $(t_{j+1},j) \in \dom \phi_s $ such that   $(t_{j+1},j+1) \in \dom \phi_s$, we have
\begin{eqnarray*}
V(\phi_s(t_{j+1},j))
&\leq&   
L_D(\phi_s(t_{j+1},j),u_{ D}(t_{j+1},j))
+ V(G(\phi_s(t_{j+1},j),u_{D}(t_{j+1},j)))\\
&=&L_D(\phi_s(t_{j+1},j),u_{ D}(t_{j+1},j))
+ V( \phi_s(t_{j+1},j+1))
\end{eqnarray*}
Summing both sides from  $j=0$ to $j=J_{\phi_s}-1$, 
we obtain
\begin{eqnarray}
\sum_{j=0}^{J_{\phi_s}-1}V(\phi_s(t_{j+1},j)) 
\leq 
 \sum_{j=0}^{J_{\phi_s}-1}L_D(\phi_s(t_{j+1},j),u_{ D}(t_{j+1},j))
+\sum_{j=0}^{J_{\phi_s}-1} V(\phi_s(t_{j+1},j+1))
\nonumber
\end{eqnarray}
Then, solving for $V$ at the first jump time $(t_1,0)$, we obtain
%
\begin{eqnarray}
V(\phi_s(t_{1},0))
\leq V( \phi_s(t_{1},1))
+ \sum_{j=0}^{J_{\phi_s}-1}L_D(\phi_s(t_{j+1},j),u_{ D}(t_{j+1},j))\nonumber
\\
+\sum_{j=1}^{J_{\phi_s}-1}\left( V(\phi_s(t_{j+1},j+1))-V(\phi_s(t_{j+1},j)) \right)
\label{DiscreteCost2goNOih}
\end{eqnarray}
In addition, given that $\phi_s(0,0)=\xi$, upper bounding $V(\phi_s(t_{1},0))$ in (\ref{ContinCost2goNOih}) by the right-hand side of (\ref{DiscreteCost2goNOih}), we obtain 
\begin{eqnarray}
V(\xi)
&\leq& \sum_{j=0}^{J_{\phi_s}}  \int_{t_{j}}^{t_{j+1}}  L_C(\phi_s(t,j),u_{C}(t,j))dt+ V(\phi_s(t_1,0)) \nonumber 
\\&&+\sum_{j=1}^{J_{\phi_s}} \left( V(\phi_s(t_{j+1},j))-V(\phi_s(t_j,j)) \right) \nonumber \\
&\leq&
\sum_{j=0}^{J_{\phi_s}}  \int_{t_{j}}^{t_{j+1}}  L_C(\phi_s(t,j),u_{C}(t,j))dt
+ \sum_{j=0}^{J_{\phi_s}-1}L_D(\phi_s(t_{j+1},j),u_{D}(t_{j+1},j))\nonumber
\\&&+\sum_{j=1}^{J_{\phi_s}-1}\left( V(\phi_s(t_{j+1},j+1))-V(\phi_s(t_{j+1},j)) \right) \nonumber \\
&&+ V(\phi_s(t_{1},1))+\sum_{j=1}^{J_{\phi_s}} \left( V(\phi_s(t_{j+1},j))-V(\phi_s(t_j,j)) \right)\nonumber 
\nonumber
\end{eqnarray}
\pnn{Since 
\begin{eqnarray}
&& V(\phi_s(t_{1},1))+\sum_{j=1}^{J_{\phi_s}-1}\left( V(\phi_s(t_{j+1},j+1))-V(\phi_s(t_{j+1},j)) \right) \nonumber 
+\sum_{j=1}^{J_{\phi_s}} \left( V(\phi_s(t_{j+1},j))-V(\phi_s(t_j,j)) \right) \nonumber
\\ &=& V(t_{{J_{\phi_s}}+1},J_{\phi_s}) +V(\phi_s(t_{1},1))+\sum_{j=1}^{J_{\phi_s}-1}\left( V(\phi_s(t_{j+1},j+1)) \right) \nonumber 
-\sum_{j=1}^{J_{\phi_s}} \left(V(\phi_s(t_j,j)) \right) \nonumber
\\ &=&V(t_{{J_{\phi_s}}+1},J_{\phi_s}) \nonumber
\end{eqnarray}
then we have}
\begin{eqnarray}
V(\xi)
&=&\sum_{j=0}^{J_{\phi_s}}  \int_{t_{j}}^{t_{j+1}}  L_C(\phi_s(t,j),u_{C}(t,j))dt
+ \sum_{j=0}^{J_{\phi_s}-1}L_D(\phi_s(t_{j+1},j),u_{D}(t_{j+1},j))\nonumber
\\&&+V(\phi_s(t_{J_{\phi_s}+1},J_{\phi_s}))
\nonumber
\end{eqnarray}
By taking the limit when $(t_{{J_{\phi_s}}+1},J_{\phi_s}) \rightarrow \sup \dom \phi_s$, and given that (\ref{TerminalCond}) holds, we have
\begin{eqnarray}
V(\xi)&\leq&\sum_{j=0}^{\sup_j \dom \phi_s} \int_{t_{j}}^{t_{j+1}}  L_C(\phi_s(t,j),u_{C}(t,j))dt
+ \sum_{j=0}^{{\sup_j \dom \phi_s}-1}L_D(\phi_s(t_{j+1},j),u_{D}(t_{j+1},j))\nonumber
\\ &&
+\underset{(t,j) \in \textup{dom}\phi_s}{\limsup_{t+j\rightarrow \infty} } V(\phi_s(t,j))
\nonumber \\ 
&\leq& \mathcal{J}(\xi,u)
\label{BoundCostNO}
\end{eqnarray}
for $u$ defined by any $\bar{\IfIh{\kappa}{\gamma}} \in \mathcal{K}$. 

\pnn{Finally, proceeding as in item 4 above, under Assumption \ref{SolExistence},
by applying the $\min$ on each side of  (\ref{BoundCostNO}) over the set $\mathcal{U}_\HS^\infty (\xi)$, we obtain
\begin{equation*}
V(\xi) \leq  \underset{ u \in \mathcal{U}_{\HS  }^\infty(\xi)
}{\text{min}}\mathcal{J}(\xi, u).
\end{equation*}
Given that $V(\xi)=\mathcal{J}(\xi, u^
*)$ from (\ref{OptimalCostih}), we have that for any $\xi \in \Pi ( \overline{C}) \cup \Pi(D)$, each 
$u^*$ defined by $\IfIh{\kappa}{\gamma}$ satisfies
\begin{equation*}
\mathcal{J}(\xi,u^*) \leq \underset{ u \in \mathcal{U}_{\HS  }^\infty(\xi)
}{\text{min}} \mathcal{J}(\xi, u)
\end{equation*}
which implies
\begin{equation}
\mathcal{J}(\xi,u^*) = \underset{ u \in \mathcal{U}_{\HS  }^\infty(\xi)
}{\text{min}} \mathcal{J}(\xi, u).
\label{LowerBoundd}
\end{equation}
Thus, from (\ref{OptimalCostih}) and (\ref{LowerBoundd}), $\kappa$ is the optimal strategy and $V(\xi)$ is the value function for $\HS$, as in Definition \ref{ValueFunction}.}}
\end{proof}
\begin{remark}{Equivalent conditions}
If there exists a continuously differentiable function $V: \mathbb{R}^n \rightarrow \mathbb{R}$ on $\Pi(C)$ that satisfies (\ref{HJBih}), (\ref{Bellmanih}) and a feedback law $\kappa:=(\kappa_C,\kappa_D): \mathbb{R}^n \rightarrow \mathbb{R}^{m_C} \times \mathbb{R}^{m_D}$ that satisfies (\ref{kHJBeqih}), (\ref{kBeqih}), then
\begin{eqnarray}
&& \IfAutom{\mathcal{L}_C(x,\kappa_C(x))}
  {L_C(x,\kappa_C(x))
+
\left\langle \nabla V(x),F(x,\kappa_C(x)) \right\rangle}=0  \hspace{1.2cm} \IfConf{\nonumber \\&& \hspace{6cm} }{} \forall x \in C_\kappa 
\label{kHJBeqiha}
\\
&&L_C(x,u_C)
+\left\langle \nabla V(x),F(x,u_C) \right\rangle \geq 0 \hspace{1.2cm}\IfConf{\nonumber \\&& \hspace{5.2cm} }{}\forall (x,u_C) \in C 
\label{kHJBeqihb}
\\
&& L_D(x,\kappa_D(x))
+ V(G(x,\kappa_D(x)))-V(x)=0  \hspace{1.2cm}\IfConf{\nonumber \\&& \hspace{6cm} }{}\forall x \in D_\kappa 
\label{kBeqihc}
\\
&& L_D(x,u_D)
+ V(G(x,u_D)) -V(x) \geq0 \hspace{1.2cm} \IfConf{\nonumber \\&& \hspace{5.2cm} }{}\forall (x,u_D) \in D
\label{kBeqihd}
\end{eqnarray}
In fact, from (\ref{kHJBeqih}), (\ref{kBeqih}) we have
\begin{multline}
\underset{u_C \in \Pi_u(x,C)}{\min} \left\{ 
 L_C(x,u_C)
+\left\langle \nabla V(x),F(x,u_C) \right\rangle \right\} \\
=L_C(x,\kappa_C(x))
+\left\langle \nabla V(x),F(x,\kappa_C(x))\right\rangle \hspace{0.9cm}\forall x \in \Pi(C) \IfConf{\\}{}
\label{kHJBeqihproof}
\end{multline}
and
\begin{multline}
\underset{u_D \in \Pi_u(x,D)}{\min}  \left\{ 
 L_D(x,u_D)
+ V(G(x,u_D)) \right\}\IfConf{\\}{}=L_D(x,\kappa_D(x))
+ V(G(x,\kappa_D(x))) \hspace{1cm} \forall x \in \Pi(D) 
\label{kBeqihproof}
\end{multline}
Thus, (\ref{kHJBeqiha}) and (\ref{kBeqihc}) are implications of (\ref{kHJBeqihproof}), (\ref{HJBih}) and (\ref{kBeqihproof}), (\ref{Bellmanih}) respectively. Similarly, (\ref{kHJBeqihb}) and (\ref{kBeqihd}) are derived straightforwardly from (\ref{HJBih}) and (\ref{Bellmanih}). 
\end{remark}
From \cite[Remark 3]{ferrante2019certifying}, if $C_\kappa=\Pi(C), D_\kappa=\Pi(D)$, then $V$ and $\kappa$ satisfying (\ref{kHJBeqiha}), (\ref{kHJBeqihb}), (\ref{kBeqihc}), and (\ref{kBeqihd}) satisfy also (\ref{HJBih}), (\ref{Bellmanih}), (\ref{kHJBeqih}), and (\ref{kBeqih}).

Next, we illustrate in some examples with linear dynamics and quadratic costs how Theorem \ref{thHJBih} provides conditions to characterize
the value function and the Nash equilibrium in a noncooperative infinite-horizon one-player game covering both continuous, discrete and hybrid systems.
\begin{example}{Linear Quadratic Differential Games}
Consider a one-player hybrid game with state $x \in \reals^n$, input $u=(u_C, u_D) \in \reals^{m_C}\times \reals^{m_D}$, and dynamics $\HS$ as in (\ref{Heq}) with 
\begin{eqnarray}
C&:=&\mathbb{R}^n \times \mathbb{R}^{m_C} \nonumber
\\
F(x,u_c)&:=&A_C x+B_C u_C \nonumber
\label{LinearGame}
\end{eqnarray}
and jump set $D=\emptyset$, where $A_C \in \mathbb{R}^{n\times n}$ and $B_C \in \mathbb{R}^{n\times m_C}$. 
Consider the stage costs $L_C(x,u_C):=x^\top Q_C x+ u_C^\top R_C u_C$ and $L_D(x,u_D):=0$, the terminal cost $q(x):=x^\top P x$ and the function $V(x):=x^\top P x$, where $Q_C \in \mathbb{S}^n_+$, $R_C \in \mathbb{S}^{m_C}_+$ and  $P\in \mathbb{S}^n_+$. We refer to this game as an infinite-horizon one-player linear quadratic differential game $(\mathcal{LQC})$. The function $V$ is continuously differentiable in $\reals^n$, and
\begin{eqnarray}
 \min_{u \in \Pi_u(x,C)}
\left\{L_C(x,u) +\left\langle \nabla V(x),F(x,u_C) \right\rangle \right\}
&=&
  \min_{u_C \in \mathbb{R}^{m_C}}
\left\{x^\top Q_C x+ u_C^\top R_C u_C +2x^\top P (A_C x+B_C u_C ) \right\}
\nonumber
\\&=&
  \min_{u_C \in \mathbb{R}^{m_C}}
\left\{ x^\top (Q_C+2 P A_C) x+ u_C^\top R_C u_C
+ 2 x^\top  P B_C u_C  \right\}
\nonumber
\\&=&0
\label{HJBLQGC}
\end{eqnarray}
is satisfied $ \forall x \in \reals^n$ (and attained by $\kappa_C(x)=-R_C^{-1}B_C^\top P x$)
if the continuous time algebraic Riccati equation (CARE)
\begin{equation}
0=Q_C +2PA_C- P B_C R_C^{-1}B_C ^\top P
\label{Riccati}
\end{equation}
is satisfied. 
Given that $V$ is continuously differentiable, and that (\ref{Bellmanih}) and (\ref{HJBih}) are satisfied thanks to  (\ref{HJBLQGC}) and $D=\emptyset$, from Theorem \ref{thHJBih}, the value function is
\begin{eqnarray}
\mathcal{J}^*(\xi) &:=& 
 \xi^\top P \xi.
\label{ContLCost}
\end{eqnarray}
for any $\xi \in \reals^n$.
%
\IfPers{Given that this is a purely continuous system, the calculation of the cost for when applying the feedback law $\kappa_C(x)=-R_C^{-1}B_C^\top P x$, which renders a solution $(\phi,u_C^*)$ with $\dom \phi \ni (t,j) \mapsto u_C^*(t,j)= \kappa_C(\phi(t,j))$, is as follows
\begin{equation*}
\mathcal{J}(\xi,u_C^*) =
 \int_{t_0}^{\infty} L_C(\phi(t,j),u_{C}^*(t,j))dt  + \underset{(t,j) \in \textup{dom}\phi}{\limsup_{t+j\rightarrow \infty} } q(\phi(t,j))
\end{equation*}
and, when (CARE) is satisfied, one has
\begin{eqnarray*}
 \int_{t_0}^{\infty} L_C(\phi(t,j),u_{C}^*(t,j))dt  &=& 
 \int_{0}^{\infty} \phi^\top(t,j) Q_C \phi(t,j) +{u_{C}^*}^\top(t,j) R_C u_{C}^*(t,j)dt  
 \\&=&  \int_{0}^{\infty} \phi^\top(t,j) (Q_C +P B_C R_C^{-1}B_C ^\top P)\phi(t,j) dt  
  \\&=&  \int_{0}^{\infty} \phi^\top(t,j) (2P B_C R_C^{-1}B_C ^\top P-2PA_C)\phi(t,j) dt  
  \\&=&-\int_{0}^{\infty} \frac{d V(\phi(t,j))}{d t} dt
    \\&=&V(\phi(0,0))-\limsup_{t\rightarrow \infty}V(\phi(t,0))
  \\&=& 
  \phi^\top(0,0) P \phi(0,0)- \limsup_{t\rightarrow \infty} \phi^\top(t,0) P \phi(t,0)
\end{eqnarray*}
which, given that $\phi(0,0)=\xi$ and $q(\phi(t,j))=\phi^\top(t,0) P \phi(t,0)$, implies
$	\mathcal{J}(\xi,u_C^*) =\xi^\top P \xi$, showing that the feedback law $\kappa_C$ attains the optimal cost.} 
From Corollary \ref{FBLrem}, the feedback law $\kappa_C$ is a pure strategy Nash equilibrium for the $(\mathcal{LQC})$ game.
\end{example}
\begin{example}{Linear Quadratic Difference Games}
Consider a one-player hybrid game with state $x \in \reals^n$, input $u=(u_C, u_D) \in \reals^{m_C}\times \reals^{m_D}$, and dynamics $\HS$ as in (\ref{Heq}) with
\begin{eqnarray}
D&=&\mathbb{R}^n \times \mathbb{R}^{m_D} \nonumber
\\
G(x,u_D)&=&A_D x+B_D u_D
\label{LinearDGame}
\end{eqnarray}
and flow set $C=\emptyset$, where $A_D \in \mathbb{R}^{n\times n}$ and $B_D \in \mathbb{R}^{n\times m_D}$. 
Consider the stage costs $L_C(x,u_C):=0$ and $L_D(x,u_D):=x^\top Q_D x+ u_D^\top R_D u_D$, the terminal cost $q(x):=x^\top P x$ and the function $V(x):=x^\top P x$, where $Q_D\in \mathbb{S}^n_+$, $R_D\in \mathbb{S}^{m_D}_+$ and  $P \in \mathbb{S}^n_+$. We refer to this game as an infinite-horizon one player linear quadratic difference game $(\mathcal{LQD})$. The function $V$ is such that
\begin{eqnarray}
\min_{u \in \Pi_u(x,D)}
&&\left\{ 
 L_D(x,u)
+ V(G(x,u)) \right\} \nonumber
\\&=&
  \min_{u_D \in \mathbb{R}^{m_D}}
\left\{x^\top Q_D x+ u_D^\top R_D u_D +(A_D x+B_D u_D )^\top P (A_D x+B_D u_D ) \right\}
\nonumber
\\&=&
  \min_{u_D \in \mathbb{R}^{m_D}}
\left\{  x^\top (Q_D+A_D^\top P A_D) x+ u_D^\top (R_D +B_D^\top P B_D) u_D
+ 2 x^\top A_D^\top  P B_D u_D\right\}
\nonumber
\\&=&x^\top P x,
\label{HJBLQGD}
\end{eqnarray}
is satisfied $ \forall x \in \reals^n$ (and attained by $\kappa_D(x)=-(R_D+B_D^\top PB_D)^{-1}B_D^\top P A_D x$)
when the discrete-time algebraic Riccati equation (DARE)
\begin{equation}
P=Q_D+A_D^\top P A_D-A_D^\top P B_D (R_D+B_D^\top P B_D)^{-1}B_D^\top P A_D
\label{DARE}
\end{equation}
is satisfied. 
Given that $V$ is continuously differentiable, and that (\ref{Bellmanih}) and (\ref{HJBih}) are satisfied thanks to $C=\emptyset$ and (\ref{HJBLQGD}), from Theorem \ref{thHJBih}, the value function is

\begin{eqnarray}
\mathcal{J}^*(\xi) &:=& 
 \xi^\top P \xi
\label{ContLCost}
\end{eqnarray}
for any $\xi \in \reals^n$.
\IfPers{Given that this is a purely discrete system, the calculation of the cost for when applying the feedback law $\kappa_D(x)=-(R_D+B_D^\top PB_D)^{-1}B_D^\top P A_D x$, which renders a solution $(\phi,u_D^*)$ with $\dom \phi \ni (t,j) \mapsto u_D^*(t,j)= \kappa_D(\phi(t,j))$, is as follows
\begin{equation*}
\mathcal{J}(\xi,u_D^*) =
  \sum_{j=0}^{\infty} L_D(\phi(t_{j+1},j),u_{D}^*(t_{j+1},j))+ \underset{(t,j) \in \textup{dom}\phi}{\limsup_{t+j\rightarrow \infty} } q(\phi(t,j))
\end{equation*}
and, when (DARE) is satisfied, one has
\begin{eqnarray*}
 \sum_{j=0}^{\infty} L_D(\phi(t_{j+1},j),u_{D}^*(t_{j+1},j)) &=& 
 \sum_{j=0}^{\infty} \phi^\top(t,j)Q_D \phi(t,j) +{u_{D}^*}^\top(t,j) R_D u_{D}^*(t,j)  
 \\&=&  \sum_{j=0}^{\infty} \phi^\top(t,j) (Q_D +A_D^\top P B_D (R_D+B_D^\top P B_D)^{-1}B_D^\top P A_D)\phi(t,j)   
  \\&=&  \sum_{j=0}^{\infty} \phi^\top(t,j) (2A_D^\top P B_D (R_D+B_D^\top P B_D)^{-1}B_D^\top P A_D+P-A_D^\top PA_D)\phi(t,j)  
  \\&=&\sum_{j=0}^{\infty} V(\phi(t,j))-V(G(\phi(t,j),u_D^*(t,j))
  \\&=&V(\phi(0,0))-\limsup_{j\rightarrow \infty}V(\phi(0,j))
  \\&=& 
  \phi^\top(0,0) P \phi(0,0)- \limsup_{j\rightarrow \infty} \phi^\top(0,j) P \phi(0,j)
\end{eqnarray*}
which, given that $\phi(0,0)=\xi$ and $q(\phi(t,j))=\phi^\top(0,j) P \phi(0,j)$, implies
$	\mathcal{J}(\xi,u_D^*) =\xi^\top P \xi$, showing that the feedback law $\kappa_D$ attains the optimal cost.} 
From Corollary \ref{FBLrem}, the feedback law $\kappa_D$ is a pure strategy Nash equilibrium for the $(\mathcal{LQD})$ game.
\end{example}
\begin{example}{Numerical  Hybrid Game with Nonunique Solutions}
Consider a one-player hybrid game with state $x \in \reals$, input $u=(u_C, u_D) \in \reals \times \reals$, and dynamics $\HS$ as in (\ref{Heq}) described by
\begin{eqnarray}
C&:=&[0,2] \nonumber
\\
F(x,u_C)&:=&a x_1+b u_C \nonumber
\\
D&:=&\{1\} \nonumber
\\
G(x,u_D)&:=&0.5
\label{nuHGame}
\end{eqnarray}
where $a=-1,b=-1$. Consider the stage costs $L_C(x,u_C):=x^2Q_C+u_C^2 R_C $ and $L_D(x,u_D):=P(x^2-0.5^2)$, the terminal cost $q(x):=P x^2$ and the function $V(x):=P x^2$, where $Q_C=1$, $R_C=1.304$, and  $P =0.429$.  The function $V$ is such that
\begin{eqnarray}
 &&\min_{u_C \in \Pi_u(x,C)}
\left(L_C(x,u_C) +\left\langle \nabla V(x),F(x,u_C) \right\rangle \right) 
\nonumber\\&=&
\min_{u_C \in \mathbb{R}}
\left( x^\top (Q+ 2P a) x+ u_C^\top R_C u_C
+  2x^\top P b u_C \right)
 \nonumber \\ &=&
\min_{u_C \in \mathbb{R}}
\left( x^2 0.142+ 1.304u_C^2 
+  0.858x u_C \right)\nonumber \\
&=&0
\label{HJBLQGnusC}
\end{eqnarray}
holds $ \forall x \in \reals^n$ (and attained by $\kappa_C(x)=-R_C^{-1}b P x$),
leading $V(x)= Px^2$ as a solution to (\ref{HJBih}). In addition, the function $V$ is such that
\begin{eqnarray}
&&\min_{u \in \Pi_u(x,D)}
\left\{
 L_D(x,u)
+ V(G(x,u)) \right\} \nonumber
\\&=&
  \min_{u_D \in \mathbb{R}^{m_D}}
\left\{P(x^2-0.5^2)+P0.5^2 \right\} 
\nonumber
\\&=& P x^2,
\label{HJBLQGnusD}
\end{eqnarray}
 which leads $V(x)= Px^2$ as a solution to (\ref{Bellmanih}).
Thus, given that $V$ is continuously on a neighborhood of $\Pi(C)$, and that (\ref{HJBih}) and (\ref{Bellmanih}) hold thanks to (\ref{HJBLQGnusC}) and (\ref{HJBLQGnusD}), from Theorem \ref{thHJBih} the value function is
\begin{eqnarray}
\mathcal{J}^*(\xi) &:=& 
 P \xi^2
\label{HyLCost}
\end{eqnarray}
for any $\xi \in \Pi ( \overline{C}) \cup \Pi(D)$.
Notice that, since $C \cap D \neq \emptyset$, the set of maximal solutions from $x=2$ is given by $\mathcal{S}_\kappa(2)=\{\phi_\kappa,\phi_h\}$, where $\phi_\kappa$ has a domain $ \dom \phi_\kappa=\realsgeq \times \{0\}$, and is given by
\begin{equation*}
  \phi_\kappa(t,0)=  2 e^{-1.3289 t} 
\end{equation*}
In simple words, $\phi_\kappa $ flows from $x=2$ up to $x=0$. The maximal solution $\phi_h$ has a domain $ \dom \phi_h=[0,2.10747)\times \{0\} \cup [2.10747, \infty)\times \{1\}$, and is given by
\begin{equation*}
 \phi_h(t,0)= 2 e^{-1.3289 t} ,\hspace{0.5cm} \phi_h(t,1)=
0.5 e^{-1.3289 (t-2.10747)} 
\end{equation*}
In simple words, the trajectory $\phi_h$ flows from $x=2$ up to $x=1$, then it jumps to $x=0.5$ and flows up to $x=0$.\\
%
By denoting the corresponding input signals as $u_\kappa=\kappa(\phi_\kappa)$ and $u_h=\kappa(\phi_h)$, and given that the cost of flowing from $x=1$ to $x=0.5$ is equal to the cost of jumping, we have that the costs of the solutions $(\phi_\kappa,u_\kappa)$ and $(\phi_h,u_h)$ are equal to $P2^2$. 
%
%
Figure \ref{fig:my_label} illustrates this behavior.
\begin{figure}[H]
    \centering
    \includegraphics[width=15cm, center]{Figures/Counterexample_[211]Fx.png}
    \caption{Nonunique solutions attaining optimal cost. Continuous solution (green). Hybrid solution (blue-red).}
    \label{fig:my_label}
\end{figure}
This aligns to the Nash equilibrium Definition \ref{NashEqCoop}, and Assumption \ref{SolExistence}, since every maximal solution that $\kappa$ renders from $\xi=2$ attains the minimum value of the cost  $\mathcal{J}$ over the set of input actions that yield maximal solutions to $\HS$ with data (\ref{nuHGame}).
\label{NumericalEx}
\end{example}
\IfIncd{\subsection{Linear Quadratic Hybrid Games with Spontaneous Jumps}}
{\subsection{Linear Quadratic Hybrid Games with Periodic Jumps}}
Given a time $\bar{T}\in \mathbb{R}$, consider a one-player hybrid game with state $x=(x_p,\tau)   \in \reals^n \times [0, \bar{T}]$, input $u=(u_C, u_D) \in \reals^{m_C}\times \reals^{m_D}$, and dynamics $\HS$ as in (\ref{Heq}) described by
\begin{eqnarray*}
C&:=&\mathbb{R}^n \times [0,\bar{T}] \times \mathbb{R}^{m_C}
\\
F&:=&(A_C x_p+B_C u_C,1)
\\
D&:=&\mathbb{R}^n \times \{\bar{T}\} \times \mathbb{R}^{m_D}
\\
G&:=&(A_D x_p+B_D u_D,0)
\end{eqnarray*}
where $A_C,A_D \in \mathbb{R}^{n\times n}$ and $B_C \in \mathbb{R}^{n\times m_C}$, $B_D \in \mathbb{R}^{n \times m_D}$.
Consider the stage costs $L_C(x,u_D):=x_p^\top Q_C x_p+ u_C^\top R_C u_C$ and $L_D(x,u_D):=x_p^\top Q_D x_p+ u_D^\top R_D u_D$, the terminal cost $q(x):=x_p^\top P(\tau) x_p$ and the function $V(x):=x_p^\top P(\tau) x_p$, where $Q_C,Q_D\in \mathbb{S}^n_+$, $R_C \in \mathbb{S}^{m_C}_+,R_D\in \mathbb{S}^{m_D}_+$ and  $P(\tau)\in \mathbb{S}^n_+ \forall \tau \in [0,\bar{T}]$. We refer to this game as an infinite-horizon one player linear quadratic game with periodic jumps $(\mathcal{LQH})$.

\begin{corollary}{Hybrid Riccati equation} Given $\bar{T}\in \mathbb{R}$, $A_C,A_D \in \mathbb{R}^{n\times n}$, $B_C \in \mathbb{R}^{n\times m_C}$, $B_D \in \mathbb{R}^{n \times m_D}$, $Q_C,Q_D\in \mathbb{S}^n_+$, $R_C \in \mathbb{S}^{m_C}_+,R_D\in \mathbb{S}^{m_D}_+$,
suppose there exists a matrix $P:[0,\bar{T}]\rightarrow \mathbb{S}^n_+ $ continuously differentiable such that
\begin{equation}
-\frac{\mathrm{d} P(\tau)}{\mathrm{d} \tau} =Q_C+P(\tau) A_C+A_C^\top P(\tau)-P(\tau) B_C R_C^{-1}B_C^\top P(\tau)  \hspace{1cm} \forall \tau  \in (0,\bar{T})
\label{DiffRiccatiih} 
\end{equation}
and 
\begin{equation}
R_D+B_D^\top P B_D \in \mathbb{S}_{0+}^{m_D}
\label{lqeqinv}
\end{equation}
with 
\begin{equation}
P (\bar{T})=Q_D+A_D^\top P(0) A_D-A_D^\top P(0) B_D (R_D+B_D^\top P(0) B_D)^{-1}B_D^\top P(0) A_D
\label{Riccatiih} 
\end{equation}

Then, the feedback law $\kappa:=(\kappa_C,\kappa_D)$, \pnn{defined as}
\begin{equation}
\kappa_C(x)=-R_C^{-1}B_C^\top P(\tau) x_p \hspace{1cm} \forall x \in \Pi(C)
\label{NashkCLQih}
\end{equation}
\begin{equation}
\kappa_D(x)=-(R_D+B_D^\top P(0) B_D)^{-1}B_D^\top P(0) A_D x_p \hspace{1cm} \forall x \in \Pi(D)
\label{NashkDLQih}
\end{equation}
is a Nash equilibrium for the $(\mathcal{LQH})$ game. In addition, given  $\xi=(\xi_p,\xi_\tau) \in \Pi(\overline{C}\cup D)$, the Nash equilibrium outcome is equal to $\xi_p^\top P(\xi_\tau) \xi_p$.
\end{corollary}
\begin{proof}
We show that when conditions (\ref{DiffRiccatiih})-(\ref{Riccatiih}) hold, by using the result in Theorem \ref{thHJBih},  the function $V$ is equal to the value function for this game. By using Corollary \ref{FBLrem}, we show that the feedback law as in (\ref{NashkCLQih}), (\ref{NashkDLQih}) attains such a cost. 
We can write (\ref{HJBih}) in Theorem \ref{thHJBih} as
\begin{eqnarray}
&&0=\min_{u \in \Pi_u(x,C)} \mathcal{L}_C(x,u_C),
 \nonumber
\\ 
&&\mathcal{L}_C(x,u_C)= 
 x_p^\top Q_C x_p+ u_C^\top R_C u_C
+2x_p^\top P(\tau) (A_C x_p+B_C u_C)  
+ x_p^\top \frac{\mathrm{d} P(\tau)}{\mathrm{d} \tau} x_p
\label{HJBihlq}
\end{eqnarray}
First, 
given that (\ref{DiffRiccatiih}) holds, and $x_p^\top(P(\tau) A_C+A_C^\top P(\tau))x_p=x_p^\top (2P(\tau)A)x_p$ for every $x\in C$, one has
\begin{equation}
\mathcal{L}_C(x,u_C)=
 x_p^\top (
 P(\tau) B_C R_C^{-1}B_C^\top P(\tau)) x_p+ u_C^\top R_C u_C
+ 2 x_p^\top  P(\tau) B_C u_C 
\end{equation}
The first order necessary condition for optimality
\begin{equation*}
\frac{\partial}{\partial u_C}\left.\left( 
 x_p^\top (
 P(\tau) B_C R_C^{-1}B_C^\top P(\tau)) x_p+ u_C^\top R_C u_C
+ 2 x_p^\top  P(\tau) B_C u_C  \right)\right|_{u_C^*} =0
\end{equation*}
is satisfied by the point
\begin{equation}
u_C^*=-R_C^{-1}B_C^\top P(\tau) x_p
\label{LQHPu*}
\end{equation}
Given that $R_C \in \mathbb{S}^{m_D}_+$, the second-order necessary condition for optimality 
\begin{equation*}
\frac{\partial^2}{\partial u_C^2}\left.\left( 
 x_p^\top (
  P(\tau) B_C R_C^{-1}B_C^\top P(\tau)) x_p+ u_C^\top R_C u_C
+ 2 x_p^\top  P(\tau) B_C u_C  \right)\right|_{u_C^*} \succeq 0,
\end{equation*}
holds, rendering $u_C^*$ in (\ref{LQHPu*}) as a minimizer of the optimization problem in (\ref{HJBihlq}).
In addition, 
it satisfies $\mathcal{L}_C(x,u_C^*)=0$, rendering $V(x)=x_p^\top P(\tau)x_p$ to be a solution of (\ref{HJBih}) in Theorem \ref{thHJBih}. 

On the other hand, we can write (\ref{Bellmanih}) in Theorem \ref{thHJBih} as
\begin{eqnarray}
&&x_p^\top P(\bar{T}) x_p =\min_{u \in \Pi_u(x,D)} \mathcal{L}_D(x,u_D),
 \nonumber
\\ 
&&\mathcal{L}_D(x,u_D)= 
x_p^\top Q_D x_p+ u_D^\top R_D u_D
+ (A_D x_p+B_D u_D)^\top P(0) (A_D x_p+B_D u_D) 
\label{Bellmanihlq}
\end{eqnarray}
which can be expanded as
\begin{equation}
\mathcal{L}_D(x,u_D)=  
 x_p^\top (Q_D+A_D^\top P(0) A_D) x_p+ u_D^\top (R_D +B_D^\top P(0) B_D) u_D
+ 2 x_p^\top A_D^\top  P(0) B_D u_D 
\end{equation}
The first order necessary condition for optimality
\begin{equation*}
\frac{\partial}{\partial u_D}\left.\left( 
 x_p^\top (Q_D+A_D^\top P(0) A_D) x_p+ u_D^\top (R_D +B_D^\top P(0) B_D) u_D
+ 2 x_p^\top A_D^\top  P(0) B_D u_D  \right)\right|_{u_D^*} =0
\end{equation*}
is satisfied by the point
\begin{equation}
u_D^*=-(R_D+B_D^\top P(0) B_D)^{-1}B_D^\top P(0) A_D x_p
\label{LQHPud*}
\end{equation}
Given that (\ref{lqeqinv}) holds, the second-order necessary condition for optimality 
\begin{equation*}
\frac{\partial^2}{\partial u_D^2}\left.\left( 
 x_p^\top (Q_D+A_D^\top P(0) A_D) x_p+ u_D^\top (R_D +B_D^\top P(0) B_D) u_D
+ 2 x_p^\top A_D^\top  P(0) B_D u_D  \right)\right|_{u_D^*} \succeq 0,
\end{equation*}
is satisfied, rendering $u_D^*$ as in (\ref{LQHPud*}) as a minimizer of the optimization problem in (\ref{Bellmanihlq}).
In addition,  $u_D^*$
satisfies $\mathcal{L}_D(x,u_D^*)=x_p^\top P(\bar{T})x_p$, with $P(\bar{T})$ as in (\ref{Riccatiih}),  rendering $V(x)=x_p^\top P(\tau)x_p$ to be a solution of (\ref{Bellmanih}) in Theorem \ref{thHJBih}. 

Then, from Theorem \ref{thHJBih} and Corollary (\ref{FBLrem}), for every $ x=(x_p,\tau) \in \Pi(\overline{C}\cup D)$ the value function is 
\begin{equation*}
\J^*(x)= x_p^\top P(\tau)x_p
\end{equation*}
and $\kappa=(\kappa_C,\kappa_D)$  as in (\ref{NashkCLQih}), (\ref{NashkDLQih}) is a pure strategy Nash equilibrium.
%
%
%
%
\end{proof}
%
%
\subsection{Lyapunov Functions for Hybrid Games}
Next, we introduce some definitions to present a result that connects optimality and stability properties for hybrid noncooperative games.
\begin{definition}{Class-$\mathcal{K}_\infty$ functions}
	A function $\rho: \mathbb{R}_{\geq 0} \rightarrow  \mathbb{R}_{\geq 0}$ is a class 
	$\mathcal{K}_\infty$ function, also written as $\rho \in \mathcal{K}_\infty$, if $\rho$ is
	 zero at zero, continuous, strictly increasing, and unbounded.
\end{definition}

\begin{definition}{Positive definite functions}
Consider a $\kappa:=(\kappa_C,\kappa_D): \mathbb{R}^n \rightarrow \mathbb{R}^{m_C} \times \mathbb{R}^{m_D}$. 
We say that a function $\rho: \reals_{\geq 0}\rightarrow \reals_{\geq 0}$ is positive definite, also written as $\rho \in \mathcal{PD}$, if $\rho(s)>0$ and $\rho(0)=0$. We say that a function $\rho: \reals^n \times \reals^{m_C} \rightarrow \reals_{\geq 0}$  is positive definite with respect to a set $\A \subset \reals^n$, in composition with $\kappa_C$, also written as $\rho \in \mathcal{PD}_{\kappa_C}(\A)$, if $\rho(s,\kappa_C(s)) > 0$ for all $s \notin \A$ and $\rho(\A, \kappa_C(\A)) = \{0\}$. Likewise, we say that a function $\rho: \reals^n \times \reals^{m_D} \rightarrow \reals_{\geq 0}$  is positive definite with respect to a set $\A \subset \reals^n$, in composition with $\kappa_D$, also written as $\rho \in \mathcal{PD}_{\kappa_D}(\A)$, if $\rho(s,\kappa_D(s)) > 0$ for all $s \notin \A$ and $\rho(\A, \kappa_D(\A)) = \{0\}$.
\end{definition}

\begin{corollary}{Nash equilibrium under the existence of a Lyapunov function} 
Consider a one-player hybrid game with closed-loop dynamics $\HS_{\kappa}$ as in (\ref{Hkeq}) with data $(C,F,D,G)$ and $\kappa:=(\kappa_C,\kappa_D): \mathbb{R}^n \rightarrow \mathbb{R}^{m_C} \times \mathbb{R}^{m_D}$ such that $C_\kappa = \Pi(C), D_\kappa = \Pi(D)$ and \pnn{suppose that Assumption \ref{AssLipsZ} holds.} Given a closed set $\mathcal{A} \subset \reals^n$, continuous stage costs $L_C:
\pnn{C} \rightarrow \mathbb{R}_{\geq 0}, L_D:
\pnn{D} \rightarrow \mathbb{R}_{\geq 0}$, such that
$L_C \in \mathcal{PD}_{\kappa_C}(\mathcal{A}\cap C_\kappa)$, $L_D \in \mathcal{PD}_{\kappa_D}(\mathcal{A}\cap D_\kappa)$, and a terminal cost $q:\reals^n \rightarrow \reals$, suppose there exist a function $V: \mathbb{R}^n \rightarrow \mathbb{R}$ that is continuously differentiable on an open set containing $\Pi(\overline{C})$,  satisfying (\ref{kHJBeqiha}), (\ref{kHJBeqihb}), (\ref{kBeqihc}) and (\ref{kBeqihd}), and such that for each $\xi \in \Pi ( \overline{C}) \cup \Pi(D)$, for each $\phi \in \mathcal{S}_{\HS_{\kappa}  }^\infty(\xi)$, (\ref{TerminalCond}) holds. If there exist $\alpha_1,\alpha_2 \in \mathcal{K}_\infty$ such that for all $x \in \Pi ( \overline{C}) \cup \Pi(D)$
\begin{equation}
\alpha_1(|x|_\A) \leq V(x) \leq \alpha_2(|x|_\A)
\label{alphabound}
\end{equation}
\pnn{and either
\begin{enumerate}
\item $\A \cap C_\kappa \neq \emptyset $ and $\A \cap D_\kappa \neq \emptyset $;
\item $\A \cap C_\kappa \neq \emptyset $, $\A \cap D_\kappa = \emptyset $, and there exists a continuous $\lambda \in \mathcal{PD}$ such that $L_D(x,\kappa_D(x)) \geq \lambda(|x|_\A)$ for all $x \in D_\kappa$;
\item $\A \cap C_\kappa = \emptyset $, $\A \cap D_\kappa \neq \emptyset $, and there exists a continuous $\lambda \in \mathcal{PD}$ such that $L_C(x,\kappa_D(x)) \geq \lambda(|x|_\A)$ for all $x \in C_\kappa$;
\end{enumerate}}
 then,   
\begin{equation}
\J^*(\xi)= V(\xi),
\end{equation}
for each $\xi \in \Pi ( \overline{C}) \cup \Pi(D)$. The feedback law
$\kappa$ is the Nash equilibrium, and it renders $\A$ uniformly globally asymptotically stable for $\HS_{\kappa}$ in $\Pi ( \overline{C}) \cup \Pi(D)$. 
\label{CorStability}
\end{corollary}
\begin{proof}
Given that $C_\kappa = \Pi(C), D_\kappa = \Pi(D)$, and $L_C,L_D,V,\kappa:=(\kappa_C,\kappa_D)$ are such that (\ref{kHJBeqiha}), (\ref{kHJBeqihb}), (\ref{kBeqihc}) and (\ref{kBeqihd}) hold, then, $V$ and $\kappa$ satisfy (\ref{HJBih}), (\ref{Bellmanih}), (\ref{kHJBeqih}), and (\ref{kBeqih}). 
Since (\ref{TerminalCond}) holds for each $\xi \in \Pi ( \overline{C}) \cup \Pi(D)$, for each $\phi \in \mathcal{S}_{\HS_{\kappa}  }^\infty(\xi)$, we have from Theorem \ref{thHJBih} that $V$ is the value function as in (\ref{cost2goih}) for $\HS_{\kappa}$ at $\Pi ( \overline{C}) \cup \Pi(D)$. 
From Corollary \ref{FBLrem} we have that the feedback law
 $\kappa:=(\kappa_C,\kappa_D)$ that satisfies (\ref{kHJBeqih}), (\ref{kBeqih}) is the pure strategy Nash equilibrium for this game. 
In addition, given that $V$ is a Lyapunov candidate for $\HS_{\kappa}$ \cite[Definition 3.16]{65nu}, 	the functions $\alpha_1, \alpha_2$ satisfy (\ref{alphabound}), and given that (\ref{kHJBeqihb}), (\ref{kBeqihd}) hold, then, from \cite[Theorem 3.18]{65nu} we have that $\A$ is uniformly globally asymptotically stable for $\HS_{\kappa}$. \pnn{Notice that the role of $\rho$ therein is played by the stage costs $L_C \in \mathcal{PD}_{\kappa_C}(\A \cap C_\kappa)$, $L_D \in \mathcal{PD}_{\kappa_D}(\A \cap D_\kappa)$ or the function $\lambda \in \mathcal{PD}$.}
 \end{proof}

\begin{example}{Linear Quadratic Differential Games}
Taking up the $(\mathcal{LQHC})$ game, recall that $\kappa:=(-R_C^{-1} B_C^\top P x,\cdot)$ and consider the case in which $(A_C-B_C R_C^{-1}B_C^\top P)$ is Hurwitz and $\A=\{0\}$. Given that $L_C \in \mathcal{PD}_{\kappa_C}(\A \cap C_\kappa)$, $V$ is continuously differentiable, and (\ref{kHJBeqiha}) and (\ref{kHJBeqihb}) hold, by setting, $\alpha_1(|x|_\A)=x^\top (P-I) x$ and $\alpha_2(|x|_\A)=x^\top (P+I) x$, from Corollary \ref{CorStability}, we have that $\kappa$ is the Nash equilibrium and renders $\A$ uniformly globally asymptotically stable for $\HS_{\kappa}$.
\end{example}
\begin{example}{Linear Quadratic Difference Games}
Taking up the $(\mathcal{LQHD})$ game, recall that $\kappa:=(\cdot,-(R_D+B_D^\top PB_D)^{-1}B_D^\top P A_D x)$ and consider the case in which $(A_D-B_D (R_D+B_D^\top PB_D)^{-1}B_D^\top P A_D )$ is Schur and $\A=\{0\}$. Given that $L_D \in \mathcal{PD}_{\kappa_D}(\A \cap D_\kappa)$, and  (\ref{kBeqihc}), and (\ref{kBeqihd}) hold, by setting $\alpha_1(|x|_\A)=x^\top (P-I) x$ and $\alpha_2(|x|_\A)=x^\top (P+I) x$, from Corollary \ref{CorStability}, we have that $\kappa$ is the Nash equilibrium and renders $\A$ uniformly globally asymptotically stable for $\HS_{\kappa}$.
\end{example}
\begin{example}{Numerical Hybrid Game with Nonunique Solutions}
Taking up Example \ref{NumericalEx}, recall that $a<0$ and $\kappa:=(\kappa_C,\kappa_D)=(-R_C^{-1}bPx,\cdot)$ and let $\A=\{0\}$. Given that $L_C \in \mathcal{PD}_{\kappa_C }(\A \cap C_\kappa)$, \pnn{$L_D \in \mathcal{PD}_{\kappa_D}(\A \cap D_\kappa)$ vacuously}, and (\ref{kHJBeqiha}),  (\ref{kHJBeqihb}), (\ref{kBeqihc}), and (\ref{kBeqihd}) hold, by setting $\alpha_1(|x|_\A)=x^\top (P-I) x$ and $\alpha_2(|x|_\A)=x^\top (P+I) x$, 
and \pnn{for $|x|_\A \mapsto \lambda(|x|_\A)=P \frac{x^2}{2}$, $L_D(x,\kappa_D(x)) \geq \lambda(|x|_\A)$ for all $x \in D_\kappa$,}  
from Corollary \ref{CorStability}, we have that $\kappa$ is the Nash equilibrium and renders $\A$ uniformly globally asymptotically stable for $\HS$ with data as in (\ref{nuHGame}). 
\end{example}}
%
%
\IfPers{Let $\mathcal{A} \subset \mathbb{R}^n$ be a closed set and consider
\begin{equation*}
	\mathcal{X}_\mathcal{A}:=\{ \phi \in \mathcal{X}: \displaystyle \lim \limits_{\substack{(t,j) \in \dom \phi\\
	(t,j)\rightarrow \sup \dom \phi}} \left | 
	\phi(t,j) \right|_\mathcal{A}=0 \}
\end{equation*}
\begin{equation*}
\mathcal{U}_\mathcal{A}(\xi) :=  \left\{ u\in \mathcal{U}: \exists  \phi \in \mathcal{R} (\xi,u) \cap \mathcal{X}_{\mathcal{A}} \right\}
\end{equation*}
To move to another section
}
%
\subsection{Problem Statement}\label{subsec: probstat}
In this section, we formulate an 
optimization problem to solve the two-player zero-sum hybrid game and provide sufficient conditions to 
characterize the solution.
\pno{\pnn{Compared to \cite{sjlAutom26}, we relax the assumption on uniqueness of solutions to the hybrid system defining the constraints of the game.
Although in this case it might not be possible to construct a saddle-point equilibrium
, we optimize the worst-case (due to nonuniqueness of solutions) value of the cost functional and refer to it as a weak solution.
Specifically,}} {following Definition \ref{elements}, consider a two-player zero-sum hybrid game with dynamics $\HS\IfIncd{}{}$ as in (\IfIncd{\ref{Hinc}}{\ref{Heq}}) \IfTp{}{with $N=2$ }
for given $(C,F,D,G)$. }
\IfIncd{}{Uniqueness of solutions for a given input implies a unique correspondence {from cost} to control {action}, which allows this type of games to be \textit{well-defined}, so 
{that an} equilibrium solution is defined \cite[Remark 5.3]{basar1999dynamic}.
This justifies the following assumption.
\begin{assumption}{}
  The flow map $F$ and the flow set $C$ are such that solutions to {$\dot{x}=F(x,u_C)$ $(x,u_C)\in C$} are unique for each input $u_C$. The jump map $G$ is single valued{, i.e., $D_1=D_2$}.
 \label{AssLipsZ}
\end{assumption}
Sufficient conditions to guarantee that Assumption \ref{AssLipsZ} holds include Lipschitz continuity of the flow map $F$, provided it is a single-valued function.
Under Assumption \ref{AssLipsZ}, the conditions in Proposition \ref{UniquenessHu} are satisfied, so solutions to $\HS$ are unique\footnote{{Under Assumption \ref{AssLipsZ}, the domain of the input $u$ specifies whether at $C\cap D$ there is a jump or flow.}} for each $u\in \mathcal{U}$.

}
Given $\xi \in C \cup D$, a joint input action $u=(u_C, u_D)\in \mathcal{U}$, 
the stage cost for flows $L_C:\mathbb{R}^n  \times \mathbb{R}^{m_C} \rightarrow \mathbb{R}_{\geq 0}$, the stage cost for jumps $L_D:\mathbb{R}^n  \times \mathbb{R}^{m_D} \rightarrow \mathbb{R}_{\geq 0}$, and the terminal cost $q: \mathbb{R}^n \rightarrow \mathbb{R}$, we conveniently define the cost {associated to the {worst-case} solution to $\HS\IfIncd{}{}$ from the initial condition $\xi$} and for the hybrid input $u$, \IfIncd{}{under Assumption \ref{AssLipsZ},}
as
\IfIncd{\begin{equation}
\label{defJTNCincinc}
\mathcal{J}(\xi,u) := 
\sup_{\phi \in \R (\xi,u)}
\widetilde{\mathcal{J}}(\phi,u)
\end{equation}
where\footnote{Notice that $\mathcal{J}$ depends on the initial condition $\xi$ and input $u$, while $\widetilde{\mathcal{J}}$ depends on the solution pair $(\phi,u)$ {provided, which in \eqref{defJTNCincinc} is such that} $\phi(0,0)=\xi$.}}{}
\begin{equation}\label{defJTNCinc}
  \begin{split}
  &\IfIncd{\widetilde{\mathcal{J}}(\phi,u)}{\mathcal{J}(\xi,u)}:=
  \sum_{j=0}^{\sup_j \dom \phi} \int_{t_{j}}^{t_{j+1}} L_C(\phi(t,j),u_{C}(t,j))dt  
  \\& \hspace{1.5cm}
  + 
  \sum_{j=0}^{\sup_j \dom \phi -1} 
  L_D(\phi(t_{j+1},j),u_{D}(t_{j+1},j))
  +\underset{(t,j) \in \textup{dom}\phi}{{\limsup_{\pnn{(t,j)\rightarrow \sup \dom \phi}} }} q(\phi(t,j)),
\end{split}
\end{equation}
 $\{t_j\}_{j=0}^{\sup_j \dom \phi}$ is a nondecreasing sequence
associated to the definition of the hybrid time domain of {$(\phi,u)$} -- see Definition \ref{htd} -- and $\R (\xi,u)$ is the set of maximal state trajectories to $\HS\IfIncd{}{}$ for the hybrid input $u$ from the initial condition $\xi$, as defined in Section \ref{Subsect HySyswIn}. 
{The last term in \eqref{defJTNCinc} penalizes the value of
the solution $(\phi,u)$ at the terminal time of 
$\dom \phi$ when $\dom \phi$ is bounded, or its limit when unbounded.
}

A solution to the two-player zero-sum game can be obtained by solving the following problem:

\textit{Problem ($\diamond$):} 
Given $\xi \in \mathbb{R}^n$, \IfIncd{}{under Assumption \ref{AssLipsZ},} 
solve
\begin{eqnarray}
\underset{u=(u_1,u_2) \in \mathcal{U}_{\HS\IfIncd{}{}}(\xi)}
{ \underset{u_{1}}{\textup{minimize}}\>\> \underset{u_{2}}{\textup{maximize}}}
&& \mathcal{J}(\xi, u)
\label{problemzsih}
\end{eqnarray}
where $\mathcal{U}_{\HS\IfIncd{}{}}$ is the set of joint input actions yielding maximal solutions to $\HS\IfIncd{}{}$, as defined in Section 2.1.
Compared to \cite{leudohygamesnu,sjlAutom26}, Problem $(\diamond)$ {considers} the cost associated to the worst-case solution from $\xi$ for $u$ as defined in \eqref{defJTNCincinc}.
  \NotACC{
\begin{remark}{\IfIh{}{Pure strategy }Feedback \pnn{weak} saddle-point equilibrium and \pnn{weak} min-max control}
A solution to Problem ($\diamond$), when it exists, can be expressed in terms of the \pnn{weak} saddle-point equilibrium strategy $\IfIh{\kappa=(\kappa_1,\kappa_2)}{\gamma}$ for the two-player zero-sum game, {as in Definition \ref{NashEqNonCoop}.} Each $u^*=(u^*_1, u^*_2) \in \mathcal{U}$ that renders a \IfIncd{worst-case cost}{} 
state trajectory $\phi^*\in \mathcal{R}(\xi,u^*)$, 
with components defined as $\dom \phi^* \ni (t,j) \mapsto u^*_i(t,j)= \IfIh{\kappa}{\gamma}_i(\IfIh{}{t,j,}\phi^*(t,j))$  for each $i\in \IfTp{\{1,2\}}{\mathcal{V}}$, satisfies
\begin{eqnarray*}
u^*\pno{\in}
\underset{u=(u_1,u_2) \in \mathcal{U}_{\HS\IfIncd{}{}}(\xi)}
{\arg \min_{u_{1}} \max_{u_{2}}}
 \mathcal{J}(\xi, u)
=
 \underset{u=(u_1,u_2) \in \mathcal{U}_{\HS\IfIncd{}{}}(\xi)}
{\arg  \max_{u_{2}}\min_{u_{1}}}
 \mathcal{J}(\xi, u),
\label{Nasheq1pzsih}
\end{eqnarray*}
and it is referred to as a \pnn{weak} min-max control at $\xi$.
{The set of maximal state trajectories $\R (\xi,u^*)$ is defined as in Section \ref{Subsect HySyswIn}.}
\label{Rmk: minmax control}
\end{remark}

\begin{definition}[Value function] Given $\xi \in \Pi(\overline{C}) \cup \Pi(D)$, 
the value function at $\xi$, {when it exists,} is given by 
\begin{equation}
  \begin{split}
\mathcal{J}^*(\xi):=
\underset{u=(u_1,u_2) \in \mathcal{U}_{\HS\IfIncd{}{}}(\xi)}
{\min_{u_{1}} \max_{u_{2}}}
\mathcal{J}(\xi,u)
=\underset{u=(u_1,u_2) \in \mathcal{U}_{\HS\IfIncd{}{}}(\xi)}
{\max_{u_{2}} \min_{u_{1}} }
\mathcal{J}(\xi,u)
  \end{split}
 \label{cost2gozsih}
\end{equation}
\label{ValueFunctionz}
\end{definition} 
}
%
%
%
%
As a preliminary step, following the approach in \cite{ferrante2018optimality}, we present {in the next section} a framework for cost evaluation for {hybrid systems with set-valued flow and jump maps} as in \eqref{Hinc}, under the special case of no inputs present.
\section{Cost Evaluation for Autonomous Hybrid Inclusions}
\subsection{Upper Bounds}
\label{sec:Cost_evaluation}
By following the ideas proposed in \cite{bernstein1993nonquadratic},  in this section we investigate how a Lyapunov-like function can be used to provide estimates of nonlinear cost functionals for a given hybrid inclusion {potentially exhibiting nonunique solutions}. 

\pnn{We define the hybrid system $\HS_0$ as in \eqref{Hinc} with no inputs, namely, with $u_C=0$ and $u_D=0$.}
For each initial condition $\xi\in\overline{C}\cup D$ to $\HS_0$, 
 the {worst-case cost associated to solutions to $\HS_0$ from $\xi$, defined} in \eqref{defJTNCincinc} and \eqref{defJTNCinc}, can be expressed as
\begin{equation}
  \label{defJTNCincincNI}
  \mathcal{J}(\xi) := 
  \sup_{\phi \in \mathcal{S}_\HS (\xi)}
  \widetilde{\mathcal{J}}(\phi)
  \end{equation}
  where
  \begin{equation}\label{defJTNCincNI}
    \begin{split}
    &\widetilde{\mathcal{J}}(\phi):=
    \sum_{j=0}^{\sup_j \dom \phi} \int_{t_{j}}^{t_{j+1}} L_C(\phi(t,j))dt  
    + 
    \sum_{j=0}^{\sup_j \dom \phi -1} 
    L_D(\phi(t_{j+1},j))
    \\& \hspace{5.7cm} 
       +\underset{(t,j) \in \textup{dom}\phi}
       {{\limsup_{\pnn{(t,j)\rightarrow \sup \dom \phi }} }} 
       q(\phi(t,j)),
  \end{split}
  \end{equation}
\IfPers{\begin{equation}
\begin{aligned}
\label{eq:cost_J}
\mathcal{J}(\xi)=&\sup_{\phi\in\mathcal{S}_{\HS} (\xi)}\left(\int_{\dom_t \phi} L_C(\phi(s,j(s)))ds\right.+\\&\left.\sum^{\sup_j\dom\phi}_{j=1}L_D(\phi(t(j), j-1))\right)
\end{aligned}
\end{equation}
where $L_C\colon C\rightarrow\R_{\geq 0}$ and $L_D\colon D\cup G(D)\rightarrow\R_{\geq 0}$.
\begin{remark}{}
Given $\phi\in\mathcal{S}_{\mathcal{H}_0}$, the definition of the cost \eqref{eq:cost_J} implies that, when $\dom_j\phi$ is bounded, the value of $\phi(t_J, \sup_j\dom\phi)$, 
with $t_J=\inf \{t\in\R_{\geq 0}\colon (t, \sup_j\dom\phi)\in\dom\phi\}$, does not contribute to the ``jump cost''. 
Obviously, this does not lead to any difference in the cost value when $\dom\phi$ is unbounded in the $j$-direction. As it will be clear in the following, such a formulation turns out to be convenient for our analysis.
\end{remark}
}
where {the second argument of the stage costs is dropped \pnn{since we consider the case of no inputs}.}
\IfPers{Throughout the paper, given a solution $\phi$ to \eqref{eq:HybridPlant}, we denote
$$
\begin{aligned}
\mathcal{J}_\phi\coloneqq &\lim_{\substack{(\tau,\iota)\in\dom\phi\\(\tau,\iota)\rightarrow\sup\dom\phi}}
\left(\int_{[0,\tau]} L_C(\phi(s,j(s)))ds\right.+\\
&\left.\sum_{j=1}^\iota L_D(\phi(t(j),j-1))\right)
\end{aligned}
$$
}
{We introduce the following definitions to establish an upper bound on the cost.
\begin{definition}[Generalized gradient]\label{def:GenGrad}
  The change of a locally Lipschitz function $V$ along flows is given by its Clarke generalized directional derivative, which, at a point $x$ in the direction of $\chi$ is given by
  \begin{equation}
 V^\circ (x,\chi) := 
\max_{\nu \in \partial V(x)} \langle \nu, \chi \rangle 
  \end{equation}
  where  $\partial V(x)$ is the Clarke
  generalized gradient of $V$ at $x$, defined as\pnn{\footnote{\pnn{Given a set of points $x_i$ with $i \in S:=\{1,2,\dots, N\}$, $N \in \nats_{\geq 0}$, their convex hull is defined as 
   $ \textup{co}\Big\{\bigcup_{i
    } 
    x_i
    \Big\} :=
    \Big\{\sum_{i
    } \lambda_i 
    x_i
     : \lambda_i \geq 0 \textup{ for all }i \in S\textup{ and } \sum_{i
    } \lambda_i =1  
      \Big\}
  $
  }}}
  \begin{equation}
    \partial V(x) := \textup{co}\Big\{ \lim_{i \rightarrow \infty} \nabla V(x_i): x_i \rightarrow x, \>\frac{d V}{dx} (x_i) \textup{  exists} 
      \Big\}
  \end{equation}
\end{definition}
Notice that\pnn{, when $F$ has compact values,} this generalized derivative leads to the infinitesimal quantity 
\begin{equation}\label{defVcirc}
\max_{\chi \in F(x)} V^\circ (x,\chi) = \max_{\chi \in F(x)}
\max_{\nu \in \partial V(x)} \langle \nu, \chi \rangle 
\end{equation}
which, on $C$, bounds the change of $V$ along flows of solutions to $\HS_0$.
When $C$ and $F$ satisfy $1)$ and $2)$ in the hybrid basic conditions \pnn{(see \IfSiam{\cite[Definition A.5]{TRHSVGArxiv}}{Definition \ref{def:hbc}}) and $F$ has compact values}, then, \pnn{for each solution $\phi$ to $\HS_0$,}
\begin{equation}
\frac{d}{dt} V(\phi(t,j)) \leq \dot{V}(\phi(t,j)) := 
  \pnn{
  \max_{\chi \in F(\phi(t,j)) \cap T_C(\phi(t,j))} V^\circ (\phi(t,j),\chi) 
  }
  \label{dotV}
  \end{equation}
 \pnn{for every $(t,j) \in \textup{int}(I_\phi^j) \times \{j\}$ with int$(I_\phi^j)$ nonempty.}
See \cite[Section 3.3.1]{220} and \IfSiam{\cite[Definition A.6]{TRHSVGArxiv}}{Definition \ref{def:TanCon}} for details. 
}

The following \pnn{proposition, by building on a suitable function $V$, provides an upper bound on the cost $\widetilde{\mathcal{J}}(\phi)$ that depends on the solution \pnn{$\phi$} chosen from $\xi$.}

\begin{proposition}[Upper bound for a given trajectory]
\label{prop:Prop1}
{Given a system $\HS_0 = (C,F,D,G)$
,} let $\xi\in \pno{\overline C}\cup D$, $L_C\colon \pno{\overline C}\rightarrow\mathbb R_{\geq 0}$, $L_D\colon D\cup G(D)\rightarrow\mathbb R_{\geq 0}$, and $q: \pno{\overline C}\cup D \cup G(D) \to \reals$. Let $V\colon\dom V\rightarrow\mathbb R$ with $\dom V\supset \pno{\overline C}\cup D\cup G(D)$ be \IfGG{locally Lipschitz}{continuously differentiable} on \NotGG{an open set containing} $\overline{C}$. {Suppose} that 
\begin{subequations} 
\label{eq:LyapCostHybr}
\begin{eqnarray}
\label{eq:LyapCostHybr_flow}
&\IfGG{\dot{V}(x)}{\displaystyle{\sup_{f\in F(x)}}\langle\nabla V(x),f\rangle}+L_C(x)\leq 0&\qquad \forall x\in C
\\
\label{eq:LyapCostHybr_jump}
&\displaystyle{\sup_{\chi\in G(x)}}V(\chi)-V(x)+L_D(x)\leq 0&\qquad \forall x\in D
\end{eqnarray}
\end{subequations}
\sj{where $\dot{V}$ is defined in \eqref{dotV}}.
\NotGG{where \pno{the second argument of the flow and jump maps is dropped for simplicity}.} 
Let $\phi\colon\dom\phi\rightarrow
\mathbb{R}^n$ be a solution to {$\HS_0$}
from $\xi$. 
{Suppose} that $(t,j)\mapsto V\circ\phi(t,j)$ is bounded
and
\begin{equation}
\underset{(t,j) \in \textup{dom}\phi}
{{\limsup_{\pnn{(t,j)\rightarrow \sup \dom \phi }} }} 
V(\phi(t,j))= \underset{(t,j) \in \textup{dom}\phi}
{{\limsup_{\pnn{(t,j)\rightarrow \sup \dom \phi }} }} 
q(\phi(t,j)).
\label{eq:TerminalCondCGinc0}
\end{equation}
Then, $
\widetilde{\mathcal{J}}(\phi)$ is a finite {number, and in particular,} 

\begin{equation}
\label{eq:item1_prop1}
\widetilde{\mathcal{J}}(\phi)
\leq V(\xi)
\end{equation} 
\end{proposition}
\medskip

\begin{proof}
  {
     {Given that $\phi$ is a solution to $\HS_0$} 
   from $\xi \in \overline{C} \cup D
    $, 
  where $\{t_j\}_{j=0}^{\sup_j \dom \phi}$ is a nondecreasing sequence
  associated to the definition of the hybrid time domain of $\phi$ as in  Definition \ref{htd},
  observe that, \pnn{for each $(T,J) \in \dom \phi$,}
  \begin{equation} \label{eq:Vattj}
  \begin{split}
      V(\phi&(T,J))-V(\xi)
      = 
      \\&
      \sum_{j=0}^{J} \int_{t_j}^{t_{j+1}}
      \frac{d}{dt} V(\phi(t,j))dt
      +
      \sum_{j=0}^{J-1} \Big(
      V(\phi(t_{j+1},j+1))-V(\phi(t_{j+1},j))\Big)
      \end{split}
  \end{equation}
  where \pnn{$t_{J+1}:=T$}, and $V \circ \phi$ is locally Lipschitz on every $I_\phi^j=[t_j,t_{j+1}]$, $j\in \nats$ with nonempty interior.  
  In particular, for each \pnn{$j\in \nats$ such that $j\leq J$ and for almost all $t\in 
  I_\phi^j $},
  \begin{equation}\label{eq:bound_LCsup}
    \frac{d}{dt}V(\phi(t,j)) \leq \IfGG{\dot V (\phi(t,j))}{\sup_{f \in F(\phi(t,j))}\left\langle\nabla V(\phi(t,j)), f \right\rangle }
 \end{equation}  
 \IfGG{where $\dot V$ is defined as in \eqref{dotV}.}{}
  Moreover, \eqref{eq:LyapCostHybr_flow} \IfGG{and \eqref{eq:bound_LCsup} imply}{implies} that for each \pnn{$j\in \nats$ such that $j\leq J$ and for almost all $t\in 
  I_\phi^j $},
  \begin{equation}\label{eq:bound_LC}
    \frac{d}{dt}V(\phi(t,j))
    \leq
    - L_C(\phi(t,j)).
  \end{equation}
  Similarly, for every $(t,j) \in \dom \phi$ such that $(t,j+1) \in \dom \phi$ \pnn{and $j \leq J-1$}, \eqref{eq:LyapCostHybr_jump} implies that
    \begin{equation}\label{eq:bound_LD}
        V(\phi(t,j+1))- V(\phi(t,j))
        \leq
        -L_D(\phi(t,j)).
  \end{equation}
Then, by combining \eqref{eq:Vattj}, \eqref{eq:bound_LC}, and \eqref{eq:bound_LD}, we obtain
  \begin{equation} \label{eq:VattjL}
  \begin{split}
      &V(\phi(T,J))-V(\xi)
      \leq
      -
      \sum_{j=0}^{J} \int_{t_j}^{t_{j+1}} L_C(\phi(t,j))dt
      -\sum_{j=0}^{J-1} L_D(\phi(t_{j+1},j))
      \end{split}
  \end{equation}
  By taking the limit when  
 \pnn{$(T,J)$ approaches $\sup \dom \phi$}, 
  thanks to \eqref{eq:TerminalCondCGinc0}, \eqref{eq:VattjL} implies 
  \begin{equation*}
    \begin{split}
      &\widetilde{\mathcal{J}}(\phi)  = 
      \sum_{j=0}^{\sup_j \dom \phi} \int_{t_{j}}^{t_{j+1}} L_C(\phi(t,j))dt  
  + 
  \sum_{j=0}^{\sup_j \dom \phi -1} 
  L_D(\phi(t_{j+1},j))
  \\& \hspace{5.6cm}
   +\underset{(t,j) \in \textup{dom}\phi}{\limsup_
   {\pnn{(t,j)\rightarrow \sup \dom \phi}} 
   } q(\phi(t,j))
      \leq  V(\xi)
    \end{split}
  \end{equation*}}
which gives \eqref{eq:item1_prop1}, concluding the proof.
\end{proof}

\begin{remark}{Upper bound for an initial condition}
  Notice that since Proposition~\ref{prop:Prop1} provides an upper bound for the cost of any solution to $\HS$ from $\xi$, it follows that under the conditions stated therein, we have 
  \begin{equation}
  \sup_{\phi \in \mathcal{S}_\HS (\xi)}
  \widetilde{\mathcal{J}}(\phi)
\leq V(\xi)
    \label{eq:item1_rmk}
    \end{equation} 
  which, {using} \eqref{defJTNCincincNI}, implies
  \begin{equation}
    \mathcal{J}(\xi)\leq V(\xi).
    \label{eq:item2_rmk}
    \end{equation} 
\end{remark}

\IfPers{Next, by relying on further assumptions, for a given initial condition $\xi\in \overline{C}\cup D$, we provide an upper bound on the cost $\mathcal{J}(\xi)$ that is solution independent.
\begin{corollary}{Upper bound for an initial condition}
\label{coro:CorollaryBound}
Let $\A\subset\mathbb R^n$ be closed, $\xi\in \CUD$, $L_C\colon C\rightarrow\mathbb R_{\geq 0}$, $L_D\colon D\cup G(D)\rightarrow\mathbb R_{\geq 0}$, and $q:C\cup D \cup G(D) \to \reals$. Suppose there exists a function $V\colon\dom V\rightarrow\reals$, $\dom V\supset C\cup D\cup G(D)$, that is continuously differentiable on an open set containing $\scalebox{0.89}{$\overline{C}$}$ and uniformly continuous on a neighborhood of $\A$ such that $V(\A\cap\dom V)=\{0\}$ and \eqref{eq:LyapCostHybr} holds. 
\IfPers{Furthermore, assume that each $\phi\in\mathcal{S}_{\mathcal{H}_0}(\xi)$ is such that
\begin{equation}
\lim_{\substack{(t, j)\in\dom\phi\\(t, j)\rightarrow\sup\dom\phi}}\vert\phi(t, j)\vert_\mathcal{A}=0
\label{eq:limit1}
\end{equation}
}
Then 
\begin{equation}
\mathcal{J}(\xi)\leq V(\xi)
\label{eq:item1_coro1}
\end{equation} 
\end{corollary}
\medskip
\begin{proof}
Pick any $\phi\in\mathcal{S}_{\HS} (\xi)$. Then, since $V$ is uniformly continuous on neighborhood of $\A$ and $V(\A\cap\dom V)=\{0\}$, thanks to Lemma~\ref{lemm:tech:2} in the Appendix, \eqref{eq:limit1} yields
$$
\lim_{\substack{(t, j)\in\dom\phi\\(t, j)\rightarrow\sup\dom\phi}}V(\phi(t, j))=0
$$ 
implying that $\dom\phi\ni (t,j)\mapsto V(\phi(t, j))$ is bounded. Therefore, thanks to Proposition~\ref{prop:Prop1} and in particular from \eqref{eq:item1_prop1}, it follows that $\mathcal{J}_{\phi}$ is finite and in particular  
$$
\mathcal{J}_{\phi}\leq V(\xi)\qquad \forall \phi\in\mathcal{S}_{\HS} (\xi)
$$
which in turn yields \eqref{eq:item1_coro1}, concluding the proof.
\end{proof}

\begin{remark}{}
To get a solution independent upper bound on the cost, in the above result we assumed $V$ to be uniformly continuous on a neighborhood of $\A$. Indeed, since $V(\dom V\cap \A)=\{0\}$, uniform continuity on a neighborhood of $\dom V\cap \A$, as shown in Lemma~\ref{lemm:tech:2} in the Appendix, ensures that for any $(t,j)\mapsto\phi(t, j)$ such that $\phi$ approaches $\A$, $V\circ\phi$ approaches zero.
Alternatively, building upon Proposition~\ref{prop:Prop1},  to relax the uniform continuity requirement of $V$, one could assume that for each $\phi\in\mathcal{S}_{\mathcal{H}_0}(\xi)$ such that $\vert\phi\vert_{\A}$ approaches zero, one has 
$$
\limsup_{\substack{(t, j)\in\dom\phi\\(t,j)\rightarrow\sup\dom\phi}}V(\phi(t,j))=0
$$
 On the other hand, observe that when $\A$ is compact (which is often the case in applications), by the \emph{Heine-Cantor Theorem}, continuity of $V$ on a neighborhood of $\A$ is enough.  
\end{remark}
\begin{remark}{}
Corollary~\ref{coro:CorollaryBound} shows that when maximal solutions from $\xi$ converge to $\A$, then an upper bound on the cost $\mathcal{J}(\xi)$ (which is solution independent) is given by $V(\xi)$. On the one hand, when $L_D, L_C$, and $V$ are positive definite with respect to $\A$, \eqref{eq:LyapCostHybr} implies for any complete solution $\phi$, that $V\circ\phi$ approaches zero.
On the other hand, for maximal solutions that are not complete, finite time convergence to $\A$ is needed. Conditions for finite-time convergence for hybrid systems are given in \cite[Theorem 3.9]{li2016results}.
\end{remark}
}
\subsection{Exact Cost Evaluation}
\label{sec:Exact}
In this section, our main objective is to obtain the exact value of the cost $\mathcal{J}$ in \eqref{defJTNCincincNI}
, without explicitly computing it. 
{The following result asserts that the exact value of the cost $\mathcal{J}(\xi)$ for a given initial condition $\xi$
is $V(\xi)$. 
It requires further assumptions on the system data and {equality in}
\eqref{eq:LyapCostHybr}.}
\begin{corollary}[Exact cost evaluation]
\label{coro:SingleValued}
{Given a system $\HS_0 = (C,F,D,G)$,
} let 
$\xi\in\CUD$, $L_C\colon \pno{ \overline C}\rightarrow \mathbb R_{\geq 0}$, $L_D\colon D\cup G(D)\rightarrow\reals_{\geq 0}$, $q:\pno{ \overline C} \cup D \cup G(D) \to \reals$, and 
\pnn{$F$ and $G$ have compact values in $C$ and $D$, respectively.}
Suppose that there exists a continuous function $V\colon\dom V\rightarrow\reals$, $\dom V\supset \pno{ \overline C}\cup D\cup G(D)$, that is \IfGG{locally Lipschitz}{continuously differentiable} on \NotGG{an open set containing }$\overline{C}$ such that
\begin{subequations} 
\label{eq:LyapCostHybr_2}
\begin{eqnarray}
\label{eq:LyapCostHybr_flow_2}
&\IfGG{\dot{V}(x)}{\displaystyle{\max_{f\in F(x)}}\langle\nabla V(x),f\rangle}+L_C(x)=0&\qquad \forall x\in C
\\
\label{eq:LyapCostHybr_jump_2}
&\displaystyle{\max_{\chi\in G(x)}}V(\chi)-V(x)+L_D(x)=0&\qquad \forall x\in D
\end{eqnarray}
\end{subequations} 
\sj{where $\dot{V}$ is defined in \eqref{dotV}}. Furthermore, {suppose} that for any solution $\phi$ to {$\HS_0$} 
 from $\xi$, $V\circ\phi$ is bounded and \eqref{eq:TerminalCondCGinc0} holds.
Pick any solution $\phi^*$ to the ``maximal'' hybrid system\footnote{{Solutions to the ``maximal system'' in \eqref{eq:argmax} exist under \IfACC{compactness of the set-valued maps, regularity of $V$ and a proper selection of the initial condition.}{the conditions in \IfSiam{\cite[Proposition A.8]{TRHSVGArxiv}}{Proposition \ref{Prop:ExistHmax}}.}}} {defined as}
\pno{
\begin{equation}
\begin{array}{lll}
\dot{x}\phantom{^+}\in 
\IfGG{\underset{\chi \in F(x) \cap T_C(x)}{\argmax} V^\circ (x,\chi)}{\underset{f\in F(x)}{\argmax}\,\langle\nabla V(x),f \rangle}&\quad x\in C\\
x^+\in\underset{\chi\in G(x)}{\argmax}\,V(\chi)&\quad x\in D
\end{array}
\label{eq:argmax}
\end{equation}
}
with $\phi^*{(0,0)}=\xi$, \sj{where $V^\circ$ is defined in \eqref{defVcirc}.} \pn{Let $\phi$ be any solution to $\HS_0$ 
from $\xi$}. Then, \pno{$\mathcal{J}({\phi})$ and $\mathcal{J}({\phi^*})$ are finite, and}
\begin{equation}
\label{eq:item0_coro2}
\begin{aligned}
&\widetilde{\mathcal{J}}({\phi})
\leq 
\widetilde{\mathcal{J}}({\phi^*})
=\mathcal{J}(\xi)
=V(\xi)
\end{aligned}
 \end{equation}
{where $\J$ and $\widetilde{\mathcal{J}}$ are defined in \eqref{defJTNCincincNI} and \eqref{defJTNCincNI}, respectively.}
\end{corollary}
\begin{proof}
{Given that $\phi^*$ is a solution} to \eqref{eq:argmax}, observe that since for each 
$x\in C$
$$
\IfGG{\underset{\chi \in F(x) \cap T_C(x)}{\argmax} V^\circ (x,\chi)}{\underset{f\in F(x)}{\argmax}\langle \nabla V(x), f\rangle}
\subset F(x)
$$
and for each $x\in D$
$$
\underset{\chi\in G(x)}{\argmax}V(\chi)\subset G(x),
$$
$\phi^*$ is a solution to {$\HS_0$
} as well. Let $\{t_j\}_{j=0}^{\sup_j \dom \phi^*}$ be a nondecreasing sequence
associated to the definition of the hybrid time domain of $\phi^*$ as in  Definition \ref{htd}.  
Moreover, for each $j \in \nats$ such that $I_{\phi^*}^j=[t_j,t_{j+1}]$ has a nonempty interior, we have, \pno{for almost all $t \in 
I_{\phi^*}^j$},
\begin{subequations}
  \label{eq:hjbe_c}
\begin{equation}
  \displaystyle{\frac{d}{dt}V(\phi^*(t,j))}=
  \IfGG{\dot V (\phi^*(t,j))}{\displaystyle{\max_{f\in F(\phi^*(t,j))}\langle\nabla V(\phi^*(t,j)), f\rangle}}
  \end{equation}
 and thanks to \eqref{eq:LyapCostHybr_flow_2}
 \begin{equation}
  \displaystyle{\frac{d}{dt}V(\phi^*(t,j))}+L_C(\phi^*(t,j))=0.
  \end{equation}
\end{subequations}
 Likewise, 
 \begin{subequations}
  \label{eq:hjbe_d}
 for each $(t,j)\in\dom\phi^*$ such that $(t,j+1)\in\dom\phi^*$
\begin{equation}
V(\phi^*(t, j+1))=\max_{\chi\in G(\phi^*(t, j))}V(\chi).
\end{equation}
and thanks to \eqref{eq:LyapCostHybr_jump_2}
\begin{equation}
  V(\phi^*(t, j+1))-V(\phi^*(t, j))+L_D(\phi^*(t, j))=0.
  \end{equation}
\end{subequations}
Following the same arguments as in the proof of Proposition~\ref{prop:Prop1}, for each $\phi\in\mathcal{S}_{\HS_0} (\xi)$, \eqref{eq:hjbe_c} and \eqref{eq:hjbe_d} yield 
\begin{equation}
\begin{array}{lll}
\label{eq:VmaxSup}
\widetilde{\mathcal{J}}({\phi})\leq
\widetilde{\mathcal{J}}({\phi^*})=V(\xi)
\end{array}
\end{equation}

{
Notice that in light of \eqref{eq:VmaxSup}, since $\phi^*\in\mathcal{S}_{\HS_0} (\xi)$ and
$\widetilde{\mathcal{J}}(\phi^*)=V(\xi)$, one has that $\phi^*\in\argmax_{\psi\in\mathcal{S}_{\HS_0} (\xi)}\widetilde{\mathcal{J}}({\psi})$. 
On the other hand, by definition, one has that
$$
\mathcal{J}(\xi)=\sup_{\phi\in\mathcal{S}_{\HS} (\xi)}\widetilde{\mathcal{J}}({\phi}).
$$
That is, $\mathcal{J}(\xi)=\widetilde{\mathcal{J}}(\phi^*)$, and this concludes the proof.}
\end{proof}
\smallskip

The results given in this section extend previous results on cost evaluation for continuous-time nonlinear systems \cite{bernstein1993nonquadratic} and constrained difference inclusions \cite{FerranteACC2018} to hybrid inclusions. Similarly to   \cite{bernstein1993nonquadratic,ferrante2018optimality}, our results have strong connections to Lyapunov {theory}. More specifically, the applicability of our results to specific examples requires the search for a suitable Lyapunov-like function, which is in general a challenging task.
In Section \ref{Sec:WealSPE}, we 
study 
games with nonautonomous hybrid systems as in \eqref{Hinc}.
%
{
\section{\pn{Weak Saddle-Point Equilibrium Solution}}
\label{Sec:WealSPE}
The cost evaluation tools employed in approaches based on dynamic programming fall short of characterizing strategies to attain a saddle-point equilibrium solution for a two-player zero-sum game with dynamics given by hybrid inclusions. 
The classical conditions involved therein do not guarantee the existence of a lower bound in the costs for a given input action. 
Nevertheless, conditions can still be established to characterize the worst-case  cost (due to {nonuniqueness of solutions}) associated to it. 
Thus,\IfPers{ to avoid setting an ill-defined game, in this section, \pnn{we 
provide sufficient conditions to find an upper value function, namely, to solve the following problem.

\textit{Problem ($\star$):} 
Given $\xi \in \mathbb{R}^n$ and $u_1^* \in \mathcal{U}_1$, 
solve
\begin{eqnarray}
\underset{(u_1^*,u_2) \in \mathcal{U}_{\HS\IfIncd{}{} }^\infty(\xi)}
{ \underset{u_{2}}{\textup{maximize}}}
&& \mathcal{J}(\xi, (u_1^*,u_2))
\label{problemzsih}
\end{eqnarray}
where $u_1^*$ is the optimal action selected by player $P_1$
and $\mathcal{U}_{\HS\IfIncd{}{} }^\infty$ is the set of joint input actions yielding maximal solutions to $\HS\IfIncd{}{} $ complete, as defined in Section 2.1.} }
in this section, we provide sufficient conditions to solve Problem $(\diamond)$ via finding a control strategy that minimizes the worst-case cost under the maximizing adversarial action, in this case, leading to a solution of a min-max problem with nonunique solutions due to $F$ or $G$ being possibly set valued, or $C\cap D$ being nonempty. In addition, such conditions allow us evaluating the value function without computing solutions. 
First, we provide pointwise conditions that allow to upper bound the cost for an initial condition and an input action, {recalling that, according to Definition \ref{def:GenGrad},
when $F$ has compact values, the generalized derivative leads to the infinitesimal quantity 
\begin{equation}\label{defVcirc2}
\max_{\chi \in F(x,u_C)} V^\circ (x,\chi) = \max_{\chi \in F(x,u_C)}
\max_{\nu \in \partial V(x)} \langle \nu, \chi \rangle 
\end{equation}
which, {for points $(x,u_C)$ in} $\overline C$, bounds the change of $V$ along flows of solutions to $\HS$.
When $C$ and $F$ satisfy $1)$ and $2)$ in the hybrid basic conditions (see \IfSiam{\cite[Definition A.5]{TRHSVGArxiv}}{Definition \ref{def:hbc}}) and $F$ has compact values, then for each solution $(\phi,u_C)$ to $\HS$
\begin{equation}
  \frac{d}{dt} V(\phi(t,j)\pn{,u_C(t,j)}) \leq \dot{V}(\phi(t,j)\pn{,u_C(t,j)}) := 
  \hspace{-0.2cm}
  \max_{\chi \in F(\phi(t,j),u_C(t,j)) \cap T_C(\phi(t,j))} V^\circ (\phi(t,j),\chi) 
  \label{dotV2}
\end{equation}
 \pnn{for every $(t,j) \in \textup{int}(I_\phi^j) \times \{j\}$ with int$(I_\phi^j)$ nonempty.}
See \cite[Section 3.3.1]{220} and \IfSiam{\cite[Definition A.6]{TRHSVGArxiv}}{Definition \ref{def:TanCon}} for details. 

}

\begin{proposition}[Upper bound for a given control action]\label{UpperBoundOpenLoop}
  Given a system $\HS\IfIncd{}{}=(C,F,D,G) $ as in (\IfIncd{\ref{Hinc}}{\ref{Heq}}), stage costs $L_C:\mathbb{R}^n  \times \mathbb{R}^{m_C} \rightarrow \mathbb{R}_{\geq 0}$ {and} $L_D:\mathbb{R}^n  \times \mathbb{R}^{m_D} \rightarrow \mathbb{R}_{\geq 0}$, terminal cost $q: \mathbb{R}^n \rightarrow \mathbb{R}$, suppose there exists a function 
  $V\colon\dom V\rightarrow\mathbb R$ with $\dom V \supset \Pi(\overline C) \cup \Pi(D) \cup G(D) $ that is \IfGG{locally Lipschitz}{continuously differentiable} on \NotGG{an open set containing\>}$\Pi(\overline{C})$ such that
    \begin{subequations} 
  \begin{eqnarray}
   L_C(x,u_C) 
  +
  \IfGG{\dot{V}(x,u_C)}{{\underset{f \in F(x,u_C)}{\sup}} \left\langle \nabla V(x),f \right\rangle }
  \leq 0 
  \hspace{0.6cm} \forall (x,u_{C}) \in C, 
\label{ubcolinc}
  \\
         L_D(x,u_D)
        + {\underset{\chi \in G(x,u_D)}{\sup}}
        V(\chi)-V(x) \leq 0
        \hspace{0.5cm} \forall (x,u_{D}) \in D .
    \label{ubdolind}
  \end{eqnarray}
\end{subequations}
\sj{where $\dot{V}$ is defined in \eqref{dotV2}}.
  Let $(\phi,u)$ be any solution to $\HS\IfIncd{}{} $ from $\xi \in \Pi ( \overline{C}) \cup \Pi(D)$.
  {Suppose that $(t,j)\mapsto V\circ\phi(t,j)$ is bounded
and
\begin{equation}
\underset{(t,j) \in \textup{dom}\phi}
{{\limsup_{\pnn{(t,j)\rightarrow \sup \dom \phi }} }} 
V(\phi(t,j))= \underset{(t,j) \in \textup{dom}\phi}
{{\limsup_{\pnn{(t,j)\rightarrow \sup \dom \phi }} }} 
q(\phi(t,j)).
\label{eq:TerminalCondCGinc5}
\end{equation}
  } 
  Then, 
  \begin{equation}
    {\widetilde{\mathcal{J}}}(\phi,u) \leq V(\xi)
  \end{equation}
  {where ${\widetilde{\mathcal{J}}}$ is defined as in \eqref{defJTNCinc}.} 
\end{proposition}
{
\begin{proof}
From \eqref{ubcolinc}, {and given that $(\phi,u)$ is a solution to $\HS\IfIncd{}{} $ from $\xi \in \Pi ( \overline{C}) \cup \Pi(D)$}, for each $j \in \nats$ such that $I_{\phi}^j=[t_j,t_{j+1}]$ has a nonempty interior, we have, for almost all $t \in 
I_{\phi}^j$,
\begin{equation}
 \begin{split}
  L_C({\phi(t,j),u_C(t,j)}) 
 + 
 \frac{d V}{d t}(\phi(t,j)\pn{,u_C(t,j)}) 
 \leq 
 \hspace{5cm}
 \\
    L_C(\phi(t,j),u_C(t,j)) 
 + 
 \IfGG{\dot{V}({\phi(t,j),u_C(t,j)})}{\hspace{5cm}
    \\
 {\underset{f \in F({\phi(t,j),u_C(t,j)})}{\sup}}\left\langle \nabla V(\phi(t,j)),f \right\rangle} 
 \leq 0
  \end{split}
 \label{eq:PrHamiltonC1}
\end{equation}
\NotAutom{by integrating over the interval $[t_j,t_{j+1}]$, 
we obtain
\begin{eqnarray*}
0\geq  \int_{t_{j}}^{t_{j+1}} \left( L_C(\phi(t,j),u_{ C}(t,j))
+ \frac{d V}{d t}(\phi(t,j)) \right )dt
\end{eqnarray*}
from where we have
%
\begin{equation*}
0\geq  \int_{t_{j}}^{t_{j+1}} L_C(\phi(t,j),u_{ C}(t,j)) dt
+  V(\phi^*(t_{j+1},j))- V(\phi(t_j,j))
\end{equation*}
%
Summing from $j=0$ to $j=J_{\phi}$, 
we obtain 
\begin{equation*}
0
\geq\sum_{j=0}^{J_{\phi}}  \int_{t_{j}}^{t_{j+1}} L_C(\phi(t,j),u_{ C}(t,j)) dt
\\+\sum_{j=0}^{J_{\phi}} \left( V(\phi(t_{j+1},j))-V(\phi(t_j,j)) \right)
\end{equation*}
Then, solving for $V$ at the initial condition $\phi(0,0)$, we obtain
\begin{equation} \label{ContinCost2goihzol}
\begin{medsize}
\begin{split} 
V(\phi(0,0))
\geq\sum_{j=0}^{J_{\phi}}  \int_{t_{j}}^{t_{j+1}} L_C(\phi(t,j),u_{ C}(t,j)) dt
\\  + V(\phi(t_1,0))
+\sum_{j=1}^{J_{\phi}} \left( V(\phi(t_{j+1},j))-V(\phi(t_j,j)) \right) 
\end{split}
\end{medsize}
\end{equation}}
\IfGG{where $\dot{V}({\phi(t,j),u_C(t,j)})$ is defined as in \eqref{dotV2}.}{} 
In addition, from \eqref{ubdolind}, for every {$(t,j) \in \dom \phi$ such that $(t,j+1) \in \dom \phi$}, we have 
\begin{equation}
\begin{split}
  L_D(\phi(t,&j),u_D(t,j))
        + V(\phi(t,j+1)) 
        - V(\phi(t,j))
        \\
        \leq& L_D(\phi(t,j),u_D(t,j))
        + {\underset{\chi \in G(\phi(t,j),u_D(t,j))}{\sup}}
        V(\chi) - V(\phi(t,j))
        \leq 0
        \label{eq:PrHamiltonD1}
\end{split}
\end{equation}
\NotAutom{Summing from $j=0$ to $j=J_{\phi}-1$,
we obtain
\begin{equation}
\begin{medsize}
\begin{split} 
\sum_{j=0}^{J_{\phi}-1}V(\phi(t_{j+1},j)) 
\geq 
\sum_{j=0}^{J_{\phi}-1}L_D(\phi(t_{j+1},j),u_{ D}(t_{j+1},j))\nonumber
\\
+\sum_{j=0}^{J_{\phi}-1} V(\phi(t_{j+1},j+1))
\nonumber
\end{split}
\end{medsize}
\end{equation}
Then, solving for $V$ at the first jump time, we obtain
\begin{eqnarray} \nonumber
&&V(\phi(t_{1},0))
\geq V( \phi(t_{1},1))
\\&&
+ \sum_{j=0}^{J_{\phi}-1}L_D(\phi(t_{j+1},j),u_{ D}(t_{j+1},j)) \label{DiscreteCost2goNOihzwincol}
\\&&
+\sum_{j=1}^{J_{\phi}-1}\left( V(\phi(t_{j+1},j+1))-V(\phi(t_{j+1},j)) \right) \nonumber
\end{eqnarray}
In addition, given that $\phi(0,0)=\xi$, lower bounding $V(\phi(t_{1},0))$ in (\ref{ContinCost2goihzol}) by the right-hand side of (\ref{DiscreteCost2goNOihzwincol}), we obtain 
\begin{equation*}
\begin{medsize}
\begin{split} 
V(\xi)
\geq& \sum_{j=0}^{J_{\phi}} \int_{t_{j}}^{t_{j+1}} L_C(\phi(t,j),u_{C}(t,j))dt+ V(\phi(t_1,0)) 
\\&
+\sum_{j=1}^{J_{\phi}} \left( V(\phi(t_{j+1},j))-V(\phi(t_j,j)) \right)  \\
\geq&
\sum_{j=0}^{J_{\phi}} \int_{t_{j}}^{t_{j+1}} L_C(\phi(t,j),u_{C}(t,j))dt
\\&
+ \sum_{j=0}^{J_{\phi}-1}L_D(\phi(t_{j+1},j),u_{D}(t_{j+1},j))
 \\
&+\sum_{j=1}^{J_{\phi}-1}\left( V(\phi(t_{j+1},j+1))-V(\phi(t_{j+1},j)) \right) 
\\&
+ V(\phi(t_{1},1))+\sum_{j=1}^{J_{\phi}} \left( V(\phi(t_{j+1},j))-V(\phi(t_j,j)) \right) 
\end{split}
\end{medsize}
\end{equation*}
Since 
\begin{eqnarray*} \nonumber
&& V(\phi(t_{1},1))
\\&&
+\sum_{j=1}^{J_{\phi}-1}\left( V(\phi(t_{j+1},j+1))-V(\phi(t_{j+1},j)) \right) \nonumber 
\\&&
+\sum_{j=1}^{J_{\phi}} \left( V(\phi(t_{j+1},j))-V(\phi(t_j,j)) \right) \nonumber
\\ &=& V(\phi(t_{{J_{\phi}}+1},J_{\phi})) +V(\phi(t_{1},1))\nonumber
\\&&
+\sum_{j=1}^{J_{\phi}-1}\left( V(\phi(t_{j+1},j+1)) \right) \nonumber 
-\sum_{j=1}^{J_{\phi}} \left(V(\phi(t_j,j)) \right) \nonumber
\\ &=&V(\phi(t_{{J_{\phi}}+1},J_{\phi})) 
\end{eqnarray*}
then we have
\begin{eqnarray*}\nonumber
V(\xi)
&\geq&\sum_{j=0}^{J_{\phi}} \int_{t_{j}}^{t_{j+1}} L_C(\phi(t,j),u_{C}(t,j))dt
\\&&
+ \sum_{j=0}^{J_{\phi}-1}L_D(\phi(t_{j+1},j),u_{D}(t_{j+1},j))
\\&&+V(\phi(t_{J_{\phi}+1},J_{\phi}))
\nonumber
\end{eqnarray*}
By taking the limit when $(t_{J_{\phi}+1},J_{\phi}) \rightarrow \sup \dom \phi$, given the definition of the cost (\ref{defJTNCinc}), we have
\begin{eqnarray} \nonumber
V(\xi)&\geq&\sum_{j=0}^{\sup_j \dom \phi} \int_{t_{j}}^{t_{j+1}} L_C(\phi(t,j),u_{C}(t,j))dt
\\&&
+ \sum_{j=0}^{{\sup_j \dom \phi}-1}L_D(\phi(t_{j+1},j),u_{D}(t_{j+1},j))
\nonumber \\
&& 
+\underset{(t,j) \in \textup{dom}\phi}{\limsup_{t+j\rightarrow \infty} } V(\phi(t,j))
\nonumber \\ 
&=&J(\phi,u)
\label{BoundCostNOzwinc}
\end{eqnarray}}
\IfAutom{Then, thanks to \eqref{eq:PrHamiltonC1} and \eqref{eq:PrHamiltonD1}, from \IfSiam{\cite[Proposition A.1]{TRHSVGArxiv}}{Proposition \ref{Pp:ProofDerivation}}, we have that 
$\widetilde{\mathcal{J}}(\phi, u) \leq V(\xi)$, with $\widetilde{\mathcal{J}}$ defined as in \eqref{defJTNCinc}.
}{}
\end{proof}
}
{
In the following result, we study a {hybrid system
that is a version with inputs of the maximal hybrid system in \eqref{eq:argmax}}.
Solutions to this system are a subset of the solutions to $\HS$ as in \eqref{Hinc} and attain the worst-case cost {over the potential nonunique} solutions to $\HS$. Following \cite{181}, we provide conditions to exactly evaluate such a cost and show that it is an upper bound for the cost of any other solution to $\HS$.  

\begin{proposition}[Cost evaluation for maximal system]\label{lemma:SolsMaxSys}
  Consider a system $\HS\IfIncd{}{}=(C,F,D,G)$ as in (\IfIncd{\ref{Hinc}}{\ref{Heq}}), \pno{where $F$ and $G$ have compact values in $C$ and $D$, respectively,}
   stage costs $L_C:\mathbb{R}^n  \times \mathbb{R}^{m_C} \rightarrow \mathbb{R}_{\geq 0}$ {and} $L_D:\mathbb{R}^n  \times \mathbb{R}^{m_D} \rightarrow \mathbb{R}_{\geq 0}$, terminal cost $q: \mathbb{R}^n \rightarrow \mathbb{R}$, 
   and suppose that there exists a {continuous} function $V: \dom V \rightarrow \mathbb{R}$ with $\dom V \supset \Pi(\overline C) \cup \Pi(D) \cup G(D) $, 
   that is \IfGG{locally Lipschitz}{continuously differentiable} on \NotGG{an open set containing\>}$\Pi(\overline{C})$.
   Given $\xi \in \Pi ( \overline{C}) \cup \Pi(D)$ and a solution\footnote{{Solutions to the ``maximal system'' in \eqref{Hmax} exist under \IfACC{compactness of the set-valued maps, regularity of $V$ and a proper selection of the initial condition.}{the conditions in \IfSiam{\cite[Proposition A.8]{TRHSVGArxiv}}{Proposition \ref{Prop:ExistHmax}}.}}} $(\phi^*,u)$ to 
   \begin{equation} \label{Hmax}
 {\HS_{\textup{max}}}\left\{ 
 \begin{matrix}
 \dot{x} &\in& 
 \IfGG{\underset{\chi \in F(x,u_{C}) \cap T_C(x)}{\argmax} V^\circ (x,\chi)}{ 
 \underset{f \in F(x,u_{C})}{\argmax}\left\langle \nabla V(x),f \right\rangle}
  \quad  (x,u_C) \in C \\ 
 x^+ &\in& 
 \underset{\chi \in G(x,u_D)}{\textup{argmax}}V(\chi)\hspace{1.2cm}   (x,u_D) \in D
 \end{matrix} \right.
   \end{equation}
   from $\xi$
   with $u=(u_C,u_D)$, 
   \sj{where $V^\circ$ is defined in \eqref{defVcirc2},}
   if 
   \begin{subequations}
    \begin{eqnarray}\label{eq:OptimalHamiltonianC}
    0=L_C(x,u_C) +
  \IfGG{\dot{V}(x,u_C)}{{\underset{f \in F(x,u_C)}{\sup}} \left\langle \nabla V(x),f \right\rangle }
   \hspace{0.8cm} \forall (x,u_C) \in C,
  \\
    \label{eq:OptimalHamiltonianD}
     V(x) = L_D(x,u_D) + {\underset{\chi \in G(x,u_D)}{\sup}} V(\chi)\hspace{0.6cm} \forall (x,u_D) \in D,
  \end{eqnarray} 
  \end{subequations} 
  \sj{where $\dot{V}$ is defined in \eqref{dotV2}},
  and {if each $(\phi,u) \in \mathcal{S}_\HS(\xi)$ satisfies}
  \begin{equation}
  \underset{(t,j) \in \textup{dom}\phi}{\limsup_{(t,j)\rightarrow \sup \dom \phi} } V(\phi(t,j))=\underset{(t,j) \in \textup{dom}\phi}{\limsup_{(t,j)\rightarrow \sup \dom \phi} } q(\phi(t,j)),
  \label{TerminalCondCGinc0}
  \end{equation}
    then
\begin{equation}\label{eq:worstcasecost}
  \J (\xi, u) = \widetilde \J (\phi^*, u) {= V(\xi)}.
\end{equation}
\end{proposition}
{
\begin{proof}
  Following \cite{ferrante2019certifying,ferrante2018optimality}, given that $\xi \in \Pi ( \overline{C}) \cup \Pi(D)$ and $(\phi^*,u^*)$ is a solution to \eqref{Hmax},
  (that thanks to \IfSiam{\cite[Lemma A.3]{TRHSVGArxiv}}{Lemma \ref{Pp:sol2Hmaxsol2Hk}} it is also a solution to $\HS\IfIncd{}{}$ as in \eqref{Hinc}), 
  and given that $V$ satisfies \eqref{eq:OptimalHamiltonianC},
  for each $j \in \nats$ such that $I_{\phi^*}^j=[t_j,t_{j+1}]$ has a nonempty interior, we have
  \begin{enumerate}[label=\arabic*)]
    \item 
   for almost all $t \in I_{\phi^*}^j$, 
  \begin{equation}
  \begin{split}
   0 
  =&  L_C(\phi^*(t,j),u_{ C}(t,j))
  +
  \IfGG{\dot{V}({\phi^*(t,j),u_C(t,j)})}
  {
    \\&
  {\underset{f \in F(\phi^*(t,j),u_C(t,j))}{\textup{max}}}\left\langle \nabla V(\phi^*(t,j)), f \right\rangle 
  }
  \end{split}
  \label{eq:PrHamiltonC0b}
  \end{equation}
  Given that $(\phi^*,u)$ is a solution to $\HS_{\textup{max}}$, and $V$ is \IfGG{locally Lipschitz on}{continuously differentiable on a neighborhood of} $\Pi(\overline C)$, we can express its total derivative 
   for \pnn{almost all $t \in I_{\phi^*}^j$ with $\textup{int} I_{\phi^*}^j$ nonempty 
  as}
  $$\frac{d V}{d t}(\phi^*(t,j)\pn{,u_C(t,j)})=
  \IfGG{\dot{V}({\phi^*(t,j),u_C(t,j)})}
  {\underset{f \in F(\phi^*(t,j),u_C(t,j))}{\max}\left\langle \nabla V(\phi^*(t,j)), f \right\rangle}$$
   \IfGG{where $\dot V$ is defined in \eqref{dotV2}}{}, and thanks to \eqref{eq:PrHamiltonC0b}
  \begin{equation}
   0 
  =  L_C(\phi^*(t,j),u_{ C}(t,j))
  +\frac{d V}{d t}(\phi^*(t,j)\pn{,u_C(t,j)})
  \label{eq:PrHamiltonC0b2}
  \end{equation}
\end{enumerate}
 Given that $V$ satisfies \eqref{eq:OptimalHamiltonianD}, we have
  \begin{enumerate}[label=\arabic*)]
  \setcounter{enumi}{1}
    \item 
  for every $(t,j) \in \dom \phi^*$ such that $(t,j+1) \in \dom \phi^*$, 
  \begin{equation}
  \begin{split} 
  V(\phi^*(t,j))
  &= 
  L_D(\phi^*(t,j),u_{ D}(t,j))
  +{\underset{\chi \in G(\phi^*(t,j), u_D(t,j))}{\textup{max}}} V(\chi) \\
  &= L_D(\phi^*(t,j),u_{ D}(t,j))
  + V(\phi^*(t,j+1)) 
  \end{split}
  \label{eq:PrHamiltonD0b}
  \end{equation}
  %
  where $(\phi^*(t,j),u(t,j)) \in D$.
\end{enumerate}
 Now, thanks to \eqref{eq:PrHamiltonC0b2} and \eqref{eq:PrHamiltonD0b}, from \IfSiam{\cite[Proposition A.1 and Corollary A.2]{TRHSVGArxiv}}{Proposition \ref{Pp:ProofDerivation} and Corollary \ref{Rmk:ProofDerivation}}, we have that 
  \begin{equation}
  V(\xi) = \widetilde{\mathcal{J}}(\phi^*, u).
    \label{WCaseCostihzinc0b}
  \end{equation}
  {Using \eqref{eq:OptimalHamiltonianC}, for all $(x,u_C) \in C$, it holds
%
  \pn{\begin{equation}
    0\geq
    L_C(x,u_C) 
   +
   \IfGG{\dot{V}(x,u_C)}
   {\left\langle \nabla V(x),f \right\rangle \quad \forall f \in F(x,u_C)} 
  \end{equation}}}
  and using \eqref{eq:OptimalHamiltonianD}, for all 
  $(x,u_D) \in D$, we have
  \begin{equation}
    V(x)\geq
    L_D(x,u_D)
  + V(\chi)
  \quad \forall \chi \in G(x,u_D).
  \end{equation}
  Thus, for any arbitrary $(\phi, u)\in \mathcal{S}_{\HS\IfIncd{}{}}(\xi)$, we have from following the same arguments as in Proposition \ref{UpperBoundOpenLoop} and \eqref{WCaseCostihzinc0b} that
  \begin{equation}
    \widetilde{\mathcal{J}}(\phi,u) \leq \pnn{\widetilde{\mathcal{J}}}(\phi^*,u)=V(\xi)  
  \end{equation}
  which also implies that for the control action $u$, the highest cost of solutions from $\xi$ satisfies 
  \begin{equation}
    \J(\xi,u) = V(\xi).  
	\label{OptimalCostihzincb}
  \end{equation}
\end{proof}
}

\pnn{In Proposition 5.2,} a solution to  \eqref{Hmax} attains the worst-case cost among the potential nonunique solutions to $\HS$ due to the set-valuedness of its data. Furthermore, the worst-case cost associated to the input action whose values satisfy \eqref{eq:OptimalHamiltonianC} and \eqref{eq:OptimalHamiltonianD} can be evaluated without computing solutions and equates $V(\xi)$.

{\begin{corollary}[Change of Signs]
 If the conditions in Proposition \ref{lemma:SolsMaxSys} hold {with} inequality, {namely, if ``$=$'' in \eqref{eq:OptimalHamiltonianC} and \eqref{eq:OptimalHamiltonianD}  is replaced with ``$\leq$'' (or ``$\geq$''),}  then \eqref{eq:worstcasecost} holds {with} ``$\leq$'' (or ``$\geq$'', respectively).
  \label{Cor:ProofDerivation0}
\end{corollary}}
}{}
 Based on Proposition \ref{UpperBoundOpenLoop} providing an upper bound on the cost ${\widetilde{\mathcal{J}}}$, \IfMS{and the exact cost evaluation in Proposition \ref{lemma:SolsMaxSys},}{} we introduce the main result of the section: sufficient conditions to 
characterize the \pnn{weak} saddle-point equilibrium strategy and evaluate the value function without computing solutions. 
\begin{ttheorem}{Sufficient conditions for Problem $(\diamond)$}\label{thHJBszihinc} 
  Given a 
  system \\$\HS\IfIncd{}{}=(C,F,D,G) $ as in (\IfIncd{\ref{Hinc}}{\ref{Heq}}),
   stage costs $L_C:\mathbb{R}^n  \times \mathbb{R}^{m_C} \rightarrow \mathbb{R}_{\geq 0}$ {and} $L_D:\mathbb{R}^n  \times \mathbb{R}^{m_D} \rightarrow \mathbb{R}_{\geq 0}$, and terminal cost $q: \mathbb{R}^n \rightarrow \mathbb{R}$, {suppose the following hold:} 
   \begin{enumerate}[label=\arabic*)]
     \item There exists a {continuous} function $V: \dom V \rightarrow \mathbb{R}$, \pnn{$\dom V \supset \Pi(\overline C) \cup \Pi(D) \cup G(D)$}, that is \IfGG{locally Lipschitz on}{continuously differentiable on a neighborhood of} $\Pi(C)$ and a feedback law
   $\kappa:=(\kappa_C,\kappa_D){=((\kappa_{C1},\kappa_{C2}),(\kappa_{D1},\kappa_{D2}))}: \mathbb{R}^n \rightarrow \mathbb{R}^{m_C} \times \mathbb{R}^{m_D}$ {such that $F(x,\kappa_C(x))$ and 
   $G(x,\kappa_D(x))$ are compact for every $x$ such that $(x,\kappa_C(x))\in C$ and $(x,\kappa_D(x))\in D$, respectively,}  
   and {such that the functions
   $\mathcal{L}_C(x,u_C):=  L_C(x,u_C)+
   \IfGG{\dot{V}(x,u_C)}{{\underset{f \in F(x,u_C)}{\sup}} \left\langle \nabla V(x),f \right\rangle }$, and
   $\mathcal{L}_D(x,u_D):=  L_D(x,u_D) + {\underset{g \in G(x,u_D)}{\sup}} V(g)$
     satisfy} 
  \begin{equation}
    0= 
    \mathcal{L}_C(x,\kappa_C(x))
  \hspace{1cm} \forall x: (x,\kappa_C(x)) \in C, 
  \label{HJBzsihinc}
  \end{equation}
  \begin{equation}
        \begin{split}
      0\leq &\>  \mathcal{L}_C(x,(u_{C1},\kappa_{C2}(x)))
    \hspace{1cm} 
    \forall (x,u_{C1}):(x,(u_{C1},\kappa_{C2}(x))) \in C, 
  \end{split}
  \label{HJBzsihinclb}
    \end{equation}
  \begin{equation}
      \begin{split}
    0\geq &\>  \mathcal{L}_C(x,(\kappa_{C1}(x),u_{C2}))
  \hspace{1cm} \forall (x,u_{C2}):(x,(\kappa_{C1}(x),u_{C2})) \in C, 
\end{split}
\label{HJBzsihincub}
  \end{equation}
  \begin{equation}
  \begin{split}
    V(x)= \mathcal{L}_D(x,\kappa_D(x)) 
   \hspace{0.6cm} \forall x: (x,\kappa_D(x)) \in D,
  \end{split}
  \label{Bellmanzsihinc}
  \end{equation}
      \begin{equation}
            \begin{split}
             V(x)\leq  &\> \mathcal{L}_D(x,(u_{D1},\kappa_{D2}(x)))
              \hspace{0.4cm} \forall (x,u_{D1}):(x,(u_{D1},\kappa_{D2}(x)))) \in D ,
            \end{split}
          \label{Bellmanzsihinclb}
        \end{equation}
  \begin{equation}
      \begin{split}
        V(x)\geq  &\> \mathcal{L}_D(x,(\kappa_{D1}(x),u_{D2}))
        \hspace{0.4cm} \forall (x,u_{D2}):(x,(\kappa_{D1}(x),u_{D2})) \in D,
      \end{split}
    \label{Bellmanzsihincub}
  \end{equation}
  \sj{where $\dot V$ is defined in \eqref{dotV2}.}

  \item {For each $\xi \in \Pi ( \overline{C}) \cup \Pi(D)$, each 
  {$\phi \in \mathcal{S}_{\HS\IfIncd{_\kappa}{}}(\xi)$ satisfies}
  \begin{equation}
  \underset{(t,j) \in \textup{dom}\phi}{\limsup_{(t,j)\rightarrow \sup \dom \phi} } V(\phi(t,j))=\underset{(t,j) \in \textup{dom}\phi}{\limsup_{(t,j)\rightarrow \sup \dom \phi} } q(\phi(t,j)),
  \label{TerminalCondCGinc}
  \end{equation}}
  \end{enumerate}
  Then, {the feedback law $\kappa$ is a \pnn{weak} saddle-point equilibrium, and}  
    \begin{equation}
      \begin{split}
\J^*(\xi)
    =V(\xi)  \hspace{0.5cm} \forall \xi \in \Pi ( \overline{C}) \cup \Pi(D). 
      \end{split}
  \label{ResultValueinc}
  \end{equation}
  \end{ttheorem}
  %
  {
  \begin{proof}
  To show the claim we proceed as follows:
  \begin{enumerate}[label=\alph*)]
  \item Pick {any} initial condition $\xi$ and evaluate the cost associated to the 
  solutions yield by $\IfIh{\kappa=(\kappa_C,\kappa_D)}{\gamma}$ 
  from $\xi$. Find an upper bound for this cost. 
  \item \pno{Lower bound the cost associated to the worst-case solutions from $\xi$ when $P_{2}$ plays $\IfIh{\kappa}{\gamma}_2:=(\IfIh{\kappa}{\gamma}_{C2},\IfIh{\kappa}{\gamma}_{D2})$ by the value of the function $V$ evaluated at $\xi$.}
  \item Upper bound the cost associated to the solutions from $\xi$ when $P_{1}$ plays $\IfIh{\kappa}{\gamma}_1:=(\IfIh{\kappa}{\gamma}_{C1},\IfIh{\kappa}{\gamma}_{D1})$ by the value of the function $V$ evaluated at $\xi$. 
  \item By showing that the cost of the solutions from $\xi$ when $P_{1}$ plays $\IfIh{\kappa}{\gamma}_1$ are not greater than the worst-case cost of the solutions yielded by $\IfIh{\kappa=(\kappa_1,\kappa_2)}{\gamma}$ from $\xi$,
  and by showing that the worst-case cost of the solutions from $\xi$ when $P_{2}$ plays $\IfIh{\kappa}{\gamma}_2$ is not smaller than the worst-case cost of the solutions yielded by $\IfIh{\kappa}{\gamma}$ from $\xi$, 
  we {establish \eqref{ResultValueinc}}.
  \end{enumerate}
  Proceeding as in item a above, \IfMS{pick any $\xi \in \Pi ( \overline{C}) \cup \Pi(D)$ and any solution $(\phi^*,u^*)$ to \eqref{Hmax} from $\xi$ with $\dom \phi^* \ni (t,j) \mapsto u^*(t,j)= \IfIh{\kappa}{\gamma}(\IfIh{}{t,j,}\phi^*(t,j))$.
  Thanks to Proposition~\ref{lemma:SolsMaxSys}, we have that
  \begin{equation}
    V(\xi) = \widetilde{\mathcal{J}}(\phi^*, u^*).
      \label{WCaseCostihzinc0}
    \end{equation} 
    and for any arbitrary $(\phi, u^*)\in \mathcal{S}_{\HS\IfIncd{}{}}(\xi)$, we have that
    \begin{equation}
     \widetilde{\mathcal{J}}(\phi,u^*) \leq \widetilde{\mathcal{J}}(\phi^*,u^*)=V(\xi)  
    \end{equation}
    which also implies that for the control action $u^*$, defined by the feedback law $\kappa$, the highest cost of solutions from $\xi$ satisfies 
    \begin{equation}
      \J(\xi,u^*) = V(\xi).  
    \label{OptimalCostihzinc}
    \end{equation}
  }{
    following \cite{ferrante2019certifying,ferrante2018optimality}, pick any $\xi \in \Pi ( \overline{C}) \cup \Pi(D)$ and any solution $(\phi^*,u^*)$ to 
  \begin{equation}
  \begin{medsize}
{\HS_{\textup{max}}}\left\{ 
\begin{matrix}
\dot{x} &\in& 
\pnn{\underset{f \in F(x,u_{C})}{\textup{argmax}}}\left\langle \nabla V(x),f \right\rangle \quad  (x,u_C) \in C \\ 
x^+ &\in& 
\pnn{\underset{g \in G(x,u_D)}{\textup{argmax}}}V(g)\hspace{1.2cm}   (x,u_D) \in D
\end{matrix} \right.
\end{medsize}
\label{Hmax}
  \end{equation}
  {(that thanks to Lemma \ref{Pp:sol2Hmaxsol2Hk} is also a solution to $\HS\IfIncd{}{}$),} 
  with $\dom \phi^* \ni (t,j) \mapsto u^*(t,j)= \IfIh{\kappa}{\gamma}(\IfIh{}{t,j,}\phi^*(t,j))$. 
  %
  Given that $V$ and $\kappa$ satisfy 
  (\ref{HJBzsihinc}), 
  for each $j \in \nats$ such that $I_{\phi^*}^j=[t_j,t_{j+1}]$ has a nonempty interior int$I_{\phi^*}^j$, 
  we have:
  \begin{enumerate}
    \item 
   for all $t \in \textup{int}I_{\phi^*}^j$, 
  \begin{equation}
  \begin{medsize}
  \begin{split}
   0
  =&  L_C(\phi^*(t,j),\IfIh{\kappa}{\gamma}_{ C}(\phi^*(t,j)))
  \\&\hspace{1cm}+{\underset{f \in F(\phi^*(t,j),\IfIh{\kappa}{\gamma}_C(\phi^*(t,j)))}{\textup{max}}}\left\langle \nabla V(\phi^*(t,j)), f \right\rangle 
  \end{split}
  \end{medsize}
  \label{eq:PrHamiltonC0}
  \end{equation}
  and $\phi^*(t,j) \in C_{\IfIh{\kappa}{\gamma}}$, as in (\ref{Hkeq}). 
\end{enumerate}
  Given that {$(\phi^*,u^*)$ is a solution to $\HS_{\textup{max}}$, and $V$ is continuously differentiable on a neighborhood of $\Pi(C)$, we can express its total derivative as}
  $$\frac{d V}{d t}(\phi^*(t,j))=\pnn{\underset{f \in F(\phi^*(t,j),\IfIh{\kappa}{\gamma}_C(\phi^*(t,j)))}{\max}}\left\langle \nabla V(\phi^*(t,j)), f \right\rangle$$
   for every $(t,j)\in \textup{int}(I_{\phi^*}^j)\times \{j\}$ with $\textup{int}(I_{\phi^*}^j)$ nonempty. 
  \NotConf{By integrating over the interval $[t_j,t_{j+1}]$, 
  we obtain
  \begin{eqnarray*}
  0=  \int_{t_{j}}^{t_{j+1}} \left( L_C(\phi^*(t,j),\IfIh{\kappa}{\gamma}_{ C}(\phi^*(t,j)))
  + \frac{d V}{d t}(\phi^*(t,j)) \right )dt
  \end{eqnarray*}
  from where we have
  \begin{multline*}
  0=  \int_{t_{j}}^{t_{j+1}} L_C(\phi^*(t,j),\IfIh{\kappa}{\gamma}_{ C}(\phi^*(t,j))) dt
  \\+  V(\phi^*(t_{j+1},j))- V(\phi^*(t_j,j))
  \end{multline*}
  %
  %
  %
  Summing from $j=0$ to $j=J_{\phi^*}$, 
  we obtain 
  \begin{multline*}
  0
  =\sum_{j=0}^{J_{\phi^*}}  \int_{t_{j}}^{t_{j+1}} L_C(\phi^*(t,j),\IfIh{\kappa}{\gamma}_{ C}(\phi^*(t,j))) dt
  \\+\sum_{j=0}^{J_{\phi^*}} \left( V(\phi^*(t_{j+1},j))-V(\phi^*(t_j,j)) \right)
  \end{multline*}
  Then, solving for $V$ at the initial condition $\phi^*(0,0)$, we obtain
  %
  %
  \begin{equation} \label{ContinCost2goihzinc}
  \begin{medsize}
  \begin{split} 
  V(\phi^*(0,0))
  =\sum_{j=0}^{J_{\phi^*}}  \int_{t_{j}}^{t_{j+1}} L_C(\phi^*(t,j),\IfIh{\kappa}{\gamma}_{ C}(\phi^*(t,j))) dt
  \\  + V(\phi^*(t_1,0))
  +\sum_{j=1}^{J_{\phi^*}} \left( V(\phi^*(t_{j+1},j))-V(\phi^*(t_j,j)) \right) 
  \end{split}
  \end{medsize}
  \end{equation}
  }
 Given that $V$ and $\IfIh{\kappa}{\gamma}_D$ satisfy 
  (\ref{Bellmanzsihinc}), we have: 
  \begin{enumerate}
    \item 
  for every $(t,j) \in \dom \phi^*$ such that $(t,j+1) \in \dom \phi^*$, 
  \begin{equation}
  \begin{medsize}
  \begin{split} 
  V(\phi^*(t,j))
  = 
  L_D(\phi^*(t,j),\IfIh{\kappa}{\gamma}_{ D}(\phi^*(t,j)))
  \hspace{1.2cm} 
  \\   \hspace{0.5cm} 
  +{\underset{g \in G(\phi^*(t,j), \IfIh{\kappa}{\gamma}_D(\phi^*(t,j)))}{\textup{max}}} V(g) \\
  {=} L_D(\phi^*(t,j),\IfIh{\kappa}{\gamma}_{ D}(\phi^*(t,j)))
  + V(\phi^*(t,j+1)) 
  \end{split}
  \end{medsize}
  \label{eq:PrHamiltonD0}
  \end{equation}
  where $\phi^*(t,j)\in D_{\IfIh{\kappa}{\gamma}}$ as in (\ref{Hkeq}).
\end{enumerate}
\NotAutom{  Summing both sides from $j=0$ to $j=J_{\phi^*}-1$, 
  we obtain
  \begin{equation*}
  \begin{medsize}
  \begin{split} 
  \sum_{j=0}^{J_{\phi^*}-1} V(\phi^*(t,j))
  =& \sum_{j=0}^{J_{\phi^*}-1} L_D(\phi^*(t,j),\IfIh{\kappa}{\gamma}_{ D}(\phi^*(t,j))) 
  \\
  &+\sum_{j=0}^{J_{\phi^*}-1} V( \phi^*(t,j+1))
  \end{split}
  \end{medsize}
  \end{equation*}
  Then, solving for $V$ at the first jump time, we obtain
  %
  \begin{eqnarray} \nonumber
  &&V(\phi^*(t_{1},0))
  = 
  \\ \nonumber
  &&V( \phi^*(t_{1},1))
  + \sum_{j=0}^{J_{\phi^*}-1} L_D(\phi^*(t,j),\IfIh{\kappa}{\gamma}_{ D}(\phi^*(t,j))) 
  \\ \label{DiscreteCost2goihzinc}
  &&+\sum_{j=1}^{J_{\phi^*}-1}\left( V( \phi^*(t,j+1))-V(\phi^*(t,j)) \right)
  \end{eqnarray}
  Given that $\phi^*(0,0)=\xi$, by substituting the right-hand side of (\ref{DiscreteCost2goihzinc}) in  (\ref{ContinCost2goihzinc}), we obtain
  \begin{equation} 
  \begin{medsize}
  \begin{split} 
  V(\xi)
  =&\sum_{j=0}^{J_{\phi^*}}  \int_{t_{j}}^{t} L_C(\phi^*(t,j),\IfIh{\kappa}{\gamma}_{ C}(\phi^*(t,j))) dt
  \nonumber \\&
  + V(\phi^*(t_1,0))\nonumber 
  +\sum_{j=1}^{J_{\phi^*}} \left( V(\phi^*(t,j))-V(\phi^*(t_j,j)) \right) \nonumber \\
  =&\sum_{j=0}^{J_{\phi^*}}  \int_{t_{j}}^{t} L_C(\phi^*(t,j),\IfIh{\kappa}{\gamma}_{ C}(\phi^*(t,j))) dt
  \nonumber \\&
  + \sum_{j=0}^{J_{\phi^*}-1}L_D(\phi^*(t,j),\IfIh{\kappa}{\gamma}_{ D}(\phi^*(t,j)))
  + V(\phi^*(t_{1},1))
  \nonumber \\&
  +\sum_{j=1}^{J_{\phi^*}-1}\left( V(\phi^*(t,j+1))-V(\phi^*(t,j)) \right) 
  \nonumber \\&
  +\sum_{j=1}^{J_{\phi^*}} \left( V(\phi^*(t,j))-V(\phi^*(t_j,j)) \right) \nonumber 
  \end{split}
  \end{medsize}
  \end{equation}
  Since 
  \begin{eqnarray}
  && V(\phi^*(t_{1},1))
  \nonumber \\&&
  +\sum_{j=1}^{J_{\phi^*}-1}\left( V(\phi^*(t,j+1))-V(\phi^*(t,j)) \right) 
  \nonumber \\&&
  +\sum_{j=1}^{J_{\phi^*}} \left( V(\phi^*(t,j))-V(\phi^*(t_j,j)) \right) \nonumber
  \\ &=& V(\phi^*(t_{J_{\phi^*}+1},J_{\phi^*}))+V(\phi^*(t_{1},1))
  \nonumber \\&&
  +\sum_{j=1}^{J_{\phi^*}-1}\left( V(\phi^*(t,j+1)) \right) \nonumber 
  -\sum_{j=1}^{J_{\phi^*}} \left(V(\phi^*(t_j,j)) \right) \nonumber
  \\ &=&V(\phi^*(t_{J_{\phi^*}+1},J_{\phi^*})) 
  \end{eqnarray}
  then we have
  \begin{eqnarray}
  V(\xi)
  &=&\sum_{j=0}^{J_{\phi^*}}  \int_{t_{j}}^{t} L_C(\phi^*(t,j),\IfIh{\kappa}{\gamma}_{ C}(\phi^*(t,j))) dt
  \nonumber \\&&
  + \sum_{j=0}^{J_{\phi^*}-1}L_D(\phi^*(t,j),\IfIh{\kappa}{\gamma}_{ D}(\phi^*(t,j))) 
  \\&&+V(\phi^*(t_{J_{\phi^*}+1},J_{\phi^*})) \nonumber
  \end{eqnarray}
  By taking the limit when $(t_{J_{\phi^*}+1},J_{\phi^*}) \rightarrow \sup \dom \phi^*$, and given that (\ref{TerminalCondCGinc}) holds, we have 
  \begin{eqnarray}
  V(\xi)&=& \sum_{j=0}^{\sup_j \dom \phi^*}   \int_{t_{j}}^{t} L_C(\phi^*(t,j),\IfIh{\kappa}{\gamma}_{ C}(\phi^*(t,j))) dt
  \nonumber \\&&
  + \sum_{j=0}^{\sup_j \dom \phi^* -1} L_D(\phi^*(t,j),\IfIh{\kappa}{\gamma}_{ D}(\phi^*(t,j))) 
  \nonumber
  \\ && 
  +\underset{(t,j) \in \textup{dom}\phi^*}{\limsup_{t+j\rightarrow \infty} } V(\phi^*(t,j))
  \nonumber \\ 
  &=& J(\phi^*,u^*)
  \label{WCaseCostihzinc}
  \end{eqnarray}}
  \IfConf{Now, 
  \pnn{thanks to \eqref{eq:PrHamiltonC0} and \eqref{eq:PrHamiltonD0},} from Proposition \ref{Pp:ProofDerivation} and Corollary \ref{Rmk:ProofDerivation}, we have that 
  \begin{equation}
  V(\xi) = \widetilde{\mathcal{J}}(\phi^*, u^*).
    \label{WCaseCostihzinc0}
  \end{equation}
  }{}
  \pnn{Using \eqref{HJBzsihinc}, for all $x$ such that $(x,\kappa_C(x)) \in C$, it holds
  \begin{equation}
    0\geq
    L_C(x,\kappa_C(x)) 
   +\left\langle \nabla V(x),f \right\rangle 
   \quad \forall f \in F(x,\kappa_C(x))
  \end{equation}
  and using \eqref{Bellmanzsihinc}, for all $x$ such that $(x,\kappa_D(x)) \in D$, we have
  \begin{equation}
    V(x)\geq 
    L_D(x,\kappa_D(x))
  + V(g)
  \quad \forall g \in G(x,\kappa_D(x)).
  \end{equation}
  Thus, for any arbitrary $(\phi, u^*)\in \mathcal{S}_{\HS\IfIncd{}{}}(\xi)$, we have from Proposition \ref{UpperBoundOpenLoop} and \eqref{WCaseCostihzinc0} that
  \begin{equation}
    \pnn{\widetilde{\mathcal{J}}}(\phi,u^*) \leq \pnn{\widetilde{\mathcal{J}}}(\phi^*,u^*)=V(\xi)  
  \end{equation}
  which also implies that for the control action $u^*$, defined by the feedback law $\kappa$, the highest cost of solutions from $\xi$ satisfies 
  \begin{equation}
    \J(\xi,u^*) = V(\xi).  
	\label{OptimalCostihzinc}
  \end{equation}
  }
  }
  Proceeding with item b as above, \IfMS{
          pick any  $(\phi_s,u^s) \in \mathcal{S}_{\HS_{\max}  }^s (\xi)$
   with ${u^s=(u_{1}^s,u_{2}^s)}$, where 
 $\mathcal{S}_{\HS_{\max}  }^s (\xi):= \{(\phi,(u_{1},u_{2}))\in \mathcal{S}_{\HS_{\max}} (\xi):
\dom \phi \ni (t,j) \mapsto u_1(t,j)= \bar{\IfIh{\kappa}{\gamma}}_1(\IfIh{}{t,j,}\phi(t,j))$, $\bar{\IfIh{\kappa}{\gamma}}_1 \in \mathcal{K}_1\}$
   and $\dom \phi \ni (t,j) \mapsto u_2(t,j)= {\IfIh{\kappa}{\gamma}}_2(\IfIh{}{t,j,}\phi(t,j))
\}$, \pnn{where $\HS_{\max}$ is given in \eqref{Hmax}.
}
      Thanks to Proposition~\ref{lemma:SolsMaxSys} and Corollary~\ref{Cor:ProofDerivation0}, we have that
      \begin{equation}
      V(\xi) \leq  \widetilde{\mathcal{J}}(\phi_s,u^s)
      =   \mathcal{J}(\xi,u^s),
          \label{BoundCostNOz}
      \end{equation}
      where the right-hand equality is established given that for any arbitrary $(\phi, u^s)\in \mathcal{S}_{\HS\IfIncd{}{}}(\xi)$, we have that
      \begin{equation}
        {\widetilde{\mathcal{J}}}(\phi,u^s) \leq {\widetilde{\mathcal{J}}}(\phi_s,u^s).
      \end{equation}
      %
      Notice that it is not possible to guarantee that $V(\xi)$ is a lower bound for every solution {obtained with} $u^s$.
  }
  {
  following \cite{ferrante2019certifying,ferrante2018optimality}, pick 
  any solution $(\phi_s,u^s)$ to 
  $\HS_{\max}$ as in \eqref{Hmax} from the initial condition $\xi$,
  (that thanks to Lemma \ref{Pp:sol2Hmaxsol2Hk} is also a solution to $\HS\IfIncd{}{}$), 
  with $u=(u_{1},u_{2})$, 
  $\dom \phi_s \ni (t,j) \mapsto u_1(t,j)= \bar{\IfIh{\kappa}{\gamma}}_1(\IfIh{}{t,j,}\phi_s(t,j))$ for some $\bar{{\IfIh{\kappa}{\gamma}}}_1:=(\bar{\kappa}_{C1},\bar{\kappa}_{D1}) \in \mathcal{K}_1$
  , and $\dom \phi_s \ni (t,j) \mapsto u_2(t,j)= {\IfIh{\kappa}{\gamma}}_2(\IfIh{}{t,j,}\phi_s(t,j))$ for ${\IfIh{\kappa}{\gamma}}_2 := (\kappa_{C2},\kappa_{D2})$. 
  %
%
  Since $\bar{\IfIh{\kappa}{\gamma}}_1$ does not necessarily attain the lower bound in 
  \eqref{HJBzsihinclb}, then, 
  for each $j \in \nats$ such that $I_{\phi_s}^j=[t_j, t_{j+1}]$ has a nonempty interior int$I_{\phi_s}^j$, we have:
  \begin{enumerate}
    \item for every $t \in \textup{int}I_{\phi_s}^j$,
  \begin{equation*}
\begin{medsize}
\begin{split}
 & 0
  \leq L_C(\phi_s(t,j),u_{ C}^s(t,j))\\
  & \hfill+\underset{f \in F(x,(\bar{\kappa}_{C1}(x),\kappa_{C2}(x)))}{\sup} \left\langle \nabla V(\phi_s(t,j)),f \right\rangle \hspace{2cm}
\end{split}
\end{medsize}
  \end{equation*}
  which implies
  \begin{equation} 
  0 \leq L_C(\phi_s(t,j),u_{ C}^s(t,j))+\frac{d V}{d t}(\phi_s(t,j))\IfAutom{.}{,}
  \label{eq:PrHamiltonC2}
  \end{equation}
\end{enumerate}
  \NotConf{ and 
  by integrating over the interval $[t_j,t_{j+1}]$, 
  we obtain
  \begin{eqnarray*}
  0 \geq \int_{t_{j}}^{t_{j+1}} \left( L_C(\phi_s(t,j),u_{ C}^s(t,j))
  + \frac{d V}{d t}(\phi_s(t,j)) \right) dt
  \end{eqnarray*}
  from which we have
  \begin{eqnarray*}
  V(\s(t_j,j))
  \geq \int_{t_{j}}^{t_{j+1}} L_C(\phi_s(t,j),u_{ C}^s(t,j))dt
  \\
  +V(\phi_s(t_{j+1},j))
  \end{eqnarray*}
  Summing both sides from $j=0$ to $j=J_{\phi_s}$, 
  we obtain 
  \begin{eqnarray}
  \sum_{j=0}^{J_{\phi_s}}V(\phi_s(t_j,j)) 
  \geq \sum_{j=0}^{J_{\phi_s}} \int_{t_{j}}^{t_{j+1}} L_C(\phi_s(t,j),u_{C}^s(t,j))dt
  \nonumber \\
  +\sum_{j=0}^{J_{\phi_s}}V(\phi_s(t_{j+1},j))
  \nonumber
  \end{eqnarray}
  Then, solving for $V$ at the initial condition $\phi_s(0,0)$, we obtain
  \begin{equation}\label{ContinCost2goNOihzwinc}
  \begin{medsize}
  \begin{split} 
  &V(\phi_s(0,0))
  \geq \sum_{j=0}^{J_{\phi_s}} \int_{t_{j}}^{t_{j+1}} L_C(\phi_s(t,j),u_{C}^s(t,j))dt
  \\&
  + V(\phi_s(t_1,0))+\sum_{j=1}^{J_{\phi_s}} \left( V(\phi_s(t_{j+1},j))-V(\phi_s(t_j,j)) \right) 
  \end{split}
  \end{medsize}
  \end{equation}
  }
In addition, since $\bar{\IfIh{\kappa}{\gamma}}_1$ does not necessarily attain the lower bound in \eqref{Bellmanzsihinclb}, we have:
\begin{enumerate}
  \setcounter{enumi}{1}
\item
for every $(t,j) \in \dom \phi_s $ such that $(t,j+1) \in \dom \phi_s$, 
   \begin{equation*}
\begin{split}
  &V(\phi_s(t,j))
  \leq 
  L_D(\phi_s(t,j),u_{ D}^s(t,j))
  \\ & \hspace{4cm}
  + \underset{g\in G(\phi_s(t,j),u_{D}^s(t,j))}{\sup} V(g)
    \end{split}
  \end{equation*}
  and 
     \begin{equation}
\begin{split}
    &V(\phi_s(t,j))
    \leq
     L_D(\phi_s(t,j),u_{ D}^s(t,j))
    + V( \phi_s(t,j+1))
      \end{split}
   \label{eq:PrHamiltonD2}
  \end{equation}
  \NotConf{Summing both sides from $j=0$ to $j=J_{\phi_s}-1$, 
  we obtain
  \begin{equation}
  \begin{medsize}
  \begin{split} 
  \sum_{j=0}^{J_{\phi_s}-1}V(\phi_s(t,j)) 
  \geq 
  \sum_{j=0}^{J_{\phi_s}-1}L_D(\phi_s(t,j),u_{ D}^s(t,j))\nonumber
  \\
  +\sum_{j=0}^{J_{\phi_s}-1} V(\phi_s(t,j+1))
  \nonumber
  \end{split}
  \end{medsize}
  \end{equation}
  Then, solving for $V$ at the first jump time, we obtain
  %
  \begin{eqnarray} \nonumber
  &&V(\phi_s(t_{1},0))
  \geq V( \phi_s(t_{1},1))
  \\&&
  + \sum_{j=0}^{J_{\phi_s}-1}L_D(\phi_s(t,j),u_{ D}^s(t,j)) \label{DiscreteCost2goNOihzwinc}
  \\&&
  +\sum_{j=1}^{J_{\phi_s}-1}\left( V(\phi_s(t,j+1))-V(\phi_s(t,j)) \right) \nonumber
  \end{eqnarray}
  In addition, given that $\phi_s(0,0)=\xi$, lower bounding $V(\phi_s(t_{1},0))$ in (\ref{ContinCost2goNOihzwinc}) by the right-hand side of (\ref{DiscreteCost2goNOihzwinc}), we obtain 
  \begin{equation*}
  \begin{medsize}
  \begin{split} 
  V(\xi)
  \geq& \sum_{j=0}^{J_{\phi_s}} \int_{t_{j}}^{t} L_C(\phi_s(t,j),u_{C}^s(t,j))dt+ V(\phi_s(t_1,0)) 
  \\&
  +\sum_{j=1}^{J_{\phi_s}} \left( V(\phi_s(t,j))-V(\phi_s(t_j,j)) \right)  \\
  \geq&
  \sum_{j=0}^{J_{\phi_s}} \int_{t_{j}}^{t} L_C(\phi_s(t,j),u_{C}^s(t,j))dt
  \\&
  + \sum_{j=0}^{J_{\phi_s}-1}L_D(\phi_s(t,j),u_{D}^s(t,j))
   \\
  &+\sum_{j=1}^{J_{\phi_s}-1}\left( V(\phi_s(t,j+1))-V(\phi_s(t,j)) \right) 
  \\&
  + V(\phi_s(t_{1},1))+\sum_{j=1}^{J_{\phi_s}} \left( V(\phi_s(t,j))-V(\phi_s(t_j,j)) \right) 
  \end{split}
  \end{medsize}
  \end{equation*}
  Since 
  \begin{eqnarray*} \nonumber
  && V(\phi_s(t_{1},1))
  \\&&
  +\sum_{j=1}^{J_{\phi_s}-1}\left( V(\phi_s(t,j+1))-V(\phi_s(t,j)) \right) \nonumber 
  \\&&
  +\sum_{j=1}^{J_{\phi_s}} \left( V(\phi_s(t,j))-V(\phi_s(t_j,j)) \right) \nonumber
  \\ &=& V(\phi_s(t_{{J_{\phi_s}}+1},J_{\phi_s})) +V(\phi_s(t_{1},1))\nonumber
  \\&&
  +\sum_{j=1}^{J_{\phi_s}-1}\left( V(\phi_s(t,j+1)) \right) \nonumber 
  -\sum_{j=1}^{J_{\phi_s}} \left(V(\phi_s(t_j,j)) \right) \nonumber
  \\ &=&V(\phi_s(t_{{J_{\phi_s}}+1},J_{\phi_s})) 
  \end{eqnarray*}
  then we have
  \begin{eqnarray*}\nonumber
  V(\xi)
  &\geq&\sum_{j=0}^{J_{\phi_s}} \int_{t_{j}}^{t} L_C(\phi_s(t,j),u_{C}^s(t,j))dt
  \\&&
  + \sum_{j=0}^{J_{\phi_s}-1}L_D(\phi_s(t,j),u_{D}^s(t,j))
  \\&&+V(\phi_s(t_{J_{\phi_s}+1},J_{\phi_s}))
  \nonumber
  \end{eqnarray*}
  By taking the limit when $(t_{J_{\phi_s}+1},J_{\phi_s}) \rightarrow \sup \dom \phi_s$, and given that (\ref{TerminalCondCGinc}) holds, we have
  \begin{eqnarray} \nonumber
  V(\xi)&\geq&\sum_{j=0}^{\sup_j \dom \phi_s} \int_{t_{j}}^{t_{j+1}} L_C(\phi_s(t,j),u_{C}^s(t,j))dt
  \\&&
  + \sum_{j=0}^{{\sup_j \dom \phi_s}-1}L_D(\phi_s(t_{j+1},j),u_{D}^s(t_{j+1},j))
  \nonumber \\
  && 
  +\underset{(t,j) \in \textup{dom}\phi_s}{\limsup_{t+j\rightarrow \infty} } V(\phi_s(t,j))
  \nonumber \\ 
  &=&J(\phi_s,u^s)
  \label{BoundCostNOzwinc}
  \end{eqnarray}
  }
\end{enumerate}
  \IfConf{
    \pno{Now, 
  thanks to \eqref{eq:PrHamiltonC2} and \eqref{eq:PrHamiltonD2}, from Proposition \ref{Pp:ProofDerivation} and Corollary \ref{Rmk:ProofDerivation}, we have that 
  \begin{equation}
  V(\xi) \leq  \widetilde{\mathcal{J}}(\phi_s,u^s).
    \label{BoundCostNOzuincp}
  \end{equation}
    Using \eqref{HJBzsihinclb}, for all $(x,u_{C1})$ such that $(x,(u_{C1},\kappa_{C2}(x))) \in C$, it holds
    \begin{equation}
      \begin{split}
      &0\leq
      L_C(x,(u_{C1},\kappa_{C2}(x))) 
     +\left\langle \nabla V(x),f \right\rangle \\
     &\hspace{3.5cm}\forall f \in F(x,(u_{C1},\kappa_{C2}(x)))
     \end{split}
    \end{equation}
    and using \eqref{Bellmanzsihinclb}, for all $(x,u_{D1})$ such that $(x,(u_{D1},\kappa_{D2}(x))) \in D$, we have
    \begin{equation}
      \begin{split}
      &V(x)\leq 
      L_D(x,(u_{D1},\kappa_{D2}(x)))
    + V(g) \\
    &\hspace{3.5cm} \forall g \in G(x,(u_{D1},\kappa_{D2}(x))).
    \end{split}
    \end{equation}
    Thus, for any arbitrary $(\phi, u^s)\in \mathcal{S}_{\HS\IfIncd{}{}}(\xi)$, we have from Proposition \ref{UpperBoundOpenLoop} and \eqref{BoundCostNOz} that
    \begin{equation}
      {\widetilde{\mathcal{J}}}(\phi,u^s) \leq {\widetilde{\mathcal{J}}}(\phi_s,u^s)
    \end{equation}
    which implies
    \begin{equation}
      V(\xi) \leq  \mathcal{J}(\xi,u^s).
        \label{BoundCostNOz}
      \end{equation}
    Notice that it is not possible to guarantee that $V(\xi)$ is a lower bound for any arbitrary $(\phi, u^s)\in \mathcal{S}_{\HS\IfIncd{}{}}(\xi)$.
    }
  }{}
  }
  \\\\
  Proceeding with item c as above, pick any $(\phi_w,u^w) \in \mathcal{S}_{\HS\IfIncd{}{}  }^w (\xi)$ {with $u^w:=(u^w_C,u^w_D)$ where $\mathcal{S}_{\HS\IfIncd{}{}  }^w (\xi):= \{(\phi,(u_{1},u_{2}))\in \mathcal{S}_{\HS\IfIncd{}{}  } (\xi):\dom \phi \ni (t,j) \mapsto u_1(t,j)= \IfIh{\kappa}{\gamma}_1(\IfIh{}{t,j,}\phi(t,j))$, 
   and $\dom \phi \ni (t,j) \mapsto u_2(t,j)= \bar{\IfIh{\kappa}{\gamma}}_2(\IfIh{}{t,j,}\phi(t,j))$,  $\bar{\IfIh{\kappa}{\gamma}}_2 \in \mathcal{K}_2\}$.}
  Since $\bar{\IfIh{\kappa}{\gamma}}_2$ does not necessarily attain the upper bound in 
  \eqref{HJBzsihincub}, then, 
  for each $j \in \nats$ such that $I_{\phi_w}^j=[t_j, t_{j+1}]$ has a nonempty interior, we have that for almost all $t \in I_{\phi_w}^j$,
  \begin{equation*}
\begin{split}
  0
  \geq L_C(\phi_w(t,j),u_{ C}^w(t,j))
  +
  \IfGG{{\dot{V}({\phi_w(t,j),u_C^w(t,j)})}}
  {\left\langle \nabla V(\phi_w(t,j)),f \right\rangle \hspace{1.8cm}
  \\ 
   \forall f \in F(\phi_w(t,j),u_C^w(t,j))}
\end{split}
  \end{equation*}
  which implies
  \begin{equation} 
  0 \geq L_C(\phi_w(t,j),u_{ C}^w(t,j))+\frac{d V}{d t}(\phi_w(t,j))\IfAutom{.}{,}
  \label{eq:PrHamiltonC3}
  \end{equation}
  \NotConf{ and 
  by integrating over the interval $[t_j,t_{j+1}]$, 
  we obtain
  \begin{eqnarray*}
  0 \geq \int_{t_{j}}^{t_{j+1}} \left( L_C(\phi_w(t,j),u_{ C}^w(t,j))
  + \frac{d V}{d t}(\phi_w(t,j)) \right) dt
  \end{eqnarray*}
  from which we have
  \begin{eqnarray*}
  V(\phi_w(t_j,j))
  \geq \int_{t_{j}}^{t_{j+1}} L_C(\phi_w(t,j),u_{ C}^w(t,j))dt
  +V(\phi_w(t_{j+1},j))
  \end{eqnarray*}
  Summing both sides from $j=0$ to $j=J_{\phi_w}$, 
  we obtain 
  \begin{eqnarray}
  \sum_{j=0}^{J_{\phi_w}}V(\phi_w(t_j,j)) 
  \geq \sum_{j=0}^{J_{\phi_w}} \int_{t_{j}}^{t_{j+1}} L_C(\phi_w(t,j),u_{C}^w(t,j))dt
  +\sum_{j=0}^{J_{\phi_w}}V(\phi_w(t_{j+1},j))
  \nonumber
  \end{eqnarray}
  Then, solving for $V$ at the initial condition $\phi_w(0,0)$, we obtain
  \begin{equation}\label{ContinCost2goNOihzwinc}
  \begin{split} 
    V(\phi_w(0,0))
  \geq 
  &\sum_{j=0}^{J_{\phi_w}} \int_{t_{j}}^{t_{j+1}} L_C(\phi_w(t,j),u_{C}^w(t,j))dt
  \\&
  + V(\phi_w(t_1,0))+\sum_{j=1}^{J_{\phi_w}} \left( V(\phi_w(t_{j+1},j))-V(\phi_w(t_j,j)) \right) 
  \end{split}
  \end{equation}
  }
In addition, since $\bar{\IfIh{\kappa}{\gamma}}_2$ does not necessarily attain the upper bound in \eqref{Bellmanzsihincub}, then, 
for every $(t,j) \in \dom \phi_w $ such that $(t,j+1) \in \dom \phi_w$, we have
   \begin{equation*}
\begin{split}
  V(\phi_w(t,j))
  \geq 
  L_D(\phi_w(t,j),u_{ D}^w(t,j))
  + V(\chi)
  \hspace{1cm} 
   \quad \forall \chi \in G(\phi_w(t,j),u_{D}^w(t,j))
    \end{split}
  \end{equation*}
  and 
     \begin{equation}
\begin{split}
    V(\phi_w(t,j))
    \geq&
     L_D(\phi_w(t,j),u_{ D}^w(t,j))
    + V( \phi_w(t,j+1))
      \end{split}
   \label{eq:PrHamiltonD3}
  \end{equation}
  \NotConf{Summing both sides from $j=0$ to $j=J_{\phi_w}-1$, 
  we obtain
  \begin{equation}
  \begin{split} 
  \sum_{j=0}^{J_{\phi_w}-1}V(\phi_w(t,j)) 
  \geq 
  \sum_{j=0}^{J_{\phi_w}-1}L_D(\phi_w(t,j),u_{ D}^w(t,j))\nonumber
  +\sum_{j=0}^{J_{\phi_w}-1} V(\phi_w(t,j+1))
  \end{split}
  \end{equation}
  %
  Then, solving for $V$ at the first jump time, we obtain
  %
  \begin{eqnarray} 
  V(\phi_w(t_{1},0))
  \geq V( \phi_w(t_{1},1))
  \label{DiscreteCost2goNOihzwinc}
  &&
  + \sum_{j=0}^{J_{\phi_w}-1}L_D(\phi_w(t,j),u_{ D}^w(t,j)) 
  \\&&
  +\sum_{j=1}^{J_{\phi_w}-1}\left( V(\phi_w(t,j+1))-V(\phi_w(t,j)) \right) \nonumber
  \end{eqnarray}
  In addition, given that $\phi_w(0,0)=\xi$, lower bounding $V(\phi_w(t_{1},0))$ in (\ref{ContinCost2goNOihzwinc}) by the right-hand side of (\ref{DiscreteCost2goNOihzwinc}), we obtain 
  \begin{equation*}
  \begin{split} 
  V(\xi)
  \geq& \sum_{j=0}^{J_{\phi_w}} \int_{t_{j}}^{t} 
  \hspace{-0.1cm}
   L_C(\phi_w(t,j),u_{C}^w(t,j))dt+ V(\phi_w(t_1,0)) 
  \hspace{-0.06cm}
 + \hspace{-0.15cm}
  \sum_{j=1}^{J_{\phi_w}} \hspace{-0.1cm}
  \left( V(\phi_w(t,j))-V(\phi_w(t_j,j)) \right)  \\
  \geq&
  \sum_{j=0}^{J_{\phi_w}} \int_{t_{j}}^{t} L_C(\phi_w(t,j),u_{C}^w(t,j))dt
   +\sum_{j=1}^{J_{\phi_w}-1}\left( V(\phi_w(t,j+1))-V(\phi_w(t,j)) \right) 
  \\&
  + \hspace{-0.2cm}
  \sum_{j=0}^{J_{\phi_w}-1}
  \hspace{-0.1cm}
  L_D(\phi_w(t,j),u_{D}^w(t,j))
  + V(\phi_w(t_{1},1))+\hspace{-0.1cm}
  \sum_{j=1}^{J_{\phi_w}} \left( V(\phi_w(t,j))-V(\phi_w(t_j,j)) \right) 
  \end{split}
  \end{equation*}
  Since 
  \begin{eqnarray*} \nonumber
  && V(\phi_w(t_{1},1))
  +\sum_{j=1}^{J_{\phi_w}-1}\left( V(\phi_w(t,j+1))-V(\phi_w(t,j)) \right) \nonumber 
  +\sum_{j=1}^{J_{\phi_w}} \left( V(\phi_w(t,j))-V(\phi_w(t_j,j)) \right) \nonumber
  \\ &&= V(\phi_w(t_{{J_{\phi_w}}+1},J_{\phi_w})) +V(\phi_w(t_{1},1))\nonumber
  +\sum_{j=1}^{J_{\phi_w}-1}\left( V(\phi_w(t,j+1)) \right) \nonumber 
  -\sum_{j=1}^{J_{\phi_w}} \left(V(\phi_w(t_j,j)) \right) \nonumber
  \\ &&=V(\phi_w(t_{{J_{\phi_w}}+1},J_{\phi_w})) 
  \end{eqnarray*}
  then we have
  \begin{eqnarray*}\nonumber
  V(\xi)
  &\geq&\sum_{j=0}^{J_{\phi_w}} 
  \hspace{-0.1cm}
  \int_{t_{j}}^{t}
  \hspace{-0.2cm}
  L_C(\phi_w(t,j),u_{C}^w(t,j))dt
   \hspace{-0.1cm}
 + \hspace{-0.3cm}
   \sum_{j=0}^{J_{\phi_w}-1}
    \hspace{-0.2cm}
   L_D(\phi_w(t,j),u_{D}^w(t,j))
  +V(\phi_w(t_{J_{\phi_w}+1},J_{\phi_w}))
  \nonumber
  \end{eqnarray*}
  By taking the limit when $(t_{J_{\phi_w}+1},J_{\phi_w}) \rightarrow \sup \dom \phi_w$, and thanks to (\ref{TerminalCondCGinc}), we have
  \begin{eqnarray} \nonumber
  V(\xi)&\geq&\sum_{j=0}^{\sup_j \dom \phi_w} \int_{t_{j}}^{t_{j+1}} L_C(\phi_w(t,j),u_{C}^w(t,j))dt
  \\&&
  + \sum_{j=0}^{{\sup_j \dom \phi_w}-1}L_D(\phi_w(t_{j+1},j),u_{D}^w(t_{j+1},j))
  +\underset{(t,j) \in \textup{dom}\phi_w}{\limsup_{t+j\rightarrow \infty} } V(\phi_w(t,j))
  \nonumber \\ 
  &=&\widetilde\J(\phi_w,u^w)
  \label{BoundCostNOzwinc}
  \end{eqnarray}
  }
  \IfConf{Now, 
  thanks to \eqref{eq:PrHamiltonC3} and \eqref{eq:PrHamiltonD3}, from Proposition \ref{Pp:ProofDerivation} and Corollary \ref{Rmk:ProofDerivation}, we have that 
  \begin{equation}
  V(\xi) \geq  {\widetilde{\mathcal{J}}}(\phi_w,u^w).
    \label{BoundCostNOzwinc}
  \end{equation}
  }{}
  Finally, 
  by proceeding as in item d above,
  \pno{by applying the infimum on (\ref{BoundCostNOz}) 
  over the set $\mathcal{S}_{\HS_{\max}}^s (\xi)$, we obtain
 \begin{equation}\label{a}
V(\xi) \leq  
\sj{    \underset{(\phi_s, u^s) \in {\mathcal{S}_{\mathcal{H}_{\max}}^s}(\xi)
}{\text{inf}}\widetilde{\mathcal{J}}(\phi_s, u^s) 
=}
\sj{    \underset{(\phi_s, \sj{u^s}) \in {\mathcal{S}_{\mathcal{H}_{\max}}^s}(\xi)
}{\text{inf}}\mathcal{J}(\xi, u^s) }
{=: \overline V(\xi)}
\end{equation}}
\sj{Notice that the infimum in \eqref{a} is attained in ${\mathcal{S}_{\mathcal{H}_{\max}}^s}(\xi)$ if there exists $(\phi_s, \sj{u^s}) \in {\mathcal{S}_{\mathcal{H}_{\max}}^s}(\xi)$ such that $\mathcal{J}(\xi, \sj{u^s}) = V(\xi)$.}
  %
  By applying the supremum on (\ref{BoundCostNOzwinc})
  over the set $
  \mathcal{S}_{\HS\IfIncd{}{} }^w (\xi)$, we obtain
\begin{equation}\label{b}
V(\xi) \geq 
\sj{\underset{(\phi_w, \sj{u^w}) \in {\mathcal{S}_{\mathcal{H}}^w}(\xi)}{\text{sup}}\widetilde{\mathcal{J}}(\phi_w, \sj{u^w}) }
{=: \underline{ V}(\xi)}.
\end{equation}
\sj{Notice that the supremum in \eqref{b} is attained in ${\mathcal{S}_{\mathcal{H}}^w}(\xi)$ if there exists $(\phi_w, \sj{u^w}) \in {\mathcal{S}_{\mathcal{H}}^w}(\xi)$ such that $\widetilde{\mathcal{J}}(\phi_w, \sj{u^w}) = V(\xi)$, at which $\underline{ V}(\xi)=V(\xi)$.}
  Given that $V(\xi)=\mathcal{J}(\xi, u^*)$ 
  from (\ref{OptimalCostihzinc}), we have that for any $\xi \in \Pi ( \overline{C}) \cup \Pi(D)$, each 
  $(\phi^*,u^*)\in \mathcal{S}_{\pnn{\HS_{\max} } }(\xi)$ with $u^*=(\kappa_1(\phi^*),\kappa_2(\phi^*))$ satisfies
  \begin{equation}
    \underline{V}(\xi)
  \leq
  \mathcal{J}(\xi,u^*) 
  \leq 
  \overline{V}(\xi)
   \label{Boundsszinc}
  \end{equation}

  Thanks to 
 $(\phi^*,u^*) \in \mathcal{S}_{\HS_{\max}}^s(\xi) \cap \mathcal{S}_{\HS}^w(\xi)
 $, we have
 \begin{subequations}
\begin{eqnarray}
\hspace{0.6cm}&&\overline{V}(\xi)= \hspace{-0.2cm}
\underset{(\phi_s, \sj{u^s}) \in {\mathcal{S}_{\mathcal{H}_{\max}}^s}
(\xi)
}{\text{inf}}\widetilde{\mathcal{J}}(\phi_s, \sj{u^s}) =
\hspace{-0.4cm}
\pn{\underset{(\phi_s, \sj{u^s}) \in {\mathcal{S}_{\mathcal{H}_{\max}}^s}(\xi)
}{\text{inf}}}\mathcal{J}(\xi, u^s)=\mathcal{J}(\xi, u^*)= V(\xi)
\label{supJk}
\\
  && \underline{V}(\xi) = \underset{(\phi_w, \sj{u^w}) \in {\mathcal{S}_{\mathcal{H}}^w}
    (\xi)
    }{\text{sup}}\widetilde{\mathcal{J}}(\phi_w, \sj{u^w})=\mathcal{J}(\xi, u^w) =\mathcal{J}(\xi, u^*) = V(\xi)
\label{infJk}
\end{eqnarray}
\end{subequations}
Since the infimum and supremum are attained in 
\eqref{supJk} and \eqref{infJk}
, respectively, \sj{by $(\phi^*,u^*)$,} \eqref{Boundsszinc}
leads to
\begin{equation}
\mathcal{J}(\xi,u^*) 
=\underset{(u_1,u_2) \in {\mathcal{U}_{\HS}}
(\xi)}
{ \min_{u_{1}} \max_{u_{2}}}
\mathcal{J}(\xi,(u_1,u_2))
\label{LowerBounddz}
\end{equation}
Thus, from (\ref{OptimalCostihzinc}) and (\ref{LowerBounddz}), $V(\xi)$ is the value function for $\HS$, as in Definition \ref{ValueFunctionz}, and it is attained by 
the strategy $\kappa=(\kappa_1,\kappa_2)$. 
\end{proof}
}
{
\begin{remark}{Weak optimality of the saddle-point equilibrium}
  When both players play the \pnn{weak} saddle-point equilibrium strategy, due to nonuniqueness of solutions,   
  there is no reason to assume that the worst-case cost is attained by the resulting state trajectory, implying that such a strategy is not necessarily optimal in the min-max sense \pnn{with respect to the cost in \eqref{defJTNCinc}.} Nevertheless, by playing the \pnn{weak} saddle-point equilibrium, the worst-case cost {as in \eqref{defJTNCincinc}} is minimized under the adversarial action that aims to maximize it. \pnn{ We illustrate this in the examples in Section 7.}
\end{remark}}
}{}
\IfIncd{}{\subsection{
Hamilton-Jacobi-Isaacs Equations for Two-player Zero-sum Infinite-horizon Hybrid Games}

The following result provides sufficient conditions to characterize the value function $\mathcal{J}^*$, and the feedback law that attains it. It addresses the solution to Problem $(\diamond)$ for each $\xi \in \Pi ( \overline{C}) \cup \Pi(D)$ showing that the optimizer is the saddle-point equilibrium.
\begin{ttheorem}{Hamilton-Jacobi-Isaacs (HJI) for Problem ($\diamond$)}
Given a two-player zero-sum hybrid game with dynamics $\HS$ as in (\ref{Heq}) \IfTp{}{with $N=2$, }with data $(C,F,D,G)$, {satisfying Assumption \ref{AssLipsZ}},
 stage costs $L_C:\mathbb{R}^n  \times \mathbb{R}^{m_C} \rightarrow \mathbb{R}_{\geq 0}$ {and} $L_D:\mathbb{R}^n  \times \mathbb{R}^{m_D} \rightarrow \mathbb{R}_{\geq 0}$, and terminal cost $q: \mathbb{R}^n \rightarrow \mathbb{R}$, {suppose the following hold:} 
 \begin{enumerate}
   \item There exists a function $V: \mathbb{R}^n \rightarrow \mathbb{R}$ that is continuously differentiable on a neighborhood of\> $\Pi(C)$ that satisfies the \textit{Hamilton-Jacobi-Isaacs hybrid equations given as} 
{
\begin{equation}
\begin{split}
0&= \hspace{-0.6cm}  \underset{u_C=(u_{C1},u_{C2}) \in \Pi_u(x,C)}
{\min_{u_{C1}} \max_{u_{C2}}}  
{\mathcal{L}_C(x,u_C)}
\\&= \hspace{-0.6cm} \underset{u_C=(u_{C1},u_{C2}) \in \Pi_u(x,C)}
{ \max_{u_{C2}}\min_{u_{C1}}}
{\mathcal{L}_C(x,u_C)}
\hspace{0.4cm}\forall x \in \Pi(C),  
\end{split}
\label{HJBzsih}
\end{equation}}
where $\mathcal{L}_C(x,u_C):=  L_C(x,u_C)+\left\langle \nabla V(x),F(x,u_C) \right\rangle$,
\begin{equation}
\begin{split}
V(x)&=\hspace{-0.6cm}\underset{u_D=(u_{D1},u_{D2}) \in \Pi_u(x,D)}
{\min_{u_{D1}} \max_{u_{D2}}}  
{\mathcal{L}_D(x,u_D)}
\\&=\hspace{-0.6cm}\underset{u_D=(u_{D1},u_{D2}) \in \Pi_u(x,D)}
{ \max_{u_{D2}}\min_{u_{D1}}}  
{\mathcal{L}_D(x,u_D)}
\hspace{0.5cm} \forall x \in \Pi(D),
\end{split}
\label{Bellmanzsih}
\end{equation}
where $\mathcal{L}_D(x,u_D):=  L_D(x,u_D) + V(G(x,u_D))$.
\item For each $\xi \in \Pi ( \overline{C}) \cup \Pi(D)$, each $(\phi, u)\in \mathcal{S}_{\HS  }^\infty(\xi)$ satisfies
\begin{equation}
\underset{(t,j) \in \textup{dom}\phi}{\limsup_{t+j\rightarrow \infty} } V(\phi(t,j))=\underset{(t,j) \in \textup{dom}\phi}{\limsup_{t+j\rightarrow \infty} } q(\phi(t,j)){.}
\label{TerminalCondCG}
\end{equation}
\end{enumerate}
Then   
\begin{equation}
\J^*(\xi)= V(\xi) \hspace{2cm} \forall \xi \in \Pi ( \overline{C}) \cup \Pi(D),
\label{ResultValue}
\end{equation}
and any 
feedback law
$\kappa:=(\kappa_C,\kappa_D)\NotConf{=((\kappa_{C1},\kappa_{C2}),(\kappa_{D1},\kappa_{D2}))}: \mathbb{R}^n \rightarrow \mathbb{R}^{m_C} \times \mathbb{R}^{m_D}$  with values 
{
\begin{equation}\begin{split}
\kappa_{C}(x)\in \IfConf{\hspace{-0.6cm}}{}  \underset{u_C=(u_{C1},u_{C2}) \in \Pi_u(x,C)}
{\arg \min_{u_{C1}} \max_{u_{C2}}} 
{\mathcal{L}_C(x,u_C)}
\hspace{1cm}
\forall x \in \Pi(C)
\end{split}
\label{kHJBeqzsihc}
\end{equation}}
and
{
\begin{equation}\begin{split}
\kappa_{D}(x)\in \IfConf{\hspace{-0.6cm}}{} \underset{u_D=(u_{D1},u_{D2}) \in \Pi_u(x,D)}
{\arg \min_{u_{D1}} \max_{u_{D2}}} 
{\mathcal{L}_D(x,u_D)}
\hspace{1cm}
\forall x \in \Pi(D) 
\end{split}
\label{kBeqzsihc}
\end{equation}
}
is a pure strategy saddle-point equilibrium for the two-player zero-sum hybrid game with {infinite horizon and} $\mathcal{J}_1=\mathcal{J}$, $\mathcal{J}_2=-\mathcal{J}$.

\label{thHJBszih} 
\end{ttheorem}

\NotConf{\begin{proofs}
To show the claim we \pno{apply cost evaluation tools built upon dynamic programming approaches} and proceed as follows:
\begin{enumerate}
\item Pick an initial condition $\xi$ and evaluate the cost associated to any 
solution yielded by $\IfIh{\kappa=(\kappa_C,\kappa_D)}{\gamma}$, with values as in (\ref{kHJBeqzsihc}) and (\ref{kBeqzsihc}), from $\xi$. Show that this cost coincides with the value of the function $V$ at $\xi$. 
\item Lower bound the cost associated to any solution from $\xi$ when $P_{2}$ plays $\IfIh{\kappa}{\gamma}_2:=(\IfIh{\kappa}{\gamma}_{C2},\IfIh{\kappa}{\gamma}_{D2})$ by the value of the function $V$ evaluated at $\xi$. 
\item Upper bound the cost associated to any solution from $\xi$ when $P_{1}$ plays $\IfIh{\kappa}{\gamma}_1:=(\IfIh{\kappa}{\gamma}_{C1},\IfIh{\kappa}{\gamma}_{D1})$ by the value of the function $V$ evaluated at $\xi$. 
\item By showing that the cost of any solution from $\xi$ when $P_{1}$ plays $\IfIh{\kappa}{\gamma}_1$ is not less than the cost of any solution yielded by $\IfIh{\kappa}{\gamma}$ from $\xi$,
and by showing that the cost of any solution from $\xi$ when $P_{2}$ plays $\IfIh{\kappa}{\gamma}_2$ is not larger than the cost of any solution yielded by $\IfIh{\kappa}{\gamma}$ from $\xi$, we show optimality of $\IfIh{\kappa}{\gamma}$ in Problem ($\diamond$) in the min-max sense. 
\end{enumerate}
\end{proofs}}
\IfIncd{The proof is presented in the Appendix.}{}
Notice that when the players select the optimal strategy, the value function equals the function $V$ evaluated at the initial condition. This makes evident the independence of the result from {needing to compute} solutions, at the price of {finding} the function $V$ satisfying the conditions therein.
{
\subsection{Proof of Theorem \ref{thHJBszih}}
To show the claim we proceed as follows:
\begin{enumerate}
\item Pick an initial condition $\xi$ and evaluate the cost associated to any 
{complete} solution yield by $\IfIh{\kappa=(\kappa_C,\kappa_D)}{\gamma}$, with values as in (\ref{kHJBeqzsihc}) and (\ref{kBeqzsihc}), from $\xi$. Show that this cost coincides with the value of the function $V$ at $\xi$. 
\item Lower bound the cost associated to any {complete} solution from $\xi$ when $P_{2}$ plays $\IfIh{\kappa}{\gamma}_2:=(\IfIh{\kappa}{\gamma}_{C2},\IfIh{\kappa}{\gamma}_{D2})$ by the value of the function $V$ evaluated at $\xi$. 
\item Upper bound the cost associated to any {complete} solution from $\xi$ when $P_{1}$ plays $\IfIh{\kappa}{\gamma}_1:=(\IfIh{\kappa}{\gamma}_{C1},\IfIh{\kappa}{\gamma}_{D1})$ by the value of the function $V$ evaluated at $\xi$. 
\item By showing that the cost of any {complete} solution from $\xi$ when $P_{1}$ plays $\IfIh{\kappa}{\gamma}_1$ is not {greater} than the cost of any solution yield by $\IfIh{\kappa}{\gamma}$ from $\xi$,
and by showing that the cost of any solution from $\xi$ when $P_{2}$ plays $\IfIh{\kappa}{\gamma}_2$ is not {less} than the cost of any solution yield by $\IfIh{\kappa}{\gamma}$ from $\xi$, we show optimality of $\IfIh{\kappa}{\gamma}$ in Problem ($\diamond$) in the min-max sense. 
\end{enumerate}
Proceeding as in item {a} above, pick any $\xi \in \Pi ( \overline{C}) \cup \Pi(D)$ and any {$(\phi^*,u^*) \in \mathcal{S}_{\HS  }
 (\xi)$} with $\dom \phi^* \ni (t,j) \mapsto u^*(t,j)= \IfIh{\kappa}{\gamma}(\IfIh{}{t,j,}\phi^*(t,j))$. 
\IfPers{and $(t_{J_{\phi^*}+1},J_{\phi^*}) \in \dom (\phi^*,u^*) $} 
We show that the cost of $(\phi^*,u^*)$ is optimal in the min-max sense.
%
Given that $V$ satisfies 
(\ref{HJBzsih}), and  $\IfIh{\kappa}{\gamma}_C$ 
is as in (\ref{kHJBeqzsihc}), for each $j \in \nats$ such that $I_{\phi^*}^j=[t_j,t_{j+1}]$ has a nonempty interior int$I_{\phi^*}^j$, we have, for all $t \in \textup{int}I_{\phi^*}^j$, 
{
\begin{equation*}
\begin{split}
0
\IfPers{=&  \underset{u_{ C}(t,j)=(u_{C1},u_{C2}) \in \Pi_u(\phi^*(t,j),C)}
{\min_{u_{C1}} \max_{u_{C2}}}  
\left\{ 
L_C(\phi^*(t,j),u_{ C}(t,j))
\right.
\\& \hspace{2.8cm}
\left.
+\left\langle \nabla V(\phi^*(t,j)),F(\phi^*(t,j),u_C(t,j)) \right\rangle \right\}\\ 
}
=&  L_C(\phi^*(t,j),\IfIh{\kappa}{\gamma}_{ C}(\phi^*(t,j)))
\\&\hspace{1.3cm}+\left\langle \nabla V(\phi^*(t,j)), F(\phi^*(t,j),\IfIh{\kappa}{\gamma}_C(\phi^*(t,j))) \right\rangle
\end{split}
\end{equation*}}
and $\phi^*(t,j) \in C_{\IfIh{\kappa}{\gamma}}$, as in (\ref{Hkeq}). Given that $V$ is continuously differentiable on a neighborhood of $\Pi(C)$, we can express its total derivative along $\phi^*$ as
\begin{equation}
\frac{d V}{d t}(\phi^*(t,j))=\left\langle \nabla V(\phi^*(t,j)), F(\phi^*(t,j),\IfIh{\kappa}{\gamma}_C(\phi^*(t,j))) \right\rangle
\label{eq:PrHamiltonC4}
\end{equation}
for every $(t,j)$ $\in \textup{int}(I^j_\phi) \times \{j\}$ with $\textup{int}(I^j_\phi)$ 
nonempty\IfAutom{.}{, and 
by integrating over the interval $[t_j,t_{j+1}]$, 
we obtain
\begin{eqnarray*}
0=  \int_{t_{j}}^{t_{j+1}} \left( L_C(\phi^*(t,j),\IfIh{\kappa}{\gamma}_{ C}(\phi^*(t,j)))
+ \frac{d V}{d t}(\phi^*(t,j)) \right )dt
\end{eqnarray*}
from where we have
\begin{multline*}
0=  \int_{t_{j}}^{t_{j+1}} L_C(\phi^*(t,j),\IfIh{\kappa}{\gamma}_{ C}(\phi^*(t,j))) dt
\\+  V(\phi^*(t_{j+1},j))- V(\phi^*(t_j,j))
\end{multline*}
%
%
%
Summing from $j=0$ to $j=J_{\phi^*}$, 
we obtain 
\begin{multline*}
0
=\sum_{j=0}^{J_{\phi^*}}  \int_{t_{j}}^{t_{j+1}} L_C(\phi^*(t,j),\IfIh{\kappa}{\gamma}_{ C}(\phi^*(t,j))) dt
\\+\sum_{j=0}^{J_{\phi^*}} \left( V(\phi^*(t_{j+1},j))-V(\phi^*(t_j,j)) \right)
\end{multline*}
Then, solving for $V$ at the initial condition $\phi^*(0,0)$, we obtain
%
%
\begin{equation} \label{ContinCost2goihz}
\begin{medsize}
\begin{split} 
V(\phi^*(0,0))
=\sum_{j=0}^{J_{\phi^*}}  \int_{t_{j}}^{t_{j+1}} L_C(\phi^*(t,j),\IfIh{\kappa}{\gamma}_{ C}(\phi^*(t,j))) dt
\\  + V(\phi^*(t_1,0))
+\sum_{j=1}^{J_{\phi^*}} \left( V(\phi^*(t_{j+1},j))-V(\phi^*(t_j,j)) \right) 
\end{split}
\end{medsize}
\end{equation}
}
Given that $V$ satisfies 
(\ref{Bellmanzsih}) and $\IfIh{\kappa}{\gamma}_D$ 
is as in (\ref{kBeqzsihc}), for every $(t_{j+1},j) \in \dom \phi^*$ such that $(t_{j+1},j+1) \in \dom \phi^*$, we have that
{
\begin{equation}
\begin{split} 
V(\phi^*(t_{j+1},j))
\IfPers{
=&\underset{u_{D}(t_{j+1},j)=(u_{D1},u_{D2}) \in \Pi_u(\phi^*(t_{j+1},j),D)}
{\min_{u_{D1}} \max_{u_{D2}}}  
\\&
\left\{ 
L_D(\phi^*(t_{j+1},j),u_{ D}(t_{j+1},j))
\right. \\ & \left. \hspace{2cm} 
+V(G(\phi^*(t_{j+1},j), u_D(t_{j+1},j) \right\}\\ 
}
=& \> 
L_D(\phi^*(t_{j+1},j),\IfIh{\kappa}{\gamma}_{ D}(\phi^*(t_{j+1},j)))
\\ &  
+ V(G(\phi^*(t_{j+1},j), \IfIh{\kappa}{\gamma}_D(\phi^*(t_{j+1},j)))) \\
=&\> L_D(\phi^*(t_{j+1},j),\IfIh{\kappa}{\gamma}_{ D}(\phi^*(t_{j+1},j)))
\\ &  \hspace{2.4cm} 
+ V(\phi^*(t_{j+1},j+1)) 
\end{split}
\label{eq:PrHamiltonD4}
\end{equation}}
where $\phi^*(t_{j+1},j)\in D_{\IfIh{\kappa}{\gamma}}$ is defined in (\ref{Hkeq}).
\NotAutom{Summing both sides from $j=0$ to $j=J_{\phi^*}-1$, 
we obtain
\begin{equation*}
\begin{medsize}
\begin{split} 
\sum_{j=0}^{J_{\phi^*}-1} V(\phi^*(t_{j+1},j))
=& \sum_{j=0}^{J_{\phi^*}-1} L_D(\phi^*(t_{j+1},j),\IfIh{\kappa}{\gamma}_{ D}(\phi^*(t_{j+1},j))) 
\\
&+\sum_{j=0}^{J_{\phi^*}-1} V( \phi^*(t_{j+1},j+1))
\end{split}
\end{medsize}
\end{equation*}
Then, solving for $V$ at the first jump time, we obtain
%
\begin{eqnarray} \nonumber
&&V(\phi^*(t_{1},0))
= 
\\ \nonumber
&&V( \phi^*(t_{1},1))
+ \sum_{j=0}^{J_{\phi^*}-1} L_D(\phi^*(t_{j+1},j),\IfIh{\kappa}{\gamma}_{ D}(\phi^*(t_{j+1},j))) 
\\ \label{DiscreteCost2goihz}
&&+\sum_{j=1}^{J_{\phi^*}-1}\left( V( \phi^*(t_{j+1},j+1))-V(\phi^*(t_{j+1},j)) \right)
\end{eqnarray}
Given that $\phi^*(0,0)=\xi$, by substituting the right-hand side of (\ref{DiscreteCost2goihz}) in  (\ref{ContinCost2goihz}), we obtain
\begin{equation} 
\begin{medsize}
\begin{split} 
V(\xi)
=&\sum_{j=0}^{J_{\phi^*}}  \int_{t_{j}}^{t_{j+1}} L_C(\phi^*(t,j),\IfIh{\kappa}{\gamma}_{ C}(\phi^*(t,j))) dt
\nonumber \\&
+ V(\phi^*(t_1,0))\nonumber 
+\sum_{j=1}^{J_{\phi^*}} \left( V(\phi^*(t_{j+1},j))-V(\phi^*(t_j,j)) \right) \nonumber \\
=&\sum_{j=0}^{J_{\phi^*}}  \int_{t_{j}}^{t_{j+1}} L_C(\phi^*(t,j),\IfIh{\kappa}{\gamma}_{ C}(\phi^*(t,j))) dt
\nonumber \\&
+ \sum_{j=0}^{J_{\phi^*}-1}L_D(\phi^*(t_{j+1},j),\IfIh{\kappa}{\gamma}_{ D}(\phi^*(t_{j+1},j)))
+ V(\phi^*(t_{1},1))
\nonumber \\&
+\sum_{j=1}^{J_{\phi^*}-1}\left( V(\phi^*(t_{j+1},j+1))-V(\phi^*(t_{j+1},j)) \right) 
\nonumber \\&
+\sum_{j=1}^{J_{\phi^*}} \left( V(\phi^*(t_{j+1},j))-V(\phi^*(t_j,j)) \right) \nonumber 
\end{split}
\end{medsize}
\end{equation}
Since 
\begin{eqnarray}
&& V(\phi^*(t_{1},1))
\nonumber \\&&
+\sum_{j=1}^{J_{\phi^*}-1}\left( V(\phi^*(t_{j+1},j+1))-V(\phi^*(t_{j+1},j)) \right) 
\nonumber \\&&
+\sum_{j=1}^{J_{\phi^*}} \left( V(\phi^*(t_{j+1},j))-V(\phi^*(t_j,j)) \right) \nonumber
\\ &=& V(\phi^*(t_{J_{\phi^*}+1},J_{\phi^*}))+V(\phi^*(t_{1},1))
\nonumber \\&&
+\sum_{j=1}^{J_{\phi^*}-1}\left( V(\phi^*(t_{j+1},j+1)) \right) \nonumber 
-\sum_{j=1}^{J_{\phi^*}} \left(V(\phi^*(t_j,j)) \right) \nonumber
\\ &=&V(\phi^*(t_{J_{\phi^*}+1},J_{\phi^*})) 
\end{eqnarray}
then we have
\begin{eqnarray}
V(\xi)
&=&\sum_{j=0}^{J_{\phi^*}}  \int_{t_{j}}^{t_{j+1}} L_C(\phi^*(t,j),\IfIh{\kappa}{\gamma}_{ C}(\phi^*(t,j))) dt
\nonumber \\&&
+ \sum_{j=0}^{J_{\phi^*}-1}L_D(\phi^*(t_{j+1},j),\IfIh{\kappa}{\gamma}_{ D}(\phi^*(t_{j+1},j))) 
\\&&+V(\phi^*(t_{J_{\phi^*}+1},J_{\phi^*})) \nonumber
\end{eqnarray}
By taking the limit when $(t_{J_{\phi^*}+1},J_{\phi^*}) \rightarrow \sup \dom \phi^*$, and given that (\ref{TerminalCondCG}) holds, we have 
\begin{eqnarray}
V(\xi)&=& \sum_{j=0}^{\sup_j \dom \phi^*}   \int_{t_{j}}^{t_{j+1}} L_C(\phi^*(t,j),\IfIh{\kappa}{\gamma}_{ C}(\phi^*(t,j))) dt
\nonumber \\&&
+ \sum_{j=0}^{\sup_j \dom \phi^* -1} L_D(\phi^*(t_{j+1},j),\IfIh{\kappa}{\gamma}_{ D}(\phi^*(t_{j+1},j))) 
\nonumber
\\ && 
+\underset{(t,j) \in \textup{dom}\phi^*}{\limsup_{t+j\rightarrow \infty} } V(\phi^*(t,j))
\nonumber \\ 
&=& \mathcal{J}(\xi,u^*)
\label{OptimalCostihz}
\end{eqnarray}
}
 \IfAutom{{Now, when $(\phi^*,u^*)$ is complete
  with $\dom \phi^* \ni (t,j) \mapsto u^*(t,j)= \IfIh{\kappa}{\gamma}(\IfIh{}{t,j,}\phi^*(t,j))$,
 thanks to \eqref{eq:PrHamiltonC4} and \eqref{eq:PrHamiltonD4}, }
 from Proposition \ref{Pp:ProofDerivation}, we have that 
  \begin{equation}
  V(\xi) = \mathcal{J}(\xi,u^*).
    \label{OptimalCostihz}
  \end{equation}
  }{}
\\\\
%
Continuing with item {b} as above, pick any 
%
$(\phi_s,u^s) \in \mathcal{S}_{\HS  }^s(\xi)$ where
$\mathcal{S}_{\HS  }^s (\xi) (\subset {\mathcal{S}_{\HS}(\xi)})$ is the set of solutions $(\phi,u)$ with $u=(u_1,u_2)$, $\dom \phi \ni (t,j) \mapsto u_1(t,j)= \bar{\IfIh{\kappa}{\gamma}}_1(\IfIh{}{t,j,}\phi(t,j))$ for some $\bar{\IfIh{\kappa}{\gamma}}_1 \in \mathcal{K}_1$, $\dom \phi \ni (t,j) \mapsto u_2(t,j)= \IfIh{\kappa}{\gamma}_2(t,j,\phi(t,j))$ for ${\IfIh{\kappa}{\gamma}}_2:=(\IfIh{\kappa}{\gamma}_{C2}, \IfIh{\kappa}{\gamma}_{D2}) $ as in (\ref{kHJBeqzsihc}) and (\ref{kBeqzsihc}). 
Since $\bar{\IfIh{\kappa}{\gamma}}_1$ does not necessarily attain the minimum in (\ref{HJBzsih}), then,  for each $j \in \nats$ such that $I_{\phi_s}^j=[t_j, t_{j+1}]$ has a nonempty interior int$I_{\phi_s}^j$, we have that for every $t \in \textup{int}I_{\phi_s}^j$,
{
\begin{equation*}
\begin{split}
0&
\leq  L_C(\phi_s(t,j),u_{ C}^s(t,j))
\\& \hspace{2cm}
+\left\langle \nabla V(\phi_s(t,j)),F(\phi_s(t,j),u_C^s(t,j)) \right\rangle.
\end{split}
\end{equation*}}
%
{Similarly to \eqref{eq:PrHamiltonC4}, we have}
\begin{equation} 
\begin{split}
\frac{d V}{d t}(\phi_s(t,j)):=\left\langle \nabla V(\phi_s(t,j)),F(\phi_s(t,j),u_C^s(t,j)) \right\rangle
\end{split}
\label{eq:PrHamiltonC5}
\end{equation}
 for every $(t,j) \in \textup{int}(I^j_{\phi_s}) \times \{j\}$ with $\textup{int}(I^j_{\phi_s})$ 
 nonempty{\IfAutom{.}{, and
by integrating over the interval $[t_j,t_{j+1}]$,  
we obtain
\begin{eqnarray*}
0 \leq \int_{t_{j}}^{t_{j+1}}  \left( L_C(\phi_s(t,j),u_{ C}^s(t,j))
+  \frac{d V}{d t}(\phi_s(t,j)) \right) dt
\end{eqnarray*}
from which we have
\begin{equation*}
\begin{medsize}
V(\phi_s(t_j,j))
\leq \int_{t_{j}}^{t_{j+1}}   L_C(\phi_s(t,j),u_{ C}^s(t,j))dt
+V(\phi_s(t_{j+1},j))
\end{medsize}
\end{equation*}
Summing both sides from  $j=0$ to $j=J_{\phi_s}$, 
we obtain 
\begin{eqnarray}
\sum_{j=0}^{J_{\phi_s}}V(\phi_s(t_j,j)) 
\leq \sum_{j=0}^{J_{\phi_s}}  \int_{t_{j}}^{t_{j+1}}  L_C(\phi_s(t,j),u_{C}^s(t,j))dt
\nonumber \\
+\sum_{j=0}^{J_{\phi_s}}V(\phi_s(t_{j+1},j))
\nonumber
\end{eqnarray}
Then, solving for $V$ at the initial condition $\phi_s(0,0)$, we obtain
%
%
\begin{equation}\label{ContinCost2goNOihz}
\begin{medsize}
\begin{split} 
&V(\phi_s(0,0))
\leq \sum_{j=0}^{J_{\phi_s}}  \int_{t_{j}}^{t_{j+1}}  L_C(\phi_s(t,j),u_{C}^s(t,j))dt
\\
&+ V(\phi_s(t_1,0))+\sum_{j=1}^{J_{\phi_s}} \left( V(\phi_s(t_{j+1},j))-V(\phi_s(t_j,j)) \right) 
\end{split}
\end{medsize}
\end{equation}
}
In addition, since 
{$\bar{\IfIh{\kappa}{\gamma}}_1$} does not necessarily attain the minimum in (\ref{Bellmanzsih}), then for every $(t_{j+1},j) \in \dom \phi_s $ such that   $(t_{j+1},j+1) \in \dom \phi_s$, we have
{
\begin{equation}
\begin{split} 
V(\phi_s(t_{j+1},j))
\leq&  \NotAutom{
  \> 
L_D(\phi_s(t_{j+1},j),u_{ D}^s(t_{j+1},j))
\\&\hspace{0.7cm}
+ V(G(\phi_s(t_{j+1},j),u_{D}^s(t_{j+1},j)))\\
{\leq}&
}
\>L_D(\phi_s(t_{j+1},j),u_{ D}^s(t_{j+1},j))
\\&\hspace{1.5cm}
+ V( \phi_s(t_{j+1},j+1)).
\end{split}
\label{eq:PrHamiltonD5}
\end{equation}}
\NotAutom{
Summing both sides from  $j=0$ to $j=J_{\phi_s}-1$, 
we obtain
\begin{eqnarray}
\sum_{j=0}^{J_{\phi_s}-1}V(\phi_s(t_{j+1},j)) 
\leq 
\sum_{j=0}^{J_{\phi_s}-1}L_D(\phi_s(t_{j+1},j),u_{ D}^s(t_{j+1},j))
\nonumber\\
+\sum_{j=0}^{J_{\phi_s}-1} V(\phi_s(t_{j+1},j+1))
\nonumber
\end{eqnarray}
Then, solving for $V$ at the first jump time, we obtain
%
\begin{equation}\label{DiscreteCost2goNOihz}
\begin{medsize}
\begin{split} 
V(\phi_s(t_{1},0))
\leq V( \phi_s(t_{1},1))
+ \sum_{j=0}^{J_{\phi_s}-1}L_D(\phi_s(t_{j+1},j),u_{ D}^s(t_{j+1},j))
\\
+\sum_{j=1}^{J_{\phi_s}-1}\left( V(\phi_s(t_{j+1},j+1))-V(\phi_s(t_{j+1},j)) \right)
\end{split}
\end{medsize}
\end{equation}
In addition, given that $\phi_s(0,0)=\xi$, upper bounding $V(\phi_s(t_{1},0))$ in (\ref{ContinCost2goNOihz}) by the right-hand side of (\ref{DiscreteCost2goNOihz}), we obtain 
\begin{equation*}
\begin{medsize}
\begin{split}
V(\xi)  
\leq& \sum_{j=0}^{J_{\phi_s}}  \int_{t_{j}}^{t_{j+1}}  L_C(\phi_s(t,j),u_{C}^s(t,j))dt+ V(\phi_s(t_1,0))
\\&
+\sum_{j=1}^{J_{\phi_s}} \left( V(\phi_s(t_{j+1},j))-V(\phi_s(t_j,j)) \right)  \\
\leq&
\sum_{j=0}^{J_{\phi_s}}  \int_{t_{j}}^{t_{j+1}}  L_C(\phi_s(t,j),u_{C}^s(t,j))dt
\\&+ \sum_{j=0}^{J_{\phi_s}-1}L_D(\phi_s(t_{j+1},j),u_{D}^s(t_{j+1},j))
 \\
&+\sum_{j=1}^{J_{\phi_s}-1}\left( V(\phi_s(t_{j+1},j+1))-V(\phi_s(t_{j+1},j)) \right)  
\\&+ V(\phi_s(t_{1},1))+\sum_{j=1}^{J_{\phi_s}} \left( V(\phi_s(t_{j+1},j))-V(\phi_s(t_j,j)) \right) 
\end{split}
\end{medsize}
\end{equation*}
Since 
\begin{eqnarray}
&& V(\phi_s(t_{1},1))\nonumber
\\&&+\sum_{j=1}^{J_{\phi_s}-1}\left( V(\phi_s(t_{j+1},j+1))-V(\phi_s(t_{j+1},j)) \right) \nonumber 
\\&&
+\sum_{j=1}^{J_{\phi_s}} \left( V(\phi_s(t_{j+1},j))-V(\phi_s(t_j,j)) \right) \nonumber
\\ &=& V(\phi_s(t_{{J_{\phi_s}}+1},J_{\phi_s})) +V(\phi_s(t_{1},1))\nonumber
\\&&
+\sum_{j=1}^{J_{\phi_s}-1}\left( V(\phi_s(t_{j+1},j+1)) \right)
-\sum_{j=1}^{J_{\phi_s}} \left(V(\phi_s(t_j,j)) \right) \nonumber
\\ &=&V(\phi_s(t_{{J_{\phi_s}}+1},J_{\phi_s}))
\end{eqnarray}
then we have
\begin{eqnarray}
V(\xi)
&\leq&\sum_{j=0}^{J_{\phi_s}}  \int_{t_{j}}^{t_{j+1}}  L_C(\phi_s(t,j),u_{C}^s(t,j))dt\nonumber
\\&&
+ \sum_{j=0}^{J_{\phi_s}-1}L_D(\phi_s(t_{j+1},j),u_{D}^s(t_{j+1},j))
\\&&+V(\phi_s(t_{J_{\phi_s}+1},J_{\phi_s}))
\nonumber
\end{eqnarray}
By taking the limit when $(t_{{J_{\phi_s}}+1},J_{\phi_s}) \rightarrow \sup \dom \phi_s$, and given that (\ref{TerminalCondCG}) holds, we have
\begin{eqnarray}
V(\xi)&\leq&\sum_{j=0}^{\sup_j \dom \phi_s} \int_{t_{j}}^{t_{j+1}}  L_C(\phi_s(t,j),u_{C}^s(t,j))dt\nonumber
\\&&\nonumber
+ \sum_{j=0}^{{\sup_j \dom \phi_s}-1}L_D(\phi_s(t_{j+1},j),u_{D}^s(t_{j+1},j))
\\&&
+\underset{(t,j) \in \textup{dom}\phi_s}{\limsup_{t+j\rightarrow \infty} } V(\phi_s(t,j))
\nonumber \\ 
&\leq& \mathcal{J}(\xi,u^s)
\label{BoundCostNOz}
\end{eqnarray}
}
\IfAutom{Now, when $(\phi_s, u_s)$ is complete, with $u^s=(u_1^s,u_2^s)$, $u_1^s$ defined by any $\bar{\IfIh{\kappa}{\gamma}}_1 \in \mathcal{K}_1$, and $u_2^s$ defined by $\IfIh{\kappa}{\gamma}_2$ as in (\ref{kHJBeqzsihc}) and (\ref{kBeqzsihc}), thanks to \eqref{eq:PrHamiltonC5} and \eqref{eq:PrHamiltonD5}, 
 from {Corollary} \ref{Rmk:ProofDerivation}, we have that 
  \begin{equation}
  V(\xi) \leq \mathcal{J}(\xi,u^s).
    \label{BoundCostNOz}
  \end{equation}
  }{}
 \\\\
%
Proceeding with item {c} as above, pick any 
$(\phi_w,u^w) \in {\mathcal{S}_{\HS}^w(\xi)}$, where  
$\mathcal{S}_{\HS}^w(\xi) (\subset \mathcal{S}_{\HS}(\xi))$ is the set of solutions $(\phi,u)$ with $u=(u_{1},u_{2})$, $\dom \phi \ni (t,j) \mapsto u_1(t,j)= \IfIh{\kappa}{\gamma}_1(\IfIh{}{t,j,}\phi(t,j))$ for ${\IfIh{\kappa}{\gamma}}_1:=(\IfIh{\kappa}{\gamma}_{C1}, \IfIh{\kappa}{\gamma}_{D1}) $ as in (\ref{kHJBeqzsihc}) and (\ref{kBeqzsihc}), $\dom \phi \ni (t,j) \mapsto u_2(t,j)= \bar{\IfIh{\kappa}{\gamma}}_2(\IfIh{}{t,j,}\phi(t,j))$ for some $\bar{\IfIh{\kappa}{\gamma}}_2 \in \mathcal{K}_2$.
Since $\bar{\IfIh{\kappa}{\gamma}}_2$ does not necessarily attain the maximum in (\ref{HJBzsih}), then, for each $j \in \nats$ such that $I_{\phi_w}^j=[t_j, t_{j+1}]$ has a nonempty interior int$I_{\phi_w}^j$, we have that for every $t \in \textup{int}I_{\phi_w}^j$,
{
\begin{equation*}
\begin{split}
0&
\geq L_C(\phi_w(t,j),u_{ C}^w(t,j))
\\& \hspace{2cm}
+\left\langle \nabla V(\phi_w(t,j)),F(\phi_w(t,j),u_C^w(t,j)) \right\rangle.
\end{split}
\end{equation*}
}
%
{Similarly to \eqref{eq:PrHamiltonC4}, we have} 
\begin{equation}
 \frac{d V}{d t}(\phi_w(t,j)):=\left\langle \nabla V(\phi_w(t,j)),F(\phi_w(t,j),u_C^w(t,j)) \right\rangle
 \label{eq:PrHamiltonC6}
 \end{equation}
  for every $(t,j) \in \textup{int}(I^j_{\phi_w}) \times \{j\}$ with $ \textup{int}(I^j_{\phi_w}) $ 
  nonempty{\IfAutom{.}{, and
by integrating over the interval $[t_j,t_{j+1}]$, 
we obtain
\begin{eqnarray*}
0 \geq \int_{t_{j}}^{t_{j+1}} \left( L_C(\phi_w(t,j),u_{ C}^w(t,j))
+ \frac{d V}{d t}(\phi_w(t,j)) \right) dt
\end{eqnarray*}
from which we have
\begin{eqnarray*}
V(\phi_w(t_j,j))
\geq \int_{t_{j}}^{t_{j+1}} L_C(\phi_w(t,j),u_{ C}^w(t,j))dt
\\
+V(\phi_w(t_{j+1},j))
\end{eqnarray*}
Summing both sides from $j=0$ to $j=J_{\phi_w}$, 
we obtain 
\begin{eqnarray}
\sum_{j=0}^{J_{\phi_w}}V(\phi_w(t_j,j)) 
\geq \sum_{j=0}^{J_{\phi_w}} \int_{t_{j}}^{t_{j+1}} L_C(\phi_w(t,j),u_{C}^w(t,j))dt
\nonumber \\
+\sum_{j=0}^{J_{\phi_w}}V(\phi_w(t_{j+1},j))
\nonumber
\end{eqnarray}
Then, solving for $V$ at the initial condition $\phi_w(0,0)$, we obtain
%
%
\begin{equation}\label{ContinCost2goNOihzw}
\begin{medsize}
\begin{split} 
&V(\phi_w(0,0))
\geq \sum_{j=0}^{J_{\phi_w}} \int_{t_{j}}^{t_{j+1}} L_C(\phi_w(t,j),u_{C}^w(t,j))dt
\\&
+ V(\phi_w(t_1,0))+\sum_{j=1}^{J_{\phi_w}} \left( V(\phi_w(t_{j+1},j))-V(\phi_w(t_j,j)) \right) 
\end{split}
\end{medsize}
\end{equation}
}
In addition, since 
$\bar{\IfIh{\kappa}{\gamma}}_2$ does not necessarily attain the maximum in (\ref{Bellmanzsih}), then for every $(t_{j+1},j) \in \dom \phi_w $ such that $(t_{j+1},j+1) \in \dom \phi_w$, we have
{
\begin{equation}
\begin{split} 
V(\phi_w(t_{j+1},j))
\geq& \NotAutom{
  \>
L_D(\phi_w(t_{j+1},j),u_{ D}^w(t_{j+1},j))
\\& \hspace{0.7cm}
+ V(G(\phi_w(t_{j+1},j),u_{D}^w(t_{j+1},j)))
\\
=&
}
\>L_D(\phi_w(t_{j+1},j),u_{ D}^w(t_{j+1},j))
\\&\hspace{1.5cm}
+ V( \phi_w(t_{j+1},j+1)).
\end{split}
\label{eq:PrHamiltonD6}
\end{equation}}
\NotAutom{
Summing both sides from $j=0$ to $j=J_{\phi_w}-1$, 
we obtain
\begin{equation}
\begin{medsize}
\begin{split} 
\sum_{j=0}^{J_{\phi_w}-1}V(\phi_w(t_{j+1},j)) 
\geq 
\sum_{j=0}^{J_{\phi_w}-1}L_D(\phi_w(t_{j+1},j),u_{ D}^w(t_{j+1},j))\nonumber
\\
+\sum_{j=0}^{J_{\phi_w}-1} V(\phi_w(t_{j+1},j+1))
\nonumber
\end{split}
\end{medsize}
\end{equation}
Then, solving for $V$ at the first jump time, we obtain
%
\begin{eqnarray} \nonumber
&&V(\phi_w(t_{1},0))
\geq V( \phi_w(t_{1},1))
\\&&
+ \sum_{j=0}^{J_{\phi_w}-1}L_D(\phi_w(t_{j+1},j),u_{ D}^w(t_{j+1},j)) \label{DiscreteCost2goNOihzw}
\\&&
+\sum_{j=1}^{J_{\phi_w}-1}\left( V(\phi_w(t_{j+1},j+1))-V(\phi_w(t_{j+1},j)) \right) \nonumber
\end{eqnarray}
In addition, given that $\phi_w(0,0)=\xi$, lower bounding $V(\phi_w(t_{1},0))$ in (\ref{ContinCost2goNOihzw}) by the right-hand side of (\ref{DiscreteCost2goNOihzw}), we obtain 
\begin{equation}
\begin{medsize}
\begin{split} 
V(\xi)
\geq& \sum_{j=0}^{J_{\phi_w}} \int_{t_{j}}^{t_{j+1}} L_C(\phi_w(t,j),u_{C}^w(t,j))dt+ V(\phi_w(t_1,0)) 
\\&
+\sum_{j=1}^{J_{\phi_w}} \left( V(\phi_w(t_{j+1},j))-V(\phi_w(t_j,j)) \right)  \\
\geq&
\sum_{j=0}^{J_{\phi_w}} \int_{t_{j}}^{t_{j+1}} L_C(\phi_w(t,j),u_{C}^w(t,j))dt
\\&
+ \sum_{j=0}^{J_{\phi_w}-1}L_D(\phi_w(t_{j+1},j),u_{D}^w(t_{j+1},j))
 \\
&+\sum_{j=1}^{J_{\phi_w}-1}\left( V(\phi_w(t_{j+1},j+1))-V(\phi_w(t_{j+1},j)) \right) 
\\&
+ V(\phi_w(t_{1},1))+\sum_{j=1}^{J_{\phi_w}} \left( V(\phi_w(t_{j+1},j))-V(\phi_w(t_j,j)) \right) 
\end{split}
\end{medsize}
\end{equation}
Since 
\begin{eqnarray} \nonumber
&& V(\phi_w(t_{1},1))
\\&&
+\sum_{j=1}^{J_{\phi_w}-1}\left( V(\phi_w(t_{j+1},j+1))-V(\phi_w(t_{j+1},j)) \right) \nonumber 
\\&&
+\sum_{j=1}^{J_{\phi_w}} \left( V(\phi_w(t_{j+1},j))-V(\phi_w(t_j,j)) \right) \nonumber
\\ &=& V(\phi_w(t_{{J_{\phi_w}}+1},J_{\phi_w})) +V(\phi_w(t_{1},1))\nonumber
\\&&
+\sum_{j=1}^{J_{\phi_w}-1}\left( V(\phi_w(t_{j+1},j+1)) \right) \nonumber 
-\sum_{j=1}^{J_{\phi_w}} \left(V(\phi_w(t_j,j)) \right) \nonumber
\\ &=&V(\phi_w(t_{{J_{\phi_w}}+1},J_{\phi_w})) 
\end{eqnarray}
then we have
\begin{eqnarray}\nonumber
V(\xi)
&\geq&\sum_{j=0}^{J_{\phi_w}} \int_{t_{j}}^{t_{j+1}} L_C(\phi_w(t,j),u_{C}^w(t,j))dt
\\&&
+ \sum_{j=0}^{J_{\phi_w}-1}L_D(\phi_w(t_{j+1},j),u_{D}^w(t_{j+1},j))
\\&&+V(\phi_w(t_{J_{\phi_w}+1},J_{\phi_w}))
\nonumber
\end{eqnarray}
By taking the limit when $(t_{J_{\phi_w}+1},J_{\phi_w}) \rightarrow \sup \dom \phi_w$, and given that (\ref{TerminalCondCG}) holds, we have
\begin{eqnarray} \nonumber
V(\xi)&\geq&\sum_{j=0}^{\sup_j \dom \phi_w} \int_{t_{j}}^{t_{j+1}} L_C(\phi_w(t,j),u_{C}^w(t,j))dt
\\&&
+ \sum_{j=0}^{{\sup_j \dom \phi_w}-1}L_D(\phi_w(t_{j+1},j),u_{D}^w(t_{j+1},j))
\nonumber \\
&& 
+\underset{(t,j) \in \textup{dom}\phi_w}{\limsup_{t+j\rightarrow \infty} } V(\phi_w(t,j))
\nonumber \\ 
&\geq& \mathcal{J}(\xi,u^w)
\label{BoundCostNOzw}
\end{eqnarray}
}
 \IfAutom{Now, when $(\phi_w,u^w)$ is complete, 
 %
 with  $u^w=(u_1^w,u_2^w)$, $u_1^w$ defined by $\IfIh{\kappa}{\gamma}_1$ as in (\ref{kHJBeqzsihc}) and (\ref{kBeqzsihc}), and $u_2^w$ defined by any $\bar{\IfIh{\kappa}{\gamma}}_2 \in \mathcal{K}_2$,}{}
 thanks to \eqref{eq:PrHamiltonC6} and \eqref{eq:PrHamiltonD6},
  from Corollary \ref{Rmk:ProofDerivation}, we have that 
  \begin{equation}
  V(\xi) \geq \mathcal{J}(\xi,u^w).
    \label{BoundCostNOzw}
  \end{equation}
  }{}
{Finally, 
by proceeding as in item {d} above,
by applying the infimum on each side of  (\ref{BoundCostNOz}) 
over the set $\mathcal{S}_{\HS}^\infty (\xi)$, we obtain
\begin{equation*}
V(\xi) \leq  \underset{{u_1:}(\phi_s, (u_1,{\kappa_2(\phi_s)})) \in \mathcal{S}_{\HS}^\infty(\xi)
}{\text{inf}}\mathcal{J}(\xi, (u_1,{\kappa_2(\phi_s)}) {=: \overline V(\xi)}
\end{equation*}
By applying the supremum on each side of (\ref{BoundCostNOzw})
over the set $\mathcal{S}_{\HS}^\infty (\xi)$, we obtain
\begin{equation*}
V(\xi) \geq \underset{{u_2:}(\phi_w, ({\kappa_1(\phi_w)},u_2)) \in \mathcal{S}_{\HS  }^\infty(\xi)
}{\text{sup}}\mathcal{J}(\xi, ({\kappa_1(\phi_w)},u_2)){=: \underline{ V}(\xi)}.
\end{equation*}
Given that $V(\xi)=\mathcal{J}(\xi, u^
*)$ from (\ref{OptimalCostihz}), we have that for any $\xi \in \Pi ( \overline{C}) \cup \Pi(D)$, {each 
$(\phi^*,u^*)\in \mathcal{S}_{\HS }^\infty(\xi)$ with $u^*=(\kappa_1(\phi^*),\kappa_2(\phi^*))$ satisfies}

\begin{multline}
{\underline{V}(\xi)}
\leq
\mathcal{J}(\xi,u^*) 
\leq 
{\overline{V}(\xi)}
\label{Boundssz}
\end{multline}
Since 
 $(\phi^*,u^*) \in \mathcal{S}_{\HS}^s(\xi) \cap \mathcal{S}_{\HS}^w(\xi)\cap \mathcal{S}_{\HS }^\infty(\xi)$,} we have
\begin{multline}
\underline{V}(\xi)
=\\
{\underset{ (\phi^*, (\kappa_1(\phi^*),\kappa_2(\phi^*))) \in \mathcal{S}_{\HS }^\infty(\xi)
}{\text{sup}} }\mathcal{J}(\xi, (\kappa_1(\phi^*),\kappa_2(\phi^*)))
\label{supJk}
\end{multline}
and
\begin{multline}
\overline{V}(\xi)
=\\ \underset{ (\phi^*, (\kappa_1(\phi^*),\kappa_2(\phi^*))) \in \mathcal{S}_{\HS}^\infty(\xi)
}{\text{inf}} \mathcal{J}(\xi, (\kappa_1(\phi^*),\kappa_2(\phi^*)).
\label{infJk}
\end{multline}}
Given that the supremum and infimum are attained in 
\eqref{supJk} and \eqref{infJk} by $\underline{V}(\xi)$ and $\overline{V}(\xi)$,  respectively, \eqref{Boundssz}
leads to
\begin{equation}
\mathcal{J}(\xi,u^*) 
=\underset{(u_1,u_2) \in \mathcal{U}_{\HS  }^\infty(\xi)}
{ \min_{u_{1}} \max_{u_{2}}}
\mathcal{J}(\xi,(u_1,u_2))
\label{LowerBounddz}
\end{equation}
Thus, from (\ref{OptimalCostihz}) and (\ref{LowerBounddz}), $V(\xi)$ is the value function for $\HS$, as in Definition \ref{ValueFunctionz} and from (\ref{Boundssz}), $\kappa$ is the saddle-point equilibrium as in Definition \ref{NashEqNonCoop}.

\hfill $\blacksquare$
}
\begin{remark}{Sufficiency of Hamilton-Jacobi-Isaacs Hybrid Equations}
Notice that satisfying the Hamilton-Jacobi-Isaacs Hybrid equations \eqref{HJBzsih} and \eqref{Bellmanzsih} is sufficient to
guarantee that any min-max strategy is also a max-min strategy. 
Thus, unlike \cite{GOEBEL2019153nu}, Theorem \ref{thHJBszih} does not explicitly require regularity properties of $\J$ (e.g., convexity/concavity)
and $\mathcal{U}_\HS^\infty$ (e.g., convexity and compactness)
to guarantee the existence of a solution to Problem $(\diamond)$.
\end{remark}
\begin{remark}{Connections between Theorem \ref{thHJBszih} and Problem $(\diamond)$} 
Given $\xi \in \Pi ( \overline{C}) \cup \Pi(D)$, 
{if there exist a function $V$ satisfying the} conditions in Theorem \ref{thHJBszih},
{then}
a solution to Problem ($\diamond$) exists, namely there exists an input action $u^*=(u_C^*,u_D^*)=((u_{C1}^*,u_{C2}^*),(u_{D1}^*,u_{D2}^*)) \in \mathcal{U}_{\HS  }^\infty (\xi)$ such that $\mathcal{J}(\xi,u^*) <\infty$, that attains the {min-max} in (\ref{problemzsih}), and, as a consequence, satisfies (\ref{SaddlePointIneq}) in Definition \ref{NashEqNonCoop}. In addition, 
{the} strategy $\IfIh{\kappa=(\kappa_C, \kappa_D)}{\gamma} \in \mathcal{K}$ {with values as in (\ref{kHJBeqzsihc}) and (\ref{kBeqzsihc}) is} such that every complete solution to the closed-loop system $\HS_{\IfIh{\kappa}{\gamma}} $ from $\xi$ has a cost that is equal to the {min-max} in (\ref{problemzsih}). 
\label{SolExistenceZ}
\end{remark}
}
\IfIncd{}{\NotAutom{
\section{Linear Quadratic Hybrid Games 
}
In this section, we consider a special case of our result that emerges {in practical scenarios with} hybrid systems with linear flow and jump maps and periodic jumps, {as in noise attenuation of cyber-physical systems,  see, e.g., \cite{ferrante2019certifying,possieri2020linear,possieri2017mathcal}}. We introduce a state variable $\tau$ that plays the role of a timer. Once $\tau$ reaches a fixed threshold $\bar{T}$, it triggers a jump in the state and resets $\tau$ to {zero. More precisely,
given $\bar{T}\in \mathbb{R}$, we consider a hybrid}
system with state $x=(x_p,\tau){=((x_{p1},x_{p2}),\tau)} \in \reals^n \times [0, \bar{T}]$, input $u=(u_C, u_D)=((u_{C1},u_{C2}),$
$(u_{D1},u_{D2})) \in \reals^{m_C}\times \reals^{m_D}$, and dynamics $\HS$ as in (\ref{Heq})\IfTp{,}{ with $N=2$ and} {defined} by
\begin{equation}\label{eq: CDFG priodic}
\begin{array}{rll}
C&:=&\mathbb{R}^n \times [0,\bar{T}] \times \mathbb{R}^{m_C}
\\
F(x,u_C)&:=&(A_C x_p+B_C u_C,1) 
\quad \forall(x,u_C) \in C
\\
D&:=&\mathbb{R}^n \times \{\bar{T}\} \times \mathbb{R}^{m_D}
\\
G(x,u_D)&:=&(A_D x_p+B_D u_D,0)
\quad \forall(x,u_D) \in D
\end{array}
\end{equation}
with $A_C=\left[\begin{smallmatrix}A_{C1} & 0 \\ 0 & A_{C2}\end{smallmatrix}\right], 
B_C=\left[ \begin{smallmatrix}B_{C1} \>B_{C2}\end{smallmatrix}\right],
A_D = \left[\begin{smallmatrix}A_{D1} & 0 \\ 0 & A_{D2}\end{smallmatrix}\right],$ and 
$ B_D = \left[ \begin{smallmatrix}B_{D1} \>B_{D2}\end{smallmatrix}\right] $.
The input $u_1=(u_{C1},u_{D1})$ is assigned by {player} $P_1$ and the input $u_2=(u_{C2},u_{D2})$ is assigned by {player} $P_2$. The problem of finding conditions for $u_1$ to minimize a cost functional $\mathcal{J}$ in the presence of the action $u_2$ that seeks to maximize it, is formulated as {the} two-player zero-sum game {as in Section \ref{Sec: Game form}.} Thus, \textcolor{black}{by solving Problem $(\diamond)$ for any $\xi \in \Pi(C\cup D)$}, 
the control objective is achieved. 
{Using Theorem \ref{thHJBszih},} the following result presents a tool for the solution of the optimal control problem for hybrid systems with linear maps and periodic jumps under an {adversarial} action.

\begin{corollary}{Hybrid Riccati equation for periodic jumps} Given {a hybrid system $\HS$ as in \eqref{Heq} defined by $(C,F,D,G)$ as in \eqref{eq: CDFG priodic}, let} $\bar{T}\in \mathbb{R}$, 
and,
{with the aim of pursuing minimum energy and distance to the origin, consider the cost functions $L_C(x,u_C):=x_p^\top Q_C x_p+ u_{C1}^\top R_{C1} u_{C1}+ u_{C2}^\top R_{C2} u_{C2}$, $L_D(x,u_D):=x_p^\top Q_D x_p+  u_{D1}^\top R_{D1} u_{D1}+ u_{D2}^\top R_{D2} u_{D2} $, and terminal cost $q(x):=x_p^\top P(\tau) x_p$ \textcolor{black}{defining $\J$ as in (\ref{defJTNCinc}), } {with} $Q_C,Q_D\in \mathbb{S}^n_+$, $R_{C1} \in \mathbb{S}^{m_{C_1}}_+$, $-R_{C2} \in \mathbb{S}^{m_{C_2}}_+$, $R_{D1} \in \mathbb{S}^{m_{D_1}}_+$, and $-R_{D2} \in \mathbb{S}^{m_{D_2}}_+$.} Suppose there exists a matrix {function} $P:[0,\bar{T}]\rightarrow \mathbb{S}^n_+ $ {that is} continuously differentiable {and} such that
{
\begin{equation}
\begin{split}
-\frac{\mathrm{d} P(\tau)}{\mathrm{d} \tau} =&-P(\tau)(  B_{C2} R_{C2}^{-1}B_{C2}^\top+B_{C1} R_{C1}^{-1}B_{C1}^\top) P(\tau)
\\&\hspace{0.1cm} +Q_C+P(\tau) A_C+A_C^\top P(\tau) \hspace{0.5cm} \forall \tau  \in (0,\bar{T}),
\end{split}
\label{DiffRiccatizsih}
\end{equation}}
\begin{equation}
\begin{array}{rrl}
-R_{D2} -B_{D2}^\top P(0) B_{D2},\>&
R_{D1} +B_{D1}^\top P(0) B_{D1}&
\in \mathbb{S}_{0+}^{m_D},
\end{array}
\label{zlqeqinv}
\end{equation}
the matrix $R_v=\left[\begin{smallmatrix}
R_{D1}+B_{D1}^\top P(0) B_{D1}
&
B_{D1}^\top P(0) B_{D2}
 \\
B_{D2}^\top P(0) B_{D1}
&
R_{D2}+B_{D2}^\top P(0) B_{D2}
\end{smallmatrix}\right] $ is invertible, and 
\begin{multline}
P (\bar{T})=Q_D+A_D^\top P(0) A_D \IfConf{\\}{}
-
\begin{bmatrix}
A_D^\top P(0) B_{D1}
 &
A_D^\top P(0) B_{D2}
\end{bmatrix}
R_v
^{-1}
\begin{bmatrix}
B_{D1}^\top P(0) A_D
\\
B_{D2}^\top P(0) A_D
\end{bmatrix}
\label{Riccatizsih} 
\end{multline}
{ where $A_C, B_{C1}, B_{C2}, A_D, B_{D1}$, and $B_{D2}$ are defined below \eqref{eq: CDFG priodic}.}
Then, the feedback law $\kappa:=(\kappa_C,\kappa_D)$, with values
\begin{eqnarray}
\kappa_C(x)=(-R_{C1}^{-1}B_{C1}^\top P(\tau) x_p,-R_{C2}^{-1}B_{C2}^\top P(\tau) x_p) \nonumber
\\ \hspace{0.4cm} \forall x \in \Pi(C),
\label{NashkCLQzsih}
\end{eqnarray}
\begin{equation}
\kappa_D(x)=-R_v
^{-1}
\begin{bmatrix}
B_{D1}^\top P(0) A_D
\\
B_{D2}^\top P(0) A_D
\end{bmatrix} x_p \hspace{1cm} \forall x \in \Pi(D)
\label{NashkDLQzsih}
\end{equation}
is the pure strategy saddle-point equilibrium for the two-player zero-sum hybrid game with \IfIncd{spontaneous}{periodic} jumps. In addition, for each $x=(x_p,\tau) \in \Pi ( \overline{C}) \cup \Pi(D)$, the value function is equal to $V(x):=x_p^\top P(\tau) x_p$.
\label{HREqPJ}
\end{corollary}
\begin{proof}
We show that when conditions (\ref{DiffRiccatizsih})-(\ref{Riccatizsih}) hold, by using Theorem \ref{thHJBszih}, the value function is equal to the function $V$, and a cost is attained with 
 the feedback law $\kappa:=(\kappa_C,\kappa_D)$ with values as in (\ref{NashkCLQzsih}) and (\ref{NashkDLQzsih}). 
We can write (\ref{HJBzsih}) in Theorem \ref{thHJBszih} as
\begin{eqnarray}
&&0=
\underset{u_C=(u_{C1},u_{C2}) \in \Pi_u(x,C)}
{\min_{u_{C1}} \max_{u_{C2}}}  \mathcal{L}_C(x,u_C),
 \nonumber
\\ 
&&\mathcal{L}_C(x,u_C)= 
 x_p^\top Q_C x_p+ u_{C1}^\top R_{C1} u_{C1}
 + u_{C2}^\top R_{C2} u_{C2}
\nonumber \\
&&\hspace{0.9cm} {
+2x_p^\top P(\tau) (A_C x_p+B_C u_C)  
+ x_p^\top \frac{\mathrm{d} P(\tau)}{\mathrm{d} \tau} x_p}
\label{HJBzsihlq}
\end{eqnarray}
First, given that (\ref{DiffRiccatizsih}) holds, and $x_p^\top(P(\tau) A_C+A_C^\top P(\tau))x_p={2x_p^\top P(\tau)A_Cx_p}$ for every $x\in \Pi(C)$, one has
\begin{equation}
\begin{split}
\mathcal{L}_C(x,u_C)= \hspace{6.5cm} \nonumber\\
 x_p^\top P(\tau)(  B_{C2} R_{C2}^{-1}B_{C2}^\top+B_{C1} R_{C1}^{-1}B_{C1}^\top)P(\tau)  x_p
 \nonumber \\
 + u_{C1}^\top R_{C1} u_{C1}+ u_{C2}^\top R_{C2} u_{C2}
+ 2 x_p^\top  P(\tau) B_C u_C 
\end{split}
\end{equation}
The {first-order} necessary conditions for optimality
{
\IfAutom{
\begin{equation*}
\frac{\partial}{\partial u_{C1}}
\mathcal{L}_C(x,u_C)  \Big|_{u_{C}^*} =0,
\quad
\frac{\partial}{\partial u_{C2}}
\mathcal{L}_C(x,u_C) \Big|_{u_{C}^*} =0
\end{equation*}
{for all $(x,u_C)\in C$}}
{\begin{multline*}
\frac{\partial}{\partial u_{C1}} \left( 
 x_p^\top P(\tau)(  B_{C2} R_{C2}^{-1}B_{C2}^\top+B_{C1} R_{C1}^{-1}B_{C1}^\top)P(\tau)  x_p\right. \\ + u_{C1}^\top R_{C1} u_{C1}+ u_{C2}^\top R_{C2} u_{C2}
\\\left.\left.+ 2 x_p^\top  P(\tau)  (B_{C1} u_{C1} + B_{C2} u_{C2}) \right)\right|_{u_{C}^*} =0
\end{multline*}
\begin{multline*}
\frac{\partial}{\partial u_{C2} }\left( 
 x_p^\top P(\tau)(  B_{C2} R_{C2}^{-1}B_{C2}^\top+B_{C1} R_{C1}^{-1}B_{C1}^\top)P(\tau)  x_p\right.+\\ u_{C1}^\top R_{C1} u_{C1}+ u_{C2}^\top R_{C2} u_{C2}
 \\ \left.\left.+ 2 x_p^\top  P(\tau) (B_{C1} u_{C1} + B_{C2} u_{C2}) \right)\right|_{u_{C}^*} =0
\end{multline*}}
}
are satisfied by the point $u_C^*=(u_{C1}^*,u_{C2}^*)$, with values
\begin{equation}
u_{C1}^*=-R_{C1}^{-1}B_{C1}^\top P(\tau) x_p, 
\quad
u_{C2}^*=-R_{C2}^{-1}B_{C2}^\top P(\tau) x_p
\label{zLQHPu2*}
\end{equation}
{for each $x=(x_p,\tau)\in \Pi(C)$.}
Given that $R_{C1},-R_{C2} \in \mathbb{S}^{m_D}_+$, the second-order sufficient conditions for optimality 
%
{
\IfAutom{
\begin{equation*}
\frac{\partial^2}{\partial u_{C1}^2} 
 \mathcal{L}_C(x,u_C) \Big|_{u_C^*} \succeq 0,
 \quad
\frac{\partial^2}{\partial u_{C2}^2}
 \mathcal{L}_C(x,u_C)\Big|_{u_C^*} \preceq 0,
\end{equation*}
{for all $(x,u_C)\in C$}
}
{
\begin{multline*}
\frac{\partial^2}{\partial u_{C1}^2} \left( 
 x_p^\top P(\tau)(  B_{C2} R_{C2}^{-1}B_{C2}^\top+B_{C1} R_{C1}^{-1}B_{C1}^\top)P(\tau)  x_p\right. \\ + u_{C1}^\top R_{C1} u_{C1}+ u_{C2}^\top R_{C2} u_{C2}
 \\\left.\left.
+ 2 x_p^\top  P(\tau)  (B_{C1} u_{C1} + B_{C2} u_{C2}) \right)\right|_{u_{C1}^*} \succeq 0
\end{multline*}
\begin{multline*}
\frac{\partial^2}{\partial u_{C2}^2 }\left( 
 x_p^\top P(\tau)(  B_{C2} R_{C2}^{-1}B_{C2}^\top+B_{C1} R_{C1}^{-1}B_{C1}^\top)P(\tau)  x_p\right. \\ + u_{C1}^\top R_{C1} u_{C1}+ u_{C2}^\top R_{C2} u_{C2}
 \\\left.\left.+ 2 x_p^\top  P(\tau) (B_{C1} u_{C1} + B_{C2} u_{C2}) \right)\right|_{u_{C2}^*} \preceq 0
\end{multline*}}
}
hold, rendering $u_{C}^*$ as in 
\eqref{zLQHPu2*} as an optimizer of the min-max problem in (\ref{HJBzsihlq}).
In addition, it satisfies $\mathcal{L}_C(x,u_C^*)=0$, making $V(x)=x_p^\top P(\tau)x_p$ a solution {to} (\ref{HJBzsih}) in Theorem \ref{thHJBszih}. 

On the other hand, we can write (\ref{Bellmanzsih}) in Theorem \ref{thHJBszih} as
\begin{eqnarray}
&&x_p^\top P(\bar{T}) x_p =\underset{u_D=(u_{D1},u_{D2}) \in \Pi_u(x,D)}
{\min_{u_{D1}} \max_{u_{D2}}}  \mathcal{L}_D(x,u_D),
 \nonumber
\\ 
&&\mathcal{L}_D(x,u_D)= 
x_p^\top Q_D x_p+ u_{D1}^\top R_{D1} u_{D1}+ u_{D2}^\top R_{D2} u_{D2} \hspace{0.4cm}
\nonumber\\
&&\hspace{1.3cm} {
+ (A_D x_p+B_D u_D)^\top P(0) (A_D x_p+B_D u_D) }
\label{Bellmanzsihlq}
\end{eqnarray}
\NotAutom{which can be expanded as
{\small
\begin{multline}
\mathcal{L}_D(x,u_D)=  
 x_p^\top (Q_D+A_D^\top P(0) A_D) x_p + 2 x_p^\top A_D^\top  P(0) B_D u_D \\
 + u_{D1}^\top (R_{D1} +B_{D1}^\top P(0) B_{D1}) u_{D1}
 + u_{D2}^\top (R_{D2} +B_{D2}^\top P(0) B_{D2}) u_{D2}\\
 + u_{D1}^\top (B_{D1}^\top P(0) B_{D2}) u_{D2}
 + u_{D2}^\top (B_{D2}^\top P(0) B_{D1}) u_{D1}
\end{multline}}
}
{Similar to the case along flows, the first-order} necessary conditions for optimality
\NotAutom{
\begin{equation*}
\frac{\partial}{\partial u_{D1}}\left.
\mathcal{L}_D(x,u_D)  \right|_{u_{D}^*} =0,
\quad
\frac{\partial}{\partial u_{D2}}\left.
\mathcal{L}_D(x,u_D) \right|_{u_{D}^*} =0
\end{equation*}
}
are satisfied by the point $u_D^*=(u_{D1}^*,u_{D2}^*)$, such that, for each $x_p \in \Pi(D)$,
%
{
\begin{equation}
{\tiny
\begin{split}
u_{D}^*= \hspace{7.5cm}
\\-\begin{bmatrix}
R_{D1}+B_{D1}^\top P(0) B_{D1}
&
B_{D1}^\top P(0) B_{D2}
 \\
B_{D2}^\top P(0) B_{D1}
&
R_{D2}+B_{D2}^\top P(0) B_{D2}
\end{bmatrix} 
^{-1}
\begin{bmatrix}
B_{D1}^\top P(0) A_D
\\
B_{D2}^\top P(0) A_D
\end{bmatrix} x_p
\end{split}}
\label{zLQHPud*}
\end{equation}}
Given that (\ref{zlqeqinv}) holds, the second-order sufficient conditions for optimality 
\NotAutom{
\begin{equation*}
\frac{\partial^2}{\partial u_{D1}^2}\left. 
 \mathcal{L}_D(x,u_D) \right|_{u_D^*} \succeq 0,
 \quad
\frac{\partial^2}{\partial u_{D2}^2}\left.
 \mathcal{L}_D(x,u_D)\right|_{u_D^*} \preceq 0,
\end{equation*}
}
are satisfied, rendering $u_D^*$ as in (\ref{zLQHPud*}) as an optimizer of the {min-max} problem in (\ref{Bellmanzsihlq}).
In addition, $u_D^*$ satisfies $\mathcal{L}_D(x,u_D^*)=x_p^\top P(\bar{T})x_p$, with $P(\bar{T})$ as in (\ref{Riccatizsih}), making $V(x)=x_p^\top P(\tau)x_p$ a solution of (\ref{Bellmanzsih}) in Theorem \ref{thHJBszih}. 

Then, given that $V$ is continuously differentiable on a neighborhood of $\Pi(C)$ and {that} Assumption \ref{AssLipsZ} holds, by applying Theorem \ref{thHJBszih}, in particular from (\ref{ResultValue}), for every $ \xi=(\xi_p,\xi_\tau) \in \Pi(\overline{C}\cup D)$ the value function is $\J^*(\xi)=\J(\xi,((u_{C1}^*,u_{D1}^*),(u_{C2}^*,u_{D2}^*))= \xi_p^\top P(\xi_\tau)\xi_p$. From (\ref{kHJBeqzsihc}) and (\ref{kBeqzsihc}), the feedback law $\kappa=(\kappa_C,\kappa_D)$ with values as in (\ref{NashkCLQzsih}) and (\ref{NashkDLQzsih}) is a pure strategy saddle-point equilibrium.
\end{proof}
}

\NotAutom{By following the same modeling approach and imposing conditions of the hybrid time domains, games for switching systems can be covered by Corollary \ref{HREqPJ}. By selecting appropriate stage costs, optimality is encoded in the satisfaction of the infinitesimal conditions instead of in the knowledge of specific solutions/trajectories. Note that for switching systems, the function $V$ might be independent of the timer state if the stage costs are independent of it as well.}

{
  As illustrated next, there are useful families of hybrid systems for which a pure strategy saddle-point equilibrium exists. 
 The following example {(with nonunique solutions)} characterizes both the pure strategy
saddle-point equilibrium and the value function in a two-player zero-sum game with a one-dimensional state, that is associated to player $P_1$. Thus, $n_1=1, n_2=0$, and the role of player $P_2$ reduces to selecting the action $u_{C2}$.
\begin{example}{Hybrid game with nonunique solutions for a given strategy}
Consider a {hybrid} system $\HS$ with state $x \in \reals$, input $u_C:=(u_{C1},u_{C2}) \in \reals^2$, and dynamics 
\begin{equation}
\begin{array}{rll}
\dot{x}=&F(x,u_C):=a x+B u_C &  x \in [0,\delta]
\\
x^+=&G(x):=\sigma & x = \mu
\end{array}
\label{nuHGamez}
\end{equation}
where $a<0,B=[b_1 \> b_2]$  and {$ \mu>\delta >\sigma>0$\footnote{Given that $ \mu>\delta$, flow from $\mu$ is not possible.}}. 
 Consider the cost functions $L_C(x,u_C):=x^2Q_C+u_C^\top R_C u_C$, $L_D(x):=P(x^2-\sigma^2)$, and terminal cost $q(x):=P x^2$, {defining $\mathcal{J}$ as in (\ref{defJTNCinc}),} with $R_C :=\left[ \begin{smallmatrix}R_{C1} & 0 \\ 0 & R_{C2}\end{smallmatrix}\right]$, {$Q_C,$ $R_{C1},$ $-R_{C2},$ $P\in \reals_{>0}$, such that} 
 \begin{equation}
 Q_C+2Pa-P^2(b_1^2 R_{C1}^{-1}+b_2^2 R_{C2}^{-1})=0.
 \label{PQ}  
 \end{equation}
 Here, $u_{C1}$ is designed by player $P_1$, which aims to minimize a cost functional $\mathcal{J}$, while player $P_2$ seeks to maximize it by means of $u_{C2}$.
  This is formulated as a two-player zero-sum hybrid game \textcolor{black}{via solving Problem ($\diamond$)} {in Section \ref{subsec: probstat}}. The function $V(x):=P x^2$ is such that
\IfConf{
 {
\begin{equation}
\begin{split}
&\underset{u_C=(u_{C1},u_{C2}) \in \reals^2}
{\underset{u_{C1}}{\min}\> \underset{u_{C2}}{\max}}  
\mathcal{L}_C(x,u_C)
\\ 
&=
\underset{u_{C1} \in \mathbb{R}}{\min} \>\underset{u_{C2} \in \mathbb{R}}{\max} 
\left\{ (Q_C+ 2P a) x^2+  R_{C1}u_{C1}^2  \right.
\\ 
&\hspace{1.2cm}
+R_{C2}u_{C2}^2+  2xP ( b_1u_{C1}+ b_2u_{C2})  \big\} 
 = 0
\end{split}
 \label{HJBLQGnusCz}
\end{equation}
}
}{
\begin{eqnarray}
&&
\underset{u_C=(u_{C1},u_{C2}) \in \reals^2}
{\min_{u_{C1}} \max_{u_{C2}}}  
\left\{L_C(x,u_C) +\left\langle \nabla V(x),F(x,u_C) \right\rangle \right\}
\nonumber\\ 
&=&
\IfPers{%
{\min_{u_{C1}\in \mathbb{R}} \max_{u_{C2}\in \mathbb{R}}}  
\left( x^\top (Q_C+ 2P a) x+ u_C^\top R_C u_C
+  2x^\top P B u_C \right)
 \nonumber \\ &=&
{\min_{u_{C1}\in \mathbb{R}} \max_{u_{C2}\in \mathbb{R}}}  
\left\{ x^2 0.142+ 1.304u_{C1}^2-2916u_{C2}^2 
+  0.858x ( b_1u_{C1}+ b_2u_{C2}) \right\}\nonumber \\&=& }%
{\min_{u_{C1}\in \mathbb{R}} \max_{u_{C2}\in \mathbb{R}}}  
\left\{ (Q_C+ 2P a) x^2+  R_{C1}u_{C1}^2
+R_{C2}u_{C2}^2+  2xP ( b_1u_{C1}+ b_2u_{C2})  \right\} \nonumber
 \\ &=&
0 \label{HJBLQGnusCz}
\end{eqnarray}}
holds for all $x \in [0,\delta]$. In fact, the min-max in (\ref{HJBLQGnusCz}) is attained by 
$\kappa_C(x)=(-R_{C1}^{-1} b_1P x,-R_{C2}^{-1}b_2P x)$.
In particular, thanks to \eqref{PQ}, {we have} \IfAutom{\pnn{$\mathcal{L}_C(x,\kappa_C(x))=0.$}}
{
\begin{equation*}\hspace{-0.2cm}
\begin{array}{l}
L_C(x,\kappa_C(x)) +\left\langle \nabla V(x),F(x,\kappa_C(x)) \right\rangle 
  \nonumber\\ \hspace{2cm}
=
 (Q_C+2Pa) x^2+  R_{C1}\kappa^2_{C1}(x)
+R_{C2}\kappa_{C2}^2(x)\nonumber
\\ \hspace{4.3cm}
+  2xP ( b_1\kappa_{C1}(x)+ b_2\kappa_{C2}(x)) \nonumber
 \\\hspace{2cm}
=
 P^2(b_1^2R_{C1}^{-1}+b_2^2 R_{C2}^{-1}) x^2+  R_{C1}\kappa^2_{C1}(x)
 \nonumber
\\ \hspace{2.3cm} 
+R_{C2}\kappa_{C2}^2(x)+  2xP ( b_1\kappa_{C1}(x)+ b_2\kappa_{C2}(x))
 \nonumber
 \\ \hspace{2cm} =
\left[P^2(b_1^2 R_{C1}^{-1}+b_2^2 R_{C2}^{-1}) +  R_{C1}^{-1} b_1^2 P^2
+R_{C2}^{-1} b_2 P^2 \right.\nonumber
\\  \hspace{3.7cm}
\left.-  2P ( b_1R_{C1}^{-1} b_1 P+ b_2R_{C2}^{-1} b_2 P)\right]x^2 \nonumber
 \\ \hspace{2cm}
 =0 
\end{array}
\label{HJBLQGnusCzk}
\end{equation*}}
Then,  $V(x)= Px^2$ is a solution to (\ref{HJBzsih}).
 In addition, the function $V$ is such that
\begin{equation}
\begin{array}{r@{}l}
\underset{(u_{D1},u_{D2}) \in \reals^2}
{\underset{u_{D1}}\min\> \underset{u_{D2}}\max}  
\left\{
 L_D(x)
+ V(G(x)) \right\} =&
P x^2
\end{array}
\label{HJBLQGnusDz}
\end{equation}
at $x=\mu$\pno{,} which makes $V(x)= Px^2$ a solution to (\ref{Bellmanzsih}) \textcolor{black}{with saddle-point equilibrium $\kappa_C$}.
Given that $V$ is continuously differentiable on $\reals$, and that (\ref{HJBzsih}) and (\ref{Bellmanzsih}) hold thanks to (\ref{HJBLQGnusCz}) and (\ref{HJBLQGnusDz}), from Theorem \ref{thHJBszih} {we have that} the value function is
$
\mathcal{J}^*(\xi) :=
 P \xi^2
\label{HyLCost}
$
for any $\xi \in [0,\delta] \cup \{\mu\}$.

{To investigate the case of nonunique solutions {yielded by the feedback law $\kappa_C$}, now let $\delta {\geq}\mu>\sigma>0$ and notice that solutions can potentially flow or jump at $x=\mu$. The set of all maximal {solutions} from $\xi=\delta$ is denoted $\mathcal{R}_\kappa(\xi)=\{\phi_\kappa,\phi_h\}$. {The continuous solution} {$\phi_\kappa$} is such that $ \dom \phi_\kappa=\realsgeq \times \{0\}$, and is given by
$ \phi_\kappa(t,0)=  \delta \exp((a-R_{C1}^{-1} b_1P-R_{C2}^{-1}b_2P) t)  $ for all $t\in[0,\infty)$.
In simple words, $\phi_\kappa $ flows from $\delta${, and converges (exponentially fast)} to $0$. 
The maximal {solution} $\phi_h$ has domain $ \dom \phi_h=([0,t^h]\times \{0\}) \cup ([t^h, \infty)\times \{1\})$, and is given by
$ \phi_h(t,0)= \delta \exp((a-R_{C1}^{-1} b_1P-R_{C2}^{-1}b_2P) t) ,\hspace{0.5cm} \phi_h(t,1)=
\sigma \exp((a-R_{C1}^{-1} b_1P-R_{C2}^{-1}b_2P) (t-t^h)).$
In simple words, $\phi_h$ flows from $\delta$ to $\mu$ in $t^h$ seconds, then it jumps to $\sigma$, and flows converging (exponentially fast) to {zero}. Figure \ref{fig:my_label} illustrates this behavior.
%
%
By denoting the corresponding input signals as $u_\kappa=\kappa(\phi_\kappa)$ and $u_h=\kappa(\phi_h)$, 
we show in the bottom of Figure \ref{fig:my_label} that the cost of {the solutions $\phi_\kappa$ and $\phi_h$, yielded by $\kappa_C$,} equal $P\delta^2$. 
\begin{figure}[h]
    \centering
    \includegraphics[width=8.5cm]{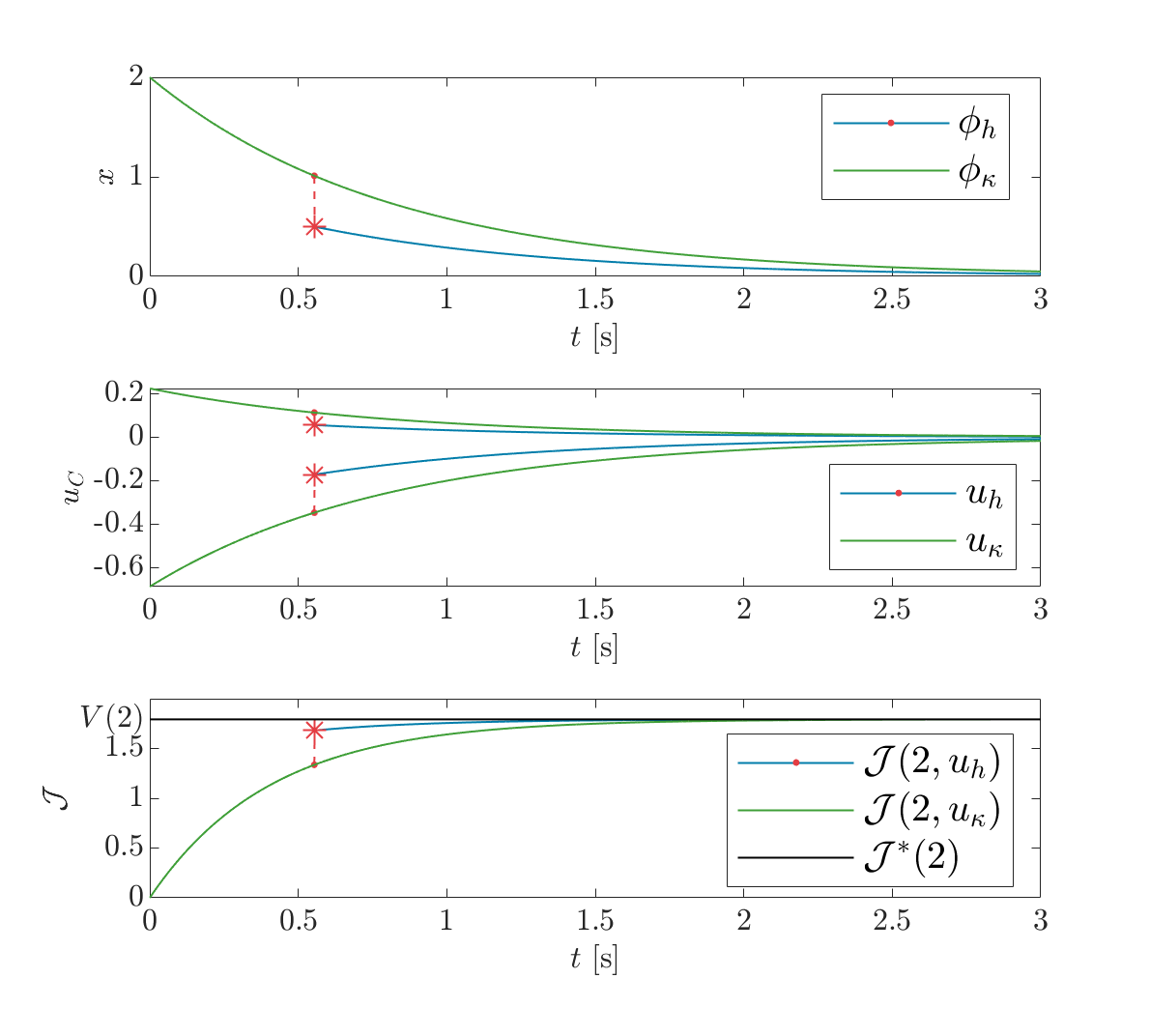}
    \vspace{-0.8cm}
    \caption{Nonunique solutions attaining min-max optimal cost for $a=-1, b_1=b_2=1, \delta=\xi=2, \mu=1, \sigma=0.5, Q_C=1,$ $R_{C1}=1.304, R_{C2}=-4$, and  $P =0.4481$. Continuous solution (green). Hybrid solution (blue and red).}
    \label{fig:my_label}
\end{figure}
This corresponds to {the optimal value with every maximal solution rendered by the equilibrium $\kappa_C$} from $\xi=2$ attaining {it}.}
\IfPers{
\begin{figure}[h]
    \centering
    \includegraphics[width=8cm]{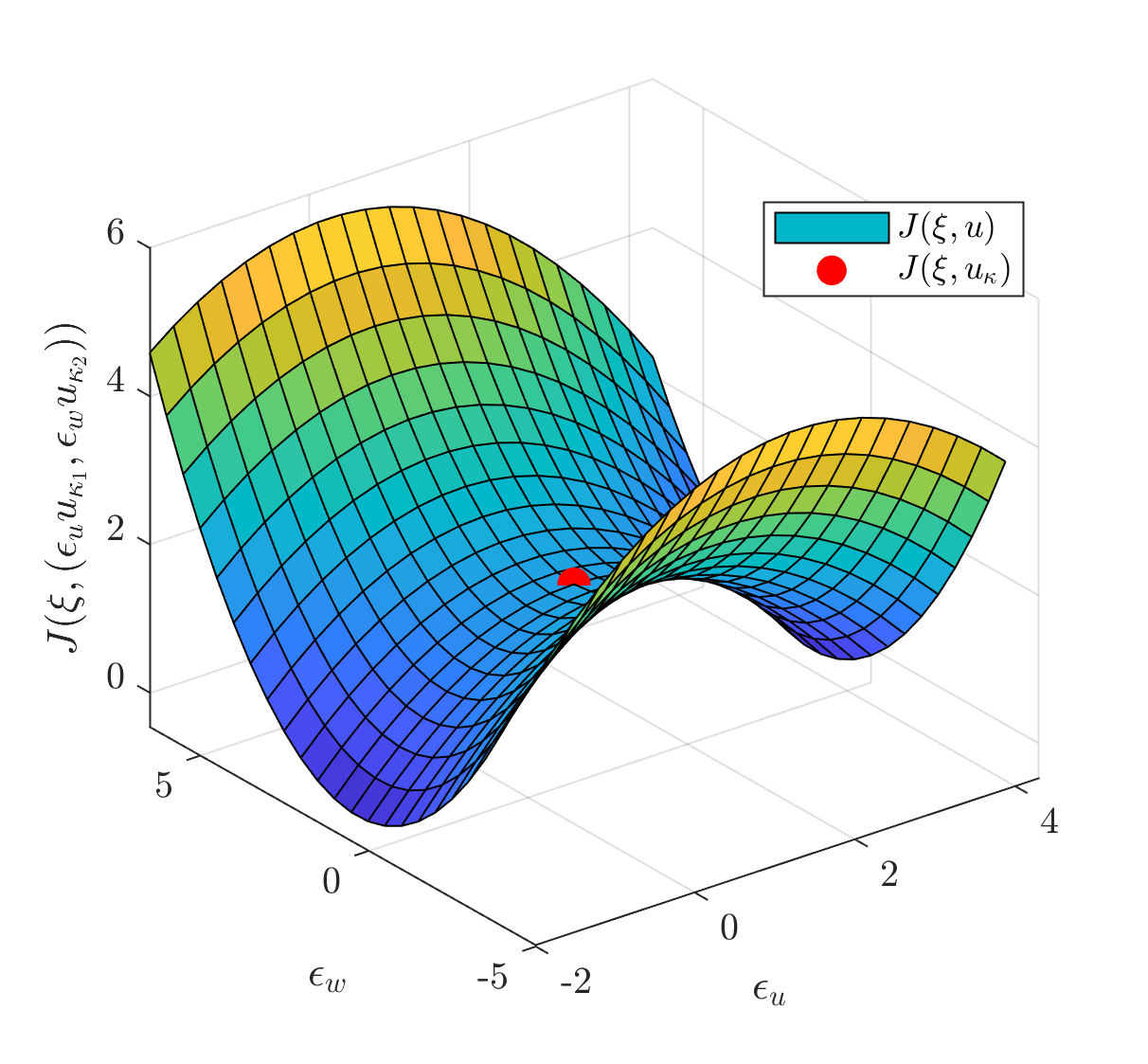}
    \caption{Saddle point behavior in the cost of continuous solutions from $\xi=2$ when varying the feedback gains around the optimal value. The cost is evaluated on solutions $(\phi,u)\in \mathcal{S}^\infty_\HS (\xi) $ with feedback law variations specified by  $\epsilon_u$ and $\epsilon_u$ in $u=(\epsilon_u \kappa_1(\phi),\epsilon_w \kappa_2(\phi))$.}
    \label{Costs}
\end{figure}
This corresponds to the saddle-point equilibrium in Definition \ref{NashEqNonCoop}, and Remark \ref{SolExistenceZ}, with every maximal solution rendered by $\kappa$ from $\xi=2$ attaining the optimal cost.}
\label{NumericalEx}
\end{example}
 }
}
{
  \section{Asymptotic Stability for Hybrid Games}
  We present a result that connects optimality and asymptotic stability for two-player zero-sum hybrid games. First, we introduce definitions of some classes of functions.
\NotACC{\begin{definition}[Class-$\mathcal{K}_\infty$ functions]
	A function $\alpha: \mathbb{R}_{\geq 0} \rightarrow  \mathbb{R}_{\geq 0}$ is a class-$\mathcal{K}_\infty$ function, also written as $\alpha \in \mathcal{K}_\infty$, if $\alpha$ is
	 zero at zero, continuous, strictly increasing, and unbounded.
\end{definition}}

\begin{definition}[Positive definite functions]
{A} function $\rho: \reals_{\geq 0}\rightarrow \reals_{\geq 0}$ is positive definite, also written as $\rho \in \mathcal{PD}$, if $\rho(s)>0$ for all $s>0$ and $\rho(0)=0$.
{A} function $\rho: \reals^n \times \reals^{m} \rightarrow \reals_{\geq 0}$  is positive definite with respect to a set $\A \subset \reals^n$, in composition with $\kappa:\mathbb{R}^n \rightarrow \mathbb{R}^{m}$, also written as $\rho \in \mathcal{PD}_{\kappa}(\A)$, if $\rho(x,\kappa(x)) > 0$ for all $x \in \reals^n\setminus \A$ and $\rho(\A, {\kappa}(\A)) = \{0\}$. 
\label{PDsets}
\end{definition}

\NotAutom{\begin{definition}[Uniform global \pnn{pre-}asymptotic stability] 
A {closed} set $\mathcal{A} \subset \mathbb{R}^n$ is uniformly globally \pnn{pre-}asymptotically stable for a {hybrid closed-loop system} $\HS_\kappa$ as in (\ref{Hkeq}) if {it is} 
\begin{itemize}
\item uniformly globally stable for $\HS_\kappa$, i.e., there exists a class-$\mathcal{K}_\infty$ function $\alpha$ such that any solution $\phi$ to $\HS_\kappa$ satisfies $|\phi(t,j)|_\mathcal{A} \leq \alpha (|\phi(0,0)|_\mathcal{A})$ for all $(t,j) \in \dom \phi$; {and}
\item uniformly globally \pnn{pre-}attractive for $\HS_\kappa$, i.e., for each $\varepsilon >0$ and $r>0$ there exists $T>0$ such that, for any solution $\phi$ to $\HS_\kappa$ with $|\phi(0,0)|_\mathcal{A} \leq r$, $(t,j) \in \dom \phi$ and $t+j \geq T$ imply   $|\phi(t,j)|_\mathcal{A} \leq \varepsilon$.
\end{itemize}
\end{definition}
}
\IfIncd{}{In the next result, we provide alternative conditions for the solution to Problem $(\diamond)$ to the conditions in Theorem \ref{thHJBszih}, {that under additional properties,} establish asymptotic stability {of a set $\mathcal{A}$}.
\begin{lemma}{Equivalent conditions}\label{Lemma:EquivCond}
 Given $\HS_{\kappa}$ as in (\ref{Hkeq})  
{with data} $(C,F,D,G)$ and {feedback} $\kappa:=(\kappa_C,\kappa_D)=((\kappa_{C1},\kappa_{C2}),(\kappa_{D1},\kappa_{D2})): \mathbb{R}^n \rightarrow \mathbb{R}^{m_C} \times \mathbb{R}^{m_D}$, if there exists a function $V: \mathbb{R}^n \rightarrow \mathbb{R}$ that is continuously differentiable on a neighborhood of $\Pi(C)$ 
 such that $C_\kappa = \Pi(C)$ {and} $D_\kappa = \Pi(D)$\footnote{The conditions  $C_\kappa = \Pi(C)$ {and} $ D_\kappa = \Pi(D)$ hold when $\kappa_C(x)\in \Pi_u(x,C)$ for all $x\in \Pi(C)$ and $\kappa_D(x)\in \Pi_u(x,D)$ for all $x\in \Pi(D)$. In words, the feedback law $\kappa$ defining the {hybrid} closed-loop system $\HS_\kappa$ does not render input actions outside $C \cup D$. }, 
then (\ref{HJBzsih}), (\ref{Bellmanzsih}), (\ref{kHJBeqzsihc}), and (\ref{kBeqzsihc}) are satisfied 
if and only if
\begin{equation}
  \IfAutom{\mathcal{L}_C(x,\kappa_C(x))}{L_C(x,\kappa_C(x))
+\left\langle \nabla V(x),F(x,\kappa_C(x)) \right\rangle}
=0 \hspace{0.5cm} \forall x \in C_\kappa ,
\label{kHJBeqihza}
\end{equation}
{
\begin{equation}
\begin{split}
  \IfAutom{&\mathcal{L}_C(x,(u_{C1},\kappa_{C2}(x)))}{
 L_C(x,(u_{C1},\kappa_{C2}(x)))
+\left\langle \nabla V(x),F(x,(u_{C1},\kappa_{C2}(x))) \right\rangle
  }
\geq 0 \hspace{0.5cm}\\ \IfAutom{&\hspace{2cm}}{}\forall (x,u_{C1}): \> (x,(u_{C1},\kappa_{C2}(x)))\in C,
\end{split}
\label{kHJBeqihzb}
\end{equation}
\begin{equation}
    \begin{split}
  \IfAutom{&\mathcal{L}_C(x,(\kappa_{C1}(x),u_{C2}))}{
 L_C(x,(\kappa_{C1}(x),u_{C2}))
+
\left\langle \nabla V(x),F(x,(\kappa_{C1}(x),u_{C2})) \right\rangle
  }
\leq 0 \hspace{0.5cm} \\ \IfAutom{&\hspace{2cm}}{} \forall (x,u_{C2}): \> (x,(\kappa_{C1}(x),u_{C2}))\in C,
\end{split}
\label{kHJBeqihzc}
\end{equation}}
\begin{equation}
  \IfAutom{\mathcal{L}_D(x,\kappa_D(x))}{
 L_D(x,\kappa_D(x))
+ V(G(x,\kappa_D(x)))}=V(x) \hspace{0.5cm} \forall x \in D_\kappa, 
\label{kBeqihzd}
\end{equation}
{
\begin{equation}
\begin{array}{r@{}l}
  \IfAutom{&\mathcal{L}_D(x,(u_{D1},\kappa_{D2}(x)))}{ L_D(x,(u_{D1},\kappa_{D2}(x)))
+ V(G(x,(u_{D1},\kappa_{D2}(x))))}  \geq V(x) \\  \IfAutom{&\hspace{1.5cm}}{} \hspace{0.4cm} \forall (x,u_{D1}):\>(x,(u_{D1},\kappa_{D2}(x))) \in D,
\end{array}
\label{kBeqihze}
\end{equation}
\begin{equation}
\begin{array}{r@{}l}
  \IfAutom{&\mathcal{L}_D(x,(\kappa_{D1}(x),u_{D2}))}{
 L_D(x,(\kappa_{D1}(x),u_{D2}))
+ V(G(x,(\kappa_{D1}(x),u_{D2})))
  }
\leq V(x) \\  \IfAutom{&\hspace{1.5cm}}{} \hspace{0.4cm} \forall (x,u_{D2}):\>(x,(\kappa_{D1}(x),u_{D2})) \in D.
\end{array}
\label{kBeqihzf}
\end{equation}}
\end{lemma}
\IfConf{{The proof is presented in the Appendix.}}{}
}
\IfIncd{}{
  \par\noindent\textbf{Proof of Lemma \ref{Lemma:EquivCond}.}
  $(\rightarrow)$
From (\ref{kHJBeqzsihc}) and (\ref{kBeqzsihc}) we have
\begin{equation}
\underset{u_C=(u_{C1},u_{C2}) \in \Pi_u(x,C)}
{ \underset{u_{C1}}{\min}\> \underset{u_{C2}}{\max}} \left\{ 
L_C(x,u_C) \right.
+  \left.\left\langle \nabla V(x),F(x,u_C) \right\rangle \right\} =L_C(x,\kappa_C(x))
+\left\langle \nabla V(x) ,F(x,\kappa_C(x)) \right\rangle  \hspace{0.49cm} \forall x \in \Pi(C)
\label{kHJBeqihzproof}
\end{equation}
and

\begin{equation}
\underset{u_D=(u_{D1},u_{D2}) \in \Pi_u(x,D)}
{ \underset{u_{D1}}{\min} \underset{u_{D2}}{\max}} \left\{ 
L_D(x,u_D)
\right. + \left. V(G(x,u_D)) \right\}=L_D(x,\kappa_D(x))
+ V(G(x,\kappa_D(x))) \hspace{1.2cm} \forall x \in \Pi(D)
\label{kBeqihzproof}
\end{equation}
Thus, (\ref{HJBzsih}) and (\ref{kHJBeqihzproof}) imply 
\begin{equation}
L_C(x,\kappa_C(x))
+\left\langle \nabla V(x) ,F(x,\kappa_C(x)) \right\rangle=0  \hspace{0.49cm} \forall x \in \Pi(C), 
\label{kHJBeqihzproofa}
\end{equation}
while, (\ref{Bellmanzsih})  and (\ref{kBeqihzproof}) imply  
\begin{equation}
L_D(x,\kappa_D(x))
+ V(G(x,\kappa_D(x))) =V(x) \hspace{0.4cm} \forall x \in \Pi(D).
\label{kBeqihzproofa}
\end{equation}
From (\ref{kHJBeqihzproofa}) and (\ref{HJBzsih}), we have 
\begin{equation}
\underset{u_{C1}:(u_{C1},\kappa_{C2}(x)) \in \Pi_u(x,C)}
{ \underset{u_{C1}}{\min} 
} \left\{ 
L_C(x,(u_{C1},\kappa_{C2}(x))) \right.
+  \left.\left\langle \nabla V(x),F(x,(u_{C1},\kappa_{C2}(x))) \right\rangle \right\}  \geq 0 \hspace{0.4cm}\forall x \in \Pi(C)
\label{kHJBeqihzproofb}
\end{equation}
and
\begin{equation}
\underset{u_{C2}:(\kappa_{C1}(x), u_{C2}) \in \Pi_u(x,C)}
{ 
\underset{u_{C2}}{\max}
} \left\{ 
L_C(x,(\kappa_{C1}(x), u_{C2})) \right.
+  \left.\left\langle \nabla V(x),F(x,(\kappa_{C1}(x), u_{C2})) \right\rangle \right\}  \leq 0 \hspace{0.4cm} \forall x \in \Pi(C)
\label{kHJBeqihzproofc}
\end{equation}
which imply (\ref{kHJBeqihzb}) and (\ref{kHJBeqihzc}), respectively. Likewise, 
From (\ref{kBeqihzproofa}) and (\ref{Bellmanzsih}), we have 
\begin{equation}
\underset{u_{D1}:(u_{D1},\kappa_{D2}(x)) \in \Pi_u(x,D)}
{ \underset{u_{D1}}{\min} 
} \left\{ 
L_D(x,(u_{D1},\kappa_{D2}(x))) \right.
+  \left.    V(G(x,(u_{D1},\kappa_{D2}(x))))   \right\} \geq V(x) \hspace{0.4cm} \forall x \in \Pi(D)
\label{kBeqihzproofb}
\end{equation}
and
\begin{equation}
\underset{u_{D2}:(\kappa_{D1}(x), u_{D2}) \in \Pi_u(x,D)}
{ 
\underset{u_{D2}}{\max}
} \left\{ 
L_D(x,(\kappa_{D1}(x), u_{D2})) \right.
+ \left.  V(G(x,(\kappa_{D1}(x), u_{D2}))) \right\} \leq V(x) \hspace{0.4cm} \forall x \in \Pi(D)
\label{kBeqihzproofc}
\end{equation}
which imply (\ref{kBeqihze}) and (\ref{kBeqihzf}), respectively.\\
$(\leftarrow)$ Given $V$ and $\kappa:=(\kappa_C,\kappa_D)=((\kappa_{C1},\kappa_{C2}),(\kappa_{D1},\kappa_{D2}))$ such that (\ref{kHJBeqihza})-(\ref{kBeqihzf}) are satisfied, and such that $C_\kappa = \Pi(C), D_\kappa = \Pi(D)$, let us prove that $V$ and $\kappa$ satisfy (\ref{HJBzsih}), (\ref{Bellmanzsih}), (\ref{kHJBeqzsihc}), and (\ref{kBeqzsihc}). From (\ref{kHJBeqihza}) and (\ref{kHJBeqihzb}) we have 
\begin{equation}
\begin{split}
\underset{u_{C1}:(u_{C1},\kappa_{C2}(x)) \in \Pi_u(x,C)}
{ \underset{u_{C1}}{\min} 
} \left\{ 
L_C(x,(u_{C1},\kappa_{C2}(x))) \right.
+  \left.\left\langle \nabla V(x),F(x,(u_{C1},\kappa_{C2}(x))) \right\rangle \right\}  
\\ =L_C(x,\kappa_C(x))
+\left\langle \nabla V(x) ,F(x,\kappa_C(x)) \right\rangle =0 \hspace{0.4cm} \forall x \in \Pi(C)
\end{split}
\label{kHJBeqihzproof1left}
\end{equation}
and from (\ref{kHJBeqihza}) and (\ref{kHJBeqihzc}) we have 
\begin{equation}
\begin{split}
\underset{u_{C2}:(\kappa_{C1}(x), u_{C2}) \in \Pi_u(x,C)}
{ 
\underset{u_{C2}}{\max}
} \left\{ 
L_C(x,(\kappa_{C1}(x), u_{C2})) \right.
+  \left.\left\langle \nabla V(x),F(x,(\kappa_{C1}(x), u_{C2})) \right\rangle \right\}  
\\=L_C(x,\kappa_C(x))
+\left\langle \nabla V(x) ,F(x,\kappa_C(x)) \right\rangle =0 \hspace{0.4cm} \forall x \in \Pi(C)
\end{split}
\label{kHJBeqihzproof2left}
\end{equation}
Thus, (\ref{kHJBeqihzproof1left}) and (\ref{kHJBeqihzproof2left}) imply (\ref{HJBzsih}) and  (\ref{kHJBeqzsihc}).
Similarly, from (\ref{kBeqihzd}) and (\ref{kBeqihze}) we have 
\begin{equation}
\begin{split}
\underset{u_{D1}:(u_{D1},\kappa_{D2}(x)) \in \Pi_u(x,D)}
{ \underset{u_{D1}}{\min} 
} \left\{ 
L_D(x,(u_{D1},\kappa_{D2}(x))) \right.
+  \left.    V(G(x,(u_{D1},\kappa_{D2}(x))))   \right\} 
\\=L_D(x,\kappa_D(x))
+ V(G(x,\kappa_D(x))) =V(x) \hspace{0.4cm} \forall x \in \Pi(D)
\end{split}
\label{kBeqihzproof1left}
\end{equation}
and from (\ref{kBeqihzd}) and (\ref{kBeqihzf}) we have 
\begin{equation}
\begin{split}
\underset{u_{D2}:(\kappa_{D1}(x), u_{D2}) \in \Pi_u(x,D)}
{ 
\underset{u_{D2}}{\max}
} \left\{ 
L_D(x,(\kappa_{D1}(x), u_{D2})) \right.
+ \left.  V(G(x,(\kappa_{D1}(x), u_{D2}))) \right\}
\\ =L_D(x,\kappa_D(x))
+ V(G(x,\kappa_D(x))) =V(x) \hspace{0.4cm} \forall x \in \Pi(D)
\end{split}
\label{kBeqihzproof2left}
\end{equation}
Thus, (\ref{kBeqihzproof1left}) and (\ref{kBeqihzproof2left}) imply (\ref{Bellmanzsih}) and  (\ref{kBeqzsihc}).
}

\begin{ttheorem}{{\pnn{Weak} saddle-point equilibrium under the existence of a Lyapunov function}}\label{zCorStability}
Consider a two-player zero-sum hybrid game with closed-loop dynamics $\HS_{\kappa}$ as in (\ref{Hkeq})  
with data $(C,F,D,G)$\IfIncd{}{ satisfying Assumption \ref{AssLipsZ}},
and $\kappa:=(\kappa_C,\kappa_D): \mathbb{R}^n \rightarrow \mathbb{R}^{m_C} \times \mathbb{R}^{m_D}$ such that $C_\kappa = \Pi(C)$ {and} $D_\kappa = \Pi(D)$. 
Given a closed set $\mathcal{A} \subset \reals^n$, continuous functions $L_C:
C \rightarrow \mathbb{R}_{\geq 0}$ and $L_D:
D \cap G(D) \rightarrow \mathbb{R}_{\geq 0}$ 
defining the stage costs for flows and jumps, respectively, {and} $q:\reals^n \rightarrow \reals$ defining the terminal cost, suppose there exists a continuous function $V: \mathbb{R}^n \rightarrow \mathbb{R}$ that is \IfGG{locally Lipschitz}{continuously differentiable} on an open set containing $\overline{C_\kappa}$,  satisfying \IfIncd{(\ref{HJBzsihinc})-(\ref{Bellmanzsihincub})}{(\ref{kHJBeqihza})-(\ref{kBeqihzf})}, and such that for each $\xi \in \overline{C_\kappa} \cup D_\kappa$, each $\phi \in \mathcal{S}_{\HS_{\kappa}  }
(\xi)$ satisfies \IfIncd{\eqref{TerminalCondCGinc}}{(\ref{TerminalCondCG})}. Furthermore, suppose that there exist $\alpha_1,\alpha_2 \in \mathcal{K}_\infty$ such that 
\begin{equation}
\alpha_1(|x|_\A) \leq V(x) \leq \alpha_2(|x|_\A) \hspace{1cm} \forall x \in \overline{C_\kappa} \cup D_\kappa \IfIncd{\cup G(D_\kappa)}{}
\label{zalphabound}
\end{equation}
and one of the following conditions\footnote{The subindexes in the sets of positive definite functions $\mathcal{PD}_*$ denote the {feedback} law that they are 
composed with as in 
Definition \ref{PDsets}.} holds
{
\begin{enumerate}[label=\arabic*)]
\item $L_C \in \mathcal{PD}_{\kappa_C}(\mathcal{A})$ and $L_D \in \mathcal{PD}_{\kappa_D}(\mathcal{A} )$;
\item $L_D \in \mathcal{PD}_{\kappa_D}(\mathcal{A})$ and there exists a continuous {function} $\eta \in \mathcal{PD}$ such that $L_C(x,\kappa_D(x)) \geq \eta(|x|_\A)$ for all $x \in C_\kappa$;
\item $L_C \in \mathcal{PD}_{\kappa_C}(\mathcal{A})$ and there exists a continuous {function} $\eta \in \mathcal{PD}$ such that $L_D(x,\kappa_D(x)) \geq \eta(|x|_\A)$ for all $x \in D_\kappa$.
\end{enumerate}}
 Then   
\begin{equation}
\J^*(\xi)= V(\xi)\hspace{2cm} {\forall \xi \in \overline{C_\kappa} \cup D_\kappa}
\end{equation}
Furthermore, the feedback law
$\kappa$ is a \pnn{weak} saddle-point equilibrium (see Definition \ref{NashEqNonCoop}) and it renders $\A$ uniformly globally \pnn{pre-}asymptotically stable 
for $\HS_{\kappa}$ as in \cite[Definition 3.6]{goebel2009hybrid}. 
\end{ttheorem}
%
%
%
{
\begin{proof}
  Since{,} by assumption{,} we have that $C_\kappa = \Pi(C), D_\kappa = \Pi(D)$, and $V,\kappa:=(\kappa_C,\kappa_D)=((\kappa_{C1},\kappa_{C2}),$ $(\kappa_{D1},\kappa_{D2}))$ are such that \IfIncd{(\ref{HJBzsihinc})-(\ref{Bellmanzsihincub})}{(\ref{kHJBeqihza})-(\ref{kBeqihzf})} hold, \IfIncd{and}{then, thanks to Lemma \ref{Lemma:EquivCond}, $V$ and $\kappa$ satisfy (\ref{HJBzsih}), (\ref{Bellmanzsih}), (\ref{kHJBeqzsihc}), and (\ref{kBeqzsihc}). 
  Since in addition,} for each $\xi \in \overline{C_\kappa} \cup D_\kappa$, each $\phi \in \mathcal{S}_{\HS_{\kappa}  }(\xi)$ satisfies \IfIncd{(\ref{TerminalCondCGinc})}{(\ref{TerminalCondCG})}, we have from Theorem \IfIncd{\ref{thHJBszihinc}}{\ref{thHJBszih}} that $V$ is the value function \NotSiam{as in (\ref{cost2gozsih}) }for $\HS_{\kappa}$ at $\overline{C_\kappa} \cup D_\kappa$ 
  and the feedback law
  $\kappa$ \IfIncd{}{with values (\ref{kHJBeqzsihc}), (\ref{kBeqzsihc})} is a \pnn{weak} saddle-point equilibrium for this game. 
Then, $V$ is a Lyapunov candidate for $\HS_{\kappa}$ \IfGG{\cite[Definition 3.17]{220}}{\cite[Definition 3.16]{65nu}} since $\overline{C_\kappa} \cup D_\kappa \cup G(D_\kappa) \subset \dom V=\reals^n$, $V$ is continuous and \IfGG{locally Lipschitz}{continuously differentiable} on an open set containing $\overline{C_\kappa}$ and \eqref{zalphabound} holds. From \IfIncd{(\ref{HJBzsihinc}) and (\ref{Bellmanzsihinc})}{(\ref{kHJBeqihza}) and (\ref{kBeqihzd})}, we have
  \IfIncd{
  \begin{subequations}    
  \begin{eqnarray}\label{kHJBeqihzabound}
  \IfGG{\dot V(x)}{\left\langle \nabla V(x),f \right\rangle}\leq- L_C(x,\kappa_C(x))\quad&\forall x \in C_\kappa\NotGG{,f \in F(x,\kappa_C(x))}
  \\
  \label{kBeqihzbound}
  V(\chi)-V(x)\leq-L_D(x,\kappa_D(x))\quad&\forall x \in D_\kappa, \chi \in G(x,\kappa_D(x)). 
  \end{eqnarray}
  \end{subequations}
  }
  {\small 
  \begin{equation}\label{kHJBeqihzabound}
  \left\langle \nabla V(x),F(x,\kappa_C(x)) \right\rangle\leq- L_C(x,\kappa_C(x))\quad\forall x \in C_\kappa,
  \end{equation}
  \begin{equation}\label{kBeqihzbound}
  V(G(x,\kappa_D(x)))-V(x)\leq-L_D(x,\kappa_D(x))\quad\forall x \in D_\kappa. 
  \end{equation}}
  Moreover, \IfSiam{
    define 
 $
    \rho(x,\kappa(x)):=
  \min\{\star_A(x), \star_B(x)\},
  x\in C_\kappa\cup D_\kappa,  
 $ 
    where if 
      \begin{enumerate}[label=\arabic*)]
  \item Item 1 above holds, then $\star_A(x)=L_C(x,\kappa_C(x))$ and $\star_B(x)=L_D(x,\kappa_D(x))$,
  \item Item 2 above holds, then $\star_A(x)=\eta(|x|_\mathcal{A})$ and $\star_B(x)=L_D(x,\kappa_D(x))$,
  \item Item 3 above holds, then $\star_A(x)=L_C(x,\kappa_C(x))$ and $\star_B(x)=\eta(|x|_\mathcal{A})$.
    \end{enumerate}
  }{{if} \pn{
  \begin{enumerate}[label=\arabic*)]
  \item 
  Item 1 above holds, {define}
   \begin{equation*}
\begin{medsize}
  \rho(x,\kappa(x)):=
  \begin{cases}
  L_C(x,\kappa_C(x)) \quad &\textup{if}\> x\in C_\kappa \setminus D_\kappa
  \\
  \min\{L_C(x,\kappa_C(x)), L_D(x,\kappa_D(x))\}\quad &\textup{if}\> x\in C_\kappa\cap D_\kappa  
  \\
    L_D(x,\kappa_D(x)) \quad &\textup{if}\> x\in D_\kappa \setminus C_\kappa
  \end{cases}
  \end{medsize}
  \end{equation*}
  \item 
  Item 2 above holds, {define}
   \begin{equation*}
\begin{medsize}
  \rho(x,\kappa(x)):=
  \begin{cases}
  \eta(|x|_\mathcal{A}) \quad &\textup{if}\> x\in C_\kappa \setminus D_\kappa
  \\
  \min\{\eta(|x|_\mathcal{A}), L_D(x,\kappa_D(x))\}\quad &\textup{if}\> x\in C_\kappa\cap D_\kappa  
  \\
    L_D(x,\kappa_D(x)) \quad &\textup{if}\> x\in D_\kappa \setminus C_\kappa
  \end{cases}
  \end{medsize}
  \end{equation*}
  \item 
 Item 3 above holds, {define}
   \begin{equation*}
\begin{medsize}
  \rho(x,\kappa(x)):=
  \begin{cases}
  L_C(x,\kappa_C(x)) \quad &\textup{if}\> x\in C_\kappa \setminus D_\kappa
 \\
  \min\{L_C(x,\kappa_C(x)), \eta(|x|_\mathcal{A})\} \quad &\textup{if}\> x\in C_\kappa\cap D_\kappa  
  \\
    \eta(|x|_\mathcal{A}) \quad &\textup{if}\> x\in D_\kappa \setminus C_\kappa
  \end{cases}
  \end{medsize}
  \end{equation*}
  \end{enumerate}}
  }
  \noindent Thus, given the functions $\alpha_1, \alpha_2$ satisfying (\ref{zalphabound}),
  and {given that from (\ref{kHJBeqihzabound}) and (\ref{kBeqihzbound}), for each case above, the continuous function $\rho$ satisfies 
   \IfIncd{
   \begin{subequations}
        \begin{eqnarray}\label{kHJBeqihzaboundrho}
    \IfGG{\dot V(x)}{\left\langle \nabla V(x),f \right\rangle}\leq- \rho(x,\kappa(x))\quad&\forall x \in C_\kappa\NotGG{, f\in F(x,\kappa_C(x))}
    \\
    \label{kBeqihzboundrho}
    V(\chi)-V(x)\leq-\rho(x,\kappa(x))\quad&\forall x \in D_\kappa,\chi \in G(x,\kappa_D(x))
    \end{eqnarray}
   \end{subequations}
    }
    {
  \begin{equation}\label{kHJBeqihzaboundrho}
  \left\langle \nabla V(x),F(x,\kappa_C(x)) \right\rangle\leq- \rho(x,\kappa(x))\quad\forall x \in C_\kappa,
  \end{equation}
  \begin{equation}\label{kBeqihzboundrho}
  V(G(x,\kappa_D(x)))-V(x)\leq-\rho(x,\kappa(x))\quad\forall x \in D_\kappa,
  \end{equation}}
  }
  thanks to \IfGG{\cite[Theorem 3.19]{220}}{\cite[Theorem 3.18]{65nu}}{, the set} $\A$ is uniformly globally \pnn{pre-}asymptotically stable for $\HS_{\kappa}$. 
  Furthermore, \IfIncd{when}{via Assumption \ref{AssLipsZ},} maximal solutions to $\HS_\kappa$ are complete, \IfIncd{we can}{which allows us to} argue uniform {global} asymptotic stability of $\A$ {as in \cite{220}}.
 \end{proof}
 }
 \IfIncd{}{
The next example illustrates Theorem \ref{zCorStability}. In addition, it shows that our results 
do not necessarily compute the value function but, in the spirit of the Lyapunov theorem, require that the conditions in Corollary \ref{zCorStability} hold. 
{
\begin{example}{Hybrid game with nonunique solutions}
From Example \ref{NumericalEx}, recall that $\kappa_C(x)=(-R_{C1}^{-1}b_1 P x,-R_{C2}^{-1}b_2 P x)$ for every $x\in \Pi(\overline{C})$. Let $\A =\{0\}$ and 
given that $L_C \in \mathcal{PD}_{\kappa_C }(\A 
)$, 
%
 (\ref{kHJBeqihza})-(\ref{kBeqihzf}) hold, and the function $s \mapsto \eta(s)=:P \frac{s^2}{2}$ is such that $L_D(x,\kappa_D(x)) \geq \eta(|x|_\A)$ for all $x \in D_\kappa$, by setting $\alpha_1(|x|_\A)=\underline{\lambda}(P)|x|^2$ and $\alpha_2(|x|_\A)=\overline{\lambda}(P)|x|^2$, 
from Corollary \ref{zCorStability} we have that $\kappa_C$ is the saddle-point equilibrium and renders $\A$ uniformly globally asymptotically stable for $\HS$ as in (\ref{nuHGamez}).
 \end{example}
 }
 }
}
{
\section{Applications and Examples}
\subsection{Robust Hybrid Linear Quadratic Problems with Spontaneous \\Jumps}
\IfIncd{
}{
\subsection{Application 1: Robust Hybrid LQR with Periodic Jumps}}
{In this section, we study a \IfIncd{hybrid game }{special case of Application \ref{RobustApp} }that emerges in practical scenarios of hybrid systems with linear flow and jump maps and \IfIncd{spontaneous}{periodic} jumps, as in noise attenuation of cyber-physical systems\IfIncd{\> with intermittent communication}{},  see, e.g., \cite{ferrante2019certifying,possieri2020linear,possieri2017mathcal}. We introduce a state variable $\tau$ that plays the role of a timer. Once $\tau$ reaches a fixed threshold $\bar{T}$, it triggers a jump in the state and resets $\tau$ to \IfIncd{a random number in $[0, \bar T]$}{zero}. More precisely,
given $\bar{T}\in \mathbb{R}$, we consider a hybrid
system with state $x=(x_p,\tau){=((x_{p1},x_{p2}),\tau)} \in \reals^n \times [0, \bar{T}]$, input $u=(u_C, u_D)=((u_{C1},u_{C2}),$
$(u_{D1},u_{D2})) \in \reals^{m_C}\times \reals^{m_D}$, and dynamics $\HS\IfIncd{}{} $ as in (\IfIncd{\ref{Hinc}}{\ref{Heq}})\IfTp{,}{ with $N=2$ and} {defined} by
\begin{equation}\label{eq: CDFG spontaneous}
\begin{array}{rll}
C&:=&\mathbb{R}^n \times [0,\bar{T}] \times \mathbb{R}^{m_C}
\\
F(x,u_C)&:=&(A_C x_p+B_C u_C,1) 
\quad \forall(x,u_C) \in C
\\
D&:=&\mathbb{R}^n \times \{\bar{T}\} \times \mathbb{R}^{m_D}
\\
G(x,u_D)&:=&(A_D x_p+B_D u_D,\IfIncd{[0,\bar T]}{0})
\quad \forall(x,u_D) \in D
\end{array}
\end{equation}
with $A_C=\left[\begin{smallmatrix}A_{C1} & 0 \\ 0 & A_{C2}\end{smallmatrix}\right], 
B_C=\left[ \begin{smallmatrix}B_{C1} \>B_{C2}\end{smallmatrix}\right],
A_D = \left[\begin{smallmatrix}A_{D1} & 0 \\ 0 & A_{D2}\end{smallmatrix}\right],$ and 
$ B_D = \left[ \begin{smallmatrix}B_{D1} \>B_{D2}\end{smallmatrix}\right] $.
\IfIncd{In this case,}{Following Application \ref{RobustApp},} the input $u_1:=(u_{C1},u_{D1})$ plays the role of the control  and is assigned by {player} $P_1$, and $u_2:=(u_{C2},u_{D2})$ is the disturbance input, which is assigned by {player} $P_2$. 
\pnn{The problem of 
 upper bounding the effect of the disturbance $u_2$ in the cost of solutions to $\HS\IfIncd{}{} $ is formulated as \pnn{a} two-player zero-sum game as in Section \ref{Sec: Game form}. Thus, by solving Problem $(\diamond)$ for every $\xi \in \Pi(C\cup D)$, 
the control objective is achieved.} 
\IfPers{Note that the control objective does not include attenuation of $u_2$.} 

The following result \pn{addresses an} 
optimal control problem for hybrid systems with linear maps and \IfIncd{spontaneous}{periodic} jumps under an {adversarial} action. \pnn{We impose a more strict condition over the regularity of the function $V$ as a special case.} 

\begin{corollary}[Hybrid Riccati equation for disturbance rejection with \IfIncd{spontaneous}{periodic} jumps] Given {a hybrid system $\HS\IfIncd{}{} $ as in (\IfIncd{\ref{Hinc}}{\ref{Heq}}) with data $(C,F,D,G)$ as in \eqref{eq: CDFG spontaneous}, let} $\bar{T}\in \mathbb{R}$, 
  and,
  $L_C(x,u_C):=x_p^\top Q_C x_p+ u_{C}^\top R_{C} u_{C}$,
  $L_D(x,u_D):=x_p^\top Q_D x_p+  u_{D}^\top R_{D} u_{D}$,
   and terminal cost $q(x):=x_p^\top P(\tau) x_p$ defining $\J$ as in (\ref{defJTNCinc}), {with} 
   $Q_C,Q_D\in \mathbb{S}^n_+$, 
   $R_C=\left[\begin{smallmatrix}R_{C1} & 0 \\ 0 & R_{C2}\end{smallmatrix}\right],
   R_D=\left[\begin{smallmatrix}R_{D1} & 0 \\ 0 & R_{D2}\end{smallmatrix}\right],$  
   $R_{C1} \in \mathbb{S}^{m_{C_1}}_+$, $-R_{C2} \in \mathbb{S}^{m_{C_2}}_+$, $R_{D1} \in \mathbb{S}^{m_{D_1}}_+$, and $-R_{D2} \in \mathbb{S}^{m_{D_2}}_+$. Suppose there exists a matrix {function} $P:[0,\bar{T}]\rightarrow \mathbb{S}^n_+ $ that is \IfIncd{nonincreasing and }{}\pno{continuously differentiable} {and} such that
  \begin{equation}
   \begin{split}
  -\frac{\mathrm{d} P(\tau)}{\mathrm{d} \tau} =
  -P(\tau)
  B_{C} R_{C}^{-1}B_{C}^\top
  P(\tau)
  +Q_C+P(\tau) A_C+A_C^\top P(\tau) \hspace{0.5cm} \forall \tau  \in (0,\bar{T}),
   \end{split}
  \label{DiffRiccatizsih}
  \end{equation}
  \begin{equation}
  \begin{split}
  -R_{D2} -B_{D2}^\top P(0) B_{D2},\>&
  R_{D1} +B_{D1}^\top P(0) B_{D1}
  \in \mathbb{S}_{0+}^{m_D},
  \end{split}
  \label{zlqeqinv}
  \end{equation}
  the matrix $R_v=
  R_D+B_D^\top P(0)B_D$
  is invertible, and 
  \begin{equation}
  P (\bar{T})=Q_D+A_D^\top P(0) A_D \IfConf{\\}{}
  -
  A_D^\top P(0)B_D
  R_v
  ^{-1}
  B_D^\top P(0)A_D
  \label{Riccatizsih} 
  \end{equation}
  { where $A_C, B_{C1}, B_{C2}, A_D, B_{D1}$, and $B_{D2}$ are defined below \eqref{eq: CDFG spontaneous}.}
  Then, the feedback law $\kappa:=(\kappa_C,\kappa_D)$,  {defined as}
  \begin{subequations}
  \begin{eqnarray}
  \kappa_C(x)=
- R_C^{-1}
B_C^\top
  P(\tau) x_p
  \hspace{1cm}
   \forall x \in \Pi(C),
  \label{NashkCLQzsih}
  \\
  \kappa_D(x)=-R_v
  ^{-1}
  B_D^\top P(0) A_D x_p \hspace{1cm} \forall x \in \Pi(D)
  \label{NashkDLQzsih}
  \end{eqnarray}
  \label{NashkLQzsih}
  \end{subequations}
  is the \pnn{weak} saddle-point equilibrium for the two-player zero-sum hybrid game with \IfIncd{spontaneous}{periodic} jumps. In addition, for each $x=(x_p,\tau) \in \Pi ( \overline{C}) \cup \Pi(D)$, the value function is equal to $V(x):=x_p^\top P(\tau) x_p$.
  \label{HREqPJ}
  \end{corollary}
  {\begin{proof}
    We show that when conditions (\ref{DiffRiccatizsih})-(\ref{Riccatizsih}) hold, by using Theorem \IfIncd{\ref{thHJBszihinc}}{\ref{thHJBszih}}, the value function is equal to the function $V$ and with 
     the feedback law $\kappa:=(\kappa_C,\kappa_D)$ with values as in (\ref{NashkCLQzsih}) and (\ref{NashkDLQzsih}), such a cost is attained. 
     \IfIncd{A sufficient condition for (\ref{HJBzsihinc})-(\ref{HJBzsihincub})}{We can write (\ref{HJBzsih})} in Theorem \IfIncd{\ref{thHJBszihinc} for the single-valued flow map  $F$ as in \eqref{eq: CDFG spontaneous} is}{\ref{thHJBszih} as}
    \begin{eqnarray}
    &&0=
    \underset{u_C=(u_{C1},u_{C2}) \in \Pi_u(x,C)}
    {\min_{u_{C1}} \max_{u_{C2}}}  \mathcal{L}_C(x,u_C),
     \nonumber
    \\ 
       &&\mathcal{L}_C(x,u_C)= 
     x_p^\top Q_C x_p+ u_{C1}^\top R_{C1} u_{C1}
     + u_{C2}^\top R_{C2} u_{C2}
    \nonumber \\
    &&
    \hspace{2cm} 
    +2x_p^\top P(\tau) (A_C x_p+B_C u_C)  
    + x_p^\top \frac{\mathrm{d} P(\tau)}{\mathrm{d} \tau} x_p
    \label{HJBzsihlq}
    \end{eqnarray}
    First, given that (\ref{DiffRiccatizsih}) holds, and $x_p^\top(P(\tau) A_C+A_C^\top P(\tau))x_p={2x_p^\top P(\tau)A_Cx_p}$ for every $x\in \Pi(C)$, one has
    \begin{equation}
    \begin{split}
    \mathcal{L}_C(x,u_C)= 
     x_p^\top P(\tau)(  B_{C2} R_{C2}^{-1}B_{C2}^\top+B_{C1} R_{C1}^{-1}B_{C1}^\top)P(\tau)  x_p
     \nonumber \\
     + u_{C1}^\top R_{C1} u_{C1}+ u_{C2}^\top R_{C2} u_{C2}
    + 2 x_p^\top  P(\tau) B_C u_C 
    \end{split}
    \end{equation}
    The {first-order} necessary conditions for optimality
    {
    \IfAutom{
    \begin{equation*}
    \frac{\partial}{\partial u_{C1}}
    \mathcal{L}_C(x,u_C)  \Big|_{u_{C}^*} =0,
    \quad
    \frac{\partial}{\partial u_{C2}}
    \mathcal{L}_C(x,u_C) \Big|_{u_{C}^*} =0
    \end{equation*}
    {for all $(x,u_C)\in C$}}
    {\begin{multline*}
    \frac{\partial}{\partial u_{C1}} \left( 
     x_p^\top P(\tau)(  B_{C2} R_{C2}^{-1}B_{C2}^\top+B_{C1} R_{C1}^{-1}B_{C1}^\top)P(\tau)  x_p\right. \\ + u_{C1}^\top R_{C1} u_{C1}+ u_{C2}^\top R_{C2} u_{C2}
    \\\left.\left.+ 2 x_p^\top  P(\tau)  (B_{C1} u_{C1} + B_{C2} u_{C2}) \right)\right|_{u_{C}^*} =0
    \end{multline*}
    \begin{multline*}
    \frac{\partial}{\partial u_{C2} }\left( 
     x_p^\top P(\tau)(  B_{C2} R_{C2}^{-1}B_{C2}^\top+B_{C1} R_{C1}^{-1}B_{C1}^\top)P(\tau)  x_p\right.+\\ u_{C1}^\top R_{C1} u_{C1}+ u_{C2}^\top R_{C2} u_{C2}
     \\ \left.\left.+ 2 x_p^\top  P(\tau) (B_{C1} u_{C1} + B_{C2} u_{C2}) \right)\right|_{u_{C}^*} =0
    \end{multline*}}
    }
    are satisfied by the point $u_C^*=(u_{C1}^*,u_{C2}^*)$, with values
    \begin{equation}
    u_{C1}^*=-R_{C1}^{-1}B_{C1}^\top P(\tau) x_p, 
    \quad
    u_{C2}^*=-R_{C2}^{-1}B_{C2}^\top P(\tau) x_p
    \label{zLQHPu2*}
    \end{equation}
    {for each $x=(x_p,\tau)\in \Pi(C)$.}
    Given that $R_{C1},-R_{C2} \in \mathbb{S}^{m_D}_+$, the second-order sufficient conditions for optimality 
    %
    {
    \IfAutom{
    \begin{equation*}
    \frac{\partial^2}{\partial u_{C1}^2} 
     \mathcal{L}_C(x,u_C) \Big|_{u_C^*} \succeq 0,
     \quad
    \frac{\partial^2}{\partial u_{C2}^2}
     \mathcal{L}_C(x,u_C)\Big|_{u_C^*} \preceq 0,
    \end{equation*}
    {for all $(x,u_C)\in C$}
    }
    {
    \begin{multline*}
    \frac{\partial^2}{\partial u_{C1}^2} \left( 
     x_p^\top P(\tau)(  B_{C2} R_{C2}^{-1}B_{C2}^\top+B_{C1} R_{C1}^{-1}B_{C1}^\top)P(\tau)  x_p\right. \\ + u_{C1}^\top R_{C1} u_{C1}+ u_{C2}^\top R_{C2} u_{C2}
     \\\left.\left.
    + 2 x_p^\top  P(\tau)  (B_{C1} u_{C1} + B_{C2} u_{C2}) \right)\right|_{u_{C1}^*} \succeq 0
    \end{multline*}
    \begin{multline*}
    \frac{\partial^2}{\partial u_{C2}^2 }\left( 
     x_p^\top P(\tau)(  B_{C2} R_{C2}^{-1}B_{C2}^\top+B_{C1} R_{C1}^{-1}B_{C1}^\top)P(\tau)  x_p\right. \\ + u_{C1}^\top R_{C1} u_{C1}+ u_{C2}^\top R_{C2} u_{C2}
     \\\left.\left.+ 2 x_p^\top  P(\tau) (B_{C1} u_{C1} + B_{C2} u_{C2}) \right)\right|_{u_{C2}^*} \preceq 0
    \end{multline*}}
    }
    hold, rendering $u_{C}^*$ as in 
    \eqref{zLQHPu2*} as an optimizer of the min-max problem in (\ref{HJBzsihlq}).
    In addition, it satisfies $\mathcal{L}_C(x,u_C^*)=0$, making $V(x)=x_p^\top P(\tau)x_p$ a solution {to} \IfIncd{(\ref{HJBzsihinc})-(\ref{HJBzsihincub})}{(\ref{HJBzsih})} in Theorem \IfIncd{\ref{thHJBszihinc}}{\ref{thHJBszih}}. 
    
    On the other hand, \IfIncd{given that $P$ is nonincreasing, then $$\underset{t \in [0, \bar T]}{\sup} (A_D x_p+B_D u_D)^\top P(t) (A_D x_p+B_D u_D) = (A_D x_p+B_D u_D)^\top P(0) (A_D x_p+B_D u_D) $$ 
    and}{} we can write \IfIncd{\eqref{Bellmanzsihinc}-\eqref{Bellmanzsihincub}}{(\ref{Bellmanzsih})} in Theorem \IfIncd{\ref{thHJBszihinc}}{\ref{thHJBszih}} as
    \IfIncd{
      \begin{equation}
       x_p^\top P(\bar{T}) x_p= \mathcal{L}_D(x,\kappa_D(x)) 
          \hspace{0.6cm} \forall x: (x,\kappa_D(x)) \in D,
        \label{Bellmanzsihincp1}
        \end{equation}
          \begin{equation}
                  \begin{split}
       x_p^\top P(\bar{T}) x_p\leq  &\> \mathcal{L}_D(x,(u_{D1},\kappa_{D2}(x)))
                    \hspace{0.4cm} \forall (x,u_{D1}):(x,(u_{D1},\kappa_{D2}(x)))) \in D ,
                  \end{split}
                 \label{Bellmanzsihinclbp1}
              \end{equation}
        \begin{equation}
            \begin{split}
       x_p^\top P(\bar{T}) x_p\geq  &\> \mathcal{L}_D(x,(\kappa_{D1}(x),u_{D2}))
               \hspace{0.4cm} \forall (x,u_{D2}):(x,(\kappa_{D1}(x),u_{D2})) \in D ,
            \end{split}
           \label{Bellmanzsihincubp1}
        \end{equation}}{}
           \begin{equation}
            \begin{split}
    \IfIncd{}{&x_p^\top P(\bar{T}) x_p =\underset{u_D=(u_{D1},u_{D2}) \in \Pi_u(x,D)}
    {\min_{u_{D1}} \max_{u_{D2}}}  \mathcal{L}_D(x,u_D),
     \nonumber
    \\
    }
    &\mathcal{L}_D(x,u_D)= 
    x_p^\top Q_D x_p+ u_{D1}^\top R_{D1} u_{D1}+ u_{D2}^\top R_{D2} u_{D2} \hspace{0.4cm}
    \nonumber\\
    &\hspace{1.3cm} +
     (A_D x_p+B_D u_D)^\top P(0) (A_D x_p+B_D u_D) 
            \end{split}
     \label{Bellmanzsihlq}
    \end{equation}
    \NotAutom{which can be expanded as
    {\small
    \begin{multline}
    \mathcal{L}_D(x,u_D)=  
     x_p^\top (Q_D+A_D^\top P(0) A_D) x_p + 2 x_p^\top A_D^\top  P(0) B_D u_D \\
     + u_{D1}^\top (R_{D1} +B_{D1}^\top P(0) B_{D1}) u_{D1}
     + u_{D2}^\top (R_{D2} +B_{D2}^\top P(0) B_{D2}) u_{D2}\\
     + u_{D1}^\top (B_{D1}^\top P(0) B_{D2}) u_{D2}
     + u_{D2}^\top (B_{D2}^\top P(0) B_{D1}) u_{D1}
    \end{multline}}
    }
    {Similar to the case along flows, \IfIncd{a sufficient condition for \eqref{Bellmanzsihincp1}-\eqref{Bellmanzsihincubp1} is 
    \begin{equation}
    x_p^\top P(\bar{T}) x_p =\underset{u_D=(u_{D1},u_{D2}) \in \Pi_u(x,D)}
    {\min_{u_{D1}} \max_{u_{D2}}}  \mathcal{L}_D(x,u_D),
    \label{minmaxsp}
  \end{equation}
    and }{}the first-order} necessary conditions for optimality
    \NotAutom{
    \begin{equation*}
    \frac{\partial}{\partial u_{D1}}\left.
    \mathcal{L}_D(x,u_D)  \right|_{u_{D}^*} =0,
    \quad
    \frac{\partial}{\partial u_{D2}}\left.
    \mathcal{L}_D(x,u_D) \right|_{u_{D}^*} =0
    \end{equation*}
    }
    are satisfied by the point $u_D^*=(u_{D1}^*,u_{D2}^*)$, such that, for each $x_p \in \Pi(D)$,
    %
    \begin{equation}
    u_{D}^*= 
    {
    -    \begin{bmatrix}
    R_{D1}+B_{D1}^\top P(0) B_{D1}
    &
    B_{D1}^\top P(0) B_{D2}
     \\
    B_{D2}^\top P(0) B_{D1}
    &
    R_{D2}+B_{D2}^\top P(0) B_{D2}
    \end{bmatrix} 
    ^{-1}
    \begin{bmatrix}
    B_{D1}^\top P(0) A_D
    \\
    B_{D2}^\top P(0) A_D
    \end{bmatrix}} x_p
    \label{zLQHPud*}
    \end{equation}
    Given that (\ref{zlqeqinv}) holds, the second-order sufficient conditions for optimality 
    \NotAutom{
    \begin{equation*}
    \frac{\partial^2}{\partial u_{D1}^2}\left. 
     \mathcal{L}_D(x,u_D) \right|_{u_D^*} \succeq 0,
     \quad
    \frac{\partial^2}{\partial u_{D2}^2}\left.
     \mathcal{L}_D(x,u_D)\right|_{u_D^*} \preceq 0,
    \end{equation*}
    }
    are satisfied, rendering $u_D^*$ as in (\ref{zLQHPud*}) as an optimizer of the {min-max} problem in (\ref{minmaxsp}).
    In addition, $u_D^*$ satisfies $\mathcal{L}_D(x,u_D^*)=x_p^\top P(\bar{T})x_p$, with $P(\bar{T})$ as in (\ref{Riccatizsih}), making $V(x)=x_p^\top P(\tau)x_p$ \IfIncd{satisfy \eqref{Bellmanzsihinc}-\eqref{Bellmanzsihincub}}{a solution of (\ref{Bellmanzsih})} in Theorem \IfIncd{\ref{thHJBszihinc}}{\ref{thHJBszih}}. 
    
    Then, given that $V$ is \pno{continuously differentiable on a neighborhood of $\Pi(C)$} {and that $q(x)=V(x)$ for all $x \in \dom V$, satisfying \eqref{TerminalCondCGinc}}\IfIncd{}{ and {that} Assumption \ref{AssLipsZ} holds}, by applying Theorem \IfIncd{\ref{thHJBszihinc}}{\ref{thHJBszih}}, in particular from \IfIncd{\eqref{ResultValueinc}}{(\ref{ResultValue})}, for every $ \xi=(\xi_p,\xi_\tau) \in \Pi(\overline{C}\cup D)$ the value function is $\J^*(\xi)=\J(\xi,((u_{C1}^*,u_{D1}^*),(u_{C2}^*,u_{D2}^*))= \xi_p^\top P(\xi_\tau)\xi_p$. \IfIncd{When}{From (\ref{kHJBeqzsihc}) and (\ref{kBeqzsihc}), when} $P_1$ plays 
    $\kappa_1=(\kappa_{C1},\kappa_{D1})$ with values as in \IfIncd{(\ref{NashkLQzsih})}{(\ref{NashkCLQzsih}) and (\ref{NashkDLQzsih})}, and $P_2$ plays any disturbance $u_2$ \pn{such that solutions to $\HS\IfIncd{}{}$ 
    are complete,} then the cost is upper bounded by $\J(\xi,u^*)$,  
satisfying (\ref{SaddlePointIneq}).
\end{proof}
  }
Notice that the \pnn{weak} saddle-point equilibrium $\kappa=(\kappa_1, \kappa_2)$ is composed of the upper bounding strategy $\kappa_1$, and the maximizing disturbance $\kappa_2$ with values as in \eqref{NashkLQzsih}.
}
\NotAutom{
We study a special case of Application \ref{RobustApp} and apply Theorem \ref{thHJBszih} in this section. Consider a hybrid system with state $x\in \reals^n $, input $u=(u_C, u_D)=((u_{C1},u_{C2}), (u_{D1},u_{D2})) \in \reals^{m_C}\times \reals^{m_D}$, and dynamics $\HS$ as in (\ref{Heq})\IfTp{,}{ with $N=2$ and} described by
\begin{equation}\label{eq: CDFG RoHyLQR}
\begin{array}{rll}
C&\subset&\mathbb{R}^n  \times \mathbb{R}^{m_C} 
\\
F(x,u_C)&:=&A_C x+B_C u_C \\ &=&\left[\begin{smallmatrix}A_{C1} & 0 \\ 0 & A_{C2}\end{smallmatrix}\right]\left[ \begin{smallmatrix}x_1\\ x_2\end{smallmatrix}\right] +\left[ \begin{smallmatrix}B_{C1} \>B_{C2}\end{smallmatrix}\right]\left[ \begin{smallmatrix}u_{C1}\\ u_{C2}\end{smallmatrix}\right]
\\
D&\subset&\mathbb{R}^n \times  \mathbb{R}^{m_D} 
\\
G(x,u_D)&:=&A_D x+B_D u_D\\ &=&\left[\begin{smallmatrix}A_{D1} & 0 \\ 0 & A_{D2}\end{smallmatrix}\right]\left[ \begin{smallmatrix}x_1\\ x_2\end{smallmatrix}\right] +\left[ \begin{smallmatrix}B_{D1} \>B_{D2}\end{smallmatrix}\right]\left[ \begin{smallmatrix}u_{D1}\\ u_{D2}\end{smallmatrix}\right]
\end{array}
\end{equation}
where 
$C\cup D$ is nonempty. Following Application \ref{RobustApp}, the input $u_1:=(u_{C1},u_{D1})$ plays the role of the control and $u_2:=(u_{C2},u_{D2})$ is the disturbance input. 
\pno{The problem of 
 upper bounding the effect of the disturbance $u_2$ in the cost of solutions to $\HS$ is formulated as the two-player zero-sum game as in Section \ref{Sec: Game form}.} Thus, \textcolor{black}{by solving Problem $(\diamond)$ for every $\xi \in \Pi(C\cup D)$,} 
the control objective is achieved. 
\IfPers{Note that the control objective does not include attenuation of $u_2$.} 

The following result presents a tool for the solution of the optimal control problem for hybrid systems with linear maps under the presence of disturbances.
\begin{corollary}{Hybrid Riccati equation for disturbance rejection} 
Given {a hybrid system $\HS$ as in \eqref{Heq} with data $(C,F,D,G)$ as in \eqref{eq: CDFG RoHyLQR}, }
and,
{with the aim of pursuing minimum energy and distance to the origin, consider the cost functions 
 $L_C(x,u_C):=x^\top Q_C x+ u_{C1}^\top R_{C1} u_{C1}+ u_{C2}^\top R_{C2} u_{C2}$, $L_D(x,u_D):=x^\top Q_D x+  u_{D1}^\top R_{D1} u_{D1}+ u_{D2}^\top R_{D2} u_{D2} $, and terminal cost $q(x):=x^\top P x$, {defining $\mathcal{J}$ as in (\ref{defJTNCinc}),}
{with} $Q_C,Q_D\in \mathbb{S}^n_+$, $R_{C1} \in \mathbb{S}^{m_{C_1}}_+$, $-R_{C2} \in \mathbb{S}^{m_{C_2}}_+$, $R_{D1} \in \mathbb{S}^{m_{D_1}}_+$, $-R_{D2} \in \mathbb{S}^{m_{D_2}}_+$   and  $P\in \mathbb{S}^n_+ $. }
Suppose there exists a matrix $P\in \mathbb{S}^n_+ $ such that
%
\begin{multline}
0 =-P(  B_{C2} R_{C2}^{-1}B_{C2}^\top+B_{C1} R_{C1}^{-1}B_{C1}^\top)P\\+Q_C+P A_C+A_C^\top P,
\label{DiffRiccatizsihrlqr} 
\end{multline}
\begin{equation}
\begin{array}{rrl}
-R_{D2} -B_{D2}^\top P B_{D2},\>&
R_{D1} +B_{D1}^\top P B_{D1}&
\in \mathbb{S}_{0+}^{m_D},
\end{array}
\label{zlqeqinvrlqr}
\end{equation}
{the matrix $R_v=\left[\begin{smallmatrix}
R_{D1}+B_{D1}^\top PB_{D1}
&
B_{D1}^\top P B_{D2}
 \\
B_{D2}^\top P B_{D1}
&
R_{D2}+B_{D2}^\top P B_{D2}
\end{smallmatrix}\right] $
is invertible,} and
\begin{equation}
\begin{split}
0 =-P+Q_D+A_D^\top P A_D \hspace{3.1cm}\IfConf{\\}{}
-
\begin{bmatrix}
A_D^\top P B_{D1}
 &
A_D^\top P B_{D2}
\end{bmatrix}
R_v
^{-1}
\begin{bmatrix}
B_{D1}^\top P A_D
\\
B_{D2}^\top P A_D
\end{bmatrix}
\end{split}
\label{Riccatizsihrlqr} 
\end{equation}
Then, {for} the feedback law $\kappa_1:=(\kappa_{C1},\kappa_{D1})$ with values\footnote{The notation $R_v^{-1}(p,q)$ denotes the $(p,q)$ entry of the matrix $R_v^{-1}$.}
\begin{equation}
\kappa_{C1}(x)=-R_{C1}^{-1}B_{C1}^\top P x
\hspace{1cm} 
\forall x \in \Pi(C),
\label{NashkCLQzsihrlqr}
\end{equation}
{\small
\begin{equation}
\kappa_{D1}(x)=-[R_v^{-1}(1,1) \> R_v^{-1}(1,2)]
\left[\begin{smallmatrix}
B_{D1}^\top P A_D
\\
B_{D2}^\top P A_D
\end{smallmatrix}\right] x \hspace{0.4cm} \forall x \in \Pi(D),
\label{NashkDLQzsihrlqr}
\end{equation}}
\noindent the cost of complete solutions to $\HS$ from $\xi$ in the presence of any disturbance $u_2$ is upper bounded by $\xi^\top P\xi$. In addition, for each $x \in \Pi ( \overline{C}) \cup \Pi(D)$, the value function is equal to $V(x):=x^\top P x$, and the maximizing disturbance is  given by $\kappa_2:=(\kappa_{C2},\kappa_{D2})$, with values
\begin{equation}
\kappa_{C2}(x)=-R_{C2}^{-1}B_{C2}^\top P x 
\hspace{1cm} 
\forall x \in \Pi(C),
\label{NashkC2LQzsihrlqr}
\end{equation}
{\small
\begin{equation}
\kappa_{D2}(x)=-[R_v^{-1}(2,1) \> R_v^{-1}(2,2)]
\left[\begin{smallmatrix}
B_{D1}^\top P A_D
\\
B_{D2}^\top P A_D
\end{smallmatrix}\right] x \hspace{0.4cm} \forall x \in \Pi(D).
\label{NashkD2LQzsihrlqr}
\end{equation}}
\label{HREqRobust}
\end{corollary}
{
  \begin{proof}
We show that when conditions (\ref{DiffRiccatizsihsec})-(\ref{Riccatizsihsec}) hold, by using the result in Theorem \ref{thHJBszih},  the value function is equal to the function $V$ and under
 the feedback law as in (\ref{NashkCLQzsihrlqr}) and (\ref{NashkDLQzsihrlqr}) such a cost is attained in the presence of the maximizing disturbance given by (\ref{NashkC2LQzsihrlqr}) and (\ref{NashkD2LQzsihrlqr}). 
We can write (\ref{HJBzsih}) in Theorem \ref{thHJBszih} as
\begin{eqnarray}
&&0=
\underset{u_C=(u_{C1},u_{C2}) \in \Pi_u(x,C)}
{\min_{u_{C1}} \max_{u_{C2}}}  \mathcal{L}_C(x,u_C),
 \nonumber
\\ 
&&\mathcal{L}_C(x,u_C)= 
 x^\top Q_C x+ u_{C1}^\top R_{C1} u_{C1}+ u_{C2}^\top R_{C2} u_{C2}
 \nonumber \\
 &&\hspace{3.7cm}+2x^\top P (A_C x+B_C u_C)  
\label{HJBzsihlqrlqr}
\end{eqnarray}
First, 
given that (\ref{DiffRiccatizsihrlqr}) holds, and $x^\top(P A_C+A_C^\top P)x=x^\top (2PA_C)x$ for every $x\in \Pi(C)$, one has
\begin{eqnarray*}
&&\mathcal{L}_C(x,u_C)=
 x^\top P(  B_{C2} R_{C2}^{-1}B_{C2}^\top +B_{C1} R_{C1}^{-1}B_{C1}^\top) P x \hspace{0.8cm}
 \nonumber \\
  &&\hspace{1.6cm}
 + u_{C1}^\top R_{C1} u_{C1}+ u_{C2}^\top R_{C2} u_{C2}
+ 2 x^\top  P B_C u_C
\end{eqnarray*}
The first order necessary conditions for optimality
\NotAutom{\begin{equation*}
\frac{\partial}{\partial u_{C1}} \mathcal{L}_C(x,u_C)
\left.\right|_{u_{C}^*} 
=0,
\quad
\frac{\partial}{\partial u_{C2} } \mathcal{L}_C(x,u_C)
\left.\right|_{u_{C}^*} 
=0
\end{equation*}}
are satisfied by the point $u_C^*=(u_{C1}^*,u_{C2}^*)$, with values at any $x \in \Pi(C)$
\begin{equation}
u_{C1}^*=-R_{C1}^{-1}B_{C1}^\top P x,
\quad
u_{C2}^*=-R_{C2}^{-1}B_{C2}^\top P x.
\label{zLQHPu2*rlqr}
\end{equation}
Given that $R_{C1},-R_{C2} \in \mathbb{S}^{m_D}_+$, the second-order sufficient conditions for optimality 
\NotAutom{
\begin{equation*}
\frac{\partial^2}{\partial u_{C1}^2}  \mathcal{L}_C(x,u_C)
\left.\right|_{u_{C}^*}
 \succeq 0,
 \quad
\frac{\partial^2}{\partial u_{C2}^2 } \mathcal{L}_C(x,u_C)
\left.\right|_{u_{C}^*} 
\preceq 0
\end{equation*}}
hold, rendering $u_{C}^*$ with values as in 
(\ref{zLQHPu2*rlqr}) for each $x\in \Pi(C)$ as an optimizer of the min-max problem in (\ref{HJBzsihlqrlqr}).
In addition, 
it satisfies $\mathcal{L}_C(x,u_C^*)=0$, making $V(x)=x^\top P x$ a solution to (\ref{HJBzsih}) in Theorem \ref{thHJBszih}. 

On the other hand, we can write (\ref{Bellmanzsih}) in Theorem \ref{thHJBszih} as
\begin{eqnarray}
&&x^\top P x =\underset{u_D=(u_{D1},u_{D2}) \in \Pi_u(x,D)}
{\min_{u_{D1}} \max_{u_{D2}}}  \mathcal{L}_D(x,u_D),
 \nonumber
\\ 
&&\mathcal{L}_D(x,u_D)= 
x^\top Q_D x+ u_{D1}^\top R_{D1} u_{D1}+ u_{D2}^\top R_{D2} u_{D2}\nonumber\\
&&\hspace{1.7cm}
+ (A_D x+B_D u_D)^\top P (A_D x+B_D u_D) 
\label{Bellmanzsihlqrlqr}
\end{eqnarray}
\NotAutom{which can be expanded as 
{\scriptsize
\begin{multline}
\mathcal{L}_D(x,u_D)=  
 x^\top (Q_D+A_D^\top P A_D) x + 2 x^\top A_D^\top  P B_D u_D \\
 + u_{D1}^\top (R_{D1} +B_{D1}^\top P B_{D1}) u_{D1}
 + u_{D2}^\top (R_{D2} +B_{D2}^\top P B_{D2}) u_{D2}\\
 + u_{D1}^\top (B_{D1}^\top P B_{D2}) u_{D2}
 + u_{D2}^\top (B_{D2}^\top P B_{D1}) u_{D1}
\end{multline}}
}
The first order necessary conditions for optimality
\NotAutom{\begin{equation*}
\frac{\partial}{\partial u_{D1}}\left.
\mathcal{L}_D(x,u_D)  \right|_{u_{D}^*} 
=0,
\quad
\frac{\partial}{\partial u_{D2}}\left.
\mathcal{L}_D(x,u_D) \right|_{u_{D}^*} 
=0
\end{equation*}
}
are satisfied by the point $u_D^*=(u_{D1}^*,u_{D2}^*)$, with values for each $x \in \Pi(D)$
\begin{equation}
{\tiny
u_{D}^*=-\begin{bmatrix}
R_{D1}+B_{D1}^\top P B_{D1}
&
B_{D1}^\top P B_{D2}
 \\
B_{D2}^\top P B_{D1}
&
R_{D2}+B_{D2}^\top P B_{D2}
\end{bmatrix} 
^{-1}
\begin{bmatrix}
B_{D1}^\top P A_D
\\
B_{D2}^\top P A_D
\end{bmatrix} x_p
}
\label{zLQHPud*rlqr}
\end{equation}

Given that (\ref{zlqeqinvrlqr}) holds, the second-order sufficient conditions for optimality of $u_D^*$\NotAutom{, namely
\begin{equation*}
\frac{\partial^2}{\partial u_{D1}^2}\left. 
 \mathcal{L}_D(x,u_D) \right|_{u_D^*} \succeq 0,
\quad
\frac{\partial^2}{\partial u_{D2}^2}\left.
 \mathcal{L}_D(x,u_D)\right|_{u_D^*} \preceq 0,
\end{equation*}}
are satisfied, rendering $u_D^*$ in (\ref{zLQHPud*rlqr}) as an optimizer of the min-max problem in (\ref{Bellmanzsihlqrlqr}).
In addition, $u_D^*$
satisfies $\mathcal{L}_D(x,u_D^*)=x^\top P x$, with $P$ as in (\ref{Riccatizsihrlqr}), making $V(x)=x^\top P x$ a solution to (\ref{Bellmanzsih}) in Theorem \ref{thHJBszih}. 

Thus, given that $V$ is continuously differentiable in $\reals^n$  and Assumption \ref{AssLipsZ} holds,
 by applying Theorem \ref{thHJBszih}, in particular from (\ref{ResultValue}), {for every $ \xi \in \Pi(\overline{C}\cup D)$ the value function is $\J^*(\xi)=\J(\xi,((u_{C1}^*,u_{D1}^*),(u_{C2}^*,u_{D2}^*))= \xi^\top P \xi$. From (\ref{kHJBeqzsihc}) and (\ref{kBeqzsihc}), when $P_1$ plays $u_1^*$ defined by $\kappa_1=(\kappa_{C1},\kappa_{D1})$ with values as in (\ref{NashkCLQzsihrlqr}), (\ref{NashkDLQzsihrlqr}), and $P_2$ plays any disturbance $u_2$ such that solutions to $\HS$ with data as in \eqref{eq: CDFG periodic} are complete, then the cost is upper bounded by $\J(\xi,u^*)$,  
satisfying (\ref{SaddlePointIneq}).
 }
\end{proof}}
Notice that the \pnn{weak} saddle-point equilibrium $\kappa=(\kappa_1, \kappa_2)$ is composed of $P_1$ playing the upper bounding strategy $\kappa_1$ with values as in (\ref{NashkCLQzsihrlqr}) and (\ref{NashkDLQzsihrlqr}), and $P_2$ playing the maximizing disturbance $\kappa_2$ with values as in (\ref{NashkC2LQzsihrlqr}) and (\ref{NashkD2LQzsihrlqr}).
}
\NotAutom{\begin{example}{Sporadic Jump Times and Linear Dynamics}
Consider the following data for \eqref{eq: CDFG spontaneous}:
\begin{equation*}
  A_C = 1.8,
  B_C = B_D = [1 ,1],
  A_D = 2,
  \bar T = 1
\end{equation*}
We solve \eqref{DiffRiccatizsih} under \eqref{Riccatizsih} and \eqref{zlqeqinv}, obtaining a nonincreasing $P(\tau)$. In Figure \ref{fig:spontaneous} we report the evolution of two solutions, with a worst-case scenario depicting periodic jumps, an arbitrary scenario with sporadic jumps. As expected, as $t+j$ goes to infinity, the costs approach $V(\xi) = \xi_p^\top P(\xi_\tau) \xi_p$, and the periodic solution attains the upper bound. 

\begin{figure}[h]
  \centering
  \includegraphics[width=8cm]{Figures/Sporadic}
  \vspace{-0.8cm}
  \caption{Nonunique solutions to \eqref{eq: CDFG spontaneous} due to set-valued dynamics for $\xi=[2,0], Q_C=Q_D=1,$ $R_{C1}=R_{D1}=1.304$, and $2R_{C2}=R_{D2}=-8$. Worst-case cost solution (with periodic jump times) in green and purple. Arbitrary solution (with spontaneous jump times) in blue and red.}
  \label{fig:spontaneous}
\end{figure}
\end{example}
}
  \subsection{Set-Valued Flow and Jump Maps}
  As illustrated next, there are useful families of hybrid systems for which a \pnn{weak} saddle-point equilibrium exists. 
  The following example characterizes both the \pnn{weak} saddle-point equilibrium and the value function in a two-player zero-sum game with a scalar state associated to player $P_1$. Thus, $n_1=1, n_2=0$, and the role of player $P_2$ reduces to selecting the action $u_{C2}$.
 \begin{example}{Scalar set-valued hybrid game}\label{NumEx}
 Consider a {hybrid} system $\HS\IfIncd{}{}$ as in \eqref{Hinc} with state $x \in \reals$, input $u_C:=(u_{C1},u_{C2}) \in \reals^2$, and dynamics 
\begin{equation}
  \begin{array}{rll}
  \dot{x}\in &F(x,u_C):=[\underline a, \overline a] x+B u_C &  x \in [0,\overline \sigma] \cup [\mu, \delta]
  \\
  x^+ \in &G(x):=[\underline \sigma, \overline \sigma] & x = \mu
  \end{array}
  \label{nuHGameznew}
  \end{equation}
  where $\underline a < \overline a<0,B=[b_1 \> b_2]$  and {
  $\delta \geq\mu>\overline \sigma>\underline \sigma > 0$\footnote{Given that $ \mu>\delta$, flow from $\mu$ is not possible.}}. 
  Consider the cost functions $L_C(x,u_C):=x^2Q_C+u_C^\top R_C u_C$, $L_D(x):=P(x^2-\overline   \sigma^2)$, and terminal cost $q(x):=P x^2$, {defining $\mathcal{J}$ as in (\ref{defJTNCinc}),} with $R_C :=\left[ \begin{smallmatrix}R_{C1} & 0 \\ 0 & R_{C2}\end{smallmatrix}\right]$, {$Q_C,$ $R_{C1},$ $-R_{C2},$ $P\in \reals_{>0}$, such that} 
  \begin{equation}
  Q_C+2P\overline a-P^2(b_1^2 R_{C1}^{-1}+b_2^2 R_{C2}^{-1})=0.
  \label{PQ}  
  \end{equation}
  Here, $u_{C1}$ is designed by player $P_1$, which aims to minimize a cost functional $\mathcal{J}$, while player $P_2$ seeks to maximize it by means of $u_{C2}$.
   This is formulated as a two-player zero-sum hybrid game via solving Problem ($\diamond$) in Section \ref{subsec: probstat}. The function $V(x):=P x^2$ satisfies the sufficient condition for (\ref{HJBzsihinc})-(\ref{HJBzsihincub}) in Theorem \ref{thHJBszihinc}, namely
 \begin{equation}
 \begin{split}
 &\underset{u_C=(u_{C1},u_{C2}) \in \reals^2}
 {\underset{u_{C1}}{\min}\> \underset{u_{C2}}{\max}}  
 \left\{  L_C(x,u_C)+{\underset{f \in F(x,u_C)}{\sup}}  \left\langle \nabla V(x),f \right\rangle\right\}
 \\ 
 =&
 \underset{u_{C1} \in \mathbb{R}}{\min} \>\underset{u_{C2} \in \mathbb{R}}{\max} 
 \left\{ (Q_C+ 2P \overline a) x^2
 \hspace{-0.1cm}
 + \hspace{-0.1cm}
   R_{C1}u_{C1}^2  \right.
  \hspace{-0.1cm}
 + \hspace{-0.1cm}
 R_{C2}u_{C2}^2
 \hspace{-0.1cm}
 + \hspace{-0.1cm}
   2xP ( b_1u_{C1}+ b_2u_{C2})  \big\} 
  = 0
 \end{split}
  \label{HJBLQGnusCz}
 \end{equation}
which holds for all $x \in [0,\overline \sigma] \cup [\mu, \delta]$. 
In fact, the min-max in (\ref{HJBLQGnusCz}) is attained by 
 $\kappa_C(x)=(-R_{C1}^{-1} b_1P x,-R_{C2}^{-1}b_2P x)$.
 In particular, thanks to \eqref{PQ}, {we have} \IfAutom{
 $-L_C(x,\kappa_C(x))={\underset{f \in F(x,\kappa_C(x))}{\sup}}  \left\langle \nabla V(x),f \right\rangle$.
 }
 {
 \begin{equation*}\hspace{-0.2cm}
 \begin{array}{l}
 L_C(x,\kappa_C(x)) +\left\langle \nabla V(x),F(x,\kappa_C(x)) \right\rangle 
   \nonumber\\ \hspace{2cm}
 =
  (Q_C+2Pa) x^2+  R_{C1}\kappa^2_{C1}(x)
 +R_{C2}\kappa_{C2}^2(x)\nonumber
 \\ \hspace{4.3cm}
 +  2xP ( b_1\kappa_{C1}(x)+ b_2\kappa_{C2}(x)) \nonumber
  \\\hspace{2cm}
 =
  P^2(b_1^2R_{C1}^{-1}+b_2^2 R_{C2}^{-1}) x^2+  R_{C1}\kappa^2_{C1}(x)
  \nonumber
 \\ \hspace{2.3cm} 
 +R_{C2}\kappa_{C2}^2(x)+  2xP ( b_1\kappa_{C1}(x)+ b_2\kappa_{C2}(x))
  \nonumber
  \\ \hspace{2cm} =
 \left[P^2(b_1^2 R_{C1}^{-1}+b_2^2 R_{C2}^{-1}) +  R_{C1}^{-1} b_1^2 P^2
 +R_{C2}^{-1} b_2 P^2 \right.\nonumber
 \\  \hspace{3.7cm}
 \left.-  2P ( b_1R_{C1}^{-1} b_1 P+ b_2R_{C2}^{-1} b_2 P)\right]x^2 \nonumber
  \\ \hspace{2cm}
  =0 
 \end{array}
 \label{HJBLQGnusCzk}
 \end{equation*}}
 Then,  $V(x)= Px^2$ is a solution to \eqref{HJBzsihinc}-\eqref{HJBzsihincub}.
  In addition, the function $V$ satisfies 
 \begin{equation}
 \begin{array}{r@{}l}
  L_D(x)
 + {\underset{g \in G(x)}{\sup}} V(g(x)) 
 =&
 P x^2
 \end{array}
 \label{HJBLQGnusDz}
 \end{equation}
 at $x=\mu$, which makes $V(x)= Px^2$ a solution to \eqref{Bellmanzsihinc}-\eqref{Bellmanzsihincub} with \pnn{weak} saddle-point equilibrium $\kappa_C$.
 Given that $V$ is \IfGG{locally Lipschitz on $\Pi(\overline C)$}{continuously differentiable on $\reals$}, and that \eqref{HJBzsihinc}-\eqref{Bellmanzsihincub} hold thanks to (\ref{HJBLQGnusCz}) and (\ref{HJBLQGnusDz}), from Theorem \ref{thHJBszihinc} we have that the value function is
 $
 \mathcal{J}^*(\xi) :=
  P \xi^2
 \label{HyLCost}
 $
 for any $\xi \in [0,\overline \sigma] \cup [\mu, \delta]$.
 
 {To study in detail the nonunique solutions {yielded by the feedback law $\kappa_C$}, notice that solutions jump at $x=\mu$ to $\sigma_s \in [\underline \sigma, \overline \sigma]$. 
  {Consider a solution} 
$\phi_h$ with domain $ \dom \phi_h=([0,t^h]\times \{0\}) \cup ([t^h, \infty)\times \{1\})$, and given by
 $ \phi_h(t,0)= \delta \exp((a_s-R_{C1}^{-1} b_1P-R_{C2}^{-1}b_2P) t), \hspace{0.5cm} \phi_h(t,1)=
 \sigma_s \exp((a_s-R_{C1}^{-1} b_1P-R_{C2}^{-1}b_2P) (t-t^h))$ with $a_s \in [\underline a, \overline a]$.
 In simple words, $\phi_h$ flows from $\delta$ to $\mu$ in $t^h$ seconds, then it jumps to $\sigma_s$, and flows converging (exponentially fast) to {zero}.
Notice that $\kappa_C$ as defined above also yields a solution {$\phi_\kappa$} with domain $ \dom \phi_\kappa=([0,t^\kappa]\times \{0\}) \cup ([t^\kappa, \infty)\times \{1\})$, and given by
$ \phi_\kappa(t,0)= \delta \exp((\overline a-R_{C1}^{-1} b_1P-R_{C2}^{-1}b_2P) t), \hspace{0.5cm} \phi_\kappa(t,1)=
\overline \sigma \exp((\overline a-R_{C1}^{-1} b_1P-R_{C2}^{-1}b_2P) (t-t^\kappa))$ attaining the worst-case cost.
Figure \ref{fig:my_label} illustrates the behavior of the solutions $\phi_h$ and $\phi_\kappa$, yielded by $\kappa_C$, with the cost of the latter equating $P\delta^2$. 
\IfSiam{
\begin{figure}
    \centering
    \subfloat[Nonunique solutions due to set-valued dynamics for $\underline a=-2, \overline a = -1, b_1=b_2=1, \delta=\xi=2, \mu=1, \underline \sigma=0.3, \overline \sigma = 0.5, Q_C=1,$ $R_{C1}=1.304, R_{C2}=-4$, and  $P =0.4481$. Worst-case cost solution (green and purple). Arbitrary solution (blue and red).]
    {
        \label{fig:my_label}
        \centering
        \includegraphics[width=0.45\linewidth]{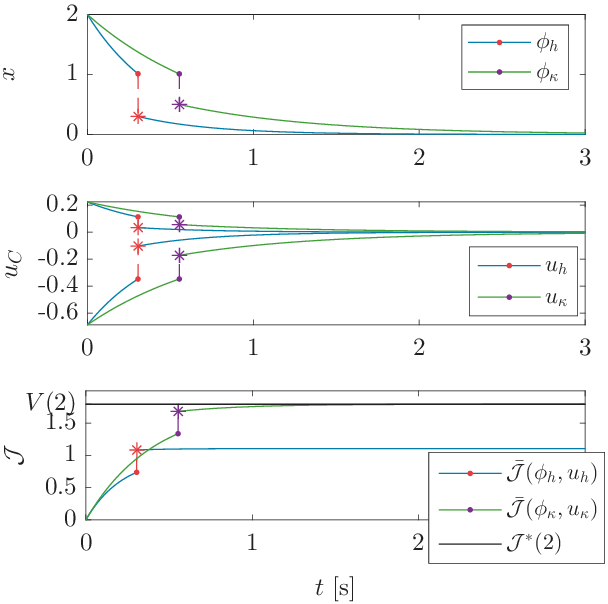}
    }
    \hfill 
    \subfloat[Nonunique costs due to set-valued dynamics when varying the feedback gains around the optimal value. The cost is evaluated on solutions to $\HS$ from $\xi$ with feedback law variations specified by  $\epsilon_u$ and $\epsilon_w$ in $u=(\epsilon_u \kappa_{C1},\epsilon_w \kappa_{C2})$. Worst-case costs (red). Costs of arbitrary solutions (blue). Value function (filled red).]
    {
        \label{fig:nonuniquecosts}
        \centering 
    \includegraphics[width=0.49\linewidth]{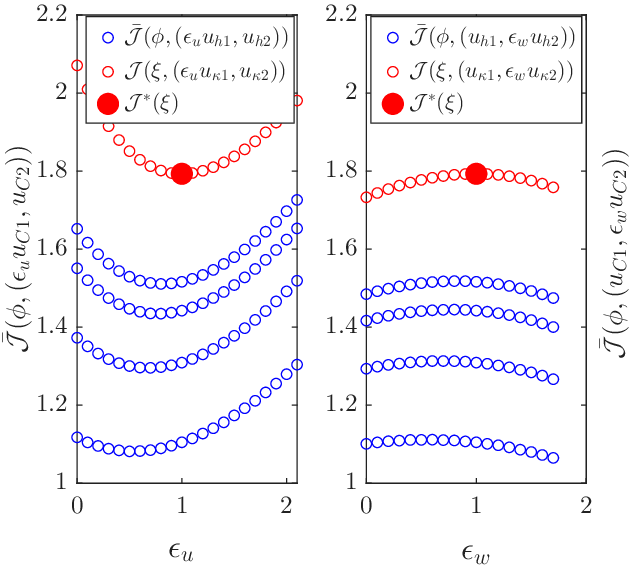}
    }
\caption{Weak saddle-point equilibrium scalar zero-sum two-player hybrid game.}
\vspace{-0.8cm}
\label{fig:csetup}
\end{figure}
}{
 \begin{figure}[h]
    \centering
     \includegraphics[width=7.5cm]{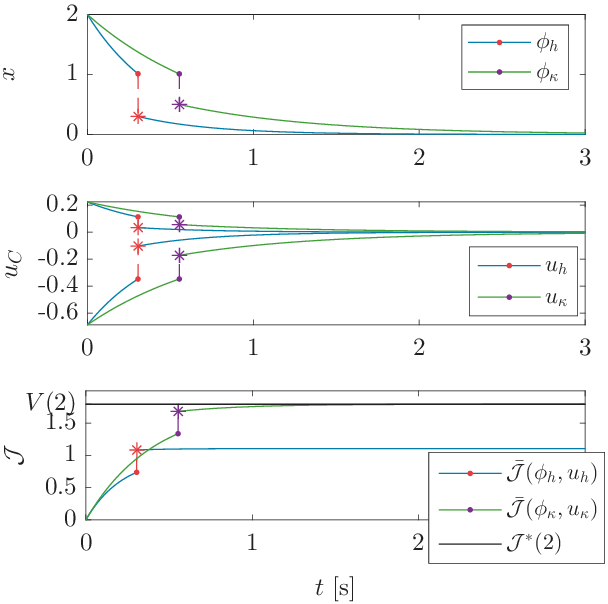}
     \caption{Nonunique solutions due to set-valued dynamics for $\underline a=-2, \overline a = -1, b_1=b_2=1, \delta=\xi=2, \mu=1, \underline \sigma=0.3, \overline \sigma = 0.5, Q_C=1,$ $R_{C1}=1.304, R_{C2}=-4$, and  $P =0.4481$. Worst-case cost solution (green and purple). Arbitrary solution (blue and red).}
     \label{fig:my_label}
 \end{figure}

 \begin{figure}[h]
  \centering
  \includegraphics[width=8cm]{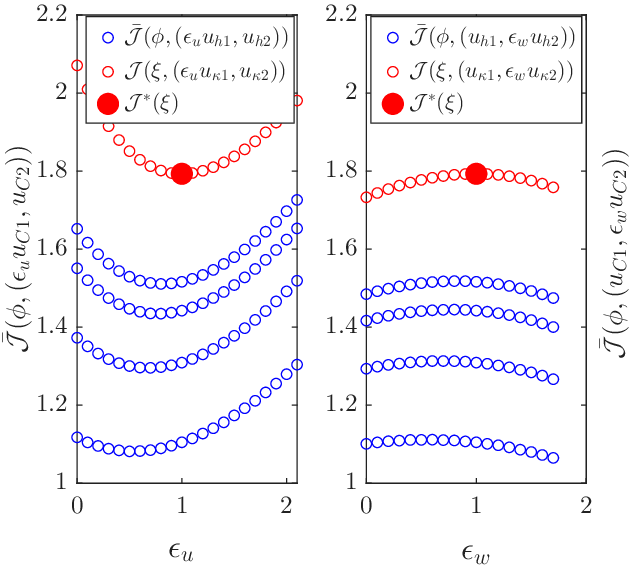}
  \caption{Nonunique costs due to set-valued dynamics when varying the feedback gains around the optimal value. The cost is evaluated on solutions to $\HS$ from $\xi$ with feedback law variations specified by  $\epsilon_u$ and $\epsilon_w$ in $u=(\epsilon_u \kappa_{C1},\epsilon_w \kappa_{C2})$.
  Worst-case costs (red). Costs of arbitrary solutions (blue). Value function (filled red).}
  \label{fig:nonuniquecosts}
\end{figure}
}

 \pnn{Notice that when both players play the weak saddle-point equilibrium strategy $\kappa_C=(\kappa_{C1},\kappa_{C2})$($\epsilon_u=1$ and $\epsilon_w=1$ in Figure \ref{fig:nonuniquecosts}),
 the worst-case cost (in red) is not necessarily attained, which implies that $\kappa_C$ is not necessarily optimal in the min-max sense. Indeed, notice that other strategies (e.g. $\epsilon_u=1.5$) 
 attain lower costs (in blue) when player $P_2$ sticks to play $\kappa_{C2}$. Nevertheless, by playing t$\kappa_{C1}$, player $P_1$ minimizes the worst-case cost (red) under the maximizing adversarial action $\kappa_{C2}$. 
 }

 \pnn{
 Let $\A=\{0\}$ and given that $L_C \in \mathcal{PD}_{\kappa_C }(\A \cap C_\kappa)$, {$L_D \in \mathcal{PD}_{\kappa_D}(\A \cap D_\kappa)$ vacuously}, and (\ref{HJBzsihinc})-(\ref{Bellmanzsihincub}) hold, by setting $\alpha_1(|x|_\A)=x^\top (P-I) x$ and $\alpha_2(|x|_\A)=x^\top (P+I) x$, 
 and \pnn{for $|x|_\A \mapsto \eta(|x|_\A)=P \frac{x^2}{2}$, $L_D(x,\kappa_D(x)) \geq \eta(|x|_\A)$ for all $x \in D_\kappa$ when $\mu^2 \geq 2 \bar \sigma ^2$,}  
 from Theorem \ref{zCorStability}, we have that $\kappa_C$ is the weak saddle-point equilibrium and renders $\A$ uniformly globally asymptotically stable for 
 (\ref{nuHGameznew}).}}
 \NotSiam
 {\begin{figure}[h]
     \centering
     \includegraphics[width=8cm]{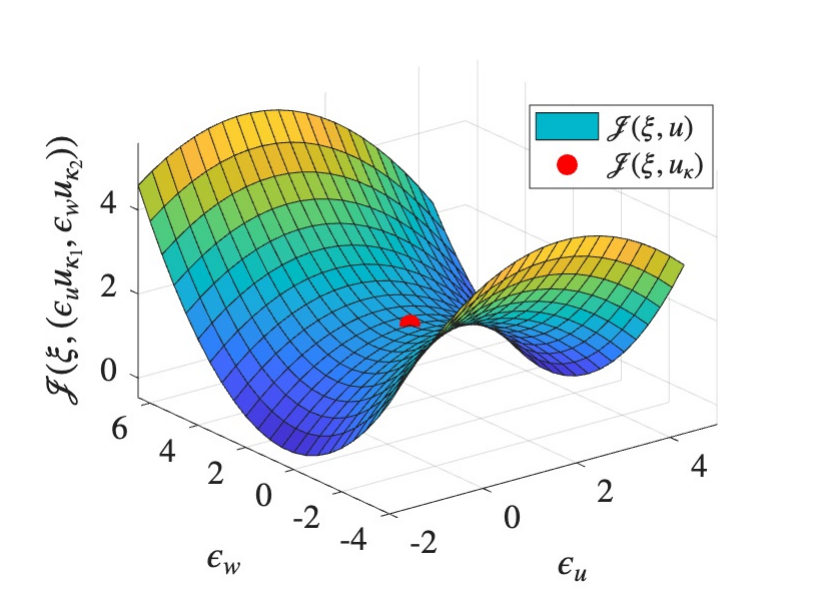}
     \caption{Saddle-point behavior in the cost of continuous solutions from $\xi=2$ when varying the feedback gains around the min-max value. The worst-case cost is evaluated on solutions $(\phi,u)\in \pnn{\mathcal{S}_\HS} (\xi) $ with feedback law variations specified by  $\epsilon_u$ and $\epsilon_u$ in $u=(\epsilon_u \kappa_1(\phi),\epsilon_w \kappa_2(\phi))$.}
     \label{Costs}
 \end{figure}
 }
 \label{NumericalEx}
 \end{example}
%
%
%
\IfIncd{}{
\subsection{Application 2: Security jumps-actuated hybrid game}
We study a special case of Application \ref{SecurityApp} and apply Theorem \ref{thHJBszih} in this section.
Consider a hybrid system with state $x\in \reals^{n} $, input $ u_D=(u_{D1},u_{D2}) \in  \reals^{m_D}$, and dynamics $\HS$ as in (\ref{Heq})\IfTp{,}{ with $N=2$ and} described by
\begin{equation}
\begin{array}{rl@{}r@{}l}
\dot{x}=&F(x)&  x\>\>\> &\in C
\\
x^+=&A_D x +\left[ \begin{smallmatrix}B_{D1} \>B_{D2}\end{smallmatrix}\right]\left[ \begin{smallmatrix}u_{D1}\\ u_{D2}\end{smallmatrix} \right]&\hspace{0.5cm} (x,u_D)& \in D
\end{array}
\label{eq: CDFG security}
\end{equation}
with Lipschitz continuous  
$F: \mathbb{R}^{n\times n}\rightarrow \mathbb{R}^{n\times n},
A_D 
\in \mathbb{R}^{n\times n}$,  and  $C \subset \mathbb{R}^{n}, D\subset\mathbb{R}^n \times  \mathbb{R}^{m_D}$, such that $C\cup D$ is nonempty. 
The input $u_{D1}$ plays the role of the control and $u_{D2}$ is the disturbance input. 
{Following Application\>\ref{SecurityApp}, the problem of 
minimizing a cost functional $\mathcal{J}$ in the presence of the maximizing attack {$u_{D2}$} is formulated as a two-player zero-sum game.} Thus, \textcolor{black}{by solving Problem $(\diamond)$ for every $\xi \in \Pi(C\cup D)$,} 
the control objective is achieved. 

The following result presents a tool for the solution of the optimal control problem for jumps-actuated hybrid systems with state-affine flow maps under a malicious input attack designed to cause as much damage as possible.
\begin{corollary}{Hybrid Riccati equation for security}
Given {a hybrid system $\HS$ as in \eqref{Heq} defined by $(C,F,D,G)$ as in \eqref{eq: CDFG security}, }
and,
{with the aim of pursuing minimum energy and distance to the origin, consider the cost functions 
$L_C(x,u_C):=0
$, $L_D(x,u_D):=x^\top Q_D x+  u_{D1}^\top R_{D1} u_{D1}+ u_{D2}^\top R_{D2} u_{D2} $, and  terminal cost $q(x):=x^\top P x$, {defining $\mathcal{J}$ as in (\ref{defJTNCinc}),} {with} $Q_D\in \mathbb{S}^n_+$, 
$R_{D1} \in \mathbb{S}^{m_{D_1}}_+$, $-R_{D2} \in \mathbb{S}^{m_{D_2}}_+$  and  $P\in \mathbb{S}^n_+ $. }
Suppose there exists a matrix $P\in \mathbb{S}^n_+ $ such that
\begin{equation}
0=
2x^\top P F(x)\hspace{1cm} \forall x \in \Pi(C),
\label{DiffRiccatizsihsec} 
\end{equation}
%
\begin{equation}
\begin{array}{rl}
-R_{D2} -B_{D2}^\top P B_{D2},R_{D1} +B_{D1}^\top P B_{D1} &\in \mathbb{S}_{0+}^{m_D},
\end{array}
\label{zlqeqinvsec}
\end{equation}
the matrix $R_v=\left[\begin{smallmatrix}
R_{D1}+B_{D1}^\top PB_{D1}
&
B_{D1}^\top P B_{D2}
 \\
B_{D2}^\top P B_{D1}
&
R_{D2}+B_{D2}^\top P B_{D2}
\end{smallmatrix}\right] $ is invertible, and 
\begin{equation}
\begin{split}
0=-P+Q_D+A_D^\top P A_D \hspace{3.1cm}\IfConf{\\}{}
-
\begin{bmatrix}
A_D^\top P B_{D1}
 &
A_D^\top P B_{D2}
\end{bmatrix}
R_v
^{-1}
\begin{bmatrix}
B_{D1}^\top P A_D
\\
B_{D2}^\top P A_D
\end{bmatrix}
\end{split}
\label{Riccatizsihsec} 
\end{equation}
Then, with the feedback law 
{\small
\begin{equation}
\kappa_{D1}(x)=-[R_v^{-1}(1,1) \> R_v^{-1}(1,2)]
\begin{bmatrix}
B_{D1}^\top P A_D
\\
B_{D2}^\top P A_D
\end{bmatrix} x  \hspace{0.4cm} \forall x \in \Pi(D)
\label{NashkDLQzsihsec}
\end{equation}}
the cost functional $\mathcal{J}$ is minimized in the presence of the maximizing attack $u_2$, given by
{\small
\begin{equation}
\kappa_{D2}(x)=-[R_v^{-1}(2,1) \> R_v^{-1}(2,2)]
\begin{bmatrix}
B_{D1}^\top P A_D
\\
B_{D2}^\top P A_D
\end{bmatrix} x  \hspace{0.4cm} \forall x \in \Pi(D)
\label{NashkD2LQzsihsec}
\end{equation}}
In addition, for each $x \in \Pi ( \overline{C}) \cup \Pi(D)$, the value function is equal to $V(x):=x^\top P x$. 
\label{HREqSec}
\end{corollary}
{\begin{proof}
We show that when conditions (\ref{DiffRiccatizsihsec})-(\ref{Riccatizsihsec}) hold, by using the result in Theorem \ref{thHJBszih},  the value function is equal to the function $V$ and under
 the feedback law as in (\ref{NashkDLQzsihsec}) such a cost is attained in the presence of the maximizing attack given by (\ref{NashkD2LQzsihsec}). 
We can write (\ref{HJBzsih}) in Theorem \ref{thHJBszih} as
\begin{eqnarray}
0=
2x^\top P F(x) \hspace{1cm} \forall x \in \Pi(C)
\label{HJBzsihlqsec}
\end{eqnarray}
which is satisfied thanks to (\ref{DiffRiccatizsihsec}).
On the other hand, we can write (\ref{Bellmanzsih}) in Theorem \ref{thHJBszih} as
\begin{eqnarray}
&&x^\top P x =\underset{u_D=(u_{D1},u_{D2}) \in \Pi_u(x,D)}
{\min_{u_{D1}} \max_{u_{D2}}}  \mathcal{L}_D(x,u_D),
 \nonumber
\\ 
&&\mathcal{L}_D(x,u_D)= 
x^\top Q_D x+ u_{D1}^\top R_{D1} u_{D1}+ u_{D2}^\top R_{D2} u_{D2}\nonumber\\
&&\hspace{1.7cm}+ (A_D x+B_D u_D)^\top P (A_D x+B_D u_D) 
\label{Bellmanzsihlqsec}
\end{eqnarray}
\NotAutom{which can be expanded as
\begin{multline}
\mathcal{L}_D(x,u_D)=  
 x^\top (Q_D+A_D^\top P A_D) x + 2 x^\top A_D^\top  P B_D u_D \\
 + u_{D1}^\top (R_{D1} +B_{D1}^\top P B_{D1}) u_{D1}
 + u_{D2}^\top (R_{D2} +B_{D2}^\top P B_{D2}) u_{D2}\\
 + u_{D1}^\top (B_{D1}^\top P B_{D2}) u_{D2}
 + u_{D2}^\top (B_{D2}^\top P B_{D1}) u_{D1}
\end{multline}}
The first order necessary conditions for optimality
\NotAutom{
\begin{equation*}
\frac{\partial}{\partial u_{D1}}\left.
\mathcal{L}_D(x,u_D)  \right|_{u_{D}^*} =0,
\quad
\frac{\partial}{\partial u_{D2}}\left.
\mathcal{L}_D(x,u_D) \right|_{u_{D}^*} =0
\end{equation*}
}
are satisfied by the point $u_D^*=(u_{D1}^*,u_{D2}^*)$, such that for each $x \in \Pi(D)$,
\begin{equation}
{\tiny
u_{D}^*=-\begin{bmatrix}
R_{D1}+B_{D1}^\top P B_{D1}
&
B_{D1}^\top P B_{D2}
 \\
B_{D2}^\top P B_{D1}
&
R_{D2}+B_{D2}^\top P B_{D2}
\end{bmatrix} 
^{-1}
\begin{bmatrix}
B_{D1}^\top P A_D
\\
B_{D2}^\top P A_D
\end{bmatrix} x
}
\label{zLQHPud*sec}
\end{equation}
Given that (\ref{zlqeqinvsec}) holds, the second-order sufficient conditions for optimality 
\NotAutom{
\begin{equation*}
\frac{\partial^2}{\partial u_{D1}^2}\left. 
 \mathcal{L}_D(x,u_D) \right|_{u_D^*} \succeq 0,
\quad
\frac{\partial^2}{\partial u_{D2}^2}\left.
 \mathcal{L}_D(x,u_D)\right|_{u_D^*} \preceq 0,
\end{equation*}}
are satisfied, rendering $u_D^*$ as in (\ref{zLQHPud*sec}) as an optimizer of the min-max problem in (\ref{Bellmanzsihlqsec}).
In addition, $u_D^*$
satisfies $\mathcal{L}_D(x,u_D^*)=x^\top P x$, with $P$ as in (\ref{Riccatizsihsec}), leading $V(x)=x^\top P x$ as a solution of (\ref{Bellmanzsih}) in Theorem \ref{thHJBszih}. 

Thus, given that $V$ is continuously differentiable in $\reals^n$ and Assumption (\ref{AssLipsZ}) holds, by applying Theorem \ref{thHJBszih}, in particular from (\ref{ResultValue}), 
{for every $ \xi \in \Pi(\overline{C}\cup D)$ the value function is $\J^*(\xi)=\J(\xi,(u_{D1}^*,u_{D2}^*))= \xi^\top P \xi$. From (\ref{kHJBeqzsihc}) and (\ref{kBeqzsihc}), when $P_2$ plays $u_2^*$ defined by $\kappa_{D2}$  as in (\ref{NashkD2LQzsihsec}), $P_1$ minimizes the cost of complete solutions to $\HS$ 
by playing $u_1^*$ defined by $\kappa_{D1}$ as in (\ref{NashkDLQzsihsec}), attaining $ \J(\xi,u^*)$, and satisfying (\ref{SaddlePointIneq}).}
\end{proof}}
\textcolor{black}{Notice that the saddle-point equilibrium $\kappa_D:=(\kappa_{D1}, \kappa_{D2})$ is composed of $P_1$ playing the minimizer strategy $\kappa_{D1}$ as in (\ref{NashkDLQzsihsec}), and $P_2$ playing the maximizing attack $\kappa_{D2}$ as in (\ref{NashkD2LQzsihsec}).}
\IfPers{
\begin{example}{Linear quadratic differential games}
Consider a two-player zero-sum hybrid game with state $x \in \reals^n$, input $u=(u_C, u_D) \in \reals^{m_C}\times \reals^{m_D}$, and dynamics $\HS$ as in (\ref{Heq}) with \IfTp{}{$N=2$ and}
\begin{eqnarray}
C&=&\mathbb{R}^n \times \mathbb{R}^{m_C} \nonumber
\\
F(x,u_C)&=&A_C x+B_C u_C \nonumber
\label{zsLinearGame}
\end{eqnarray}
and jump set $D=\emptyset$, where $A_C:=\textup{blkdiag}_{i \in \IfTp{\{1,2\}}{\mathcal{V}}}\{A_i\} \in \mathbb{R}^{n\times n}$ and $B_C:=\textup{blkdiag}_{i \in \IfTp{\{1,2\}}{\mathcal{V}}} \{B_{Ci}\} \in \mathbb{R}^{n\times m_C}$. 
Notice that $B_C u_C = B_{C1} u_{C1} + B_{C2} u_{C2}$, where $B_{C1}=(B_{C1},\textbf{0}_{n_2\times m_{C1}})$ and $B_{C2}=(\textbf{0}_{n_1\times m_{C2}},B_{C2})$. The jump map is arbitrary, as $D$ is empty.
Consider the stage costs $L_C(x,u_C):=x^\top Q_C x+ u_{C1}^\top R_{C1} u_{C1}+ u_{C2}^\top R_{C2} u_{C2}$ and $L_C(x,u_D):=0$, the terminal cost $q(x):=x^\top P x$ and the function $V(x):=x^\top P x$, where $Q_C \in \mathbb{S}^n_+$, $R_{C1} \in \mathbb{S}^{m_{C_1}}_+$, $-R_{C2} \in \mathbb{S}^{m_{C_2}}_+$  and  $P\in \mathbb{S}^n_+$. We refer to this game as a zero-sum infinite-horizon linear quadratic differential game $(\mathcal{ZLQC})$. The function $V$ is continuously differentiable in $\reals^n$, and
\begin{eqnarray}
 &&\nonumber 
\underset{u_C=(u_{C1},u_{C2}) \in \Pi_u(x,C)}
{ \min_{u_{C1}} \max_{u_{C2}}}
\left\{L_C(x,u) +\left\langle \nabla V(x),F(x,u_C) \right\rangle \right\}
\\&=&
{  \min_{u_{C1} \in \reals^{m_{C_1}}} \max_{u_{C2} \in \reals^{m_{C_2}}}  }
\left\{x^\top Q_C x+ u_{C1}^\top R_{C1} u_{C1}+ u_{C2}^\top R_{C2} u_{C2}+2x^\top P (A_C x+B_C u_C ) \right\}
\nonumber
\\
&=&
{  \min_{u_{C1} \in \reals^{m_{C_1}}} \max_{u_{C2} \in \reals^{m_{C_2}}}  }
\left\{ x^\top (Q_C+2 P A_C) x   
+u_{C1}^\top R_{C1} u_{C1} +  2 x^\top  P B_{C1} u_{C1} 
\right. 
\nonumber \\  &&
\hspace{9cm} 
\left. + u_{C2}^\top R_{C2} u_{C2} + 2 x^\top  P  B_{C2} u_{C2} \right\}
\nonumber
\\  
&=& 0
\label{zsHJBLQGC}
\end{eqnarray}
is satisfied $ \forall x \in \reals^n$, and attained by $\kappa_C(x)=(-R_{C1}^{-1}B_{C1}^\top P x, -R_{C2}^{-1}B_{C2}^\top P x)$, if the continuous time $\HS_\infty$-like Riccati equation 
\begin{equation}
0
=Q_C+PA_C+A_C^\top P-P(  B_{C2} R_{C2}^{-1}B_{C2}^\top+B_{C1} R_{C1}^{-1}B_{C1}^\top) P
\label{zsRiccati}
\end{equation}
is satisfied. 
Given that $V$ is continuously differentiable, and that (\ref{Bellmanzsih}) and (\ref{HJBzsih}) are satisfied thanks to  (\ref{zsHJBLQGC}) and $D=\emptyset$, from Theorem \ref{thHJBszih}, the value function is
\begin{eqnarray}
\mathcal{J}^*(\xi) &:=& 
 \xi^\top P \xi,
\label{zsContLCost}
\end{eqnarray}
for any $\xi \in \reals^n$ 
%
{Given that this is a purely continuous system, the calculation of the cost for when applying the feedback law $\kappa_C(x)=-R_C^{-1}B_C^\top P x$, which renders a solution $(\phi,u_C^*)$ with $\dom \phi \ni (t,j) \mapsto u_C^*(t,j)= \kappa_C(\phi(t,j))$, is as follows
\begin{equation*}
\mathcal{J}(\xi,u_C^*) =
 \int_{t_0}^{\infty} L_C(\phi(t,j),u_{C}^*(t,j))dt  + \underset{(t,j) \in \textup{dom}\phi}{\limsup_{t+j\rightarrow \infty} } q(\phi(t,j))
\end{equation*}
and, when (CARE) is satisfied, one has
\begin{eqnarray*}
 \int_{t_0}^{\infty} L_C(\phi(t,j),u_{C}^*(t,j))dt  &=& 
 \int_{0}^{\infty} \phi^\top(t,j) Q_C \phi(t,j) +{u_{C}^*}^\top(t,j) R_C u_{C}^*(t,j)dt  
 \\&=&  \int_{0}^{\infty} \phi^\top(t,j) (Q_C +P B_C R_C^{-1}B_C ^\top P)\phi(t,j) dt  
  \\&=&  \int_{0}^{\infty} \phi^\top(t,j) (2P B_C R_C^{-1}B_C ^\top P-2PA_C)\phi(t,j) dt  
  \\&=&-\int_{0}^{\infty} \frac{d V(\phi(t,j))}{d t} dt
    \\&=&V(\phi(0,0))-\limsup_{t\rightarrow \infty}V(\phi(t,0))
  \\&=& 
  \phi^\top(0,0) P \phi(0,0)- \limsup_{t\rightarrow \infty} \phi^\top(t,0) P \phi(t,0)
\end{eqnarray*}
which, given that $\phi(0,0)=\xi$ and $q(\phi(t,j))=\phi^\top(t,0) P \phi(t,0)$, implies
$	\mathcal{J}(\xi,u_C^*) =\xi^\top P \xi$, showing that the feedback law $\kappa_C$ attains the optimal cost.} 
and the feedback law $\kappa_C$ is a pure strategy saddle-point equilibrium for the $(\mathcal{ZLQC})$ game.
\end{example}
\begin{example}{Linear quadratic difference games}
Consider a two-player zero-sum hybrid game with state $x \in \reals^n$, input $u=(u_C, u_D) \in \reals^{m_C}\times \reals^{m_D}$, and dynamics $\HS$ as in (\ref{Heq}) with \IfTp{}{$N=2$ and}
\begin{eqnarray}
D&=&\mathbb{R}^n \times \mathbb{R}^{m_D} \nonumber
\\
G(x,u_D)&=&A_D x+B_D u_D
\label{zsLinearDGame}
\end{eqnarray}
and flow set $C=\emptyset$, where $A_D:=\textup{blkdiag}_{i \in \IfTp{\{1,2\}}{\mathcal{V}}}\{A_i\}  \in \mathbb{R}^{n\times n}$ and $B_D:=\textup{blkdiag}_{i \in \IfTp{\{1,2\}}{\mathcal{V}}} \{B_{Di}\} \in \mathbb{R}^{n\times m_D}$. 
Notice that $B_D u_D = B_{D1} u_{D1} + B_{D2} u_{D2}$, where $B_{D1}=(B_{D1},\textbf{0}_{n_2\times m_{D1}})$ and $B_{D2}=(\textbf{0}_{n_1\times m_{D2}},B_{D2})$. The flow map is arbitrary, as $C$ is empty.
Consider the stage costs $L_C(x,u_C):=0$ and $L_D(x,u_D):=x^\top Q_D x+u_{D1}^\top R_{D1} u_{D1}+ u_{D2}^\top R_{D2} u_{D2}$, the terminal cost $q(x):=x^\top P x$ and the function $V(x):=x^\top P x$, where $Q_D\in \mathbb{S}^n_+$, $R_{D1}\in \mathbb{S}^{m_{D_1}}_+, -R_{D2}\in \mathbb{S}^{m_{D_2}}_+$ and  $P \in \mathbb{S}^n_+$. We refer to this game as an infinite-horizon one player linear quadratic difference game $(\mathcal{ZLQD})$. The function $V$ is such that
\begin{eqnarray}
&&
\underset{u_D=(u_{D1},u_{D2}) \in \Pi_u(x,D)}
{\min_{u_{D1}} \max_{u_{D2}}}  
\left\{ 
 L_D(x,u_D)
+ V(G(x,u_D)) \right\}  \nonumber
\\&=&
{  \min_{u_{D1} \in \reals^{m_{D_1}}} \max_{u_{D2} \in \reals^{m_{D_2}}}  }
\left\{x^\top Q_D x+ u_{D1}^\top R_{D1} u_{D1}+ u_{D2}^\top R_{D2} u_{D2} \right. \nonumber \\
&&\hspace{8cm}\left. +(A_D x+B_D u_D )^\top P (A_D x+B_D u_D ) \right\}
\nonumber
\\&=&
{  \min_{u_{D1} \in \reals^{m_{D_1}}} \max_{u_{D2} \in \reals^{m_{D_2}}}  }
\left\{  x^\top (Q_D+A_D^\top P A_D) x 
+ u_{D1}^\top (R_{D1}+B_{D1}^\top P B_{D1}) u_{D1} 
\right. 
 \nonumber \\ && 
 \left. \hspace{2cm}
+ u_{D2}^\top (R_{D2}+B_{D2}^\top P B_{D2}) u_{D2} 
+2u_{D1}^\top B_{D1}^\top P B_{D2} u_{D2} 
+ 2 x^\top A_D^\top  P B_D u_D\right\} 
\nonumber
\\
&=&x^\top P x,
\label{zsHJBLQGD}
\end{eqnarray}
is satisfied $ \forall x \in \reals^n$, and attained by $\kappa_D(x)=-R^{-1}
(B_{D1}^\top P A_D, 
B_{D2}^\top P A_D)x$, with $R=\begin{bmatrix}
R_{D1}+B_{D1}^\top P B_{D1}
&
B_{D1}^\top P B_{D2}
 \\
B_{D2}^\top P B_{D1}
&
R_{D2}+B_{D2}^\top P B_{D2}
\end{bmatrix}$,
when the discrete-time $\HS_\infty$-like algebraic Riccati equation
\begin{multline}
P =Q_D+A_D^\top P A_D \IfConf{\\}{}
-
\begin{bmatrix}
A_D^\top P B_{D1}
 &
A_D^\top P B_{D2}
\end{bmatrix}
R 
^{-1}
\begin{bmatrix}
B_{D1}^\top P A_D
\\
B_{D2}^\top P A_D
\end{bmatrix}
\label{zsDARE} 
\end{multline}
is satisfied. 
Given that $V$ is continuously differentiable, and that (\ref{Bellmanzsih}) and (\ref{HJBzsih}) are satisfied thanks to $C=\emptyset$ and (\ref{zsHJBLQGC}), from Theorem \ref{thHJBszih}, the value function is
\begin{eqnarray}
\mathcal{J}^*(\xi) &:=& 
 \xi^\top P \xi,
\label{zsContLCost}
\end{eqnarray}
{Given that this is a purely discrete system, the calculation of the cost for when applying the feedback law $\kappa_D(x)=-(R_D+B_D^\top PB_D)^{-1}B_D^\top P A_D x$, which renders a solution $(\phi,u_D^*)$ with $\dom \phi \ni (t,j) \mapsto u_D^*(t,j)= \kappa_D(\phi(t,j))$, is as follows
\begin{equation*}
\mathcal{J}(\xi,u_D^*) =
  \sum_{j=0}^{\infty} L_D(\phi(t_{j+1},j),u_{D}^*(t_{j+1},j))+ \underset{(t,j) \in \textup{dom}\phi}{\limsup_{t+j\rightarrow \infty} } q(\phi(t,j))
\end{equation*}
and, when (DARE) is satisfied, one has
\begin{eqnarray*}
 \sum_{j=0}^{\infty} L_D(\phi(t_{j+1},j),u_{D}^*(t_{j+1},j)) &=& 
 \sum_{j=0}^{\infty} \phi^\top(t,j)Q_D \phi(t,j) +{u_{D}^*}^\top(t,j) R_D u_{D}^*(t,j)  
 \\&=&  \sum_{j=0}^{\infty} \phi^\top(t,j) (Q_D +A_D^\top P B_D (R_D+B_D^\top P B_D)^{-1}B_D^\top P A_D)\phi(t,j)   
  \\&=&  \sum_{j=0}^{\infty} \phi^\top(t,j) (2A_D^\top P B_D (R_D+B_D^\top P B_D)^{-1}B_D^\top P A_D+P-A_D^\top PA_D)\phi(t,j)  
  \\&=&\sum_{j=0}^{\infty} V(\phi(t,j))-V(G(\phi(t,j),u_D^*(t,j))
  \\&=&V(\phi(0,0))-\limsup_{j\rightarrow \infty}V(\phi(0,j))
  \\&=& 
  \phi^\top(0,0) P \phi(0,0)- \limsup_{j\rightarrow \infty} \phi^\top(0,j) P \phi(0,j)
\end{eqnarray*}
which, given that $\phi(0,0)=\xi$ and $q(\phi(t,j))=\phi^\top(0,j) P \phi(0,j)$, implies
$	\mathcal{J}(\xi,u_D^*) =\xi^\top P \xi$, showing that the feedback law $\kappa_D$ attains the optimal cost.} 
and the feedback law $\kappa_D$ is a pure strategy saddle-point equilibrium for the $(\mathcal{ZLQD})$ game.
\end{example}}
\begin{example}{Bouncing ball} {Inspired by the problem in \cite{79c}, consider a simplified model of a juggling system as in \cite{69c},  }
with state $x \in \reals^2$, input $ u_D:=(u_{D1},u_{D2}) \in \reals^2$, and dynamics $\HS$ as in (\ref{Heq})\pno{,} with data
\begin{equation}
\begin{array}{rll}
C&=&\realsgeq \times \reals 
\\
F(x)&=&(x_2,-1)
\\
D&=&\{0\} \times \reals_{\leq 0} \times \reals^2 
\\
G(x,u_D)&=&(0, -\lambda  x_2+u_{D1}+u_{D2})
\end{array}
\label{BouncingBallDynamics}
\end{equation}
{where $u_{D1}$ is the control input, $u_{D2}$ is the action of an attacker, and $\lambda\in (0,1)$ is the coefficient of restitution of the ball. {Let $\A=\{0\}$, encoding the goal of stabilizing the ball to a rest state under an optimization objective.} As an instance of  Application \ref{SecurityApp}, the scenario in which $u_{D1}$ is designed to minimize a cost functional $\J$ under the presence of the maximizing attack $u_{D2}$ is formulated as a two-player zero-sum game.
With the aim of pursuing minimum energy and distance to the origin at jumps,} consider the cost functions $L_C(x,u_C):=0$, $L_D(x,u_D):=x_2^2Q_D+u_D^\top R_D u_D$, and terminal cost $q(x):=\frac{1}{2} x_2^2+x_1$ \textcolor{black}{defining $\J$ as in (\ref{defJTNCinc}),} with $R_D :=\left[ \begin{smallmatrix}R_{D1} & 0 \\ 0 & R_{D2}\end{smallmatrix}\right]$ and {$Q_D,$ $R_{D1},$ $-R_{D2} \in \reals_{>0}$}.  Here, $u_{D1}$ is designed by player $P_1$ which aims to minimize $\mathcal{J}$ while player $P_2$ seeks to maximize it by means of choosing $u_{D2}$. 
This is formulated as a two-player zero-sum hybrid game \textcolor{black}{via solving Problem ($\diamond$)}. 
The function $V(x):=\frac{1}{2} x_2^2+x_1$ is such that $\left\langle \nabla V(x),F(x) \right\rangle=0$
for all $x \in C$, making $V$ a solution to (\ref{HJBzsih}). In addition, the function $V$ is such that
%
{
\begin{equation}
  \begin{split}
 &\underset{u_D=(u_{D1},u_{D2}) \in \reals^2}
{  \underset{u_{D1}}{\min}\>  \underset{u_{D2}}{\max}  }  
\IfAutom{ \mathcal{L}_D(x,u_D)}
{
\left\{
 L_D(x,u_D)
+ V(G(x,u_D)) \right\} }
\\
&  
= \underset{u_{D1}\in \mathbb{R}}{\min} \> \underset{u_{D2}\in \mathbb{R}}{\max}  
\Big\{ x_2^2Q_D+u_D^\top R_D u_D
\\ & \hspace{2cm} +\frac{(-\lambda x_2 + u_{D1}+u_{D2})^2}{2} \Big\} 
= \frac{1}{2} x_2^2
\end{split}
\label{IBB}
\end{equation}}
for all $(x,u_D)\in D$, and attained by $\kappa_D(x)=(
\kappa_{D1}(x),\kappa_{D2}(x))$ with 
$\kappa_{D1}(x)=\frac{R_{D2}\lambda }{R_{D1} + R_{D2} + 2 R_{D1} R_{D2}}x_2$ and $\kappa_{D2}(x)=\frac{R_{D1}\lambda }{R_{D1} + R_{D2} + 2 R_{D1} R_{D2}}x_2$ when 
\begin{equation}
Q_D=\frac{- 2R_{D1} R_{D2} \lambda^2 + R_{D1} + R_{D2} + 2 R_{D1}R_{D2}}{2 R_{D1} + 2 R_{D2} + 4 R_{D1} R_{D2}},
\label{QdBB}
\end{equation}
which makes $V$ a solution to (\ref{Bellmanzsih}) \textcolor{black}{with saddle-point equilibrium $\kappa_D$}.
Thus, given that $V$ is continuously differentiable on $\reals^2$, and that (\ref{HJBzsih}) and (\ref{Bellmanzsih}) hold thanks to (\ref{IBB}) and (\ref{QdBB}), from Theorem \ref{thHJBszih}\pno{,} the value function is
$
\mathcal{J}^*(\xi_1,\xi_2) :=
  \frac{\xi_2^2}{2}+\xi_1.
\label{HyLCostBB}
$
Figure \ref{fig:bb} displays this behavior. 
\vspace{-0.2cm}
\begin{figure}[h]
    \centering
    \includegraphics[width=9cm]{Figures/bbfp2}
        \vspace{-0.5cm}
    \caption{Bouncing ball solutions attaining minimum cost under maximizing $u_{2}$, with $\lambda=0.8, R_{D1}=10, R_{D2}=-20$, and $Q_D=0.189$. 
    }
    \label{fig:bb}
\end{figure}

{Furthermore, given that 
$L_D \in \mathcal{PD}_{\kappa_D}(\A 
)$, and (\ref{kHJBeqihza})-(\ref{kBeqihzf}) hold, by setting $\alpha_1(s)=\min \left\{\frac{1}{2} \left(\frac{s}{ \sqrt 2}\right)^2, \frac{s}{ \sqrt 2}\right\}$ and $\alpha_2 (s)=\frac{1}{2}s^2+s$, 
from Corollary \ref{zCorStability}, we have that $\kappa_D$ is the saddle-point equilibrium and renders $\A=\{0\}$ uniformly globally asymptotically stable for $\HS$}.  

{In Figure \ref{CostsBB}, we let the players select feedback laws close to the Nash equilibrium and calculate the cost associated to the new laws. The variation of the cost along the changes in the feedback laws makes evident the saddle-point geometry. This example illustrates how our results apply to Zeno systems.}
\begin{figure}[h]
    \centering
    \includegraphics[width=9cm]{Figures/SaddleBB}
    \caption{Saddle point behavior in the cost of solutions to bouncing ball from $\xi=(1,1)$ when varying the feedback gains around the optimal value. The cost is evaluated on solutions $(\phi,u)\in \mathcal{S}^\infty_\HS (\xi) $ with feedback law variations specified by  $\epsilon_u$ and $\epsilon_w$ in $u=(\epsilon_u \kappa_1(\phi),\epsilon_w \kappa_2(\phi))$.}
    \label{CostsBB}
\end{figure}
\label{BouncingBallEx}
\end{example}
}
\IfPers{
\begin{example}{Linear quadratic differential games}
Taking up the $(\mathcal{ZLQC})$ game, recall that $\kappa_C:=(-R_{C1}^{-1}B_{C1}^\top P x,-R_{C2}^{-1}B_{C2}^\top P x)$ and consider the case in which $A_C-B_{C1}R_{C1}^{-1}B_{C1}^\top P-B_{C2}R_{C2}^{-1}B_{C2}^\top P$ is Hurwitz and $\A=\{0\}$. Given that $L_C \in \mathcal{PD}_{\kappa_C}(\A)$, $V$ is continuously differentiable, and (\ref{kHJBeqihza}), (\ref{kHJBeqihzb}), and (\ref{kHJBeqihzc}) hold, by setting $\alpha_1(|x|_\A)=x^\top (P-I) x$ and $\alpha_1(|x|_\A)=x^\top (P+I) x$, from Corollary \ref{zCorStability}, we have that $\kappa:=(\kappa_C, \cdot )$ is the pure strategy saddle-point equilibrium and renders $\A$ uniformly globally asymptotically stable for $\HS_{\kappa}$.
\end{example}

\begin{example}{Linear quadratic difference games}
Taking up the $(\mathcal{ZLQD})$ game, recall that $\kappa_D:=-
(R^{-1}B_{D1}^\top P A_D x, 
R^{-1}B_{D2}^\top P A_D x)$ and consider the case in which $A_D-[B_{D1},B_{D2}] R^{-1}(B_{D1}^\top P(0) A_D,
B_{D2}^\top P(0) A_D)$ is Schur and $\A=\{0\}$. Given that $L_D \in \mathcal{PD}_{\kappa_D}(\A)$, and  (\ref{kBeqihzd}), (\ref{kBeqihze}), and (\ref{kBeqihzf}) hold, by setting $\alpha_1(|x|_\A)=x^\top (P-I) x$ and $\alpha_1(|x|_\A)=x^\top (P+I) x$, from Corollary \ref{zCorStability}, we have that $\kappa:=(\cdot, \kappa_D)$ is the pure strategy saddle-point equilibrium and renders $\A$ uniformly globally asymptotically stable for $\HS_{\kappa}$. 
\end{example}}

}
{}\section{Conclusion and Future Work}
%
%
%

In this paper, {we formulate a two-player zero-sum game under dynamic constraints given in terms of hybrid dynamical systems, as in \cite{65nu}.}
Scenarios in which the control action is selected by a player $P_1$ to accomplish an objective and countereffect the damage by an adversarial player $P_2$ are studied. By encoding the objectives of the players in the optimization of a cost functional, sufficient conditions in Hamilton-Jacobi-Isaacs form are provided to bound and exactly evaluate it. 
The main result characterizes the {strategy of $P_1$} that minimizes the worst-case cost under the maximizing {adversarial action}. {Additional conditions are proposed to allow the \pnn{weak} saddle-point strategy to render a set of interest asymptotically stable by letting the value function take the role of a Lyapunov function.}
Future work includes \IfIncd{studying structural}{the extension of the results to 
 zero-sum hybrid games {with the relaxation of the uniqueness of solutions assumption. We expect the results \pnn{in this paper} can be generalized to randomized strategies, in particular, through the connection between set-valued dynamics and nonuniqueness of solutions, which captures nondeterminism.} 
Structural} conditions on the system that do not involve \pnn{the use of} $V$ and guarantee the existence of a solution to Problem $(\diamond)$ based on the smoothness and regularity of the data of the system, similar to those in \cite{GOEBEL2019153nu}\IfIncd{.}{, will also be studied.}
%
\bibliographystyle{siamplain}
\bibliography{GT,RGSOnline, biblio}

\NotSiam{
\appendix
\section{Appendix}
\begin{proposition}[Time-dependent conditions for upper bound]
  Consider a solution $(\phi,u)$ to $\HS\IfIncd{}{} $ as in (\IfIncd{\ref{Hinc}}{\ref{Heq}}) from $\xi$ with $u=(u_C, u_D)$, such that
%
 \begin{itemize}
\item for each $j \in \nats$ such that $I_{\phi}^j=[t_j,t_{j+1}]$ has a nonempty interior
\begin{equation}
  \begin{split}
    L_C(\phi(t,j),u_C(t,j)) 
   +\frac{d V}{d t}(\phi(t,j)\pn{,u_C(t,j)})  \leq
 0 \quad \>\>
  \textup{for almost all } t \in I_{\phi}^j
\end{split}
\label{eq:upperboundhamiltonian}
\end{equation}
  and 
\item  for every $(t_{j+1},j) \in \dom \phi$ such that $(t_{j+1},j+1) \in \dom \phi$, 
\begin{equation}
\begin{split}
  L_D(\phi(t_{j+1},j),u_D(t_{j+1},j))
        + V(\phi(t_{j+1},j+1))-V(\phi(t_{j+1},j)) 
         \leq 0.
\end{split}
         \label{upperbounddiscrete}
\end{equation}
\end{itemize}
Then
\begin{equation}
\begin{split}
   \sum_{j=0}^{\sup_j \dom \phi} \int_{t_{j}}^{t_{j+1}} L_C(\phi(t,j),u_{C}(t,j))dt
     &
 + \hspace{-0.3cm}
 \sum_{j=0}^{{\sup_j \dom \phi}-1}L_D(\phi(t_{j+1},j),u_{D}(t_{j+1},j))
  \\
 &+\underset{(t,j) \in \textup{dom}\phi}
 {\pnn{\limsup_{(t,j)\rightarrow \sup \dom \phi}}}
   V(\phi(t,j)) 
 \leq V(\xi).
 \end{split}
 \label{eq:costboundVxi}
 \end{equation}
 \label{Pp:ProofDerivation}
\end{proposition}
\begin{proof}
Given a solution $(\phi,u)$ to $\HS$ from $\xi$, where $\{t_j\}_{j=0}^{\sup_j \dom \phi}$ is a nondecreasing sequence associated with the hybrid time domain of $(\phi,u)$ as in Definition \ref{htd},  
{
for each $j \in \nats$ such that $I_{\phi}^j=[t_j, t_{j+1}]$ has a nonempty interior,
}
by integrating \eqref{eq:upperboundhamiltonian} over 
$I_{\phi}^j$, 
we obtain
\begin{equation*}
\begin{split}
0\geq  \int_{t_{j}}^{t_{j+1}} \left( L_C(\phi(t,j),u_{ C}(t,j))
+ \frac{d V}{d t}(\phi(t,j)\pn{,u_C(t,j)}) \right)dt
\end{split}
\end{equation*}
from where we have
\begin{equation*}
\begin{split}
0\geq  \int_{t_{j}}^{t_{j+1}} L_C(\phi(t,j),u_{ C}(t,j)) dt
+  V(\phi(t_{j+1},j))- V(\phi(t_j,j))
\end{split}
\end{equation*}
{Pick $(t^*,j^*) \in \dom (\phi,u)$.} Summing from $j=0$ to $j={j^*}$ 
we obtain 
\begin{equation*}
0
\geq\sum_{j=0}^{{j^*}}  \int_{t_{j}}^{t_{j+1}} L_C(\phi(t,j),u_{ C}(t,j)) dt
+\sum_{j=0}^{{j^*}} \left( V(\phi(t_{j+1},j))-V(\phi(t_j,j)) \right)
\end{equation*}
Then, solving for $V$ at the initial condition $\phi(0,0)$, we obtain
\begin{equation} \label{ContinCost2goihzol}
\begin{split} 
V(\phi(0,0))
\geq&\sum_{j=0}^{{j^*}}  \int_{t_{j}}^{t_{j+1}} L_C(\phi(t,j),u_{ C}(t,j)) dt
\\  &+ V(\phi(t_1,0))
+\sum_{j=1}^{{j^*}} \left( V(\phi(t_{j+1},j))-V(\phi(t_j,j)) \right) 
\end{split}
\end{equation}`
In addition, {if $j^*>0$,} 
adding \eqref{upperbounddiscrete} from $j=0$ to $j={j^*}-1$, 
we obtain
\begin{equation}
\begin{split} 
\sum_{j=0}^{{j^*}-1}V(\phi(t_{j+1},j)) 
\geq 
\sum_{j=0}^{{j^*}-1}L_D(\phi(t_{j+1},j),u_{ D}(t_{j+1},j))\nonumber
+\sum_{j=0}^{{j^*}-1} V(\phi(t_{j+1},j+1))
\nonumber
\end{split}
\end{equation}

Then, solving for $V$ at the first jump time, we obtain
\begin{equation}
\begin{split} 
   V(\phi(t_{1},0))
&\geq V( \phi(t_{1},1))
+ \sum_{j=0}^{{j^*}-1}L_D(\phi(t_{j+1},j),u_{ D}(t_{j+1},j)) 
\\
&\hspace{0.5cm}+\sum_{j=1}^{{j^*}-1}\left( V(\phi(t_{j+1},j+1))-V(\phi(t_{j+1},j)) \right) 
\end{split}
\label{DiscreteCost2goNOihzwincol}
\end{equation}
In addition, given that $\phi(0,0)=\xi$, lower bounding $V(\phi(t_{1},0))$ in (\ref{ContinCost2goihzol}) by the right-hand side of (\ref{DiscreteCost2goNOihzwincol}), we obtain 
\begin{equation*}
\begin{split} 
V(\xi)
\geq& 
  \sum_{j=0}^{{j^*}} \int_{t_{j}}^{t_{j+1}} 
  \hspace{-0.3cm}
  L_C(\phi(t,j),u_{C}(t,j))dt+ V(\phi(t_1,0)) 
+
\hspace{-0.1cm}
\sum_{j=1}^{{j^*}} \left( V(\phi(t_{j+1},j))-V(\phi(t_j,j)) \right)  
\\
\geq&
\sum_{j=0}^{{j^*}} \int_{t_{j}}^{t_{j+1}} 
\hspace{-0.3cm}
L_C(\phi(t,j),u_{C}(t,j))dt
+ V(\phi(t_{1},1))
+ \sum_{j=0}^{{j^*}-1}L_D(\phi(t_{j+1},j),u_{D}(t_{j+1},j))
 \\&
 +\sum_{j=1}^{{j^*}-1}\left( V(\phi(t_{j+1},j+1))-V(\phi(t_{j+1},j)) \right) 
+\sum_{j=1}^{{j^*}} \left( V(\phi(t_{j+1},j))-V(\phi(t_j,j)) \right) 
\end{split}
\end{equation*}
Since 
\begin{eqnarray*} \nonumber
&& V(\phi(t_{1},1))
+\sum_{j=1}^{{j^*}-1}\left( V(\phi(t_{j+1},j+1))-V(\phi(t_{j+1},j)) \right) \nonumber 
 \\ && \hspace{5.6cm} 
 +\sum_{j=1}^{{j^*}} \left( V(\phi(t_{j+1},j))-V(\phi(t_j,j)) \right) \nonumber
  \\
 &=& V(\phi(t_{{{j^*}}+1},{j^*})) +V(\phi(t_{1},1))\nonumber
+\sum_{j=1}^{{j^*}-1}\left( V(\phi(t_{j+1},j+1)) \right) \nonumber 
-\sum_{j=1}^{{j^*}} \left(V(\phi(t_j,j)) \right) \nonumber
\\ &=&
V(\phi(t_{{{j^*}}+1},{j^*})) 
\end{eqnarray*}
then we have
{
\begin{eqnarray*}\nonumber
V(\xi)
&\geq&\sum_{j=0}^{{j^*}} \int_{t_{j}}^{t_{j+1}} L_C(\phi(t,j),u_{C}(t,j))dt
+ \sum_{j=0}^{{j^*}-1}L_D(\phi(t_{j+1},j),u_{D}(t_{j+1},j))
 \\&&
+V(\phi(t_{{j^*}+1},{j^*}))
\nonumber
\end{eqnarray*}}

By taking the limit when \pnn{
  $(t_{{j^*}+1},{j^*})\rightarrow  \sup \dom \phi $},
we establish \eqref{eq:costboundVxi}.
Notice that if $j^*=0$, the solution $(\phi,u)$ is continuous and \eqref{eq:costboundVxi} reduces to 
{
\begin{equation*}\nonumber
  \begin{split}
    &
V(\xi)
\geq
\underset{t^* \rightarrow \sup_t \dom \phi}{\limsup}
\bigg\{
\int_{t_{0}}^{t^*} L_C(\phi(t,0),u_{C}(t,0))dt
+V(\phi(t^*,0))
\bigg\}.
  \end{split}
\end{equation*}}
On the other hand, if $t_{j^*+1}=0$ for all $j^*$, the solution $(\phi,u)$ is discrete and 
\eqref{eq:costboundVxi} reduces to 
{
\begin{equation*}
  \begin{split}
 &   V(\xi)
\geq
\underset{j^* \rightarrow \sup_j \dom \phi}{\limsup}
\bigg\{
 \sum_{j=0}^{{j^*}-1}L_D(\phi(0,j),u_{D}(0,j))
+V(\phi(0,{j^*}))
\bigg\}.
\end{split}
\end{equation*}}
%
\end{proof}
\begin{corollary}[Change of signs]
  If the inequalities in the conditions in Proposition \ref{Pp:ProofDerivation} are inverted, {namely, if ``$\leq$'' in \eqref{eq:upperboundhamiltonian} and \eqref{upperbounddiscrete} is replaced with ``$\geq$'',} then \eqref{eq:costboundVxi} holds with the inequality inverted. Likewise, if the conditions in Proposition \ref{Pp:ProofDerivation} hold {with} equalities, then \eqref{eq:costboundVxi} holds {with} equality.
  \label{Rmk:ProofDerivation}
\end{corollary}
{
{
\begin{lemma}[Solutions to $\HS_{\textup{max}}$]\label{Pp:sol2Hmaxsol2Hk}
  Any solution to $\HS_{\textup{max}}$ as in \eqref{Hmax} is a solution to $\HS\IfIncd{}{} $ as in \IfIncd{\eqref{Hinc}}{\eqref{Heq}}.
\end{lemma}
\begin{proof}
  Pick any solution $(\phi,u)$ to $\HS_{\textup{max}}$ from $\xi \in \Pi(C \cup D)$. Then,
  $\phi \in \mathcal{X}$, $u=(u_C,u_D) \in \mathcal{U}$,   $\textup{dom}\phi=\textup{dom}u$, and 
  \begin{itemize}
    \item $(\phi(0,0),u_C(0,0)) \in \overline{C}$ or $(\phi(0,0),u_D(0,0)) \in D$,
    \item For each $j\in \mathbb{N}$ such that $I^j_{\phi}$ has a nonempty interior $\textup{int} I^j_{\phi}$, we have, for all
    $t \in \textup{int} I^j_{\phi}$,
    \begin{equation*}
        (\phi(t,j),u_C(t,j)) \in C
    \end{equation*}  
    and, for almost all $t \in I^j_{\phi}$,
    \begin{eqnarray*}
        \frac{d}{dt}\phi(t,j)&=& 
        \underset{\chi \in F(\phi(t,j), u_C(t,j))\IfGG{\cap T_C(\phi(t,j))}{}}{\textup{argmax}}\IfGG{V^\circ(\phi(t,j),\chi)}{\left\langle \nabla V(\phi(t,j)),\chi \right\rangle}     \\ 
        &\subset& F(\phi(t,j), u_C(t,j))
    \end{eqnarray*}
    \item For all $(t,j)\in \dom \phi$ such that $(t,j+1)\in \dom \phi$,  
    \begin{eqnarray*}
	(\phi(t,j),u_D(t,j)) &\in& D, 
  \quad
	\phi(t,j+1) =
  \hspace{-0.5cm}
  \underset{\chi \in G(\phi(t,j),u_D(t,j))}{\textup{argmax}}
  \hspace{-0.2cm}
  V(\chi) 
   \subset G(\phi(t,j),u_D(t,j)) 
\end{eqnarray*}
\end{itemize}
which, according to Definition \ref{SolutiontocalH}, defines a solution pair to \eqref{Hinc}.
\end{proof}
\begin{lemma}[Existence of solutions to optimization problems in maximal system]\label{Lemma:sol2OptinHmax}
Consider the maximal system $\HS_{\max}$ as in \eqref{Hmax} with $(C,F,D,G)$ and $u=(u_1,u_2)$. 
If the maps $F(x,u_C)$ and $G(x,u_D)$ are compact for each  $(x,u_C) \in C$ and each $(x,u_D) \in D$, respectively, and
the function $V: \dom V \rightarrow \mathbb{R}$, $\dom V \supset \Pi(\overline C) \cup \Pi(D) \cup G(D) $ is continuous on $\dom V$  and \IfGG{locally Lipschitz on}{continuously differentiable on a neighborhood of} $\Pi(C)$, then there exists a solution to 
\begin{equation}\label{Hmaxflow}
 \IfGG{\underset{\chi \in F(x,u_{C}) \cap T_C(x)}{\max} V^\circ (x,\chi)}{\underset{f \in F(x,u_{C})}{\textup{max}}\left\langle \nabla V(x),f \right\rangle} 
\end{equation}
for each $(x,u_C) \in C$, 
and to 
\begin{equation}
 \underset{\chi \in G(x,u_D)}{\textup{max}}V(\chi)
\end{equation}
for each $(x,u_D) \in D$.
\end{lemma}
\begin{proof}
\pno{Given that $V$ is \IfGG{locally Lipschitz on}{continuously differentiable on a neighborhood of} $\Pi(C)$, \IfGG{the Clarke generalized directional derivative $V^\circ$}{the gradient $\nabla V(x)$} exists for \pnn{each} $x \in \Pi(C)$ \IfGG{and \pnn{each} $\chi \in F(x,u_C)$ with $(x,u_C) \in C$}{}.
Since $F(x,u_C)\cap T_C(x)$ is compact for each  $(x,u_C) \in C$ 
given that the map $F(x,u_C)$ is compact for each  $(x,u_C) \in C$, 
and \pnn{by} the function $\chi \mapsto \IfGG{V^\circ(x,\chi)}{\left\langle \nabla V(x),f \right\rangle }$ \pnn{being upper semicontinuous} for any selection $\chi \in F(x,u_C)$, then
\eqref{Hmaxflow} attains its maximum in $F(x,u_C) \cap T_C(x)$.}

Given that $V$ is continuous in $\Pi(\overline C) \cup \Pi(D) \cup G(D) $, 
and the map $G(x,u_D)$ is compact for each $(x,u_D) \in D$, 
then $g \mapsto V(g)$  attains its maximum in $G(x,u_D)$.
\end{proof}
\begin{definition}[Hybrid basic conditions]\label{def:hbc} Given a system $\HS$ and a feedback law $\kappa:=(\kappa_C,\kappa_D)$ with $\kappa:\reals^n \rightarrow \reals^{m_C} \times \reals^{m_D}$, the resulting closed-loop system from their composition,  denoted $\HS_\kappa$ and defined as in \eqref{Hkeq}, is said to satisfy the {hybrid basic conditions} if its data $(C_\kappa,F,D_\kappa,G)$ satisfies the following properties:
  \begin{enumerate}[label=\arabic*)]
    \item $C_\kappa$ and $D_\kappa$ are closed subsets of $\reals^n$;
    \item $x \mapsto F(x,\kappa_C(x))$ is outer semicontinuous and locally bounded relative to $C_\kappa$, every $x\in C_\kappa$ is such that $(x,\kappa(x)) \in  \dom F$, and $x \mapsto F(x,\kappa_C(x))$ is convex  for each $x\in C_\kappa$.
    \item $x \mapsto G(x,\kappa_D(x))$ is outer semicontinuous and locally bounded relative to $D_\kappa$, and every $x\in D_\kappa$ is such that $(x,\kappa_D(x)) \in  \dom G$.
  \end{enumerate}
\end{definition}
\begin{definition}[Tangent cone]\label{def:TanCon}
  The tangent cone to a set $S\subset \reals^n$ at a point $x \in \reals^n$, denoted $T_S(x)$, is the set of all vectors $w \in \reals^n$ for which there exist sequences $x_i \in S, \tau_i>0$ with $x_0 \rightarrow x, \tau_i \searrow 0 $, and $w = \lim_{i\rightarrow \infty} \frac{x_i-x}{\tau_i}$.
\end{definition}
The following result provides conditions that guarantee that solutions to a hybrid closed-loop system exist. Below, $T_C (x)$ denotes the tangent cone of $C_\kappa$ at $x$
\begin{proposition}[Existence of solutions to $\HS$] \label{Prop:ExtSolH}
  Consider a system $\HS$ and a feedback law $\kappa:=(\kappa_C,\kappa_D)$ with $\kappa:\reals^n \rightarrow \reals^{m_C} \times \reals^{m_D}$, and suppose the resulting closed-loop system from their composition, denoted $\HS_\kappa$ and defined as in \eqref{Hkeq} with data $(C_\kappa, F, D_\kappa, G)$, satisfies the hybrid basic conditions in Definition \ref{def:hbc}.
  Let $\xi \in C_\kappa \cup D_\kappa$ be arbitrary. 
  If $\xi \in D_\kappa$ or
  \begin{enumerate}[label=(VC)]
  \item there exists a neighborhood $\mathcal{O}$ of $\xi$ such that for every $x \in \mathcal{O} \cap C_\kappa$,
  \begin{equation*}
    F(x,\kappa_C(x))\cap T_C(x) \neq \emptyset
  \end{equation*}
  \end{enumerate}
  then there exists a nontrivial solution $\phi$ to $\HS_\kappa$ with $\phi(0,0)=\xi$. If $(VC)$ holds for every $\xi \in C_\kappa \setminus D_\kappa$, then there exists a nontrivial solution to $\HS_\kappa$ for every initial point in $C_\kappa \cup D_\kappa$.
\end{proposition}
 \begin{proposition}[Existence of solutions to $\HS_{\max}$]\label{Prop:ExistHmax}
   Consider 
   a feedback law $\kappa:=(\kappa_C,\kappa_D)$ with $\kappa:\reals^n \rightarrow \reals^{m_C} \times \reals^{m_D}$,
   a function $V: \dom V \rightarrow \mathbb{R}$, $\dom V \supset \Pi(\overline C) \cup \Pi(D) \cup G(D) $, that is \IfGG{continuous on $\dom V$ and locally Lipschitz on}{continuously differentiable on a neighborhood of} $\Pi(C)$,
  and the closed-loop system
  \begin{equation}
  \HS_{\max_\kappa}\left\{
  \begin{matrix}
  \dot{x} &\in&   
  \IfGG{\underset{\chi \in F(x,\kappa_{C}(x)) \cap T_C(x)}{\argmax} V^\circ (x,\chi)}
  {\underset{\chi \in F(x,\kappa_{C}(x))}{\textup{argmax}}\left\langle \nabla V(x),\chi \right\rangle}
   \quad  (x,\kappa_C(x)) \in C_\kappa 
  \\ 
  x^+ &\in& 
  \underset{\chi \in G(x,\kappa_D(x))}{\textup{argmax}}
  V(\chi)\hspace{1.2cm}   (x,\kappa_D(x)) \in D_\kappa
  \end{matrix} \right.
  %
  %
  \label{MaxSystCL}
    \end{equation}
    Suppose \pn{$\HS_{\max_\kappa}$
  satisfies the hybrid basic conditions in Definition \ref{def:hbc} }
  and 
    $F(x,\kappa_C(x))$ and $G(x,\kappa_D(x))$ are compact for each  $x \in C_\kappa$ and each $x \in D_\kappa$, respectively.
  Let $\xi \in C_\kappa \cup D_\kappa$ be arbitrary. 
  If $\xi \in D_\kappa$ or
  \begin{enumerate}[label=(VCM)]
  \item there exists a neighborhood $\mathcal{O}$ of $\xi$ such that for every $x \in \mathcal{O} \cap C_\kappa$,
  \begin{equation*}
    F(x,\kappa_C(x))\cap T_C(x) \neq \emptyset
  \end{equation*}
  \end{enumerate}
  then there exists a nontrivial solution $\phi$ to $\HS_{\max_\kappa}$ with $\phi(0,0)=\xi$. If $(VC)$ holds for every $\xi \in C_\kappa \setminus D_\kappa$, then there exists a nontrivial solution to $\HS_{\max_\kappa}$ for every initial point in $C_\kappa \cup D_\kappa$.
 \end{proposition}
 %
 %
  The proof, similar to the proof of Proposition \ref{Prop:ExtSolH}, follows from the derivations in \cite[Proposition 2.10, 6.10]{65nu}.
\pnn{ Notice that we can get the flow and jump map to satisfy the conditions in the hybrid basic conditions via Berge's maximum theorem.
 }
}{}
}
}
\end{document}